\documentclass[11pt]{book}
\usepackage{epic,eepic}
\usepackage{amscd,amsfonts,amsmath,amsxtra,amssymb}
\usepackage[all]{xy}
\input xypic
\newfont{\gothic}{eufm10 scaled 1100}
\newfont{\mediumgothic}{eufm10 scaled 1000}
\newfont{\smallgothic}{eufm10 scaled 900}
\newfont{\verysmallgothic}{eufm10 scaled 700}

\def \PP {{\mathbb P}}
\def \QQ {{\mathbb Q}}
\def \ZZ {{\mathbb Z}}

\def \NN {{\mathbb N}}
\def \cald {{\cal D}}
\def \calm {{\cal M}}
\newcommand{\mg}{{M_g}}
\newcommand{\mgbar}{{\overline{M}_g}}
\newcommand{\cmg}{{{\cal M}_g}}
\newcommand{\cmgbar}{{\overline{{\cal M}}_g}}

\newcommand{\mgbarquot}{{\overline{M}_{g^+,\mgamma}}}

\newcommand{\cmgbarquot}{{\overline{{\cal M}}_{g^+,\mgamma}}}
\newcommand{\cmgplus}{{\overline{{\cal M}}_{g^+,m}}}
\newcommand{\mgprimequot}{{M_{g^+,\mgamma}'}}
\newcommand{\mgprimebarquot}{{\overline{M}_{g^+,\mgamma}'}}

\newcommand{\cmgprimequot}{{{\cal M}_{g^+,\mgamma}'}}
\newcommand{\cmgprimebarquot}{{\overline{{\cal M}}_{g^+,\mgamma}'}}

\newcommand{\dpo}{{D({C_0;P_1})}}
\newcommand{\dpobar}{{\overline{D}({C_0;P_1})}}

\newcommand{\dpohat}{{\hat{D}({C_0;P_1})}}
\newcommand{\dpp}{{D({P_1})}}

\newcommand{\dppbar}{{\overline{D}({P_1})}}

\newcommand{\cdpo}{{{\cal D}({C_0;P_1})}}

\newcommand{\cdpobar}{{\overline{{\cal D}}({C_0;P_1})}}

\newcommand{\cdpohat}{{\hat{{\cal D}}({C_0;P_1})}}
\newcommand{\mgprime}{{M_{g^+,m}'}}

\newcommand{\mgprimebar}{{\overline{M}_ {g^+,m}'}}
\newcommand{\hgprime}{{H_{g^+,n,m}'}}

\newcommand{\hgprimequot}{{H_{g^+,n,\mgamma}'}}
\newcommand{\hgprimebar}{{\overline{H}_ {g^+,n,m}'}}

\newcommand{\hgprimebarquot}{{\overline{H}_ {g^+,n,\mgamma}'}}

\newcommand{\hgbarquot}{{\overline{H}_ {g^+,n,\mgamma}}}

\newcommand{\cmgprimebar}{{\overline{{\cal M}}_ {g^+,m}'}}

\newcommand{\thx}{{\Theta^\times}}
\newcommand{\thxbar}{{\overline{\Theta}{}^\times}}
\newcommand{\thxbarGamma}{{\overline{\Theta}{}^\times_\Gamma}}
\newcommand{\hx}{{H_{g^+,n,m}^\times}}
\newcommand{\hxbar}{{\overline{H}{}^\times_{g^+,n,m}}}
\newcommand{\hgammax}{{H_{g^+,n,\mgamma}^\times}}
\newcommand{\hgammaxbar}{{\overline{H}{}^\times_{g^+,n,\mgamma}}}
\newcommand{\hgammadetailx}{{{H}{}^\times_{g^+,n,m/\Gamma(C_0^-;{\cal Q})}}}
\newcommand{\hgammadetailxbar}{{\overline{H}{}^\times_{g^+,n,m/\Gamma(C_0^-;{\cal Q})}}}
\newcommand{\mx}{{M_{g^+,m}^\times}}
\newcommand{\cmx}{{{\cal M}_{g^+,m}^\times}}
\newcommand{\mxbar}{{\overline{M}{}^\times_{g^+,m}}}
\newcommand{\cmxbar}{{\overline{\cal M}{}^\times_{g^+,m}}}
\newcommand{\mgammax}{{M_{g^+,\mgamma}^\times}}
\newcommand{\mgammaxbar}{{\overline{M}{}^\times_{g^+,\mgamma}}}
\newcommand{\kx}{{K^\times(C_0;P_1)}}
\newcommand{\kxbar}{{\overline{K}{}^\times(C_0;P_1)}}
\newcommand{\kpx}{{K^\times(P_1)}}
\newcommand{\kpxbar}{{\overline{K}{}^\times(P_1)}}
\newcommand{\dx}{{D^\times(C_0;P_1)}}
\newcommand{\dxbar}{{\overline{D}{}^\times(C_0;P_1)}}
\newcommand{\dpx}{{D^\times(P_1)}}
\newcommand{\dpxbar}{{\overline{D}{}^\times(P_1)}}
\newcommand{\cdx}{{{\cal D}^\times(C_0;P_1)}}
\newcommand{\cdxbar}{{\overline{{\cal D}}{}^\times(C_0;P_1)}}
\newcommand{\cmgammax}{{{\cal M}_{g^+,\mgamma}^\times}}
\newcommand{\cmgammaxbar}{{\overline{{\cal M}}{}^\times_{g^+,\mgamma}}}
\newcommand{\cmgammadetailxbar}{{\overline{\cal M}{}^\times_{g^+,n,m/\Gamma(C_0^-;{\cal Q})}}}
\newcommand{\cmgammadetailx}{{{\cal M}{}^\times_{g^+,n,m/\Gamma(C_0^-;{\cal Q})}}}
\newcommand{\mgammadetailx}{{{M}{}^\times_{g^+,n,m/\Gamma(C_0^-;{\cal Q})}}}

\newcommand{\thxx}{{\hat{\Theta}^{\times}}}

\newcommand{\hxxbar}{{\hat{H}{}^\times_{g^{++},n,\mu}}}

\newcommand{\hgammaxxbar}{{\hat{H}{}^\times_{g^{++},n,\mu/\hat{\Gamma}}}}
\newcommand{\hgammadetailxxbar}{{\hat{H}{}^\times_{g^{++},n,\mu/\Gamma(C_0^{--};{\cal R})}}}
\newcommand{\dxxbar}{{\hat{D}{}^\times(C_0;P_1)}}
\newcommand{\kxxbar}{{\hat{K}{}^\times(C_0;P_1)}}
\newcommand{\mxxbar}{{\hat{M}{}^\times_{g^{++},\mu}}}
\newcommand{\cmxxbar}{{\hat{\cal M}{}^\times_{g^{++},\mu}}}
\newcommand{\mgammaxxbar}{{\hat{M}{}^\times_{g^{++},\mu/\hat{\Gamma}}}}
\newcommand{\cmgammaxxbar}{{\hat{{\cal M}}{}^\times_{g^{++},\mu/\hat{\Gamma}}}}
\newcommand{\cmgammadetailxxbar}{{\hat{\cal M}{}^\times_{g^{++},n,\mu/\Gamma(C_0^{--};{\cal R})}}}
\newcommand{\FDPxx}{{\hat{{\cal D}}{}^\times(C_0;P_1)}}

\newcommand{\doo}{{D^\circ(C_0;P_1)}}
\newcommand{\cdoo}{{{\cal D}^\circ(C_0;P_1)}}
\newcommand{\koo}{{K^\circ(C_0;P_1)}}
\newcommand{\hoo}{{H_{g^+,n,m}^\circ}}
\newcommand{\hoohat}{\hat{H}_{g^{++},n,\mu}^\circ}
\newcommand{\hoogamma}{{H_{g^+,n,\mgamma}^\circ}}
\newcommand{\hoogammahat}{{\hat{H}_{g^{++},n,\mu/\hat{\Gamma}}^\circ}}
\newcommand{\hoogammadetail}{{{H}{}^\circ_{g^{+},n,m/\Gamma(C_0^{-};{\cal Q})}}}
\newcommand{\hoogammahatdetail}{{\hat{H}{}^\circ_{g^{++},n,\mu/\Gamma(C_0^{--};{\cal R})}}}

\newcommand{\dii}{{D^\diamond(C_0;P_1)}}
\newcommand{\cdii}{{{\cal D}^\diamond(C_0;P_1)}}
\newcommand{\kii}{{K^\diamond(C_0;P_1)}}
\newcommand{\hii}{{H_{g^+,n,m}^\diamond}}
\newcommand{\hiihat}{\hat{H}_{g^{++},n,\mu}^\diamond}
\newcommand{\hiigamma}{{H_{g^+,n,\mgamma}^\diamond}}
\newcommand{\hiigammahat}{{\hat{H}_{g^{++},n,\mu/\hat{\Gamma}}^\diamond}}
\newcommand{\mii}{{M_{g^+,m}^\diamond}}
\newcommand{\miihat}{\hat{M}_{g^{++},\mu}^\diamond}
\newcommand{\miigamma}{{M_{g^+,\mgamma}^\diamond}}
\newcommand{\miigammahat}{{\hat{M}_{g^{++},\mu/\hat{\Gamma}}^\diamond}}

\newcommand{\cms}{{\cal M}}
\newcommand{\cmsbar}{{\overline{\cal M}}}
\newcommand{\mgamma}{{m/\Gamma}}

\newcommand{\cat}[1]{{\hbox{\gothic #1}}}

\newcommand{\st}[1]{{\hbox{\gothic #1}}}
\newcommand{\scat}[1]{{\hbox{\smallgothic #1}}}
\newcommand{\sscat}[1]{{\hbox{\verysmallgothic #1}}}

\newcommand{\cata}{\cat{A}}
\newcommand{\catb}{\cat{B}}
\newcommand{\catc}{\cat{C}}
\newcommand{\cate}{\cat{E}}
\newcommand{\catf}{\cat{F}}
\newcommand{\catg}{\cat{G}}
\newcommand{\cath}{\cat{H}}
\newcommand{\catm}{\cat{M}}

\newcommand{\catq}{\cat{Q}}

\newcommand{\catz}{\cat{Z}}

\newcommand{\stackf}{{\mathbb F}}

\newcommand{\Sch}{\text{Sch}}

\newcommand{\Cat}{{\text{Cat}}}
\newcommand{\TwoCat}{{\text{2-Cat}}}
\newcommand{\Af}{\text{Aff}}
\newcommand{\GrC}{\text{(Gp/$\catc$)}}
\newcommand{\Gp}{\text{(Gp)}}
\newcommand{\Gpc}{\text{(Gp/$\catc$)}}

\newcommand{\St}{\text{St}}
\newcommand{\Set}{\text{(Set)}}

\newcommand{\bq}{\mathbf q}
\newcommand{\calr}{{\cal R}}
\newcommand{\calq}{{\cal Q}}
\newcommand{\ms}{{{M}}}     					
\newcommand{\msbar}{\overline{{M}}}				
\newcommand{\mgstack}[1]{{\overline{\cal M}_{{#1}} }}		
\newcommand{\mgstackopen}[1]{{{\cal M}_{{#1}}}}			
\newcommand{\Hgn}{\overline{H}_{g,n,m}}
\newcommand{\Hgnopen}{{H}_{g,n,m}}
\newcommand{\oka}{{\cal O}}
\newcommand{\Hgmp}{\overline{H}_{g,n,\mgamma}}			
\newcommand{\FDP}{{\cdpobar}}					
\newcommand{\FDPtilde}{{\tilde{\cald}(C_0;P_1)}}	
\newcommand{\FDPopen}{{\cdpo}}					
\newcommand{\FDPhat}{{\cdpohat}}				
\newcommand{\FDPc}{{\overline{\cald}{}^\times(C_0;P_1)}}	
\newcommand{\FDPcopen}{{{\cald}{}^\times(C_0;P_1)}}	
\newcommand{\hilb}[1]{\overline{H}_{{{#1}}}}
\newcommand{\CQ}{{\cal Q}}
\newcommand{\CR}{{\cal R}}
\newcommand{\CN}{{\cal N}}
\renewcommand{\epsilon}{\varepsilon}
\renewcommand{\rho}{\varrho}
\renewcommand{\theta}{\vartheta}
\newcommand{\ob}{\text{\rm Ob}\,} 
\newcommand{\simto}{\stackrel{\sim}{\longrightarrow}}
\newcommand{\mor}{\text{\rm Mor}\,}
\newcommand{\Isom}{\text{\rm Isom}\,}
\newcommand{\ISOM}{{\cal I}som\,}
\newcommand{\Hilb}{\text{\rm Hilb}\,}
\newcommand{\Lax}{\text{\rm Lax}\,}

\newcommand{\Aut}{\text{\rm Aut}  }
\newcommand{\Stab}{\text{\rm Stab}}
\newcommand{\id}{\text{\rm id}}
\newcommand{\pr}{\text{\rm pr}}
\newcommand{\bpr}{\text{\rm \bf pr}}
\newcommand{\GL}{\text{\rm GL}}
\newcommand{\PGL}{\text{\rm PGL}}
\newcommand{\Spec}{\text{\rm Spec\,}}
\newcommand{\im}{\text{\rm im}}
\newcommand{\prlim}{\lim\limits_{\stackrel{\longleftarrow}{\mathbf R}}}
\newcommand{\Prlim}[2]{\lim\limits_{\stackrel{\longleftarrow}{{#1}}}{#2}}
\newcommand{\ebew}{\hfill$\Box$ \par}
\renewcommand{\Text}[1]{\mbox{\rm #1}}

\newcommand{\proof}{{\em \underline{Proof.}\quad}}
\newcommand{\proofof}[1]{{\em \underline{Proof of #1.}\quad}}
\newcommand{\function}[5]{\begin{array}{llcl}
#1 : & #2 & \rightarrow & #3\\
     & #4 & \mapsto     & #5 
\end{array} }

\newcommand{\funktor}[7]{\begin{array}{lccc}
#1 : & #2 & \rightarrow & #3\\[1mm]
     & #4 & \mapsto     & #5\\
     & #6 & \mapsto     & #7 
\end{array} }

\newtheorem{thm}{Theorem}[chapter]
\newtheorem{defi}[thm]{Definition}
\newtheorem{prop}[thm]{Proposition}
\newtheorem{ex}[thm]{Example}
\newtheorem{rk}[thm]{Remark}
\newtheorem{remark}[thm]{Remark}
\newtheorem{cor}[thm]{Corollary}
\newtheorem{lemma}[thm]{Lemma}
\newtheorem{notation}[thm]{Notation}
\newtheorem{abschnitt}[thm]{}
\newtheorem{construction}[thm]{Construction}

\renewcommand{\footnote}[1]{{}}
\begin{document}

%
%

%
%
%
\thispagestyle{empty}
\vspace*{3mm}
\begin{center}
{\huge \bf
One-dimensional substacks \\[4mm] 
of the moduli stack\\[7mm]
of Deligne-Mumford stable curves
}\\[25mm]

{\Large
Dr. rer. nat. J\"org Zintl}\\[80mm]

{\Large 
Habilitationsschrift (rev.)\\[4mm]
Kaiserslautern \\[4mm]
2005
}

\end{center}

\pagebreak
\thispagestyle{empty}
\hfill
\pagebreak

%
%


\section*{Abstract}\label{abstract}

{\em
It is well known that there is a stratification of the moduli space of  Deligne-Mumford stable curves of genus $g$ by subschemes 
\[ \mgbar =  \mg^{(0)} \cup \ldots \cup \mg^{(3g-3)} ,\]
where each $\mg^{(i)}$ is of codimension $i$, and defined as the locus of stable curves with precisely $i$ nodes. 
A connected  component $D$ of  $\mg^{(3g-4)}$ can be identified with a scheme $M'$, which is the moduli space of a moduli stack $\calm'$  of certain pointed stable curves of some genus $g'\le g$. 

Our objects of study are those one-dimensional substacks $\cald$ of the moduli stack $\cmgbar$  of Deligne-Mumford stable curves of genus $g$, which are defined as reductions of the preimages of components $D$ in the above stratification. Even though there is an isomorphism $M'\cong D$, we find that the stacks $\cald$ and $\calm'$ are in general not isomorphic. We construct explicitely a surjective morphism on dense and open substacks 
\[ \calm^\times \: \longrightarrow \: \cald^\times ,\]
and show that is representable, unramified and finite of degree equal to the order of a certain group of automorphisms. Using this, we can give a new description of the substacks $\cald$ generically as global quotients. We explain how the above morphism can be extended to the stacks corresponding to the natural compactifications of $M'$ and $D$, and we discuss when and how the quotient description of $\, \cald$ extends. 

Finally, it is indicated how the results can be generalized  to substacks defined by boundary strata of codimension less than $3g-4$. The main theorem, which relates the stack $\cald$ to  a quotient of a substack of a moduli stack of pointed stable curves of genus less or equal to $g$, carries over immediately.

In an  appendix, an extensive compendium on stacks is provided, with a particular focus on  quotient stacks, including a few as yet unpublished observations. 
}

\pagebreak

\tableofcontents

%
%

\setcounter{chapter}{-1}
\chapter{Introduction}

\thispagestyle{plain}

The study of moduli spaces often turns out to be technically demanding. At the same time, these challenges may lead to the developement of intricate and intriguing new theories. The theory of geometric invariants is just one example of many.

Algebraic stacks are another instance of highly sophisticated techniques to approach moduli problems. They prove to be increasingly effective, as many questions on moduli can  naturally be formulated in the language of categories and functors. Often answers turn out to be much simpler, when phrased in terms of stacks, than they were in terms of schemes. At least, this is true  formally. The simplification gained by the use of these tools is paid for with an increased  effort necessary  to built up the formal framework.  As J. Harris  and I. Morrison put it: ``the initial learning curve is steep...'', \cite[p. 140]{HM}. What one gets in the end is a powerful instrument indeed, but also a concept which is less intuitively manageable than the concept of schemes, say. The starting point of our work is analyse one instance, where stacks do not behave in the way one would expect them to from looking at the corresponding moduli schemes only. 

This is a technical paper, but I am starting with a straightforward geometrical question. Consider the moduli space $\mgbar$ of Deligne-Mumford stable curves of a given genus $g$. Contained in it is a subscheme $D$ as the locus of stable curves sharing a certain property. Suppose that this property can be formulated as a second moduli problem of its own right, with a moduli space $M'$ isomorphic to $D$. Now compare the two stacks, one defined as the reduced substack $\cald$ of the moduli stack $\cmgbar$ of Deligne-Mumford stable curves, which lies   over $D$, the other as the stack $\calm'$ corresponding to the second moduli problem. What is the relation of $\cald$  and $\calm'$? 

In general, the two stacks will not be isomorphic. The moduli stack $\cald$  still  ``parametrizes''  those features of objects represented by points of $D$, which do not vary within $D$. In our case, it will turn out that $\cald$  and $\calm'$ differ generically by the action of a finite group $A$, which acts trivially on the moduli space  $M'$. While there is obviously an isomorphism between the schemes $M'$ and $M'/A$, the stacks $\calm'$, $\calm'/A $ and $ \cald$ may be different. Here, the intuition  aquired from working with schemes can no longer be used as a guide. 

This work is organized as follows. The first chapter illustrates in an example a question asked by A. Hirschowitz. Fix a stable curve $C_0$ of some given genus $g$, which has the maximal number of nodes that is possible, and fix one of its nodes $P_1$. Now look at such deformations of $C_0$, which preserve all nodes except $P_1$. This distinguishes an one-dimensional substack $\cdpobar$ of the moduli stack $\cmgbar$ of all stable curves of genus $g$. What can be said about $\cdpobar$?

In the second chapter, I recall some basic definitions and facts about pointed stable curves, as far as I need them later on. This includes a slight variation of Edidin's theorem \cite{Ed} about the quotient representation of the moduli stack of Deligne-Mumford stable curves. 
I also present the notion of an $\mgamma$-pointed stable curve of genus $g$, where $\Gamma$ is a subgroup of the permutation group on $m$ elements. Such a curve should be thought of as an $m$-pointed stable curve, but where the labels of its marked points are determined only up to permutations in $\Gamma$. The actual definiton \ref{mgdef} is somewhat more technical to get the behaviour of curves in families right.  

In chapter 3, the basic objects of study are introduced. The locus $\mgbar^{(3g-4)}$ of stable curves of genus $g$ with at least $3g-4$ nodes is a one-dimensional subscheme of $\mgbar$, which decomposes into several irreducible components. To specify one such component, I am fixing  one stable curve $C_0$ with $3g-3$ nodes, together with one node $P_1$ on it. Deformations of $C_0$, preserving all nodes but $P_1$, distinguish a unique irreducible component $\dpobar$ of $\mgbar^{(3g-4)}$. After a characterization of those stable curves, which are represented by points of $\dpobar$, a description of $\dpobar$ in terms of a moduli space $\mgprimebar$ of certain pointed stable curves is given. This characterization leads also to the definition of a substack $\cdpobar$ of the stack $\cmgbar$, which has $\dpobar$ as its moduli space. This substack $\cdpobar$ is the object I am mainly interested in. 

The main theorem \ref{MainThm} is stated in chapter 4. My result concerns the open substack $\cdpo$ of $\cdpobar$, which relates to those stable curves with exactly $3g-4$ nodes, and describes it generically as a quotient stack. More importantly, this exhibits a dense open substack of $\cdpo$ as a quotient of a very simple moduli stack. There may be other ways to obtain this result. But on reading Edidin's paper \cite{Ed} I liked his approach, which made use of quotient stacks of Hilbert schemes, very much, and I decided that I wanted to apply this technique to my problem as well.  

Chapter 5 deals with a variation of this approach in order to include stable curves with $3g-3$ nodes as well. It turns out that one has to take care of each of the ``extremal'' curves separately. Let $\dpohat$ denote the one-point partial compactification of $\dpo$ at the point representing $C_0$, and let $\cdpohat$ denote the corresponding substack of $\cdpobar$. Then theorem \ref{49}, which describes the stack $\cdpohat$ by a quotient stack, can be thought of as an equivariant version of the main theorem with respect to the action of the group of automorphisms of $C_0$. 

The question of the compatibility and the glueing of the partial compactifications is discussed in chapter 6. This amounts mainly to the study of various groups of automorphisms associated to pointed stable curves. 

Finally, in chapter 7 all local information is collected to obtain some global results about the closed substack $\cdpobar$. In particular, in proposition \ref{P62} a criterion is given to decide when the quotient description of $\cdpo$ extends to the partial compactification $\cdpohat$. Essentially this is determined by the group of automorphisms of $C_0$. I also take a closer look at the natural morphisms from the moduli stack $\cmgprimebar$ of certain pointed stable curves to $\cdpobar$, and some of its properties. A short summary of the results is given in remark \ref{R722}. 

Chapter 8 is intended to illustrate some of the ideas presented here. The case of genus $g=3$ is studied in more detail. I also explain how my ideas generalize to a description of boundary stata of $\cmgbar$ of dimensions greater than one. Finally relations to intersection theory are briefly indicated. 

I am including  an extensive appendix on algebraic stacks in general. The purpose of the first and general part is to introduce the reader to the notations and conventions I am using. It seem that still everyone writing about stacks has to do this, although the situation is improving,  thanks to the efforts of people like Laumon and Bailly \cite{LM} and Vistoli \cite[appendix]{Vi}, \cite{Vi3}. 
The second part on quotient stacks contains also a number of observations, which I have not seen mentioned  elsewhere. Still, I can hardly claim much originality. Certainly the experts will find nothing here that surprises them, and I could as well have left much of it as an exercise to the reader. However, some of the facts there have never been writen down explicitely, or no proof has been given.  It seems important to my work to have them easily accessible for reference. 

My studies on these issues were initiated after a number of discussions with A. Hirschowitz on the subject of moduli stacks, for which I am grateful to him. Without D. Edidin's paper \cite{Ed} on the quotient construction of the moduli stack of Deligne-Mumford stable curves it could not have been carried out the way it was. And most of all I am indebted to G. Trautmann for his support while I was completing my work. This is a revised version of my original thesis, correcting a few mistakes, and clarifying one or two unprecise statements. Following the suggestion of some early readers, I included more guiding prose, so that it will be easier to keep track of what is going on between the technical bits. One major change is the modernization of my notation. 

\thispagestyle{plain}

\hfill\begin{tabular}{l}
J\"org Zintl\\
 Kaiserslautern\\
 2 February, 2005\\
 30 May,  2006
\end{tabular}


%
%

\chapter{Motivation}\label{intro}

\begin{abschnitt}\em
Our objects of study are stable curves of some genus $g$ over an
algebraically closed ground field $k$, as defined in \cite{DM}. Unless stated otherwise we will always assume $g \ge 3$.  Let us recall some well known facts first, most of which can be found for example in \cite{HM}. 

There exists a coarse moduli space of stable curves of genus $g$, which is usually denoted by $\mgbar$\label{n1}. It is an irreducible projective
variety of dimension $3g-3$. The open subscheme corresponding to
smooth stable curves is denoted by $\mg$. 

A $k$-valued point of $\mgbar$ represents a stable curve of genus $g$
over $k$, with nodes as its only singularities. The number of its
nodes is at most $3g-3$, and there are only finitely many curves
attaining this number of nodes. 

If $C_0$ is a stable curve over $k$ with $3g-3$ nodes, then all of
its irreducible components are rational curves, with at most three
singularities lying on each of them. The normalization of an
irreducible component is always a rational curve with 3 distiguished points. There
are exactly $2g-2$  irreducible components of such a curve $C_0$. 

Let $C_0$ be a stable curve of arithmetic genus $g$, with $\delta$
nodes. If $C_1,\ldots,C_n$ are the irreducible components of $C_0$, each with
geometric genus $g_i$ for $i=1,\ldots ,n$, then
\[ g = \sum_{i=1}^n (g_i-1) +\delta+1 .\]

The number of irreducible components of $C_0$ is at most $2g-2$, and the maximal number of components  is attained if and only if $C_0$ has also the maximal number of singularities.

The locus of stable curves with at least $k$ nodes forms a closed 
subscheme $\mgbar^{(k)}$\label{n2}
of codimension $k$ in $\mgbar$. In particular, the boundary $\mgbar
\setminus \mg$ decomposes into divisors $\Delta_0$ and $\Delta_1,
\ldots , \Delta_{[\frac{g}{2}]}$, where $\Delta_0$ is the closure of
the locus of  stable curves consisting of one singular irreducible  component,
and $\Delta_i$ is the closure of the locus of  stable curves consisting of
two non-singular irreducible  components  of
genus $i$, and genus $g-i$ respectively, meeting transversally in one point.
\end{abschnitt} 

\begin{abschnitt}\em
Fix a stable  curve $C_0$ with $3g-3$ nodes, represented by a
point $[C_0] \in \mgbar$, together with an enumeration of its
singularities $P_1,\ldots ,P_{3g-3}$. Choose one of the nodes $P_1$,
say. 

The objects we want to study in this paper are stable curves $C\rightarrow S$,
which fibrewise have the same singularities as $C_0$, except that the node $P_1$
is allowed to ``smoothen out''. Before we make this rigorous in chapter \ref{sect3}, we will illustrate the geometric idea in the case of curves of
genus $3$.  Essentially, what we want to consider are  stable curves $C\rightarrow S$ of genus $g$, such that for all closed points $s\in S$ the fibre $C_s$ has at least $3g-4$
nodes, and such that  there exists a 
subcurve $C_s^{-} \subset C_s$ with 
\[ C_s \supset C_s^{-} \quad \cong \quad C_0^{-} \subset C_0 ,\]
where  $C_0^{-}$ is the reduced subscheme of $C_0$ consisting of all
irreducible components of $C_0$ which do not contain $P_1$. Similarly, we define $C_0^{+} \subseteq  C_0$ as the union of all
irreducible components of $C_0$ containing $P_1$. 

Define $\dpobar \subset \mgbar$ as the (closed, reduced) subscheme
determined by all such stable curves $C \rightarrow S$.  In
fact, as we will see below, 
$\dpobar$ is an irreducible component of $\mgbar^{(3g-4)}$, so it is 
1-dimensional,  and even a smooth curve. Let $\dpo$ denote the open subscheme of $\dpobar$ representing stable curves with exactly $3g-4$ nodes. 
\end{abschnitt}

\begin{ex}\em\label{E11}
Consider stable curves of genus $g=3$. So the maximal number of nodes
is $6$, and the picture below gives a schematic drawing of an example of a
curve $C_0$, where this 
number is attained.

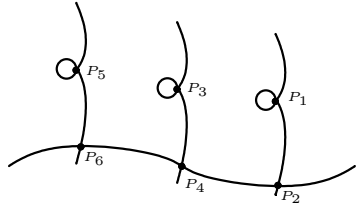
\begin{center}
\setlength{\unitlength}{0.00016667in}
\begingroup\makeatletter\ifx\SetFigFont\undefined%
\gdef\SetFigFont#1#2#3#4#5{%
  \reset@font\fontsize{#1}{#2pt}%
  \fontfamily{#3}\fontseries{#4}\fontshape{#5}%
  \selectfont}%
\fi\endgroup%
{\renewcommand{\dashlinestretch}{30}

}

\end{center}

We fix an enumeration of the nodes as indicated. The subcurve
$C_0^{+}$ of $C_0$ is a stable curve of arithmetic genus $g^{+} = 1$, with
one node, and 
with one marked point $P_2$ different from the node, as shown on the right hand side of  the next 
picture.

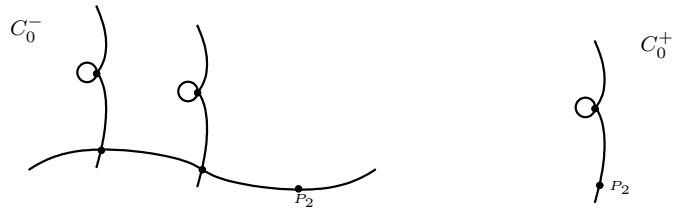
\begin{center}
\setlength{\unitlength}{0.00016667in}
\begingroup\makeatletter\ifx\SetFigFont\undefined%
\gdef\SetFigFont#1#2#3#4#5{%
  \reset@font\fontsize{#1}{#2pt}%
  \fontfamily{#3}\fontseries{#4}\fontshape{#5}%
  \selectfont}%
\fi\endgroup%
{\renewcommand{\dashlinestretch}{30}
%
}
\end{center}

The subcurve $C_0^{-}$ consists of three rational curves, with
a total of four nodes, and one marked point, as on the left hand side of the previous picture. 

A general closed point
of $\dpobar$ will represent a stable curve with $3g-4=5$ nodes, as shown 
 in the next picture. The numbers given in the picture indicate the
geometric genus of the irreducible component next to them. 

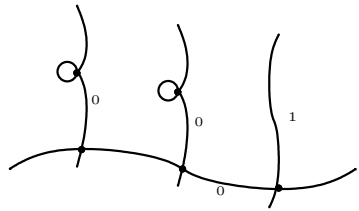
\begin{center}
\setlength{\unitlength}{0.00016667in}
\begingroup\makeatletter\ifx\SetFigFont\undefined%
\gdef\SetFigFont#1#2#3#4#5{%
  \reset@font\fontsize{#1}{#2pt}%
  \fontfamily{#3}\fontseries{#4}\fontshape{#5}%
  \selectfont}%
\fi\endgroup%
{\renewcommand{\dashlinestretch}{30}

}
\end{center}

There is no other point  
on $\dpobar$ representing a curve with the maximal number of nodes apart
from $[C_0]$.

In fact, since there is an invariant part $C^{-} \cong C_0^{-}$ of 
each curve $C$ with $[C] \in \dpobar$,  all of the information about $[C]$ is contained in the 
subcurve $C^{+}$, which is the closure of the complement of $C^{-}$ in $C$. 
In other words, a curve in $\dpobar$, or rather a point of the moduli space representing
that curve, is uniquely determined by a curve of genus $g^{+}=1$, with
one  marked point. But such curves  are in one-to-one correspondence with
the points of the moduli space $\msbar_{1,1}$  of stable
curves of genus $g^{+}=1$ with 1 marked point. So there is an
isomorphism
\[ \msbar_{1,1} \: \cong \: \dpobar.\]
The unique point of $\msbar_{1,1}$ representing a singular curve corresponds to the point $[C_0] \in \dpobar$. Obviously, this restricts to an isomorphism $\ms_{1,1} \cong \dpo$. 

The scheme $\msbar_3$ is the moduli space of an algebraic stack
$\mgstack{3}$, the stack of Deligne-Mumford stable curves of genus $3$. This is known
to be a smooth Deligne-Mumford stack, which is projective,
irreducible, 
and a quotient stack.

Similarly, $\msbar_{1,1}$ is the moduli space of the stack
$\mgstack{1,1}$ of stable
curves of genus 1 with one marked point.  Define $\cdpobar$
as the reduction of the preimage substack of $\dpobar$ in $\mgstack{3}$. We obtain a commutative diagram
\[ \diagram
\mgstack{1,1} \dto \xdotted[r]|<{\rotate\tip}|>\tip^{?} & \FDP 
\rto|<\hole|<<\ahook
\dto  &
\mgstack{3} \dto \\ 
\msbar_{1,1} \rto^\cong &\dpobar  \rto|<\hole|<<\ahook &  \msbar_3.
\enddiagram\]

It is a natural question to ask whether the dotted arrow exists, in either direction, and
if it exists, whether it is an isomorphism of stacks.

The answer to the second question is no. There is however a morphism 
$\mgstack{1,1} \rightarrow  \FDP$ making the above
diagram commutative, which is surjective and of finite degree. Let $\cdpo$ denote the reduced preimage substack of $\dpo$ in $\mgstack{3}$.  Then there is a morphism from $\cdpo$ to $\cms_{1,1}$, such that the composition with the restriction of the above morphism is the identity on $\cms_{1,1}$, but  it is not an isomorphism. To discuss relations of stacks like these  in general is  our main goal. 

Loosely speaking, encoded in  the stack $\cdpobar$ are not only
isomorphism classes of stable curves of genus 3, as it is the case for the moduli spaces,  but also information
about their automorphisms. There are usually less
automorphisms on a 
one-pointed curve of genus 1 than on an ambient curve of genus 3 which has  the
maximal number of nodes, and this is
exactly the piece of information we 
are missing. However, as we will see, there is a way to add some extra data to the
stack $\cms_{1,1}$ to recover $\cdpo$. 
\end{ex}

\begin{ex}\em\label{E12}
In general however the situation will be more complicated. We can see
this in the same example as above, if we consider $\overline{D}(C_0;P_2)$ instead of
$\dpobar$, i.e. stable curves which maintain  the singularities $P_1$
and $P_3,\ldots,P_6$.  

Now the subcurve
$C_0^{+}$ of $C_0$ is a stable  curve of genus $g^{+} = 1$, with
two nodes, and 
with two marked points $P_4$ and $P_6$ different from the nodes, as shown in the 
picture below.
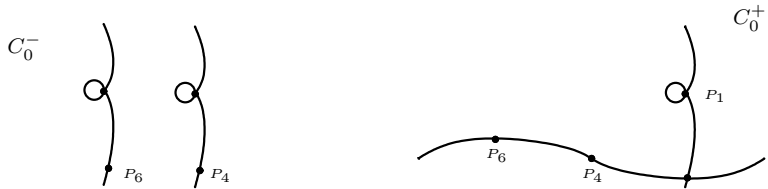
\begin{center}
\setlength{\unitlength}{0.00016667in}
\begingroup\makeatletter\ifx\SetFigFont\undefined%
\gdef\SetFigFont#1#2#3#4#5{%
  \reset@font\fontsize{#1}{#2pt}%
  \fontfamily{#3}\fontseries{#4}\fontshape{#5}%
  \selectfont}%
\fi\endgroup%
{\renewcommand{\dashlinestretch}{30}

}
\end{center}

The curve $C_0^{-}$ consists of two disjoint rational curves, each with
exactly one node, and with one marked point different from the node. 

A general closed point
of $\overline{D}(C_0;P_2)$ will represent a stable   curve $C$ with $3g-4$ nodes, as 
shown on the left hand side of the next picture.  There are two special
points on $\overline{D}(C_0;P_2)$. One of them represents $C_0$, the other one the
curve $C_1$ shown on the right hand side below. 

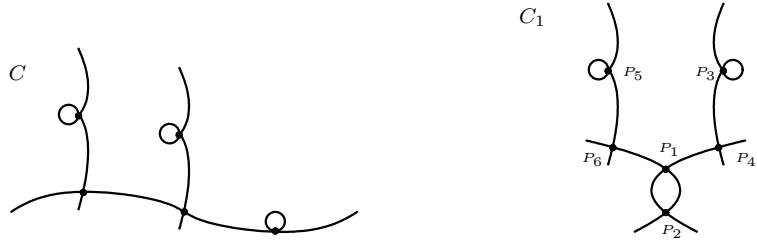
\begin{center}
\setlength{\unitlength}{0.00016667in}
\begingroup\makeatletter\ifx\SetFigFont\undefined%
\gdef\SetFigFont#1#2#3#4#5{%
  \reset@font\fontsize{#1}{#2pt}%
  \fontfamily{#3}\fontseries{#4}\fontshape{#5}%
  \selectfont}%
\fi\endgroup%
{\renewcommand{\dashlinestretch}{30}

}
\end{center}

A curve in $\overline{D}(C_0;P_2)$  is determined by deformations of  $C_0^{+}$, preserving the node at $P_1$, and is hence determined by a singular curve of genus
$g^{+}=1$, with 
at least one node and two marked points. These points form a 
divisor $\Delta_0$ on the moduli space $\msbar_{1,2}$  of stable
curves of genus $g^{+}=1$ with 2 marked points. In fact, there is the
forgetful morphism $f : \msbar_{1,2} \rightarrow \msbar_{1,1}$,
essentially forgetting the second marked point,  and
$\Delta_0$ is the fibre over the unique point of $\msbar_{1,1}$ which
represents a singular curve.  

Note, however, that there is no
distinguished ordering of the nodes of $C_0$. Actually, there is an
automorphism of $C_0$, which interchanges the points $P_4$ and
$P_6$. Therefore the marking of the two points on the curve $C_0^{+}$
is only determined up to permutation.  This is what we will later on call an $\mgamma$-pointed stable curve, for $m=2$, and $\Gamma$ the full permutation group on two elements. 
So what we obtain is an
isomorphism 
\[ \msbar_{1,2}^{\, \diamond} \supset \Delta_0^\diamond \quad \cong \quad  \overline{D}(C_0;P_2) \subset \msbar_3, \]
where $ \Delta_0^\diamond$ is the image of $ \Delta_0$ in $\msbar_{1,2}^{\, \diamond} := 
\msbar_{1,2}/ {\mathbb Z}_2$. Again, the question arises how the associated moduli stacks relate to each other.
\end{ex}

\chapter{Moduli of pointed stable curves}\begin{abschnitt}\em
In this chapter we want to collect some basics about pointed stable curves, so that we can refer to them later on. We also introduce the notion of an $\mgamma$-pointed curve, which will play an important role in our considerations.  
Let us first recall the definition of pointed stable curves as it is given in the paper of Knudsen  \cite[def. 1.1]{Kn}.
\end{abschnitt}

\begin{defi}\em
An {\em $m$-pointed stable curve of genus $g$} is a flat and proper
morphism $f : C \rightarrow S$ of schemes, together with $m$ sections
$\sigma_i : S \rightarrow C$, for $i=1,\ldots,m$, such that for all closed points $s\in
S$ holds
\begin{itemize}
\item[$(i)$] the fibre $C_s$ is a reduced connected algebraic curve with at most
ordinary double points as singularities; 

\item[$(ii)$] the arithmetic genus of $C_s$ is $\dim
H^1(C_s,\oka_{{C_s}}) = g$;

\item[$(iii)$] for $1\le i \le m$, the point $\sigma_i(s)$ is a smooth
point of $C_s$;

\item[$(iv)$] for all $ 1\le i,j \le m$ holds $\sigma_i(s) \neq
\sigma_j(s)$ if $i \neq j$;

\item[$(v)$] The number of points where a nonsingular rational
component $C_s'$ of $C_s$ meets the rest of $C_s$ plus the number of
points $\sigma_i(s)$ which lie on $C_s'$ is at least $3$.
\end{itemize}
\end{defi}
 
\begin{rk}\em
Condition $(v)$ guarantees that the group $\Aut(C_s)$ of automorphisms of a fibre $C_s$ is finite. Because of condition $(iv)$ the permutation group $\Sigma_m$\label{n3} acts faithfully as well on the set of sections $\{\sigma_1,\ldots,\sigma_m\}$, as on the set of marked points $\{\sigma_1(s),\ldots,\sigma_m(s)\}$ of each fibre $C_s$. 
\end{rk} 

\begin{rk}\em\label{M23}
In \cite[def. 3.1]{Kn} there is also the notion of {\em prestable curves}, where the condition on the connectedness of the fibres $C_s$ and condition $(v)$ are dropped. However, for us prestable curves will always satisfy condition $(v)$, so that all connected components of a prestable curve over $\Spec(k)$ are stable curves.   
\end{rk}

\begin{abschnitt}\em
Let $f: C\rightarrow S$ be an $m$-pointed  stable curve of genus $g$ such that  $2g-2+m > 0$. The sections $\sigma_1,\ldots, \sigma_{m}$ determine effective
divisors $S_1,\ldots ,S_{m}$ on $C$. 

There is a  canonical invertible 
dualizing sheaf on $C$, which is  denoted by $\omega_{C/S}$. In \cite[cor. 1.9, cor 1.11]{Kn} it is  shown that for
$n\ge 3$ the sheaf $\left(\omega_{C/S}(S_1+\ldots+S_{m})\right)^{\otimes n}$ is
relatively very ample,   and furthermore that $f_\ast\left(\left(\omega_{C/S}(S_1+\ldots+S_{m})\right)^{\otimes n}\right)$ is locally free of
rank $(2g-2+m)n -g+1$. 

Consider the Hilbert scheme $\Hilb_{{\PP^N}}^{{P_{g,n,m}}}$ of curves $ C \rightarrow \Spec(k)$ embedded in $\PP^N$ with Hilbert polynomial $P_{g,n,m} := (2g-2+m)nt-g+1$, where  $N :=(2g-2+m)n -g$. 

The Hilbert scheme of a simple point in $\PP^N$ is of course $\PP^N$ itself. The incidence condition of $m$ points lying on a curve $C$ defines a closed subscheme 
\[ I \subset 
\Hilb_{{\PP^N}}^{{P_{g,n,m}}} \times (\PP^N)^m .\]
There is an open subset $U \subset I$ parametrizing curves, where the $m$ points are pairwise different, and smooth points of $C$. Finally, there is a closed subscheme  of $U$, which represents such embedded $m$-pointed stable curves $C \rightarrow \Spec(k)$ of genus $g$, where the embedding is determined by an isomorphism 
\[  \left(\omega_{C/\Spec(k)}(Q_1+\ldots+Q_{m})\right)^{\otimes n} \cong \oka_{{\PP^N}}(1) | C. \]
Here, $Q_1,\ldots,Q_m$ denote the marked points on $C$. 
The scheme constructed in this way 
is in fact a quasi-projective subscheme
\[ \overline{H}_{g,n,m} \subset \Hilb_{{\PP^N}}^{{P_{g,n,m}}} \times (\PP^N)^m \label{n3a}\]
of the Hilbert scheme. There is also a subscheme ${H}_{g,n,m }\subset\overline{H}_{g,n,m}$ corresponding to non-singular stable curves.  Compare \cite[section 2.3]{FP} for the more general situation of moduli of stable maps.  

There is a one-to-one correspondence  
between morphisms $\theta: S \rightarrow \Hgn$ and $m$-pointed stable curves $f :
C\rightarrow S$ together with trivializations of Grothendieck's associated projective bundle 
\[ \PP f_\ast\left(\left(\omega_{C/S}(S_1+\ldots+S_{m})\right)^{\otimes n}\right) \: 
\simto \: \PP^N\times S. \]
In particular, the curve $f: C\rightarrow S$ can be considered as embedded into $\PP^N\times S$, such that the diagram 
\[\diagram
C \rrto|<\hole|<<\ahook \drto_f &&\PP^N \dlto^{\pr_2}\times S \\
&S 
\enddiagram \]
commutes. 

An element of the group $\PGL(N+1)$ acts on $\PP^N$, changing embeddings of $C$ by an isomorphism. Hence there is a natural action of $\PGL(N+1)$  on $\Hilb_{{\PP^N}}^{{P_{g,n,m}}}$, and thus on $\Hilb_{{\PP^N}}^{{P_{g,n,m}}} \times (\PP^N)^m$, which restricts to an action on 
$\Hgn$. Furthermore, the symmetric group $\Sigma_m$ acts on $(\PP^N)^m$ by permutation of the coordinates, which is equivalent to the permutation of the $m$ sections of marked points. Hence $\Sigma_m$ acts freely on $\Hgn\times (\PP^N)^m$. The combined action of  $\PGL(N+1) \times \Sigma_m$ on $\Hgn$ shall be written as an  action from the right. 
\end{abschnitt}

\begin{abschnitt}\em
Let $f: C_0 \rightarrow \Spec(k)$ be an $m$-pointed  stable curve over an algebraically closed field $k$, so in
particular an algebraic curve of genus $g$ with nodes as its only
singularities,  and $m$ marked points $P_1, \ldots , P_{m}$.  If we
fix one (arbitrary)  embedding of $C_0$ into $\PP^N$ via an isomorphism 
\[ \PP f_\ast(\omega_{C_0/k}(P_1+\ldots+P_m))^{\otimes n} \cong \PP^N ,\]
then this
distinguishes a $k$-valued point 
$[C_0]\in \Hgn$\label{n3b}. For the group of automorphisms of $C_0$ respecting the marked points, i.e. automorphisms which map each marked point to itself, there is a
natural isomorphism
\[ \Aut(C_0) \cong \Stab_{\PGL(N+1)}([C_0]),\]
that is,   an isomorphism with the subgroup of those  elements in $\PGL(N+1)$, which stabilize  the
point $[C_0]$. Recall that $\Sigma_m$ acts on $\Hgn$ without fixed points. The stabilizer of $[C_0]$ in $\PGL(N+1)$ is in a natural way a subgroup of  the stabilizer of $[C_0]$ in $\PGL(N+1)\times \Sigma_m$. 

In general, elements $\gamma \in \PGL(N+1)$ are in
one-to-one correspondence with isomorphisms
\[ \gamma : \: \: C_0 \rightarrow C_\gamma\]
of embedded $m$-pointed stable curves, where 
 $C_\gamma$ is the curve
represented by the  point $[C_\gamma] := [C_0] \cdot \gamma$ in
$\Hgn$.
\end{abschnitt}

\begin{abschnitt}\em
Generalizing Gieseker's construction, the moduli space $\msbar_{g,m}$\label{n4}
of $m$-pointed stable
curves of genus $g$ is constructed as the GIT-quotient of the action
of $\PGL(N+1)$ on $\Hgn$, see \cite{Gs}, \cite{DM} and \cite{HM} for
the case of $m=0$, and compare \cite[Remark 2.4]{FP} for the more general case of stable maps. Note that the notion of stability of embedded $m$-pointed stable curves is compatible with the notion of stability of the corresponding points in $\Hgn$  with respect to the group action. 
In particular, there is a categorical quotient 
\[ \Hgn \rightarrow \msbar_{g,m} .\]
The construction of the moduli space is independent of the choice of $n$, provided it is large enough, and
hence independent of $N$.
For a stable algebraic curve $C_0 \rightarrow \Spec(k)$, the fibre over the point
in $\msbar_{g,m}$ representing it is isomorphic to the quotient $ \Aut(C_0) \setminus
\PGL(N+1)$.

The scheme $\msbar_{g,m}$ is at the same time a moduli space for the stack
$\mgstack{g,m}$\label{n5} of $m$-pointed stable curves of genus $g$.
\end{abschnitt}

\begin{defi}\em
The {\em moduli stack $\, \mgstack{g,m}$ of $\,m$-pointed stable curves of genus $g$} is the stack defined as the category fibred in groupoids over  the category of schemes, where for a scheme $S$ the objects in the fibre category $\mgstack{g,m}(S)$ are the $m$-pointed stable curves of genus $g$ over $S$.

 Morphisms in $\mgstack{g,m}$ are given as follows. Let $f: C \rightarrow S$ and $f': C' \rightarrow S'$ be $m$-pointed stable curves of genus $g$, i.e. objects of $\mgstack{g,m}(S)$ and $\mgstack{g,m}(S')$, respectively, with sections $\sigma_i:S \rightarrow C$ and $\sigma_i':S' \rightarrow C'$ for $i=1,\ldots,m$. A morphism from $f: C \rightarrow S$ to $f': C' \rightarrow S'$ in $\mgstack{g,m}$ is given by a pair of morphisms of schemes
\[ (g: S' \rightarrow S, \; \overline{g}: C' \rightarrow C),\]
such that the diagram
\[ \diagram
C' \rrto^{\overline{g}} \dto_{f'}&& C \dto^{f}\\
S' \rrto_g& &S
\enddiagram \]
is Cartesian, and for all $i=1,\ldots,m$ holds
\[ \overline{g}\circ \sigma_i' = \sigma_i \circ g .\]
\end{defi}

\begin{rk}\em
It was shown by Knudsen \cite{Kn} that the moduli stack $\mgstack{g,m}$ is a smooth, irreducible Deligne-Mumford stack, which is  proper over $\Spec(\ZZ)$. Compare also \cite{DM} and \cite{Ed}. 
\end{rk}

The following proposition was proven for the case $m=0$ in \cite{Ed}, but also holds true in the case of arbitrary $m$. Here, as usual, square brackets are used to denote quotient stacks. 

\begin{prop}\label{28}
There are  isomorphisms of stacks
\[ \mgstack{g,m} \cong \left[ \, \Hgn/\PGL(N+1) \right], \]
and
\[ \mgstackopen{g,m} \cong \left[ \, \Hgnopen/\PGL(N+1) \right]. \]
\end{prop}

\proof
We will only give the proof in the case of the closed moduli stack, the second 
 case being  completely analogous. 
Recall that for any given scheme $S$ an object of  the fibre category 
$\left[ \, \Hgn/\PGL(N+1) \right](S)$ is a 
triple $(E,p,\phi)$, where $p: E\rightarrow S$ is a principal
$\PGL(N+1)$-bundle, and $\phi: E\rightarrow \Hgn$ is a
$\PGL(N+1)$-equivariant morphism. 

$(i)$
An isomorphism of stacks from $\mgstack{g,m}$ to $ [\Hgn/\PGL(N+1) ]$ is
constructed as follows. Let an $m$-pointed  stable curve $f : C
\rightarrow S \in 
\ob(\mgstack{g,m})$ be given, 
with sections $\sigma_1,\ldots,\sigma_{m}: S \rightarrow C$ 
defining divisors $S_1,\ldots,S_{m}$ on $C$.   We denote  by $p :
E\rightarrow S$ the 
principal $\PGL(N+1)$-bundle associated to the projective bundle 
$\PP f_\ast(\omega_{C/S}(S_1+\ldots+S_{m}))^{\otimes
n}$ over $S$.  Consider the
Cartesian diagram
\[\diagram
C\times_S E \rto^{\overline{f}} \dto_{\overline{p}} & E \dto^p\\
C \rto_f & S.
\enddiagram \]
The pullback of $\PP f_\ast(\omega_{C/S}(S_1+\ldots+S_{m}))^{\otimes
n} $ to $E$ has a natural trivialization as 
a  projective bundle. Denote by $\tilde{S}_1,\ldots,\tilde{S}_{m}$ the divisors on $C\times_S E$ defined by the sections $\tilde{\sigma}_i : E \rightarrow C\times_S E$, with $\tilde{\sigma}_i(e) := [\sigma_i\circ p(e),e]$ for $e\in E$ and $i=1,\ldots,m$. Then by the universal property of the pullback there is a unique isomorphism  
\[ p^\ast  \PP f_\ast\left(\omega_{C/S}(S_1+\ldots+S_{m})\right)^{\otimes
n}  \cong 
\PP \overline{f}_\ast\left(\omega_{C\times_S E/E}   (\tilde{S}_1+\ldots+\tilde{S}_{m})\right)^{\otimes
n}  \]
of projective bundles over $E$. 
So we have an $m$-pointed  stable curve $ \overline{f}: C\times_S E
\rightarrow E$, together with a natural trivialization of the projective bundle $\PP
\overline{f}_\ast(\omega_{C\times_S E/E} (\tilde{S}_1+\ldots+\tilde{S}_{m}))^{\otimes
n}$. This is equivalent to specifying a morphism 
\[\phi : \: \: E \rightarrow \Hgn \]
by the universal property of the Hilbert scheme. By construction, the morphism  $\phi$
is $\PGL(N+1)$-equivariant, so we obtain an object $(E,p,\phi)\in
[\Hgn/\PGL(N+1) ](S)$.

To define a functor from $ \mgstack{g,m}$ to $ [\Hgn/\PGL(N+1) ]$ we need also to consider morphisms $(g:S'\rightarrow S, \overline{g}:C'\rightarrow C)$ between  stable curves $f: C\rightarrow S$ and $f':C'\rightarrow S'$. Because for the respective sections $ \sigma_1,\ldots,\sigma_m: S\rightarrow C$ and $\sigma_1',\ldots,\sigma_m':S'\rightarrow C'$ holds $\sigma_i\circ g = \overline{g}\circ \sigma_i'$ for all $i=1,\ldots,m$ by definition,  there is an induced morphism
\[ \PP f_\ast'(\omega_{C'/S'}(S_1'+\ldots+S_m'))^{\otimes n} \rightarrow  
\PP f_\ast(\omega_{C/S}(S_1+\ldots+S_m))^{\otimes n} ,\]
where $S_1',\ldots,S_m'$ denote the divisors of the marked points on $C'$. Hence there is an induced morphism of the associated principal $\PGL(N+1)$-bundles $\tilde{g} : E' \rightarrow E$, which fits into a commutative diagram
\[ \diagram
&& \Hgn\\
E' \rto_{\tilde{g}} \dto_{p'} \urrto^{\phi'} & E \dto^{p}\urto_{\phi}\\
S' \rto_g & S.
\enddiagram \]
One easily verifies that this assignment is functorial. In this way we obtain a functor from the fibred category $\mgstack{g,m}$ to the quotient fibred category $ [\Hgn/\PGL(N+1) ]$ over the category of schemes, and thus a morphism of stacks. 

$(ii)$ Conversely, consider a triple $(E,p,\phi)\in
[\Hgn/\PGL(N+1) ](S)$. The morphism $\phi : E \rightarrow \Hgn $
determines an $m$-pointed  stable curve $f' : C'\rightarrow E$ of genus $g$,
together with a trivialization $\PP f'_\ast(\omega_{C'/E}(S_1'+\ldots+S_{m}'))^{\otimes
n} \cong \PP^N \times E$, where $S_1',\ldots,S_{m}'$ denote the divisors on $C'$ determined by the $m$ sections $\sigma_1',\ldots,\sigma_{m}' : E \rightarrow C'$. 
The group  $\PGL(N+1)$ acts  diagonally on  $\PP^N \times E$. By the $\PGL(N+1)$-equivariance of $\phi$, we have for all $e\in E$ and all $\gamma \in \PGL(N+1)$ the identity of fibres $C'_e \cdot \gamma = C'_{\gamma(e)}$ as embedded $m$-pointed curves. Therefore there is an induced action of $\PGL(N+1)$ on the embedded curve $C'\subset \PP^N \times E$, which respects the sections $\sigma_1',\ldots,\sigma_{m}'$. 
Taking quotients we obtain 
\[ C := C' / \PGL(N+1) \longrightarrow E/ \PGL(N+1) \cong S, \]
which defines an $m$-pointed  stable curve $C \rightarrow S\in \mgstack{g,m}(S)$.

The same construction assigns to a morphism in $[\Hgn/\PGL(N+1) ]$ a morphism between $m$-pointed stable curves. This defines a functor between fibred categories, and thus
a morphism of stacks from  $[\Hgn/\PGL(N+1) ]$ to  $\mgstack{g,m}$. 

$(iii)$ It remains to show that the two functors defined above form an equivalence of categories. In fact, the composition
\[\mgstack{g,m} \rightarrow [\Hgn/\PGL(N+1) ]\rightarrow\mgstack{g,m} \]
is just the identity on $\mgstack{g,m}$.

Conversely, let $(E,p,\phi) \in [\Hgn/\PGL(N+1) ](S)$ for some scheme $S$. This defines an $m$-pointed  stable curve $f:C\rightarrow S$ in $\mgstack{g,m}(S)$, which in turn defines an object $(E',p',\phi') \in [\Hgn/\PGL(N+1) ](S)$. We claim that the triples $(E,p,\phi)$ and $(E',p',\phi')$ are isomorphic. 

By construction we have a commutative diagram
\[ \diagram
C' \rto^{f'} \dto_{\overline{p}} & E\dto^p\\
C := C'/\PGL(N+1) \rto_f & E/\PGL(N+1) = S, 
\enddiagram \]
where $f':C'\rightarrow E$ is the $m$-pointed stable curve over $E$ defined by the morphism $\phi:E\rightarrow \Hgn$, and the vertical arrows represent the natural quotient morphisms.  Hence there is an isomorphism of projective bundles
\[ p^\ast\PP f_\ast(\omega_{C/S}(S_1+\ldots +S_m))^{\otimes n} \cong 
\PP f'_\ast(\omega_{C'/E}(S'_1+\ldots +S'_m))^{\otimes n}, \]
together with a trivialization of the bundle $
\PP f'_\ast(\omega_{C'/E}(S'_1+\ldots +S'_m))^{\otimes n}$ determined by $\phi: E \rightarrow \Hgn$. 
Thus the pullback of the projective bundle $\PP f_\ast(\omega_{C/S}(S_1+\ldots+S_m))^{\otimes n}$ to the principal $\PGL(N+1)$-bundle $p: E \rightarrow S$ is trivial, and hence $E$ must be isomorphic to the principal $\PGL(N+1)$-bundle associated to $\PP f_\ast(\omega_{C/S}(S_1+\ldots+S_m))^{\otimes n}$.

It is straightforward to give the corresponding arguments for the 
morphisms in both categories to show that  the two constructions are inverse to each other, up to isomorphism. 
\ebew

\begin{rk}\em
Since $\msbar_{g,m}$ is a categorical 
quotient of $\Hgn$, there is a natural morphism of stacks 
$  \left[\Hgn/\PGL(N+1) \right]\rightarrow \msbar_{g,m}$. 
In fact, the diagram  
\[\diagram
\mgstack{g,m} \rrto^{\sim}\drto&& \left[\Hgn/\PGL(N+1) \right]\dlto \\
&  \msbar_{g,m} 
\enddiagram \]
commutes.
\end{rk}

\begin{rk}\em
Recall that the symmetric group $\Sigma_m$ acts naturally on the set of $m$-pointed stable curves by permutation of the labels of the marked points. Since the $m$ sections of marked points are disjoint, the action of $\Sigma_m$, and of any subgroup $\Gamma \subset \Sigma_m$,  on $\overline{H}_{g,n,m}$ is free. 
\end{rk}

\begin{abschnitt}\em\label{gammaequiv}
It will turn out that our results can best be described using the notion of $\mgamma$-pointed stable curves, which we are going to introduce formally in definition \ref{mgdef}. Let us give the idea behind the concept first. An $m$-pointed stable curve of genus $g$ over a point is a tuple $(f:C\rightarrow\Spec(k);\sigma_1,\ldots,\sigma_m)$, where $\sigma_i: \Spec(k) \rightarrow C$, for $i=1,\ldots,m$, define the $m$ marked points on $C$. If $\Gamma$ is a subgroup of the permutation group $\Sigma_m$, then an $\mgamma$-pointed stable curve over a point is simply an equivalence class of such tuples, where  $(f:C\rightarrow\Spec(k);\sigma_1,\ldots,\sigma_m)$ and $(f':C'\rightarrow\Spec(k);\sigma_1',\ldots,\sigma_m')$ are called equivalent if $C=C'$  and if  there exists a permutation $\gamma\in\Gamma$, such that $\sigma_i = \sigma_{\gamma(i)}'$ for all $i=1,\ldots,m$. 

Over an arbitrary base scheme $S$ the situation is more subtle. Here, an $\mgamma$-pointed stable curve cannot simply be defined as an equivalence class of $m$-pointed stable curves with respect to the analogous equivalence relation. In fact, for the curves we are interested in, there usually exist no global sections distinguishing points on the fibres. However, the definition of $\mgamma$-pointed stable curves is made in such a way that locally, with respect to the \'etale topology,    
an $\mgamma$-pointed stable curve $f:C\rightarrow S$ over an arbitrary scheme $S$ looks like an equivalence class  of $m$-pointed stable curves over $S$. The  sections which exist on an \'etale cover are only determined up to permutations in $\Gamma$, and sections on different pieces need to satisfy certain compatibility conditions. 
\end{abschnitt}

\begin{defi}\em\label{mgdef}\label{n6}
Let $\Gamma$ be a subgroup of the symmetric group $\Sigma_m$. \\
$(i)$ A {\em charted $\mgamma$-pointed stable curve of genus $g$} is a tuple
\[ (f: C\rightarrow S, \:\: u: S' \rightarrow S, \:\: \sigma_1,\ldots,\sigma_m: S'\rightarrow C'),\]
where $f: C\rightarrow S$ is a flat and proper morphism of schemes, with reduced and connected algebraic curves as its fibres,  $u: S'\rightarrow S$ is an \'etale covering, defining  the fibre product $C' := C\times_S S'$, such that the induced morphism $f' : C' \rightarrow S'$, together with the sections $\sigma_1,\ldots,\sigma_m: S'\rightarrow C'$, is an $m$-pointed stable curve of genus $g$. Furthermore, we require that for all closed points $s,s'\in S'$ with $u(s)=u(s')$, there exists a permutation $\gamma_{s,s'}\in \Gamma$, such that 
\[ (\ast) \qquad \qquad \overline{u} \circ \sigma_i(s) = \overline{u}\circ \sigma_{\gamma_{s,s'}(i)}(s') \]
holds for all $i=1,\ldots,m$, where $\overline{u}: C'\rightarrow C$ is the morphism induced by $u: S' \rightarrow S$.\\
$(ii)$ Two charted $\mgamma$-pointed stable curves are called {\em equivalent}, if there exists a third charted $\mgamma$-pointed stable curve dominating both of them, in the sense of remark \ref{explain} below.\\
$(iii)$ An {\em $\mgamma$-pointed stable curve of genus $g$} is an equivalence class of  charted $\mgamma$-pointed stable curves of genus $g$.
\end{defi}

\begin{rk}\em\label{explain}
Let a charted $\mgamma$-pointed stable curve $(f: C\rightarrow S, \:\: u': S' \rightarrow S, \:\: \sigma_1,\ldots,\sigma_m: S'\rightarrow C')$ be given. A second charted $\mgamma$-pointed stable curve $(f: C\rightarrow S, \:\: u'': S'' \rightarrow S, \:\: \tau_1,\ldots,\tau_m: S''\rightarrow C'')$ with the same underlying curve $f:C\rightarrow S$ is said to {\em dominate} the first one, if the morphism $u'': S''\rightarrow S$ factors as $u'' = u' \circ v$, where $v: S''\rightarrow S'$ is an \'etale covering, which induces an isomorphism
\[ C'' \cong C'\times_{S'}S'', \]
and such that for all closed points $s\in S''$ there exists a permutation $\gamma_s\in \Gamma$, with
\[ (\ast') \qquad \qquad \overline{v}\circ \tau_i(s) = \sigma_{\gamma_s(i)}\circ v(s) \]
holding for all $i=1,\ldots, m$. Here $\overline{v}: C''\rightarrow C'$ denotes the morphism induced by $v: S'' \rightarrow S'$.
\end{rk}

\begin{rk}\em\label{class}
Let us analyse what the definition of an $\mgamma$-pointed stable curve means on the fibres of a curve, i.e. let us consider curves $f: C\rightarrow S$ for $S= \Spec(k)$.

$(i)$  Let $(f: C\rightarrow \Spec(k), \:\: u: S' \rightarrow \Spec(k), \:\: \sigma_1,\ldots,\sigma_m: S'\rightarrow C')$ be a charted  $\mgamma$-pointed stable curve of genus $g$ over a simple point. Then, by the compatibility condition $(\ast)$ of definition \ref{mgdef}, the sections  $\sigma_1,\ldots,\sigma_m$ distinguish $m$ distinct points on $C$, which are necessarily different from any nodes $C$ might have. 

Let $[i]$ denote the equivalence class of an element $i\in\{1,\ldots,m\}$ with respect to the action of $\Gamma$. By the compatibility condition $(\ast)$, the sections of $f': C'\rightarrow S'$ associate to each distinguished point of $C$ a unique label $[i]$, for some $i\in\{1,\ldots,m\}$. Formally we define for a distinguished point $q\in C$
\[ \mbox{class}(q) := [i], \]
if $q=\overline{u}\circ \sigma_i(s)$ for some $i\in \{1,\ldots,m\}$ and some $s\in S'$. In other words,  a charted  $\mgamma$-pointed stable curve of genus $g$ defines an $m$-pointed stable curve of genus $g$, but where the labels of the marked points are only determined up to a permutation in $\Gamma$. We will make this precise in lemma \ref{onetoone04} below. 

$(ii)$ Let $(f: C\rightarrow S, \:\: u'': S'' \rightarrow S, \:\: \tau_1,\ldots,\tau_m: S''\rightarrow C'')$ be a  second charted $\mgamma$-pointed stable curve with the same underlying curve $f:C\rightarrow S$, which dominates the first one. The compatibility condition $(\ast')$ of remark \ref{explain}, together with condition $(\ast)$ of definition \ref{mgdef}, ensures that the distinguished points on $C$ are the same in both cases. Even more, the classes of the distinguished points as defined above remain the same. 

$(iii)$ The notion of a charted $\mgamma$-pointed stable curve makes it necessary to specify one \'etale cover of $f:C\rightarrow S$ by an $m$-pointed stable curve. The definition of an $\mgamma$-pointed stable curve is independent of such a choice. 
\end{rk}

\begin{lemma}\label{onetoone04}
There is one-to-one correspondence between $\mgamma$-pointed stable curves of genus $g$ over $\Spec(k)$ and equivalence classes of $m$-pointed stable curves of the same genus over $\Spec(k)$ in the sense of remark \ref{gammaequiv}. 
\end{lemma}

\proof
Let $(f: C\rightarrow \Spec(k), \sigma_1,\ldots,\sigma_m)$ represent an equivalence class of $m$-pointed stable curves with respect to the action of $\Gamma$. Then $(f: C\rightarrow \Spec(k), \id_{\Spec(k)}, \sigma_1,\ldots,\sigma_m)$ is a charted  $\mgamma$-pointed stable curve. If $(f: C\rightarrow \Spec(k), \tau_1,\ldots,\tau_m)$ is a different representative of the above class, then by definition there exists a permutation $\gamma\in\Gamma$, such that $\sigma_i= \tau_{\gamma(i)}$ for all $i=1,\ldots,m$. The disjoint union $f \coprod f': C \coprod C' \rightarrow \Spec(k) \coprod \Spec(k)$ defines a charted $\mgamma$-pointed stable curve, dominating both $(f: C\rightarrow \Spec(k), \id_{\Spec(k)}, \sigma_1,\ldots,\sigma_m)$ and $(f: C\rightarrow \Spec(k), \id_{\Spec(k)}, \tau_1,\ldots,\tau_m)$. Hence to each equivalence class of $m$-pointed stable curves there is associated a well-defined $\mgamma$-pointed stable curve.   

Conversely, let an $\mgamma$-pointed stable curve be given by a representative 
 $(f: C\rightarrow \Spec(k), \:\: u: S' \rightarrow \Spec(k), \:\: \sigma_1,\ldots,\sigma_m: S'\rightarrow C')$. Choose some closed point $s\in S'$. Then the fibre $f_s': C_s'\rightarrow \Spec(k)$ of $f': C'\rightarrow S'$ over $s$, 
together with the sections $\tau_i := \sigma_i|\{s\}$, for $i=1,\ldots,m$,  defines an $m$-pointed stable curve $(f'_s: C_s'\rightarrow \Spec(k), \tau_1,\ldots,\tau_m)$. Let $(f'_{s'}: C_{s'}'\rightarrow \Spec(k), \tau_1',\ldots,\tau_m')$  denote the $m$-pointed curve obtained from the fibre over a different point $s'\in S'$. We clearly have $C_{s'}' \cong C \cong C_s'$. Because of condition $(\ast)$ of definition \ref{mgdef}  there exists a permutation $\gamma\in \Gamma$, such that $\tau_i= \tau_{\gamma(i)}$ for all $i=1,\ldots,m$. Hence the equivalence class of the $m$-pointed stable curve is well defined.

The two constructions are inverse to each other, which concludes the proof of  the lemma.
\ebew

\begin{defi}\em 
For $k=1,2$, let $(f_k: C_k\rightarrow S_k, u_k: S_k'\rightarrow S_k, \sigma_1^{(k)},\ldots,\sigma_m^{(k)})$ be two charted $\mgamma$-pointed stable curves of genus $g$. A {\em morphism} from the second tuple to the first is a morphism of schemes 
\[ h : \:\: S_1 \rightarrow S_2,\]
which induces an isomorphism $C_1\cong C_2\times_{S_2} S_1$, and which satisfies the condition $(\diamond)$ below. 

Put $S_1'' := S_1' \times_{S_2} S_2'$. Then the natural composed morphism $u''_1: S_1''\rightarrow S_1$ is an \'etale covering of $S_1$, and for $C_1'':= C_1\times_{S_1} S_1''$ the tuple $(f_1: C_1\rightarrow S_1, u'':S_1''\rightarrow S_1, \overline{\sigma}_1^{(1)},\ldots,\overline{\sigma}_m^{(1)} : S_1''\rightarrow C_1'')$ is a charted $\mgamma$-pointed stable curve dominating $(f_1: C_1\rightarrow S_1, u':S_1'\rightarrow S_1, {\sigma}_1^{(1)},\ldots,{\sigma}_m^{(1)} : S_1'\rightarrow C_1')$. Here $\overline{\sigma}_i^{(1)}$ denotes the section induced by $\sigma_i^{(1)}$, for $i=1,\ldots,m$. There are induced morphisms $\overline{h}: S_1'' \rightarrow S_2'$ and $h'':C_1''\rightarrow C_2'$, such that the diagram
\[ \diagram
C_1'' \xto[rrr]^{h''} \xto[ddd] \drto &&& C_2' \xto[ddd] \dlto\\
& C_1 \rto \dto & C_2 \dto \\
& S_1 \rto^h & S_2 \\
S_1'' \xto[rrr]^{\overline{h}} \urto  &&& S_2' \ulto
\enddiagram \]
commutes, and all quadrangles are Cartesian. We require that there exists a permutation $\gamma\in \Gamma$, such that
\[ (\diamond) \qquad \qquad h''\circ \overline{\sigma}_i^{(1)} = \sigma_{\gamma(i)}^{(2)} \circ \overline{h} \]
holds for all $i=1,\ldots,m$. 
\end{defi}

\begin{rk}\em
$(i)$ If $h: S_1 \rightarrow S_2$ is a morphism between a pair of charted $\mgamma$-pointed stable curves, then it is also a morphism for any other pair equivalent to it. This can easily be seen by drawing the appropriate extensions of the above commutative diagram. In particular, there is a well-defined notion of a {\em morphism between $\mgamma$-pointed stable curves}.\\
$(ii)$ A morphism $h: S_1\rightarrow S_2$ between $\mgamma$-pointed stable curves is an isomorphism, if and only if it is an isomorphism of schemes $S_1\cong S_2$. Note however that not any isomorphism of schemes defines an isomorphism of $\mgamma$-pointed stable curves, it might not even define a morphism between them. \\
$(iii)$ Two $\mgamma$-pointed stable curves over the same underlying curve $f:C\rightarrow S$ are isomorphic if and only if they admit \'etale coverings by $m$-pointed stable curves, which are equivalent up to a permutation $\gamma$ of the labels of their sections, with $\gamma$ contained in  $\Gamma$, compare remark \ref{gammaequiv}. 
\end{rk}

\begin{rk}\em
A morphism $h: S_1\rightarrow S_2$ is a morphism of $\mgamma$-pointed stable curves if and only if there is a Cartesian diagram
\[ \diagram
C_1 \rto^{h'} \dto & C_2 \dto\\
S_1\rto^h &S_2
\enddiagram \]
such that for all closed points $s\in S_1$ and all distinguished points $q\in C_{1,s}$ in the fibre over $s$ holds
\[ \mbox{class}(q) = \mbox{class}(h'(q)) .\]
See remark \ref{class} for the notation. 
\end{rk}

\begin{rk}\em
$(i)$ It is straightforward to generalize definition \ref{mgdef}  to a notion of $\mgamma$-pointed prestable curves of genus $g$. \\
$(ii)$ Note that the permutations $\gamma_{s,s'} \in \Gamma$ in definition \ref{mgdef} are necessarily uniquely determined, and the same is true for the permutations $\gamma_s$ in remark \ref{explain}. 
\end{rk}

\begin{rk}\em
A subgroup $\Gamma$ of $\Sigma_m$ acts on $\msbar_{g,m}$. The quotient 
\[ \msbar_{g,\mgamma} := \msbar_{g,m} / \Gamma \label{n7}\]
is a coarse moduli space for $\mgamma$-pointed curves.  Let $f: C\rightarrow S$ be an $\mgamma$-pointed stable curve of genus $g$. By lemma \ref{onetoone04}, the closed points of $\msbar_{g,\mgamma}$ are in one-to-one correspondence with isomorphism classes of $ \mgamma$-pointed stable curves of genus $g$ over $\Spec(k)$.

Let $(f:C\rightarrow S, u:S'\rightarrow S, \sigma_1,\ldots,\sigma_m: S'\rightarrow C')$ be a charted $\mgamma$-pointed stable curve of genus $g$. Since $f': C'\rightarrow S'$ is an $m$-pointed stable curve, there is an induced morphism $\theta_{f'}:S'\rightarrow  \msbar_{g,m}$, and by composition a morphism $\overline{\theta}_{f'}:S'\rightarrow  \msbar_{g,\mgamma}$. By definition, for any closed point $s\in S$, and for any pair of points $s',s''\in S''$ with $u(s')=u(s'')=s$, holds the identity $\overline{\theta}_{f'}(s') =\overline{\theta}_{f'}(s'')$. Thus there is an induced morphism $\theta: S\rightarrow \msbar_{g,\mgamma}$, such that $\theta\circ u =  \overline{\theta}_{f'}$.  Clearly any other charted $\mgamma$-pointed stable curve dominating this one induces the same morphism from $S$ into $\msbar_{g,\mgamma}$. Hence there is a well defined morphism $\theta_f: S\rightarrow \msbar_{g,\mgamma}$ for any $\mgamma$-pointed stable curve $f: C\rightarrow S$. Furthermore, by referring back to $\msbar_{g,m}$ again, one easily sees that $\msbar_{g,\mgamma}$ dominates any other scheme with this universal property, so it is indeed a coarse moduli space for $\mgamma$-pointed stable curves.
\end{rk}

\begin{rk}\em
If $\Gamma = \{\id\}$ is the trivial subgroup of of $\Sigma_m$, then the glueing condition $(\ast)$ for $\mgamma$-pointed stable curves implies that there are  $m$ induced  global sections of marked points of $f:C\rightarrow S$.  Thus in this case an $\mgamma$-pointed stable curve is just an $m$-pointed stable curve in the classical sense.
\end{rk}

\begin{rk}\em\label{rk224}
Let $\Gamma \subset \Sigma_m$ be a subgroup. There is a free action of $\Gamma$ on $\Hgn$. The quotient shall be denoted by
\[ \Hgmp := \Hgn / \Gamma \label{n8}.\]
Let $u_{g,n,m} : {\cal C}_{g,n,m} \rightarrow \Hgn$ denote the universal curve over $\Hgn$. Note that fibres of $u_{g,n,,m}$ over such points of $\Hgn$, which correspond to each other under the action of $\Gamma$, are identical as embedded curves in $\PP^N$, with the same distinguished points on them. The only difference is in the labels of the marked points, which are permuted by elements of $\Gamma$. Thus the action of $\Gamma$ on $\Hgn$ extends to an action on ${\cal C}_{g,n,\mgamma}$, by identifying fibres over corresponding points. 
The quotient
\[ u_{g,n,\mgamma}: \quad {\cal C}_{g,n,\mgamma} := {\cal C}_{g,n,m}/ \Gamma \: \longrightarrow \: \Hgn/\Gamma = \Hgmp\]
exists as an $\mgamma$-pointed stable curve over $\Hgmp$. The original universal curve $u_{g,n,m} : {\cal C}_{g,n,m} \rightarrow \Hgn$, together with the \'etale morphism $\Hgn\rightarrow \Hgmp$ defines a charted $\mgamma$-pointed stable curve representing it. 
\end{rk}

\begin{abschnitt}\em
Let $f : C\rightarrow S$ be an $\mgamma$-pointed stable curve of genus $g$, which is given together with an embedding of the underlying curve $f: C\rightarrow S$ into $\pr_2: \PP^N\times S \rightarrow  S$. For any \'etale covering $u:S'\rightarrow S$, there is an induced embedding of $f' : C' = C\times_S S' \rightarrow S'$ into $\pr_2: \PP^N\times S' \rightarrow S'$. Hence, for a representing charted $\mgamma$-pointed curve, the curve $f' : C' \rightarrow S'$ is an embedded $m$-pointed stable curve, and thus induces a morphism 
\[ \theta_{f'}: \: \: S' \rightarrow \Hgn.\]
 For any closed point $s\in S$, the embedding of the fibre $C'_{s'}$ into $\PP^N$ is the same for  all points $s'\in S'$ with $u(s')=s$, as it is nothing else but the embedding of $C_s$ into $\PP^N$. Note that the points in $\Hgn$ representing different fibres $C'_{s'}$ need not be the same, as the labels of the marked points may differ by a permutation in $\Gamma$. However, the morphism $\overline{\theta}_{f'}: S' \rightarrow \Hgmp$ obtained by composition with the quotient map, factors through a morphism $\theta_f : S\rightarrow \Hgmp$. This morphism is even independent of the chosen charted $\mgamma$-pointed stable curve representing $f:C\rightarrow S$. 

By the universal property of the Hilbert scheme, there is an isomorphism of $m$-pointed stable curves between $\theta_{f'}^\ast {\cal C}_{g,n,m}$ and $C'$ over $S'$. As a scheme over $S'$, the curve $\theta_{f'}^\ast {\cal C}_{g,n,m}$ is isomorphic to  $\overline{\theta}_{f'}^\ast {\cal C}_{g,n,\mgamma}$, even though the latter has no natural structure as an $m$-pointed curve. Therefore, by the definition of $C' = u^\ast C$, and since $\overline{\theta}_{f'}= \theta_f\circ u$, one can construct an isomorphism of schemes over $S$ between $C$ and $\theta_f^\ast {\cal C}_{g,n,\mgamma}$. Together with the isomorphisms between the covering $m$-pointed stable curves, this shows that $f:C\rightarrow S$ is isomorphic to $\theta_f^\ast {\cal C}_{g,n,\mgamma} \rightarrow S$ as an $\mgamma$-pointed stable curve. 

In other words,  $\Hgmp$ is in fact a
fine moduli space for embedded $\mgamma$-pointed stable curves of genus $g$, with universal curve $u_{g,n,\mgamma}: {\cal C}_{g,n,\mgamma} \rightarrow \Hgmp$. 
\end{abschnitt}

\begin{rk}\em
The action of $\PGL(N+1)$ on $\Hgn$ induces an action of $\PGL(N+1)$ on $\Hgmp$. Indeed, the embedding of an $m$-pointed stable curve of genus $g$ does not depend on the ordering of the labels of its marked points. So the action of $\PGL(N+1)$ on $\Hgn$ commutes with the action of any subgroup $\Gamma\subset \Sigma_m$ on $\Hgn$.  The coarse moduli space $\msbar_{g,\mgamma}$ is a GIT-quotient of $\Hgmp$ by the action of $\PGL(N+1)$. 
\end{rk}

\begin{defi}\em
The {\em moduli stack $\, \mgstack{g,\mgamma}$\label{n9} of $\mgamma$-pointed stable curves of genus $g$} is the stack defined as the category fibred in groupoids over  the category of schemes, where for a scheme $S$ the objects in the fibre category $\mgstack{g,\mgamma}(S)$ are the $\mgamma$-pointed stable curves of genus $g$ over $S$.
Morphisms in $\mgstack{g,\mgamma}$ are morphisms between $\mgamma$-pointed stable curves.
\end{defi}

\begin{prop}\label{P217}
There is an isomorphism of stacks
\[  \mgstack{g,\mgamma} \: \cong \: \left[ \Hgmp / \PGL(N+1) \right] .\]
\end{prop}

\begin{rk}\em
Note that by the general theory of quotient stacks, and the fact that the action of $\Gamma$ on $\Hgn$ is free, there is an isomorpism of stacks
\[   \left[ \, \Hgmp / \PGL(N+1) \, \right]  \: \cong \: \left[ \, \Hgn / \Gamma \times \PGL(N+1) \right]
 .\]
\end{rk}

\proofof{the proposition}
The proof is analogous to that of proposition \ref{28}, so we may be brief here. \\
$(i)$ Let $f:C\rightarrow S \in \ob(\mgstack{g,\mgamma})$ be an $\mgamma$-pointed stable curve. For each 
\'etale cover $u: S'\rightarrow S$, the $m$ disjoint sections $\sigma_1,\ldots,\sigma_m: S'\rightarrow C'$ define hypersurfaces $S_1,\ldots,S_m$ in $C'= C\times_S S'$. Their images in $C'$ define a divisor in $C$, which will be denoted by $S_{\mgamma}$. Note that this divisor is independent of the chosen \'etale cover. 

 Let $p:E\rightarrow S$ be the principal $\PGL(N+1)$-bundle associated to the projective bundle $\PP f_\ast( \omega_{C/S}(S_{\mgamma}))^{\otimes n}$. Let $E' := E\times_S S'$ be the pullback of $E$ to some \'etale covering $u:S'\rightarrow S$ of $S$. In the same way as in the proof of proposition \ref{28}, there is a natural embedding of the $m$-pointed stable curve $C'\times_{S'} E' \rightarrow E'$, and this induces a morphism
\[ \phi' : \: E' \rightarrow \Hgn .\]
By composition with the quotient map there is a morphism
\[\overline{\phi} : \:\: E' \rightarrow \Hgn/\Gamma = \Hgmp , \]
which is $\PGL(N+1)$-equivariant. Note that the action of $\Gamma$ commutes with the action of $\PGL(N+1)$ on $\Hgn$.  

For all closed points $e,e' \in E'$ which project to the same point in $E$, we have $\overline{\phi}(e)=\overline{\phi}(e')$. Therefore  $\overline{\phi}$ factors through a morphism $\phi: E\rightarrow \Hgmp$, which is also $\PGL(N+1)$-equivariant. Thus we obtain an object $(E,p,\phi)\in [\Hgmp/\PGL(N+1](S)$. 

The construction of the functor from $\mgstack{g,\mgamma}$ to $[\Hgmp/\PGL(N+1]$ on morphisms is straightforward. 

$(ii)$ Conversely, consider a triple  $(E,p,\phi)\in [\Hgmp/\PGL(N+1](S)$ for some scheme $S$. The morphism $\phi: E \rightarrow \Hgmp$ determines an embedded $\mgamma$-pointed stable curve $f': C' \rightarrow E$ of genus $g$, together with an isomorphism between $C'$ and the pullback of the universal curve ${\cal C}_{g,n,\mgamma}$ as schemes over $E$. The quotient 
\[ f: C := C' / \PGL(N+1) \: \rightarrow  \: E / \PGL(N+1) \cong S \]
exists, and it is by construction an $\mgamma$-pointed stable curve. 

One verifies that the two functors, whose  construction we outlined  above,  are inverse to each other as morphisms of stacks.
\ebew

\begin{rk}\em
From proposition \ref{P217} it follows in particular that $\mgstack{g,\mgamma}$ is a smooth Deligne-Mumford stack, with $\msbar_{g,\mgamma}$ as its moduli space.
 There is a canonical commutative diagram
\[ \diagram
\mgstack{g,m} \dto^\cong \rto &\mgstack{g,\mgamma} \dto^\cong\\
\left[\Hgn/\PGL(N+1) \right]\dto \rto & \left[\Hgn/\Gamma\times \PGL(N+1) \right]\dto \\
 \msbar_{g,m} \rto & \msbar_{g,\mgamma}.
\enddiagram \]

All horizontal arrows represent finite morphisms of degree equal to the order of $\, \Gamma$. The morphisms between the stacks are unramified, and even finite \'etale morphisms,  while in general the morphisms between the moduli spaces are not.  
\end{rk}

\begin{rk}\em
Note that for everything we said above analogous statements hold true if we replace $\Hgn$ 
and $\msbar_{g,m}$ by the open subschemes $H_{g,n,m}$ and $\ms_{g,m}$, as well as  $\Hgmp$ 
and $\msbar_{g,\mgamma}$ by the open subschemes $H_{g,n,\mgamma}$ and 
$\ms_{g,\mgamma}$, respectively. 
\end{rk}

\chapter{The moduli problem}\label{sect3}

\begin{abschnitt}\em
We are now going to introduce the main objects of study, and some basic concepts. Recall that there is a stratification of the coarse moduli space 
\[ \msbar_g =  \ms_g^{(0)} \cup \ldots \cup \ms_g^{(3g-3)} ,\]
where each $\ms_g^{(i)}$ is of codimension $i$ in $\msbar_{g}$, and defined as the reduced locus of stable curves with precisely $i$ nodes. The closure of $\ms_g^{(i)}$ in $\msbar_g$ is denoted by $\msbar_g^{(i)}$.

First, we want to characterize those stable curves $f:C\rightarrow \Spec(k)$, which are represented by points of one fixed irreducible component ${D}$ of $ \ms_g^{(3g-4)}$. Essentially, these are curves which can be derived from deformations of one stable curve $f_0: C_0\rightarrow \Spec(k)$, which has $3g-3$ nodes, and which is represented by a point in the closure $\overline{D}$  of  ${D}$ in $\overline{\ms}_g^{(3g-4)}$. In fact, any irreducible component of $ \ms_g^{(3g-4)}$ is uniquely determined by specifying one stable curve $f_0: C_0\rightarrow \Spec(k)$ with $3g-3$ nodes, and one node $P_1$ on it, see lemma \ref{L31} below. We then write ${D} = \dpo$. 

Lemma \ref{315} gives a criterion to decide whether a given stable curve is represented by a point in $\dpo$. This description leads to the definition \ref{n26} of a  stack $\cdpo$, which is nothing else but the reduction of the preimage substack of $\dpo$ in the moduli stack $\cmgbar$. 

For any stable curve $f:C\rightarrow S$ in $\cdpo$, there exists a distinguished subcurve, see corollary \ref{K324}. This implies that the subscheme $\dpo$ is isomorphic to (a subscheme of) a (finite quotient of a) moduli space $\ms_{g^+,m}$ of $m$-pointed stable curves of some genus $g^+$, see proposition \ref{312}. There will be more to be said  about a morphism from $\ms_{g^+,m}$ to $\mg$ later on. Our main purpose will be to compare the moduli stacks corresponding to these two moduli spaces.
\end{abschnitt}

\begin{abschnitt}\em
Fix now and for all a stable curve $f_0: C_0 \rightarrow \Spec(k)$\label{n10} of genus $g\ge 3$ with $3g-3$ nodes. We also  choose an enumeration $P_1,\ldots,P_{3g-3}$ of the nodes. 
\end{abschnitt}

\begin{lemma}\label{L31}
The node $P_1$ distinguishes an irreducible component $\overline{D}$ of the stratum $\msbar_g^{(3g-4)}$ in the following way.
Let $f: C \rightarrow S$ be any stable curve of genus $g$ with irreducible and reduced base $S$, which is a deformation of $C_0$, with $C_0 \cong C_{{s_0}}$ for some $s_0 \in S$. Assume that there exist disjoint sections $\varrho_2,\ldots,\varrho_{3g-3} : S \rightarrow C$, such that for all $s \in S$ and all $i=2,\ldots,3g-3$ holds\\
\qquad $(i) \qquad \varrho_i(s)$ is a node of $C_s$, and \\
\qquad $(ii) \qquad \varrho_i(s_0) = P_i$. \\
Then for all closed points $s\in S$, the point $[C_s]\in \msbar_{g}$ representing the fibre $C_s$ is contained in $\overline{D}$, or equivalently, the induced morphism $\theta_f:S\rightarrow \msbar_g$ factors through $\overline{D}$.
\end{lemma}

\begin{defi}\em
$(i)$ The unique irreducible component of $\msbar_g^{(3g-4)}$, which contains the point $[C_0]$ representing the curve $f_0: C_0\rightarrow \Spec(k)$, and all points representing fibres of a deformation of $C_0$, which preserve all nodes except $P_1$, as in lemma \ref{L31}, is denoted by $\dpobar$.\label{n11}\\
$(ii)$ 
The intersection of $\dpobar$ with the locus  $\ms_g^{(3g-4)}$ of points parametrizing curves with exactly $3g-4$ nodes is denoted by $\dpo$.  
\end{defi}

\proofof{the lemma}
A  neighbourhood of $[C_0]\in \msbar_g$ can be constructed as
a quotient $B / \Aut(C_0)$, where $B$ is a smooth base of a universal
deformation $C \rightarrow B$ of $C_0$ of dimension $3g-3$. If we fix one projective
embedding of $C_0$, then 
 $B$ can  be  choosen to be  a subscheme of  the Hilbert scheme $\hilb{g,n,0}$ containing the point representing the fixed embedding of $C_0$. Moreover, if
$B$ is sufficiently  small, then  each node
$P_i$ of $C_0$ determines a smooth hypersurface $D_i \subset B$ 
parametrizing those fibres, which preserve the singularity $P_i$. All of these 
hypersurfaces together form a normal crossing divisor. In particular,
the intersection $D_2 \cap \ldots \cap D_{3g-3}$ determines a smooth one-dimensional subscheme  $D$ in $ B / \Aut(C_0) \hookrightarrow \msbar_g$ containing $[C_0]$. 

We now define $\overline{D}$ as the irreducible component of $\msbar_g^{(3g-4)}$ which contains  $D$.  Since $f: C \rightarrow B$ is a universal deformation, the
map $\theta_f: S \rightarrow \msbar_g$ factors locally near $s_0$ through
$D$, and since $S$ is irreducible it factors globally through 
$\overline{D}$ as well.
\ebew

\begin{notation}\em
Whenever there is no danger of misunderstanding, we will denote the distinguished irreducible component $\dpobar$ of lemma \ref{L31} simply by $\dppbar$. The open subscheme, which is the intersection of  $\dpobar$ with $\ms_g^{(3g-4)}$, i.e. the subcurve excluding points lying in $\msbar_g^{(3g-3)}$, shall be abbreviated as $\dpp$. 

Note that while an irreducible component of $\msbar_g^{(3g-4)}$ is uniquely determined by the choice of $C_0$ and a node $P_1$ on it, the converse is not true. There may exist another node $P'$ on $C_0$, or even  another stable curve $f_1:C_1\rightarrow \Spec(k)$ with $3g-3$ nodes $Q_1,\ldots,Q_{3g-3}$ on it,  such that $\overline{D}(C_0;P')= \dpobar$ or $\overline{D}(C_1;Q_1)= \dpobar$, respectively.  
\end{notation}

\begin{rk}\em\label{R32}
The irreducible components of $\msbar_g^{(3g-4)}$ are smooth curves. Together with the above description of the boundary of $\msbar_g$  this is equivalent to the observation that $\dpobar = \overline{D}(C_0;P_i)$ for some $i\in\{2,\ldots,3g-3\}$ if and only if there exists an automorphism $\gamma\in\Aut(C_0)$  such that $\gamma(P_1)=P_i$. Compare also lemma \ref{L326} below. 
\end{rk}
 
By definiton, the subscheme $\dpobar$ contains the point $[C_0]$ and the points representing deformations of $C_0$, which preserve the nodes $P_2,\ldots,P_{3g-3}$. The following lemma shows that $\dpobar$ is exactly the locus of such deformations, and finite iterations of such deformations, preserving the ``right'' $3g-3$ nodes. 

\begin{lemma}\label{33}
Let $f: C \rightarrow S$ be a stable curve of genus $g$, with
irreducible and reduced  base $S$. The induced morphism $\theta_f: S \rightarrow \msbar_g$
factors through $\dppbar$ if and only if there exists a finite family
$\{f_i: C^{(i)}\rightarrow S^{(i)}\}_{i=0,\ldots ,k}$ of stable curves
of genus $g$, with irreducible and reduced  bases, such that 
\begin{itemize}
\item[$(i)$] for all $0 \le i \le k$ there exist disjoint sections
$\rho_j^{(i)} : S^{(i)}\rightarrow C^{(i)}$, with $2 \le j \le 3g-3$,
such that for all $s\in S^{(i)}$ the point $\rho_j^{(i)}(s)$ is a
node of the fibre $C_s^{(i)}$;

\item[$(ii)$] for all $1 \le i \le k$ there exist points $s_+^{(i)}
\in S^{(i)}$ and $s_-^{(i)}
\in S^{(i-1)}$ such that $f_i: C^{(i)}\rightarrow S^{(i)}$ is a
deformation of 
$ C^{(i-1)}_{s_-^{(i)}}$ with 
\[ C^{(i)}_{s_+^{(i)}} =
C^{(i-1)}_{s_-^{(i)}} ;\]

\item[$(iii)$] under the above identification, for all $1 \le i \le
k$ and all $2 \le j \le 3g-3$ there is the identity 
\[ \rho_j^{(i)}(s_+^{(i)}) = \rho_j^{(i-1)}(s_-^{(i)}) ;\]

\item[$(iv)$] the base $S^{(k)}$ is an \'etale open subscheme of $S$, and the curve $f_k: C^{(k)}\rightarrow S^{(k)}$ is obtained from the curve
$f : C \rightarrow S$ by restriction to $S^{(k)}$;

\item[$(v)$] the curve $f_0: C^{(0)}\rightarrow S^{(0)}$ is a
deformation of $C_0$, with $C_0 \cong C_{{s^{(0)}}}$ for some point
$s^{(0)} \in  S^{(0)}$, and such that for all $2 \le j \le 3g-3$
holds
\[ \rho_j^{(0)}(s^{(0)}) = P_j .\]
\end{itemize}
\end{lemma}

To prove the above, we need an auxiliary lemma.

\begin{lemma}\label{sectnod}
Let $f: C\rightarrow S$ be a stable curve of genus $g$, with reduced and irreducible base $S$. Suppose that the number of nodes of a fibre over a general point is exactly $k$, for some $1\le k \le 3g-3$. Then there exists an \'etale covering $S'\rightarrow S$, such that the pullback $f': C'\rightarrow S'$ of $f:C\rightarrow S$ to $S'$ admits $k$ disjoint sections of nodes. 
\end{lemma}

\proof
By assumption there exists a dense open subset $S_0\subset S$, such that for all closed points $s\in S_0$ the fibre $C_s$ over $s$ has exactly $k$ nodes. For an arbitrary closed point $s\in S$, the fibre $C_s$ has at least $k$ nodes. Without loss of generality we may assume that the curve $f:C\rightarrow S$ is embedded into some projective space. 
Consider the subscheme $N \subset C$, which is defined as the reduced locus of all points of $C$ which are nodes of fibres of $f:C\rightarrow S$. Formally, it can be constructed as the locus of those points, where the sheaf $\Omega_{C/S}^1$ is not free. Because of the flatness of $f:C\rightarrow S$, these points are precisely the nodes of the fibres. Define $R$ as the reduced union of all those irreducible components of $N$, which surject onto $S$.

Let $s\in S$ be a general point, and let $n\in C_s$ be one of the $k$  nodes over $s$. Since all neighbouring fibres have $k$ nodes as well, we must have $n\in R$. Therefore, the fibre of the restricted morphism $f|R: R\rightarrow S$ over $s$ consists of exactly $k$ necessarily simple points, which are precisely all the nodes over $s$. Clearly no fibre of $f|R$ may consist of more than $k$ points. Recall that closed points of $R$ are nodes of fibres of $f:C\rightarrow S$, so they are ordinary double points. There is no deformation, which splits one ordinary double point into two, so there can also be no fibre of $f|R$ consisting of less than $k$ points. In other words, each fibre of $f|R: R\rightarrow S$ consists of  exactly $k$ simple points, and  therefore the morphism $f|R:R\rightarrow S$ is flat, and an \'etale covering of $S$. 

Incidentally, the last argument shows also that 
 the decomposition of $N$ into $R$ and $R^\#$ is disjoint, where  $R^\#$ denotes the reduced union of all those irreducible components of $N$, which do not surject onto $S$. 

We now put $S_1:= R$, and define $f_1:C_1\rightarrow S_1$ as the pullback of $f:C\rightarrow S$ to $S_1$. In particular, the stable curve $f_1:C_1\rightarrow S_1$ admits a tautological section $\sigma_1: S_1\rightarrow C_1$, such that for all closed points $s\in S_1$, the point $\sigma_1(s)$ is a node of the fibre $C_{1,s}$. 

Note that $\sigma(S_1)$ is a connected component of the pullback $R_1\subset C_1$ of $R$. Repeating the above construction with the complement of $\sigma(S_1)$ in $R_1$, and proceeding inductively, we obtain a sequence of \'etale coverings $S_i \rightarrow S_{i-1}$, for $k=2,\ldots,k$, such that the combined pullback $f': C' \rightarrow S'$ of $f:C\rightarrow S$ to $S':= S_k$ admits $k$ sections of nodes.
\ebew

Using this result, we can now give the proof of lemma \ref{33}. 

\proofof{lemma \ref{33}}
Let $f : C \rightarrow S$ be a stable curve satisfying
$(i)$ - $(v)$. Because of $(i)$ each fibre of each $f_i$ has at least $3g-4$
nodes, so the induced morphism $ \theta_{f_i} : S \rightarrow \msbar_g$ factors
through an irreducible component of $\msbar_g^{(3g-4)}$. By the
definition of $\dppbar$,  the induced morphism $\theta_{{f_0}} : S^{(0)}
\rightarrow \msbar_g$ factors through $\dppbar$. Because of properties $(ii)$ and 
$(iii)$,  the deformations $f_i: C^{(i)}\rightarrow
S^{(i)}$ preserve $3g-4$ singularities of $ C^{(i-1)}_{s_-^{(i)}}$,
hence $S^{(i)}$ maps under $\theta_{{f_i}}$ to the same irreducible
component of $\msbar_g^{(3g-4)}$ as $S^{(i-1)}$. So, by induction, $S$
maps to   $\dppbar$.

Conversely, let $f : C \rightarrow S$ be a stable curve of genus $g$, for some irreducible and reduced scheme $S$, 
such that the induced morphism $\theta_f : S \rightarrow \msbar_g$
factors through $\dppbar$. 
Choose some closed point $s\in
S$. This point maps to a point $\theta_f(s)\in \dppbar$, and $f: C
\rightarrow S$ is a deformation of the curve represented by the point
$\theta_f(s)\in  \msbar_g$. Since each fibre of $f$ has at least $3g-4$ nodes, in an \'etale  neighbourhood $S'$ of $s\in S$ there exist sections $ \rho_2,\ldots,\rho_{3g-3}: S'
\rightarrow C'= C\times_S S'$ satisfying $(i)$. 

Note that no universal family exists over $\msbar_g$. However, for
each point $x \in \dppbar$ representing a stable curve $C_x
\rightarrow \Spec(k)$ there exists an \'etale  neighbourhood $S_x$ of $x$ in $\dppbar$,
which is the base of a deformation of $C_x$, and such that there are
$3g-4$ sections of nodes. Since $\dppbar$ is proper, a finite number of
such neighbourhoods suffices to cover all of $\dppbar$. One of these neighbourhoods contains the point $\theta_f(s)$. 

From this, one constructs inductively a family of stable curves together with  an enumeration of the sections of nodes compatible
with conditions $(ii)$ - $(v)$.

Note that there is one special case: if $f : C\rightarrow S$
admits $3g-3$ sections of nodes over $S'$,  then $\theta_f : S \rightarrow \msbar_g$ is necessarily constant because $\msbar_g^{(3g-3)}$ is discrete in $\msbar_g$. The above conditions then determine the choice of the appropriate $3g-4$ sections. 
\ebew

\begin{rk}\em
Note that the enumeration of the nodes of $C_0$ does not distinguish via this construction an enumeration of the nodes of the fibres of $f : C
\rightarrow S$. 
\end{rk}

\begin{notation}\em\label{N35}
Let $C_0^{+}$\label{n12} denote the subscheme of the fixed curve $C_0$ which is the reduced union of all
irreducible components of $C_0$ containing the node $P_1$, together
with marked points $Q_1,\ldots,Q_m \in \{ P_2,\ldots,P_{3g-3}\}$, such
that $Q_i$ is not a singularity on $C_0^{+}$. Similarly, let $C_0^{-}$ denote the closure of the
complement of  $C_0^{+}$ in $C_0$, together with the same marked points $Q_1,\ldots,Q_m$. Note that the chosen enumeration of the nodes of $C_0$ distinguishes an ordering of the marked points of the subcurves.  

If we need to emphasize the defining node we will sometimes write $C_0^+(P_1)$ instead of $C_0^+$, and $C_0^-(P_1)$ instead of $C_0^-$. 
\end{notation}

\begin{rk}\em\label{R36}
The subcurve $f_0^{+}: C_0^{+} \rightarrow \Spec(k)$ is in a natural way an $m$-pointed stable curve of some genus
$g^{+}$. The subcurves $C_0^{+}$ and $C_0^{-}$
meet  exactly in the points $Q_1,\ldots,Q_m$.  The curve $C_0^{-}$
may have more than one connected  
component, so it is an $m$-pointed prestable curve of some genus
$g^{-}$ as in remark \ref{M23}.  Of course, $g^{+}+g^{-}+m-1 =
g$. The
connected components of $C_0^{-}$ are again pointed stable curves. 

Furthermore, $C_0^{+}$ has the maximal number of sigularities
possible for such a curve, and the same is true for
$C_0^{-}$. Indeed, if there were an $m$-pointed stable 
curve of genus $g^{+}$ with more singularities, glueing with
$C_0^{-}$ at the marked points would produce a stable curve of genus
$g$ with more that $3g-3$ nodes, which is impossible. 
\end{rk}

\begin{rk}\em
Of course, the enumeration of the
marked points of $C_0^{+}$ and $C_0^{-}$ is not intrinsically determined by $C_0$, but by our initial  choice of the
enumeration of the nodes of $C_0$.
\end{rk}

\begin{rk}\em\label{37}
Only very few types of $m$-pointed stable curves $C_0^{+}$ occur in
this way, and here is  a
complete list. Since each irreducible component of $C_0$ is a rational curve with at most three nodes of $C_0$ lying on it, we must have $m \le 4$. We still assume $g \ge 3$. 
\begin{description}
\item[$m=4$ :] There exist exactly three types of non-isomorphic nodal $4$-pointed
stable curves of genus $g^{+}=0$, each consisting of two rational
curves meeting each other transversally in one point. 

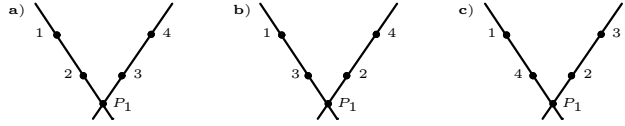
\begin{center}
\setlength{\unitlength}{0.00016667in}
\begingroup\makeatletter\ifx\SetFigFont\undefined%
\gdef\SetFigFont#1#2#3#4#5{%
  \reset@font\fontsize{#1}{#2pt}%
  \fontfamily{#3}\fontseries{#4}\fontshape{#5}%
  \selectfont}%
\fi\endgroup%
{\renewcommand{\dashlinestretch}{30}

}
\end{center}

The curves pictured above are represented by the three points which lie in the boundary of the
moduli space $\msbar_{0,4}$. If the condition on the existence of a
node at $P_1$ is dropped, we 
obtain the full moduli space $\msbar_{0,4}$.

\item[$m=2$ :] There exist two types of non-isomorphic nodal $2$-pointed
stable curves of genus $g^{+}=1$, both with two nodes and two rational irreducible components, as pictured below.\\

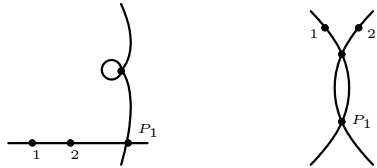
\begin{center}
\setlength{\unitlength}{0.00016667in}
\begingroup\makeatletter\ifx\SetFigFont\undefined%
\gdef\SetFigFont#1#2#3#4#5{%
  \reset@font\fontsize{#1}{#2pt}%
  \fontfamily{#3}\fontseries{#4}\fontshape{#5}%
  \selectfont}%
\fi\endgroup%
{\renewcommand{\dashlinestretch}{30}

}
\end{center}

Both of them are represented by isolated points in the
moduli space $\msbar_{1,2}$. The locus of curves in $\msbar_{1,2}$,
which preserve the node which is not $P_1$, is the boundary divisor
$\Delta_0$, representing irreducible singular stable curves of genus $1$, and  the only two degenerations thereof that are possible.

\item[$m=1$ :] The only instance here is a curve of arithmetic genus 1 with
one node at $P_1$, i.e. the one point in the boundary of $\msbar_{1,1}$.\\

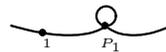
\begin{center}
\setlength{\unitlength}{0.00016667in}
\begingroup\makeatletter\ifx\SetFigFont\undefined%
\gdef\SetFigFont#1#2#3#4#5{%
  \reset@font\fontsize{#1}{#2pt}%
  \fontfamily{#3}\fontseries{#4}\fontshape{#5}%
  \selectfont}%
\fi\endgroup%
{\renewcommand{\dashlinestretch}{30}
\begin{picture}(5091,1314)(0,-10)
\thicklines
\put(2996,975){\ellipse{604}{604}}
\path(3033,675)(3032,675)(3031,674)
	(3028,672)(3023,669)(3017,665)
	(3008,659)(2996,652)(2983,644)
	(2967,635)(2948,624)(2927,612)
	(2903,600)(2878,586)(2850,572)
	(2821,558)(2789,544)(2756,529)
	(2721,515)(2684,501)(2644,487)
	(2603,473)(2559,460)(2513,448)
	(2463,436)(2410,425)(2353,415)
	(2293,405)(2227,397)(2157,389)
	(2083,383)(2004,379)(1920,376)
	(1833,375)(1761,376)(1689,378)
	(1617,381)(1545,385)(1475,390)
	(1407,396)(1340,402)(1275,410)
	(1211,418)(1149,426)(1088,435)
	(1028,444)(970,454)(913,465)
	(856,475)(801,486)(747,497)
	(693,509)(641,520)(590,532)
	(539,544)(490,555)(443,567)
	(397,578)(353,589)(311,600)
	(272,610)(235,620)(201,629)
	(170,637)(142,645)(118,651)
	(97,657)(79,662)(65,666)
	(54,669)(45,672)(39,673)
	(36,674)(34,675)(33,675)
\path(2958,675)(2959,674)(2961,673)
	(2964,670)(2970,666)(2977,660)
	(2987,653)(3000,643)(3016,633)
	(3033,620)(3054,606)(3076,592)
	(3100,576)(3127,560)(3155,544)
	(3185,527)(3216,511)(3249,495)
	(3284,480)(3321,465)(3361,450)
	(3402,437)(3447,424)(3494,413)
	(3546,402)(3601,393)(3660,386)
	(3723,380)(3789,376)(3858,375)
	(3920,376)(3982,379)(4043,383)
	(4103,389)(4160,397)(4216,405)
	(4270,415)(4322,425)(4372,436)
	(4421,448)(4469,460)(4515,473)
	(4560,487)(4604,501)(4647,515)
	(4689,529)(4730,544)(4769,558)
	(4807,572)(4843,586)(4877,600)
	(4908,612)(4937,624)(4963,635)
	(4986,644)(5005,652)(5021,659)
	(5034,665)(5044,669)(5050,672)
	(5055,674)(5057,675)(5058,675)
\put(2958,675){\blacken\ellipse{150}{150}}
\put(2958,675){\ellipse{150}{150}}
\put(1008,450){\blacken\ellipse{150}{150}}
\put(1008,450){\ellipse{150}{150}}
\put(2808,0){\makebox(0,0)[lb]{\smash{{{\SetFigFont{5}{6.0}{\rmdefault}{\mddefault}{\updefault}$P_1$}}}}}
\put(1008,0){\makebox(0,0)[lb]{\smash{{{\SetFigFont{5}{6.0}{\rmdefault}{\mddefault}{\updefault}$1$}}}}}
\end{picture}
}
\end{center}

If the condition on the existence of a node at $P_1$ is dropped, we
obtain the full moduli space $\msbar_{1,1}$. 
\end{description}

The following two cases only occur if $g=2$. Here the maximal number
of nodes possible is $3$, and this is realized precisely in the
following two curves.

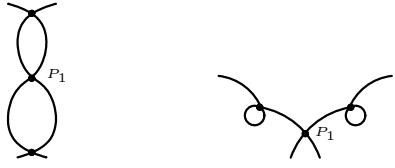
\begin{center}
\setlength{\unitlength}{0.00016667in}
\begingroup\makeatletter\ifx\SetFigFont\undefined%
\gdef\SetFigFont#1#2#3#4#5{%
  \reset@font\fontsize{#1}{#2pt}%
  \fontfamily{#3}\fontseries{#4}\fontshape{#5}%
  \selectfont}%
\fi\endgroup%
{\renewcommand{\dashlinestretch}{30}

}
\end{center}

In both cases $C_0 = C_0^{+}$ holds. Curves of genus $2$ which are deformations preserving the two nodes which are different from  $P_1$ are irreducible curves of geometric genus $0$ with two nodes, and their locus 
is a curve in $\msbar_2$. 
\end{rk}

\begin{notation}\em\label{38}
Let $C_1^{-},\ldots ,C_r^{-}$\label{n13} denote the connected components of
$C_0^{-}$. They are represented by points of not necessarily pairwise different 
moduli spaces  $\msbar_{{g_1^{-},m_1}}, \ldots ,
\msbar_{{g_r^{-},m_r}}$, with $g_1^{-} + \ldots + g_r^{-} = 
g^{-}$ and $m_1 + \ldots + m_r = m$.

From now on we assume that $n$ is chosen large enough to work
simultaneously for all of the finitely many Hilbert schemes $\hilb{g,n,0}$,
$\hilb{g^{+},n,m}$, and $\hilb{g_i^{-},n,m_i}$ for $i= 1,\ldots,r$. Let $N := (2g-2)n-g+1$, and put $N^{+} := (2g^{+}-2+m)n-g^{+}+1$. Note that always $N^{+} \le N$ holds.  
\end{notation}

\begin{construction}\label{new39}\em
Essential for our studies is the existence of a morphism
\[  \thx \: : \quad \hx \: \longrightarrow \: \hilb{g,n,0}, \]
where $\hx$ is a subscheme of the Hilbert scheme  $\hilb{g^{+},n,m}$, which will be defined formally in \ref{def_split}. The morphism $\thx$ shall have some special properties, like for example its equivariance with respect to certain group actions. To construct this map  we will proceed in several steps, which are described in paragraphs \ref{const_part_1}, \ref{const_part_2}, \ref{const_part_3}, \ref{const_part_4} and proposition \ref{const_prop}. 
\end{construction}

\begin{abschnitt}\em\label{const_part_1}
{\em Construction step 1.} 
Consider the universal embedded $m$-pointed stable curve $u^{+}: {\cal C}_{{g^{+},n,m}} \rightarrow
\hilb{g^{+},n,m}$, which is the restriction of the universal curve on the full Hilbert scheme. We glue this along the $m$ given sections to the
trivial $m$-pointed prestable curve $C_0^{-}\times \hilb{g^{+},n,m}
\rightarrow \hilb{g^{+},n,m}$. This produces a stable curve $u_0: {\cal
C}_0\rightarrow  \hilb{g^{+},n,m}$ of genus
$g$, compare \cite[Thm. 3.4]{Kn}. Actually, Knudsen calls this procedure of glueing along a pair of disjoint sections ``clutching''. 

Let  $i : {\cal C}_{{g^{+},n,m}} \hookrightarrow {\cal C}_0 $ denote the inclusion morphism from above. 
Consider the restriction of the canonical invertible dualizing sheaf $\omega_{{\cal C}_0/\hilb{g^{+},n,m}}$ to a fibre $C_h \rightarrow \Spec(k)$ of $u_0: {\cal
C}_0\rightarrow  \hilb{g^{+},n,m}$, for some closed point $h\in \hilb{g^{+},n,m}$. Let $n : \hat{C}_h \rightarrow C_h$ be the normalization of $C_h$, and for a node $P_i\in C_h$ let $A_i$ and $B_i$ denote its preimages in $\hat{C}_h$. 
The restriction of $\omega_{{\cal C}_0/\hilb{g^{+},n,m}}$ to $C_h$ can be described as the sheaf of those $1$-forms on $\hat{C}_h$ which have at most simple poles in $A_i$ and $B_i$, and such that the sum of the residues in $A_i$ and $B_i$ is zero. Compare for example \cite[section 1]{Kn}. 

By construction, the curve $u_0: {\cal
C}_0\rightarrow  \hilb{g^{+},n,m}$ comes with $m$ disjoint sections of nodes, which are given by the sections of glueing. Let $n_0 : \hat{C}_h' \rightarrow C_h$ denote the partial normalization of $C_h$, in the sense that it is the normalization map near those  nodes resulting from the glueing, and an isomorphism elsewhere. In particular, $\hat{C}_h'$ preserves exactly those nodes, which are not glueing nodes.  Thus $\hat{C}_h'$ decomposes into two disjoint subcurves, one of them isomorphic to $C_0^-$, the other one isomorphic to $C^+_h$, where $C_h^+$ is the fibre of the universal curve $u^{+}: {\cal C}_{{g^{+},n,m}} \rightarrow
\hilb{g^{+},n,m}$ over the point $h$. This induces a morphism $\iota: C_h^+ \rightarrow \hat{C}_h'$, such that $i= n_0\circ \iota$, and an isomorphism
$\iota^\ast n_0^\ast \omega_{C_h/k} \cong \omega_{C_h^+/k}(Q_1+\ldots+Q_m)$, where $Q_1,\ldots,Q_m$ denote the marked points on $C_h^+$. Thus there is a natural isomorphism of sheaves
\[  i^\ast \omega_{{\cal C}_0/\hilb{g^{+},n,m}} \:\: \cong \:\: \omega_{{\cal C}_{{g^{+},n,m}}/\hilb{g^{+},n,m}}({\cal S}_1+\ldots+{\cal S}_m) .\]
There is a natural surjection of sheaves $\omega_{{\cal C}_0 / \hilb{g^{+},n,m}} \rightarrow i_\ast i^\ast \omega_{{\cal C}_0/\hilb{g^{+},n,m}}$. After applying $n$-th tensor powers, and push-forward, we obtain a surjection 
\[ (u_0)_\ast(\omega_{{\cal C}_0/\hilb{g^{+},n,m}})^{\otimes n} \:\: \longrightarrow \:\: u_\ast^+ i^\ast (\omega_{{\cal C}_0/\hilb{g^{+},n,m}})^{\otimes n} \]
of locally free sheaves, 
using the factorization $u^+ = u_0\circ i$. Together with the above isomorphism, and applying Grothendieck's construction of the associated projective bundle, we obtain an inclusion of projective bundles
\[ \PP u^{+}_\ast(\omega_{{\cal C}_{{g^{+},n,m}}/\hilb{g^{+},n,m}}({\cal S}_1+\ldots+{\cal S}_m))^{\otimes n} \hookrightarrow  \PP(u_0)_\ast(\omega_{{\cal C}_0/\hilb{g^{+},n,m}})^{\otimes n} .\]

Recall that $f_0 : C_0 \rightarrow \Spec(k)$ is the stable curve fixed throughout this section, and $f_0^{+} : C_0^{+}  \rightarrow \Spec(k)$ is the subcurve defined in \ref{N35}. We fix once and for all an embedding of $C_0$ into $\PP^N$, which is equivalent to distinguishing an isomorphism 
\[ \PP(f_0)_\ast(\omega_{C_0/k})^{\otimes n} \cong \PP^N ,\]
together with an embedding of $C_0^{+}$ into $\PP^{N^{+}}$, corresponding to an isomorphism
\[ \PP(f_0^{+})_\ast(\omega_{C_0^{+}/k}(Q_1+\ldots+ Q_m))^{\otimes n} \cong \PP^{N^{+}} ,\]
in such a way that the diagram
\[  \diagram 
&C_0^{+}  \ddto|<\hole|<<\ahook \rto|<\hole|<<\ahook & \PP^{{N^{+}}} \ddto|<\hole|<<\ahook & \cong & \PP(f_0^{+})_\ast(\omega_{C_0^{+}/k}(Q_1+\ldots+ Q_m))^{\otimes n}\ddto|<\hole|<<\ahook \\
(\diamond)\\
&C_0 \rto|<\hole|<<\ahook & \PP^N &\cong &  \PP(f_0)_\ast(\omega_{C_0/k})^{\otimes n} 
\enddiagram \]
commutes. 

The universal curve $u^{+}: 
{\cal C}_{{g^{+},n,m}} \rightarrow
\hilb{g^{+},n,m}$ determines a trivialization 
\[ \begin{array}{c}
\PP u^{+}_\ast (\omega_{{\cal C}_{{g^{+},n,m}}/\hilb{g^{+},n,m}}({\cal S}_1+\ldots+{\cal S}_m))^{\otimes n} \cong \PP^{N^{+}} \times \hilb{g^{+},n,m} \hfill\\[4mm]
\hfill \cong \PP (f_0^{+})_\ast (\omega_{{C_0^{+}/k}}(Q_1+\ldots+Q_m))^{\otimes n}  \times \hilb{g^{+},n,m}.
\end{array}
\]
The tautological section shall be denoted by
\[ \tau^+ \: : \quad \hilb{g^{+},n,m} \: \longrightarrow \: \PP u^{+}_\ast (\omega_{{\cal C}_{{g^{+},n,m}}/\hilb{g^{+},n,m}}({\cal S}_1+\ldots+{\cal S}_m))^{\otimes n}. \] 
If $[C^+] \in \hilb{g^{+},n,m}$ is a point representing an embedded  curve $C^+$, then by definition of $\tau^+$, the point $\tau^+([C^+]) $ is fixed under the action of the group
\[ \Aut(C^+) \cong \Stab_{\PGL(N^++1)}([C^+]) \]
on the fibre of the projective bundle over $[C^+]$. Here, as always, the notation $\Aut(C^+)$ denotes the group of automorphisms of $C^+$ considered as an $m$-pointed stable curve, i.e. the group consisting of those automorphisms of the scheme $C^+$, which fix each one of its $m$ marked points.  
\end{abschnitt}

\begin{rk}\em\label{R311}\label{pgl_plus}
Via diagram $(\diamond)$ of the above construction, we may consider $\PGL(N^{+}+1)$ as a subgroup of $\PGL(N+1)$.  To do this, we view $\,\PP^{N^{+}}$ as a subspace of $\,\PP^N$ under the above embedding. Then $\PGL(N^{+}+1)$ is in a natural way  contained in the stabilizer subgroup of $\,\PP^{N^{+}}$. 

The group of automorphisms of $C_0$ can be identified with a subgroup of $\PGL(N+1)$ via
$ \Aut(C_0) \cong \Stab_{\PGL(N+1)}([C_0])$, and the diagram
\[  \diagram 
 \Aut(C_0^+) \dto|<\hole|<<\ahook & \cong  &\Stab_{\PGL(N^+ +1)}([C_0^+]) \dto|<\hole|<<\ahook & \subset& \PGL(N^+ +1) \dto|<\hole|<<\ahook \\
\Aut(C_0) & \cong &\Stab_{\PGL(N+1)}([C_0])&\subset & \PGL(N +1)
\enddiagram \]
commutes. 
 
Furthermore, for any embedded curve $C^+$, which is represented by a point $[C^+]\in \hilb{g^{+},n,m}$, there is a natural  embedding of $\Aut(C^+)$ first into the group $\PGL(N^+ +1)$, and from there into $\PGL(N+1)$. 
\end{rk}

\begin{abschnitt}\em\label{const_part_2}
{\em Construction step 2.} 
As in \ref{const_part_1}, there is also a natural surjection  of the sheaves
\[ j^\ast \omega_{{\cal C}_0/\hilb{g^+,n,m}} \:\: \rightarrow \:\: \omega_{C_0^-\times\hilb{g^+,n,m}/\hilb{g^+,n,m}}({\cal S}_1+\ldots+{\cal S}_m)  , \]
where $j : C_0^-\times\hilb{g^+,n,m} \rightarrow {\cal C}_0$ denotes the embedding as a subcurve over $\hilb{g^+,n,m}$ as  constructed above. Clearly there is a trivialization
\[ \PP^{N^-}\times \hilb{g^+,n,m} \cong \PP ( (\pr_2)_\ast (\omega_{C_0^-\times\hilb{g^+,n,m}/\hilb{g^+,n,m}}({\cal S}_1+\ldots+{\cal S}_m))^{\otimes n}) \]
as a projective bundle over $\hilb{g^+,n,m}$, where $N^- = (2g^- -2+m)n-g^- +1$, provided $n$ is chosen sufficiently large. The tautological section shall be denoted by $\tau^-$. 

We may assume that the corresponding  embedding of $C_0^-$ into $\PP^{N^-}$ is generic in the sense that the $m$ marked points on $C_0^-$ are contained in no subspace of dimension $m-2$. 

The above projective bundle is again in a natural way a subbundle of the  projective bundle
$\PP(u_0)_\ast(\omega_{{\cal C}_0/\hilb{g^+,n,m}})^{\otimes n} $. 
\end{abschnitt}

\begin{abschnitt}\em \label{const_part_3}
{\em Construction step 3.} 
Combining the constructions of \ref{const_part_1} and \ref{const_part_2}, we obtain two subbundles of $\PP(u_0)_\ast(\omega_{{\cal C}_0/\hilb{g^+,n,m}})^{\otimes n} $, and a surjective morphism $\pi$:
 \[ 
\begin{array}{ccc}
\PP u^{+}_\ast (\omega_{{\cal C}_{{g^{+},n,m}}/\hilb{g^{+},n,m}}(\sum_{i=1}^m{\cal S}_i))^{\otimes n} \\[3mm]
\oplus &\stackrel{\pi}{\longrightarrow} & \PP(u_0)_\ast( \omega_{{\cal C}_0/\hilb{g^+,n,m}})^{\otimes n}  \\
\PP (\pr_2)_\ast (\omega_{C_0^-\times\hilb{g^+,n,m}/\hilb{g^+,n,m}}(\sum_{i=1}^m{\cal S}_i))^{\otimes n}
\end{array}  \]
of projective bundles over $\hilb{g^+,n,m}$. Note that the intersection of the two subbundles can be described as the projectivization  of the cokernel of an injective  homomorphism 
\[ 
\begin{array}{ccc}
&& u^{+}_\ast (\omega_{{\cal C}_{{g^{+},n,m}}/\hilb{g^{+},n,m}}({\cal S}_1+\ldots+{\cal S}_m))^{\otimes n} \\[3mm]
 (u_0)_\ast( \omega_{{\cal C}_0/\hilb{g^+,n,m}})^{\otimes n}&\rightarrow & \oplus   \\
&& (\pr_2)_\ast (\omega_{C_0^-\times\hilb{g^+,n,m}/\hilb{g^+,n,m}}({\cal S}_1+\ldots+{\cal S}_m))^{\otimes n}
\end{array}  \]
of locally free sheaves over the base $\hilb{g^+,n,m}$. 

Thus, for a fibre of $u_0: {\cal C}_0\rightarrow  \hilb{g^{+},n,m}$, the subcurve $C_0^-$ is embedded into a subspace $\PP^{N^-} \subset \PP^N$, while the closure of its complement is embedded into a subspace $\PP^{N^+} \subset \PP^N$, and the intersection of $\PP^{N^-}$ and $\PP^{N^+}$ is a subspace of dimension $m-1$. This is just the subspace spanned by the $m$ marked points of $C_0^-$ embedded into $\PP^N$. 

By construction, for the two sections of the subbundles holds
\[ \pi \circ \tau^+ \: = \: \pi\circ \tau^- \] 
as a  section from $\hilb{g^+,n,m}$ into the bundle $\PP(u_0)_\ast(\omega_{{\cal C}_0/\hilb{g^+,n,m}})^{\otimes n} $. In \ref{const_part_4} we will use this section to construct a trivialization of the ambient $\PGL(N+1)$-bundle, extending the trivializations of the subbundles. Before we can do that, we have to take a closer look at the automorphism groups of stable curves. 
\end{abschnitt}

\begin{rk}\em \label{R311b}
In exactly the same way as in remark \ref{pgl_plus} for the group $\PGL(N^+ +1)$, the construction \ref{const_part_2} induces an embedding of $\PGL(N^- +1)$ into $\PGL(N+1)$. The group of automorphisms of $C_0^-$, considered as an $m$-pointed stable curve, can be identified with the group 
\[ \Aut(C_0^-) \cong \Stab_{\PGL(N^- +1)}([C_0^-]) ,\]
i.e. the subgroup of those elements in $\PGL(N^- +1)$, which stabilize the  embedded tuple
\[ (C_0^-,Q_1,\ldots,Q_m) \in \PP^{N^-} \times (\PP^{N^-})^m \]
with respect to the diagonal action. More formally one should thus write 
\[ \Aut(C_0^-) = \Aut(C_0^-;Q_1,\ldots,Q_m) \label{n15},\]
to set it apart from the following group. Denote by
\[ \Aut(C_0^-;\CQ) = \Aut(C_0^-;\{Q_1,\ldots,Q_m\}) \label{n16} \]
the group of all automorphisms of $C_0^-$ as a scheme over $\Spec(k)$, which map the set $\CQ := \{Q_1,\ldots,Q_m\}$\label{n16a} to itself. We can think of this group as the group of automorphisms of $C_0^-$ considered as curve with $m$ distinguished, but not ordered points. The notation $\Aut(C_0^-)$ will always stand for the group of those automorphisms of $C_0^-$, which respect the ordering of the marked points. 

Obviously, there is an inclusion
\[ \Aut(C_0^-) \: \subset \Aut(C_0^-;\CQ) .\]
We can also identify this group with a subgroup of $\PGL(N^- +1)$ via
\[ \Aut(C_0^-;\CQ) \cong \Stab_{\PGL(N^- +1)}(([C_0^-],\CQ)),\]
i.e. the subgroup of those elements in $\PGL(N^- +1)$, which stabilize the embedded pair
\[ (C_0^-,\CQ) \subset \PP^{N^-} \times \PP^{N^-} \]
with respect to the diagonal action. The group $ \Aut(C_0^-;\CQ)$ can therefore in particular  be considered as a subgroup of $\PGL(N+1)$, too. 
\end{rk}

\begin{defi}\em
Let $\Sigma_m$ be the permutation group acting on the labels of the $m$ marked points of $C_0^{-}$ contained in $\CQ$. There is an induced  action of $\Sigma_m$ on the set of $m$-pointed curves. Let $\Gamma(C_0^-;\CQ)$\label{n17} denote the subgroup of $\Sigma_m$ such that for all $\gamma\in\Gamma(C_0^-;\CQ)$ the curve $\gamma(C_0^{-})$ is isomorphic to $C_0^{-}$ as an $m$-pointed curve, i.e. such that there exists an automorphism $\phi: C_0^- \rightarrow C_0^-$ as a scheme over $\Spec(k)$ with
\[ \phi(Q_i) = Q_{\gamma(i)} \]
for $i=1,\ldots,m$. 
Define $\Gamma^{\circ}(C_0^-;\CQ)$\label{n18} to be the subgroup of $\Sigma_m$ consisting of all permutations $\gamma \in \Sigma_m$, such that $\gamma(C_0^{-})$ is isomorphic to $C_0^{-}$ and $\gamma(C_0^{+})$ is isomorphic to $C_0^{+}$, both considered as $m$-pointed curves.
\end{defi}

\begin{rk}\em
 There is obviously an inclusion
\[ \Gamma^{\circ}(C_0^-;\CQ) \: \subset \: \Gamma(C_0^-;\CQ) .\]
Note that in the cases $m=1$ and $m=2$ any automorphism of $C_0^{-}$ extends to all of $C_0$, so the groups $\Gamma^{\circ}(P_1;\CQ)$ and $\Gamma(P_1;\CQ)$ are equal in these cases. 
\end{rk}

\begin{rk}\em\label{R315b}
Let $\Aut(C_0;\{P_1\})$ denote the subgroup of automorphisms of
$C_0$ fixing the node $P_1$. Clearly, this subgroup stabilizes the
subcurves $C_0^{+}$ and $C_0^{-}$ in $C_0$.  An automorphism $\phi\in \Aut(C_0;\{P_1\})$ gives rise to a permutation of the $m$ marked points of $C_0^{+}$ and $C_0^{-}$, so there is a group homomorphism 
\[ \Aut(C_0;\{P_1\}) \: \rightarrow \: \Sigma_m .\]
The image of this homomorphism is equal to $\Gamma^{\circ}(C_0^-;\CQ)$.
\end{rk}

\begin{defi}\em\label{D318a}
Let $\Gamma \subset \Sigma_m$ be a subgroup. We define the {\em group of $\mgamma$-automorphisms of $C_0^-$} as the subgroup $\Aut_\Gamma(C_0^-;\calq)$\label{n19}  of those elements  $\phi\in\Aut(C_0^-;\calq)$, for which  there exists a permutation $\gamma \in \Gamma$ such that 
\[ \phi(Q_i) = Q_{\gamma(i)} \]
holds for all $i=1,\ldots,m$, where $\calq=\{Q_1,\ldots,Q_m\}$. 
\end{defi}

\begin{rk}\em \label{R319}
$(i)$ The group $\Aut_\Gamma(C_0^-;\calq)$ is the group of automorphisms of the curve $C_0^-$ when considered an an $\mgamma$-pointed curve. One clearly has
\[ \Aut_{\{\id\}}(C_0^-;\calq) = \Aut(C_0^-), \] 
and the identity
\[ \Aut_{\Sigma_m}(C_0^-;\calq) = \Aut(C_0^-;\CQ) .\] 
Note also that by definition
\[  \Aut_{\Gamma(C_0^-;\CQ)}(C_0^-;\calq) = \Aut(C_0^-;\CQ).\]
$(ii)$ For any $\phi \in  \Aut(C_0^-;\{Q_1,\ldots,Q_m\})$ there is a unique permutation $\gamma \in \Sigma_m$ such that $\phi(Q_i) = Q_{\gamma(i)}$
holds for all $i=1,\ldots,m$. Thus there is a natural morphism of groups
\[ \Aut(C_0^-;\{Q_1,\ldots,Q_m\}) \rightarrow \Sigma_m .\]
By definition, the image of this group is the subgroup $\Gamma(C_0^-;\CQ)$. 
For any subgroup $\Gamma$ in $\Gamma(C_0^-;\CQ)$  the sequence of homomorphisms 
\[ \id \rightarrow \Aut(C_0^{-})\rightarrow \Aut_\Gamma(C_0^-;\calq)  \rightarrow \Gamma\rightarrow \id \]
is exact. 
\end{rk}

\begin{notation}\em
Let $\Gamma= \Gamma(C_0^-;\calq)$. For a stable curve $f^+:C^+ \rightarrow \Spec(k)$ of genus $g^+$ with $m$ marked points, which is represented by a point $[C^+]\in \hilb{g^+,n,m}$, we define the group 
\[ \Aut_\Gamma(C^+;\calq) \label{m1}\]
as the group of automorphisms of $C^+$, when considered as an $\mgamma$-pointed stable curve, with $\calq$ as its set of distinguished points, analogously to definition \ref{D318a}.
\end{notation}

\begin{lemma}\label{newsplit}
Let $f^+: C^+\rightarrow \Spec(k)$ be a stable curve of genus $g^+$ with $m$ marked points, which is represented by an inner point $[C^+]\in H_{g^+,n,m}$. Let $f: C\rightarrow \Spec(k)$ be the stable curve of genus $g$, which is obtained by glueing $C^+$ with $C_0^-$ in the $m$ marked points. Suppose that the curve $C$ has exactly $3g-4$ nodes. Then there is an exact sequence of groups
\[ \id \rightarrow \Aut(C_0^-) \rightarrow \Aut(C) \rightarrow \Aut_\Gamma(C^+;\calq) \rightarrow \id ,\]
where  $\Gamma= \Gamma(C_0^-;\calq)$. 
\end{lemma}

\proof
By definition, any automorphism in $\Aut(C_0^-)$ fixes the marked points of $C_0^-$, and these marked points are the glueing points of $C_0^-$ with $C^+$. By trivial extension, the group $\Aut(C_0^-)$ becomes a subgroup of $\Aut(C)$. Conversely, since $[C^+]\in H_{g^+,n,m}$ is an inner point, the curve $C^+$  is the unique connected subcurve of $C$ of genus $g^+$, which has not the maximal number of nodes as an $m$-pointed stable curve. Hence any automorphism of $C$ restricts to an automorphism on $C_0^-$ and to an automorphism on $C^+$, which both are not necessarily fixing the marked points. However, clearly only such permutation of marked points occur in this way on $C^+$, which are also effected by automorphisms on $C_0^-$. Thus by definition, an automorphism of $C$ restricted to $C^+$ is an automorphism of $C^+$ considered as an $\mgamma$-pointed stable curve, i.e. an element of $\Aut_\Gamma(C^+;\calq)$. 

All automorphisms in $\Aut_\Gamma(C^+;\calq)$ occur in this way. To see this, let $\rho^+ \in \Aut_\Gamma(C^+;\calq)$ be an automorphism of $C^+$ as an  $\mgamma$-pointed stable curve. It induces a permutation $\gamma$ of the marked points, with $\gamma\in \Gamma$. By the definition of $\Gamma$, there exists an automorphism $\rho^-\in\Aut(C_0^-;\calq)$ inducing $\gamma$. Glueing the automorphisms $\rho^+$ and $\rho^-$ defines an automorphism $\rho$ on $C$, which restricts to $\rho^+$. Thus we have shown that there exists an exact sequence of groups, as claimed.
\ebew

\begin{rk}\em\label{rk00328}
The exact sequence of lemma \ref{newsplit} does not split for all stable curves $f^+: C^+\rightarrow \Spec(k)$ with $[C^+]\in H_{g^+,n,m}$. The smallest counterexample occurs for $m=2$. Consider the curve $C_0^-$ consisting of two disjoint $1$-pointed nodal rational curves, so that $\Aut(C_0^-) = {\ZZ_2}^2$. Let $C^+$ be a $2$-pointed stable curve of genus $g^+ =1$, with one node, and such that there exists an automorphism interchanging the two marked points, so that  $\Aut_\Gamma(C^+;\calq) \cong \ZZ_2 \times \ZZ_2$.  Glueing $C_0^-$ and $C^+$ in their marked points gives a stable curve $C$ of genus $g=3$, whose group of automorphisms $\Aut(C)$ is isomorphic to the direct product of $\ZZ_2$ and the non-trivial extension of ${\ZZ_2}^2$ by $\ZZ_2$.     
\end{rk}

\begin{rk}\em\label{newsplitrk}
For a general point $[C^+]$ as in lemma \ref{newsplit}, representing a stable curve $f^+: C^+\rightarrow \Spec(k)$, which glued to $C_0^-$ produces a stable curve represented by a point in $M_{g}^{(3g-4)}$, the exact sequence of the lemma  splits. This can be seen as follows.

By the definition of $\Aut_\Gamma(C^+;\calq)$, there is an exact sequence of groups
\[ \id \rightarrow \Aut(C^+)\rightarrow \Aut_\Gamma(C^+;\calq) \rightarrow \Gamma .\]
The last arrow might not be surjective, as not any permutation in $\Gamma$ needs to be realizeable by an automorphism in $\Aut(C^+;\calq)$. Let $\Gamma'$ denote the image of $\Aut_\Gamma(C^+;\calq)$ in $\Gamma$. 

If $\Gamma' = \{\id\}$, then we have the equality $\Aut(C^+) =  \Aut_\Gamma(C^+;\calq)$, and a splitting homomorphism $\Aut_\Gamma(C^+;\calq) \rightarrow \Aut(C)$ can be defined by trivial extension. In particular, this settles the case $m=1$. 

Suppose that $m=2$. Then, by remark \ref{37}, the subcurve $C^+$ is a stable curve of genus $g^+ =1$, with $2$ marked points and one node. An automorphism of $C^+$, which interchanges the two marked points, exists only if the two marked points are in a very special position. Or, equivalently, if $[C^+]\in M_{1,2}'$ is a general point, then $\Aut_\Gamma(C^+;\calq) = \Aut(C^+)$, and therefore $\Gamma' = \{\id\}$. This settles the case $m=2$. 

Suppose that $m=4$. Then $C^+$ is a smooth rational curve, with $4$ marked points. For a general point $[C^+]\in M_{0,4}$, one sees  that $\Aut_\Gamma(C^+;\calq) = \Aut(C^+) = \{\id\}$, and thus $\Gamma' = \{ \id \}$ again.
\end{rk}

\begin{lemma}
Suppose that  $m=1$, or that the fixed curve $C_0^-$ has the property that its group of automorphisms $\Aut(C_0^-;\calq)$ splits as 
\[ \Aut(C_0^-;\calq) \:\: \cong \:\: \Aut(C_0^-)\times \Gamma(C_0^-;\calq) .\]
Then for all stable curves $f: C\rightarrow \Spec(k)$ of genus $g$ with exactly $3g-4$ nodes, which are obtained from glueing $C_0^-$ with a curve $C^+$, which is represented by a point $[C^+]\in H_{g^+,n,m}$, its group of automorphisms splits as
\[\Aut(C) = \Aut_\Gamma(C^+;\calq)\times \Aut(C_0^-),
 \]
where $\Gamma= \Gamma(C_0^-;\calq)$. 
\end{lemma}

\proof
If $m=1$, then the splitting is trivial, as any automorphism of $C^+$, fixing the one marked point, can be  extended to all of $C$. 

If $m=4$, then the curve $C^+$ is a smooth stable curve of  genus $g^+=0$ with $4$  marked points, see remark \ref{37}. In particular, we have $\Aut(C^+) = \{\id\}$, and therefore an isomorphism $\Aut_\Gamma(C^+;\calq) = \Gamma'$, using the exact sequence of remark \ref{newsplitrk}. By the assumption on the splitting of $\Aut(C_0^-;\calq)$ there exists an injective homomorphism $j: \Gamma'\rightarrow \Aut(C_0^-;\calq)$, such that for all $\gamma\in\Gamma'$ the image $j(\gamma)$ induces the permutation $\gamma$ on the marked points of $C_0^-$. Let $\rho^+ \in \Aut_\Gamma(C^+;\calq)$ be given, inducing a permutation $\gamma\in \Gamma'$ on the marked points of $C^+$. Glueing $\rho^+$ with $j(\gamma)$ defines an automorphism $\rho\in \Aut(C)$. Note that $\rho$ commutes with all elements of $\Aut(C_0^-)$, so that this construction gives indeed a splitting homomorphism.

If $m=2$, then the curve $C^+$ is a stable curve of genus $g^+=1$ with $2$ marked points and one node. If $\Gamma' = \{\id\}$, we are done. So suppose that  $\Gamma' = \ZZ_2$. We then have a splitting $ \Aut_\Gamma(C^+;\calq) = \Aut(C^+)\times \ZZ_2$. Using trivial extension on the first factor, and the splitting of  $\Aut(C_0^-;\calq)$ for the embedding of the second factor, we can construct a splitting homomorphism $\Aut_\Gamma(C^+;\calq) \rightarrow \Aut(C)$.
\ebew

\begin{notation}\em\label{def_split}
We define the nonempty reduced subscheme
\[ \hxbar \: \subset \: \hilb{g^{+},n,m} \label{m2} \]
of the Hilbert scheme $\hilb{g^{+},n,m}$ as  the locus of points $[C^+]$ parametrizing curves  $f^+:C^+\rightarrow \Spec(k)$ with the following property: if $f: C\rightarrow \Spec(k)$ denotes the stable curve of genus $g$, which is obtained by glueing $C^+$ with $C_0^-$ in the $m$ marked points, then  $C$ has exactly $3g-4$ nodes, and its group of automorphisms splits naturally as
\[ \Aut(C) = \Aut_\Gamma(C^+;\calq)\times \Aut(C_0^-) .\]
The intersection of $ \hxbar$ with $ H_{g^{+},n,m}$ is denoted by $\hx$.
\end{notation} 

\begin{abschnitt}\em\label{const_part_4}
{\em Construction step 4.}
Consider the morphism $\pi$ of projective bundles over $\hilb{g^{+},n,m} $ as in \ref{const_part_3}, restricted to the subscheme $\hxbar$:
\[ 
\begin{array}{ccc}
\PP u^{+}_\ast (\omega_{{\cal C}_{{g^{+},n,m}}/\hilb{g^{+},n,m}}(\sum_{i=1}^m{\cal S}_i))^{\otimes n} \\[3mm]
\oplus &\stackrel{\pi}{\longrightarrow} & \PP(u_0)_\ast( \omega_{{\cal C}_0/\hilb{g^+,n,m}})^{\otimes n}.  \\
\PP (\pr_2)_\ast (\omega_{C_0^-\times\hilb{g^+,n,m}/\hilb{g^+,n,m}}(\sum_{i=1}^m{\cal S}_i))^{\otimes n}
\end{array}  \]
By composition with $\pi$, both tautological sections $\tau^+$ and $\tau^-$ of the subbundles define the same  section 
\[\tau \: : \quad \hilb{g^{+},n,m}^\times  \: \longrightarrow \: \PP(u_0)_\ast( \omega_{{\cal C}_0/\hilb{g^+,n,m}})^{\otimes n}. \]

By \ref{const_part_1} and \ref{const_part_2}, there exist trivializations
\[ \theta^+: \quad \PP u^{+}_\ast (\omega_{{\cal C}_{{g^{+},n,m}}/\hilb{g^{+},n,m}}(\sum_{i=1}^m{\cal S}_i))^{\otimes n}|\hxbar \: \rightarrow \: \PP^{N^+} \times  \hxbar \]
and 
\[ \theta^-:  \PP (\pr_2)_\ast (\omega_{C_0^-\times\hilb{g^+,n,m}/\hilb{g^+,n,m}}(\sum_{i=1}^m{\cal S}_i))^{\otimes n}|\hxbar  \rightarrow \PP^{N^+} \times  \hxbar \]
as $\PGL(N^+ + 1)$-bundles and $\PGL(N^- + 1)$-bundles, respectively. We want to extend both of these to a common trivialization of the ambient $\PGL(N+1)$-bundle $\PP(u_0)_\ast( \omega_{{\cal C}_0/\hilb{g^+,n,m}})^{\otimes n}$. If $e$ is a point in the fibre of $\PP(u_0)_\ast( \omega_{{\cal C}_0/\hilb{g^+,n,m}})^{\otimes n}$ over a point $[C^+]\in \hxbar$, representing an embedded $m$-pointed stable curve $C^+$, then $e$ can be written as
\[ e = \gamma \cdot \tau([C^+]) \]
for some $\gamma \in \PGL(N+1)$. We now define 
\[ \theta : \quad  \PP(u_0)_\ast( \omega_{{\cal C}_0/\hilb{g^+,n,m}})^{\otimes n}|\hxbar \: \longrightarrow \: \PP^{N} \times  \hxbar\]
by
\[ \theta(e) := \gamma\cdot \theta^+(\tau^+([C^+])) ,\]
using the fact that the point $\tau([C^+]) = \pi\circ \tau^+([C^+])$ is contained in the subbundle $\PP u^{+}_\ast (\omega_{{\cal C}_{{g^{+},n,m}}/\hilb{g^{+},n,m}}(\sum_{i=1}^m{\cal S}_i))^{\otimes n}$, and the fixed embedding of the subspace  $\PP^{N^+} \subset \PP^N$ from \ref{const_part_1}. 
Note that in general $\gamma$ is not uniquely determined by the point $e$. To show that the above map is indeed well-defined, it suffices to show that all elements of $\PGL(N+1)$, which stabilize the point $\tau([C^+])$ also stabilize the point $\theta^+(\tau^+([C^+]))$. By construction, the elements stabilizing the point $\tau([C^+])$ in the fibre are exactly the elements of the group
\[ \Aut(C) \cong \Stab_{\PGL(N+1)}([C]), \]
where $C$ is the embedded curve in $\PP^N$  obtained by glueing $C^+$ and $C_0^-$ in their marked points.  Recall that $\tau = \pi\circ\tau^+ = \pi\circ\tau^- $, and thus $\theta^+(\tau^+([C^+]))=\theta^-(\tau^-([C^+])) \in \PP^N$, with respect to the fixed  embeddings of \ref{const_part_1} and \ref{const_part_2}. 

By the definition of the subscheme $\hilb{g^{+},n,m}^\times$ we have a splitting of $ \Aut(C) $ as a product $ \Aut_\Gamma(C^+;\calq)\times \Aut(C_0^-)$. Because $\tau^+$ and $\tau^-$ are the tautological sections of the subbundles, the points $\theta^+(\tau^+([C^+]))$ and  $\theta^-( \tau^-([C^+]))$ are fixed by  $ \Aut(C^+) $ and $ \Aut(C_0^-)$, respectively. Since the process of glueing makes the ordering of the $m$ glueing points contained in $\calq$ disappear, the point $\theta^+ (\tau^+([C^+]))$ is also fixed under the action of $ \Aut_\Gamma(C^+;\calq) $.
Therefore, the point $\theta^+(\tau^+([C^+])) = \theta^-( 
\tau^-([C^+]))$ is fixed unter the action of  $  \Aut(C) $, and we are done. 

We have thus constructed a trivialization of the restricted projective bundle
\[ \begin{array}{rcl}
\PP(u_0)_\ast(\omega_{{\cal C}_0/\hilb{g^{+},n,m}})^{\otimes n} \: | \: \hxbar
  & \cong &\PP(f_0)_\ast(\omega_{{C}_0/k})^{\otimes n}
\times 
\hxbar\\[3mm]
&\cong&
\PP^N \times 
\hxbar,
\end{array} \]
which is fibrewise compatible with both of the natural inclusions 
\[ \PP (f_0^{+})_\ast (\omega_{{C_0^{+}/k}}(Q_1+\ldots+Q_m))^{\otimes n} \: \: \hookrightarrow \: \: \PP(f_0)_\ast(\omega_{{C}_0/k})^{\otimes n}
\]
and 
\[ \PP (f_0^{-})_\ast (\omega_{{C_0^{-}/k}}(Q_1+\ldots+Q_m))^{\otimes n} \: \: \hookrightarrow \: \: \PP(f_0)_\ast(\omega_{{C}_0/k})^{\otimes n}
\]
from above. Globally, there is a commutative diagram of projective bundles over $\hilb{g^{+},n,m}^\times$: 
\[  \diagram 
\PP^{{N^{+}}} \times \hxbar \dto|<\hole|<<\ahook  \quad \cong & 
\PP u^{+}_\ast (\omega_{{\cal C}_{{g^{+},n,m}}/\hilb{g^{+},n,m}}(\sum_{i=1}^m{\cal S}_i))^{\otimes n} |\hxbar\dto|<\hole|<<\ahook\\
\PP^{{N}} \times \hxbar \quad \: \: \cong  & 
\PP(u_0)_\ast(\omega_{{\cal C}_0/\hilb{g^{+},n,m}})^{\otimes n}|\hxbar\\
\PP^{{N^{-}}} \times \hxbar \uto|<\hole|<<\ahook \quad\cong &
 \PP ( (\pr_2)_\ast (\omega_{C_0^-\times\hxbar/\hxbar}(\sum_{i=1}^m{\cal S}_i))^{\otimes n}).\uto|<\hole|<<\ahook 
\enddiagram \]
By the universal property of the Hilbert scheme $ \hilb{g^{+},n,m}$, the trivialization of the projective bundle $\PP(f_0)_\ast(\omega_{{C}_0/k})^{\otimes n}$ restricted to $ \hxbar$ 
defines a morphism
\[ \thx : \: \: \hxbar \rightarrow \hilb{g,n,0}  \label{n14}\]
which in general is not an embedding. By construction, this morphism is equvariant with respect to the action of the group $\PGL(N^+ +1)$, considered as a subgroup of $\PGL(N+1)$ as in remark \ref{pgl_plus}.

We summarize our construction in the following proposition. Using the universal property of Hilbert schemes, it is formulated geometrically in terms of embedded curves rather than those of projective bundles, while the content is of course equivalent to our construction.
\end{abschnitt}

\begin{prop}\label{const_prop}
Let $u^{+,\times}: {\cal C}^\times_{{g^{+},n,m}} \rightarrow
\hxbar$ denote the restriction of the universal embedded $m$-pointed stable curve  on the  Hilbert scheme $\hilb{g^{+},n,m}$ to $\hxbar$, and let $u_0^\times : {\cal
C}_0^\times \rightarrow  \hxbar$ denote the stable curve of genus $g$ obtained by glueing it along the $m$ given sections to the
trivial $m$-pointed prestable curve $C_0^{-}\times \hxbar
\rightarrow \hxbar$. 

Then there exists a global  embedding of the curve $u_0^\times : {\cal
C}_0^\times \rightarrow  \hxbar$ into $\pr_2: \PP^N\times \hxbar \rightarrow \hxbar$, which makes the following diagram commutative:
\[ \diagram 
{\cal C}_{g^{+},n,m}^\times\rto|<\hole|<<\ahook \dto|<\hole|<<\ahook 
& {\cal C}_0 \dto|<\hole|<<\ahook
&C_0^-\times  \hxbar\lto|<\hole|<<\ahook \dto|<\hole|<<\ahook 
\\
\PP^{N^{+}} \times \hxbar \drto \rto|<\hole|<<\ahook 
&\PP^N \times \hxbar \dto
&\PP^{N^{-}} \times \hxbar \dlto \lto|<\hole|<<\ahook \\
& \hxbar. 
\enddiagram \]
In particular, there exists a $\PGL(N^+ +1)$-equivariant morphism
\[ \thx : \: \: \hxbar \rightarrow \hilb{g,n,0} . \]
\end{prop}

\proof The construction was given in four steps in  the  paragraphs \ref{const_part_1}, \ref{const_part_2}, \ref{const_part_3},  and \ref{const_part_4}.
\ebew

\begin{rk}\em\label{R317}
The distinguished embeddings of $C_0$ into $\PP^N$ and of $C_0^{+}$
into $\PP^{N^{+}}$ from construction \ref{const_part_1} define points $[C_0]
\in  \hilb{g,n,0}$ and $[C_0^{+}] \in \hxbar$ such that
$\thx( [C_0^{+}]) = [C_0]$. At the same time, the embedding of $C_0^-$ into $\PP^{N^-}$ distinguishes a point $[C_0^-] \in \hilb{g^-,n,m}$. Clearly the embedding of $C_0^-\times \hxbar$ into $\PP^{N^-}\times\hxbar$ is independent of the fibre. Because we have chosen compatible trivializations in the construction of $\thx$ above, the  embedding of the trivial curve $C_0^{-}\times \hxbar \rightarrow \hxbar$ into $\PP^N \times \hxbar \rightarrow \hxbar$  is independent of the fibre. 
\end{rk}

\begin{prop}\label{P318}
The  morphism $\thx: \hxbar\rightarrow \hilb{g,n,0}$ from proposition \ref{const_prop} is invariant under the action of $\, \Gamma(C_0^-;\CQ)$. 
\end{prop}

\proof
Recall that $\Gamma(C_0^-;\CQ)$ acts freely on $\hilb{g^+,n,m}$. The projective bundle $\PP u^{+}_\ast (\omega_{{\cal C}_{{g^{+},n,m}}/\hilb{g^{+},n,m}}({\cal S}_1+\ldots+{\cal S}_m))^{\otimes n}$ is invariant under the action of $\Sigma_m$, and hence in particular under the action of $\Gamma(C_0^-;\CQ)$. This is just saying that the embedding of a fibre of $u^{+} : {\cal C}_{{g^{+},n,m}} \rightarrow \hilb{g^{+},n,m}$  into $\PP^{N^{+}}$ does not depend on the order of the marked points. By construction \ref{const_part_1}, the embedding of this projective bundle into $\PP (u_0)_\ast (\omega_{{\cal C}_{0}/\hilb{g^{+},n,m}})^{\otimes n}$ is also independent of the ordering of the marked points. Thus for any point $[C']\in \hxbar$ representing an embedded curve $C'\subset \PP^{N^+}$, the embedding of $C'$ into $\PP^N$ is the same as the embedding corresponding to the point $\gamma([C'])$, for any  $\gamma\in \Gamma(C_0^-;\CQ)$. So it remains to consider the complement of $C'$ in the curve $C$ representing the point $\thx([C'])$. The closure of the complement of $C'$ in $C$ is just the distinguished embedding of the curve $C_0^-$, by construction of the morphism $\thx$. By definition of $\Gamma(C_0^-;\CQ)$, and by the choice of the embedding of $\Aut(C_0^-;\CQ)$ into $\PGL(N+1)$ this complement is also invariant under the action of $\Gamma(C_0^-;\CQ)$. Therefore the whole curve $C\subset \PP^N$, and hence the point $\thx([C']) \in \hilb{g,n,0}$,  is fixed under the action of $\Gamma(C_0^-;\CQ)$. 
\ebew

\begin{cor}\label{C318}
The  morphism $\thx: \hxbar\rightarrow \hilb{g,n,0}$ is equivariant with respect to the action of the group $\PGL(N^+ +1)\times \Gamma(C_0^-;\CQ)$.
\end{cor}

\proof The equivariance with respect to the individual groups follows from proposition \ref{const_prop} and proposition \ref{P318}. Note that the actions of the groups $\PGL(N^+ +1)$ and $\Gamma(C_0^-;\CQ)$ on $\hxbar$ commute since the embedding of an $m$-pointed stable curve into $\PP^{N^+}$ is independent of the ordering of its marked points. 
\ebew 

\begin{notation}\em\label{D322}
Let $\pi : \hilb{g^+,n,m}\rightarrow \msbar_{g^+,m}$ denote the canonical quotient morphism, and let $\mxbar$ \label{m3} denote the image of $\hxbar$ in $\msbar_{g^+,m}$ as a reduced subscheme.  

The $\PGL(N^++1)$-equivariance of the morphism $\thx:\hxbar\rightarrow \hilb{g,n,0}$ from above  induces a morphism
\[ \thxbar : \: \: \mxbar \rightarrow \msbar_g , \]
which by construction maps into the boundary stratum $\msbar_g^{(3g-4)}$.

Let $\mgprimebar$\label{n21} denote the reduced closure of $\mxbar$ in $\msbar_{g^+,m}$. As we will see in remark \ref{R0323} below, the subscheme $\mgprimebar$ is smooth and irreducible, and therefore the morphism $\thxbar$ extends to a surjective morphism
\[ \overline{\Theta}' : \quad \mgprimebar \rightarrow \dpobar .\]
Note that $\overline{\Theta}'$ is the restriction of a natural morphism
\[\overline{\Theta} : \quad \msbar_{g^+,m} \rightarrow \msbar_g \label{n20},\]
which sends the point representing an $m$-pointed stable curve $C^+$ of genus $g$ to the point represening the stable curve obtained by glueing $C^+$ and $C^-_0$ in their marked points. The subscheme $\mgprimebar$ is the reduced preimage of $\dpobar$ in $\msbar_{g^+,m}$ under $\overline{\Theta}$. 

We define the subscheme $\hgprimebar$  \label{n22} as the reduction of the preimage of $\mgprimebar$ in $ \hilb{g^+,n,m}$ under the canonical quotient morphism $\pi$. It is the closure of $\hxbar$ inside $ \hilb{g^+,n,m}$. 
  
 For the interior parts we define the intersections  $\mx := \mxbar \cap \ms_{g^{+},m}$ and $\hx  := \hxbar \cap {H}_{g^{+},n,m}$, as well as 
$\mgprime := \mgprimebar \cap \ms_{g^{+},m}$ and $\hgprime  := \hgprimebar \cap {H}_{g^{+},n,m}$. 

For $\mgamma$-pointed stable curves, as in definition \ref{n6}, the notations $\hgammaxbar$ and $\mgammaxbar$, as well as $\hgprimebarquot$ and $\mgprimebarquot$ together with their open subschemes $\hgammax$, $\mgammax$, $\hgprimequot$, and $\mgprimequot$  \label{m4}, are defined ana\-logously, using the appropriate canonical quotient morphisms.  

We define the subscheme 
\[ \dxbar \subset \dpobar \label{m5}\]
as the reduced subscheme of $\dpobar$, which is the locus of points para\-metrizing stable curves $C$, for which the group of automorphisms splits naturally as $ \Aut(C) = \Aut_\Gamma(C^+;\calq)\times \Aut(C_0^-)$. Its intersection with the open subscheme $\dpo$ is denoted by $\dx$. 

Schematically, the subschemes defined above fit into the commutative diagram below. Note that there is an anlogous diagram for their open partners. 
\[ \diagram 
\hxbar \dto  & \subset   & \hgprimebar \dto& \subset & \hilb{g^+,n,m} \dto^\pi \\
\mxbar \dto & \subset   & \mgprimebar\dto & \subset & \msbar_{g^+,m}\dto^{\overline{\Theta}}\\
\dxbar&\subset & \dpobar &\subset & \msbar_g.
\enddiagram \]
The superscript ``$\times$'' is intended as a reminder of splitting automorphisms groups. All vertical arrows except $\overline{\Theta}$ are surjections. 
\end{notation}

\begin{rk}\em\label{R0323}
Let $\Delta_0 \subset \msbar_{g^{+},m}$ denote the divisor defined by the closure of the locus of irreducible singular stable curves. From remark \ref{37} we obtain immediately the following identifications. 
\[ \msbar_{g^{+},m}' = \left\{ \begin{array}{ll}
\msbar_{g^{+},m} & \Text{for } m = 1,4;\\
\Delta_0& \Text{for } m=2.
\end{array} \right. \]
Note  that  
in each case there is an isomorphism $\msbar_{g^{+},m}' \cong \PP^1$.
\end{rk}

Before we can proceed we need the following lemma.

\begin{lemma}\label{314a}
Let $f : C \rightarrow \Spec(k)$ be a stable curve of genus $g$, with at least $3g-4$ nodes, and which contains $C_0^{-}$ as a subcurve, such that the $m$ marked points of $C_0^{-}$ are nodes on $C$. Then either the closure of the complement of 
$C_0^{-}$ in $C$ is irreducible, or $C$ has $3g-3$ nodes.
\end{lemma}

\proof
The curve $C_0$ has $3g-3$ nodes and $2g-2$ irreducible components. Thus by
 our discussion in remark \ref{37} the number of irreducible components of $C_0^{-}$ is $2g-3$ if $m=1$, or  $2g-4$ if $m=2$ or $m=4$. No stable curve of genus $g$ can have more that $2g-2$ irreducible components, so the number of irreducible components of the closure of the complement of $C_0^{-}$ is at most $2$. If there is more than one irreducible component, then $C$ has the maximal number of irreducible components, and hence also the maximal number of nodes, which is $3g-3$. 
\ebew

\begin{prop}\label{312}
The morphism $\thx : \hxbar \rightarrow \hilb{g,n,0}$
induces a morphism 
\[ \Theta : \quad  \hxbar /\Gamma(C_0^-;\CQ)   \: \rightarrow \:
\hilb{g,n,0} , \label{n23}\]
which is an embedding, 
and an
isomorphism of schemes
\[ \msbar_{g^{+},m}' / \Gamma(C_0^-;\CQ) \quad \cong \quad \dpobar .\]
\end{prop}

\begin{rk}\em
The surjectivity of the morphism $\msbar_{g^{+},m}' \rightarrow \dpobar$ implies in particular that for all stable curves $C\rightarrow \Spec(k)$ with $3g-3$ nodes, which are represented by a point on $\dpobar$, there exists a node $P\in C$, such that for the induced decomposition $C = C^+(P) \cup C^-(P)$ as in \ref{N35}  holds $C^-(P) \cong C_0^-$, and thus $\overline{D}(C;P) = \dpobar$. 
\end{rk}

 \proofof{the proposition} 
It follows from proposition \ref{P318} that the morphism $\thx : \hilb{g^+,n,m}^\times \rightarrow \hilb{g,n,0}$ factors through $\hxbar /\Gamma(C_0^-;\CQ)$. 

Let $C_1$ and $C_2$ denote two embedded $m$-pointed stable curves over $\Spec(k)$, represented by points $[C_1], [C_2]\in \hilb{g^+,n,m}^\circ$. Suppose first that $\thx([C_1]) = \thx([C_2]) \in H_{g,n,0}$, and let $C$ denote the embedded curve representing this point. By the assumption on $[C_1]$ and $[C_2]$, the curve $C$ has $3g-4$ nodes, so there is a unique irreducible component on $C$, which has not the maximal number of nodes when considered as a pointed stable curve itself. This component is of course equal to the embedded curves $C_1$ and $C_2$ in $\PP^N$, considered as curves with $m$ distinguished points, but without ordered labels on these points. The closure of its complement in $C$ is isomorphic to $C_0^-$ by definition of the morphism $\thx$. Since $\thx([C_1])$ and $\thx([C_2])$ represent the same embedded curve, the order of the marked points of $C_1$ and $C_2$ can differ at most by a reordering induced by an automorphism of $C_0^-$. Hence by definition of $\Gamma(C_0^-;\CQ)$, the points $[C_1]$ and $[C_2]$ are the same modulo the action of $\Gamma(C_0^-;\CQ)$. 

Suppose now that $\thx([C_1]) = \thx([C_2])$ lies in the boundary of $ \hilb{g,n,0}$. Then there is no such distinguished irreducible component. In fact, if $C$ denotes the embedded curve representing the  point $\thx([C_1]) = \thx([C_2])$, then there may be even more than one subcurve of $C$ isomorphic to $C_1$ or $C_2$. However, by the construction of $\thx$, the embedding of $C_0^-$ into $\PP^N$ is fixed. The embedding of the curve $C_1$ into $\PP^N$ is just the closure of the complement of $C_0^-$ in $C$, and analogously for $C_2$. Thus as embedded curves, $C_1$ and $C_2$ are equal up to a permutation of the order of their marked points, which is induced by  an automorphism of $C_0^-$. Therefore in this case also the points $[C_1]$ and $[C_2]$ are the same up to the action of $\Gamma(C_0^-;\CQ)$. 
This proves the first part of the proposition.

For the second part note that by corollary \ref{C318} the morphism $\thx$ induces a morphism
\[ \thxbarGamma: \: \:  \mxbar/ \Gamma(C_0^-;\CQ) \rightarrow \msbar_g .\]
This morphism extends to a morphism from $\msbar_{g^+,m}' / \Gamma(C_0^-;\CQ)$ to $ \msbar_g$, as the quotient is smooth.  
By the definition of $\msbar_{g^+,m}'$, the image of this morphism is contained in $\dppbar$. Since $\dppbar$ is irreducible, this non-constant morphism is necessarily surjective.

The restriction of this morphism to the open part $\ms_{g^+,m}'/\Gamma(C_0^-;\CQ)$ is injective for analogous reasons as in the first part of the proof. For each point in the image of $\thxbarGamma$,  the preimage points represent the unique irreducible component which has not the maximal number of nodes possible, up to a permutation of the labels of the marked points.  Two such curves, representing two points in the preimage of a point,  can only differ by a reordering  of the distinguished points, which is induced by an automorphism of $C_0^-$. Hence they are isomorphic as $\mgamma$-pointed curves, with $\Gamma = \Gamma(C_0^-;\CQ)$.   

However, if a curve $C \rightarrow \Spec(k)$ in the image of the morphism $\thxbarGamma$ 
lies in
the boundary of $\dppbar$, then there may be more than one
subcurve isomorphic to $C_0^{-}$. Hence to prove that the induced morphism from $\msbar_{g^{+},m}' / \Gamma(C_0^-;\CQ)$ to $\dppbar$ is injective,  one needs to
show that for any two embeddings of $C_0^{-}$ into $C$ the closures of
the complements in $C$ are isomorphic as $m/\Gamma(C_0^-;\CQ)$-pointed
curves. More precisely, suppose that there are two decompositions of $C$ as 
\[ C \: = \: C_1^- + C_1^+ \: = \: C_2^- + C_2^+ ,\]
where $C_1^-$ and $C_2^-$ are isomorphic to $C_0^-$ as $m$-pointed curves, and $C_1^+$ and $C_2^+$ are represented by points  $[C_1^+], [C_2^+] \in \msbar_{g^+,m}'$. Here ``$+$'' denotes glueing of $m$-pointed curves in their marked points, according to their labels. The morphism is injective if we can show the following claim: the curves $C_1^+$ and $C_2^+$ are isomorphic as $m/\Gamma(C_0^-;\CQ)$-pointed curves. 

If $m=1$ then this is trivial. Let $m=2$. By counting the number of ``pigtails'' of $C$, i.e. the number of irreducible components of $C$, which are nodal rational curves  with one marked point, we find that $C_1^+$ and $C_2^+$ must be of the same topological type, compare remark \ref{37}. Hence  $[C_1^+]=[C_2^+]$. For $m=4$ the claim follows from the following technical lemma, if one puts $\mu=0$ and $k=4$.
\ebew

\begin{lemma}\label{L326a}\label{L326A}
Let $C \rightarrow \Spec(k)$ be a $\mu$-pointed prestable curve of some genus with stable connected components. Suppose that $C$ has the maximal number of nodes possible. Let $C_1^+$ and $C_2^+$ be subcurves of $C$, represented by points in $\msbar_{0,4}$, such that the four marked points are nodes or marked points of $C$. Denote the closures of their complements in $C$ by $C_1^-$ and $C_2^-$, respectively. Let $k_1$ and $k_2$ denote the number of points of intersection between $C_1^+$ and $C_1^-$, and of $C_2^+$ and $C_2^-$, respectively. Suppose that $k_1=k_2 =: k$. Both $C_1^-$ and $C_2^-$ shall be considered as $\mu+2k-4$-pointed curves, and suppose that there is an isomorphism
\[ \phi : \: \: C_2^- \rightarrow C_1^- \]
of $\mu+2k-4$-pointed curves. Suppose that $\phi$ maps glueing points of $C_2^-$ and $C_2^+$ to glueing points of $C_1^-$ and $C_1^+$. \\
For $i=1,2$ let $P_1^{(i)},\ldots,P_{\mu+k-4}^{(i)},Q_1^{(i)},\ldots, Q_k^{(i)}$ denote the the marked points of $C_i^-$, where $Q_1^{(i)},\ldots, Q_k^{(i)}$ denote the points of intersection between $C_i^-$ and $C_i^+$. Let $\Aut(C_i^-;P_1^{(i)},\ldots,P_{\mu+k-4}^{(i)},\{Q_1^{(i)},\ldots, Q_k^{(i)}\} )$ denote the group of those automorphisms of $C_i^-$ as a scheme over $\Spec(k)$, which stabilize the points $P_1^{(i)},\ldots,P_{\mu+k-4}^{(i)}$ and the set $\{Q_1^{(i)},\ldots, Q_k^{(i)}\}$. 
There is a natural  homomorphism 
\[ \Aut(C_i^-;P_1^{(i)},\ldots,P_{\mu+k-4}^{(i)},\{Q_1^{(i)},\ldots, Q_k^{(i)}\} ) \: \rightarrow \: \Sigma_k \]
into the symmetric group on $k$ elements. The image of this homomorphism shall be denoted by $\Gamma(C_i^-;P_1^{(i)},\ldots,P_{\mu+k-4}^{(i)},\{Q_1^{(i)},\ldots, Q_k^{(i)}\})$. Let $Q_1^{(i)},\ldots, Q_k^{(i)}, R_1^{(i)},\ldots,R_{4-k}^{(i)} $ be the four  marked points on the subcurve  $C_i^+$, where $Q_1^{(i)},\ldots, Q_k^{(i)}$ denote those  points where  $C_i^+$ is glued to $C_i^-$ in the points of the same name. 

Then there exists a permutation
\[ \sigma \in \Gamma(C_2^-;P_1^{(2)},\ldots,P_{\mu+k-4}^{(2)},\{Q_1^{(2)},\ldots, Q_k^{(2)}\})\]
 and an isomorphism $\gamma: C_1^+ \rightarrow C_2^+ $  of curves over $\Spec(k)$,  such that 
\[ \gamma(Q_j^{(1)}) = Q_{\sigma(j)}^{(2)} \]
for $j= 1,\ldots, k$, and 
\[ \gamma(R_j^{(1)}) = R_{j}^{(2)} \]
for $j= 1,\ldots, 4-k$.
\end {lemma}

\proof
Because of its length, and its  unenlightening combinatorial structure, the proof  is relocated from the main text to appendix \ref{AppA}.
\ebew 
 
\begin{rk}\em
Note that the content of the above lemma is symmetric with respect to $C_1^+$ and $C_2^+$. Using the isomorphism $\phi : C_2^- \rightarrow C_1^-$ there is an induced action of $\Gamma(C_2^-;P_1^{(2)},\ldots,P_{\mu+k-4}^{(2)},\{Q_1^{(2)},\ldots, Q_k^{(2)}\})$ on the labels of the marked points $Q_1^{(1)},\ldots , Q_k^{(1)}$ of $C_1^+$. Analogously, in the statement the group $\Gamma(C_2^-;P_1^{(2)},\ldots,P_{\mu+k-4}^{(2)},\{Q_1^{(2)},\ldots, Q_k^{(2)}\})$ could be replaced by the group $ \Gamma(C_1^-;P_1^{(1)},\ldots,P_{\mu+k-4}^{(1)},\{Q_1^{(1)},\ldots, Q_k^{(1)}\})$.
\end{rk}

The proof of proposition \ref{312} implies the following observation. Compare also remark \ref{R32}.

\begin{lemma}\label{L326}
Let $P_2$ be a second node on $C_0$. Then $\dppbar =\overline{D}(P_2)$ if and only if there exists an automorphism $\rho\in \Aut(C_0)$ such that $\rho(P_2) = P_1$. In particular the curve $\dppbar$ is smooth. 
\end{lemma}

\proof
If there exists such an automorphism $\rho\in \Aut(C_0)$ with $\rho(P_2) = P_1$ then clearly $\dppbar =\overline{D}(P_2)$, compare the proof of lemma \ref{L31}. Conversely, suppose that $\dppbar =\overline{D}(P_2)$ holds. Let $C_0^+(P_2)$ denote the union of those irreducible comonents of $C_0$ which contain $P_2$, and let $C_0^-(P_2)$ denote the closure of its complement in $C_0$. Both shall be considered as pointed curves. Now consider two deformations of $C_0$, which preserve all nodes except $P_1$, and all nodes except $P_2$, respectively. In these deformations the subcurves $C_0^-= C_0^-(P_1)$ and $C_0^-(P_2)$, respectively, are preserved. Because $\dppbar =\overline{D}(P_2)$ there must be an isomorphism $\varphi: C_0^-(P_2) \rightarrow C_0^-(P_1)$. By proposition \ref{312} there exists an automorphism $\rho \in \Aut(C_0^-(P_1);\CQ)$, inducing a permutation $\gamma\in \Gamma(C_0^-(P_1);\CQ)$, which reorders the marked points on $C_0^+(P_1)$ in  such a way that the reordered pointed curve is isomorphic to $C_0^+(P_2)$. Then the isomorphism $\rho\circ \varphi: C_0^-(P_2) \rightarrow C_0^-(P_1)$ extends to an isomorphism between $C_0^-(P_2) + C_0^+(P_2)$ and  $C_0^-(P_1) + C_0^+(P_1)$, i.e. to an automorphism of $C_0$, mapping the node $P_2$ to $P_1$.
\ebew

\begin{cor}\label{C333}
Let $C_0=C_1^+ +C_1^-= C_2^+ +C_2^-$ be two decompositions of $C_0$ into $4$-pointed prestable curves. Suppose that $C_1^- \cong C_2^-$ as $4$-pointed prestable curves, and that  both  $C_1^+$ and $C_2^+$ consist of two lines which meet in a node $P_1\in C_1^+$ and $P_2\in C_2^+$, respectively. Then there exists an automorphism $\gamma\in\Aut(C_0)$ such that $\gamma(P_1) = P_2$.
\end{cor}

\proof
By assumption $C_1^- \cong C_2^-$, so $\dppbar = \overline{D}(P_2)$. Therefore the claim follows immediately from lemma \ref{L326}.
\ebew

\begin{rk}\em
$(i)$ If $m=1$ then $\Gamma(C_0^-;\CQ)$ is necessarily trivial. If $m=2$, then $\Gamma(C_0^-;\CQ)\subset \ZZ_2$, and it leaves the point $[C_0^{+}]\in \msbar_{g^{+},m}'$ invariant. In case  $m=4$, one has $\Gamma(C_0^-;\CQ) \subset \Sigma_4$, and $\Gamma(C_0^-;\CQ)$ may or may not stabilize $[C_0^{+}]$, depending on the structure of $C_0^{-}$.  

$(ii)$ The identity $\sigma([C_0^{+}]) = [C_0^{+}]$ holds in $\hgprimebar$ for some $\sigma \in \Gamma(C_0^-;\CQ)$ if and only if $\sigma \in \Gamma^{\circ}(C_0^-;\CQ)$.

$(iii)$ If the groups $\Aut(C_0^-;\CQ)$  and $\Aut(C_0^+;\CQ)$ are considered as as subgroups of $\PGL(N+1)$, and if 
\[ \rho^- : \:\: \Aut(C_0^-;\CQ) \rightarrow \Sigma_m \text{ \:  and \: } \rho^+ : \:\: \Aut(C_0^+;\CQ) \rightarrow \Sigma_m\]
denote the natural homomorphisms of groups, such that  $\rho^-(\Aut(C_0^-;\CQ)) = \Gamma(C_0^-;\CQ)$, then
\[ \rho^-( \Aut(C_0^-;\CQ)) \cap \rho^+( \Aut(C_0^-;\CQ))  = \Gamma^\circ(C_0^-;\CQ) .\]
\end{rk}

\begin{notation}\em
Let $\pi : \hilb{g,n,0} \rightarrow \msbar_g$ be the canonical
quotient morphism. We denote by 
\[ \overline{K}(C_0;P_1) \label{n24} \]
the reduction of the preimage of $\dpobar$ in $\hilb{g,n,0}$. Furthermore, we denote  by $\kxbar$\label{m6} the reduction of the preimage of $\dxbar$. Their intersections with the open subscheme $H_{g,n,0}$ are denote by ${K}(C_0;P_1)$ and $\kx$, respectively. 
\end{notation}

\begin{cor}\label{314}
As a subset, $\kxbar$ is equal to the orbit 
\[ \kxbar = \thx(\hxbar) \cdot \PGL(N+1) \]
of the image $ \thx(\hxbar)$ under the action of $\, \PGL(N+1)$. 
Furthermore,  it holds that 
\[ \kx = \thx(\hx) \cdot \PGL(N+1) .\]
\end{cor}

\proof
This follows immediately from the surjectivity of the morphism $\mxbar \rightarrow \dxbar$. 
\ebew

\begin{rk}\em\label{R336}
As before, we will usually write for example $\kpxbar$ instead of $\kxbar$, and $\kpx$ instead of $\kx$,  when the curve $C_0$ is understood. Note that in general the preimage of $\dxbar$ in $ \hilb{g,n,0}$ will not be reduced. 

The scheme $\kpx$  is a smooth and irreducible subscheme  of $\overline{H}_{g,n,0}$, which is a smooth scheme itself. This follows from the definition of $\dxbar$ as an irreducible component of $\msbar_g^{(3g-4)}$, and the fact that the boundary $\hilb{g,n,0} \setminus {H}_{g,n,0}$ is a normal crossing divisor, compare \cite[cor. 1.7, cor. 1.9]{DM}. 
\end{rk}

\begin{lemma}\label{315}
Let $f : C \rightarrow S$ be a stable curve of genus $g$, with reduced base $S$, such that
for all $s\in S$ the fibre $C_s$ has at least $3g-4$ nodes. Then the
induced morphism $\theta_f : S \rightarrow \msbar_g$ factors through
$\dppbar$ if and only if 
there exists an \'etale  covering $\{S_\alpha\}_{\alpha\in A}$ of $S$, and for each $\alpha \in A$  a commutative diagram
\[\diagram
C_0^{-} \times S_\alpha \drto_{\pr_2}\rrto|<\hole|<<\ahook && C_\alpha \rrto\dlto&& C \dto^f\\
&S_\alpha \xto[rrr]&&& S 
\enddiagram\]
with 
\[ \rho_i(s) = (Q_i,s) \]
for $i=1,\ldots,m$ and for 
all $s \in S_\alpha$, where $\rho_1,\ldots ,\rho_{m} : S_\alpha \rightarrow
C_\alpha:= C\times_S S_\alpha$ denote  sections of nodes as in lemma \ref{sectnod}, and $Q_1,\ldots,Q_m$ are the marked points on $C_0^-$. 
\end{lemma}

\proof
Before we start with the proof 
let us consider a stable curve $C \rightarrow \Spec(k)$ of genus $g$, which contains a subcurve isomorphic to $C_0^{-}$ in such a way that the marked points of $C_0^{-}$ are nodes of $C$. Let $\tilde{C} \rightarrow B$ be a deformation of $C$, which preserves the nodes of $C$. Let $b_0\in B$ such that the fibre $C_{b_0}$ is isomorphic to $C$. Because the number of nodes does not decrease within the family, and $C_0^{-}$ has the maximal number of nodes possible for an $m$-pointed prestable curve of genus $g^{-}$, each fibre of $\tilde{C} \rightarrow B$ contains a subcurve isomorphic to $C_0^{-}$. In an \'etale  neighbourhood $B'$ of $b_0$ there exist sections $\rho_1,\ldots,\rho_m$ of nodes, such that for each $b\in B'$ and each $i=1,\ldots,m$ the point $\rho_i(b)$ is a marked point of $C_0^{-}$. Cutting along these sections we obtain a subcurve $C'\rightarrow B'$ of $C\rightarrow B'$, such that each fibre of $C'\rightarrow B'$ is isomorphic to $C_0^{-}$. 
The locus of curves with the maximal number of nodes is discrete in
the corresponding moduli space. Therefore, using another \'etale covering if necessary,  the family $C'\rightarrow
B'$ must be locally  isomorphic to the trivial family $C_0^{-}\times B' \rightarrow B'$. 

Now consider a stable curve $f : C \rightarrow S$ of genus $g$, with reduced base $S$, such that $\theta_f : S \rightarrow \msbar_g$ factors through $\dppbar$. Lemma \ref{33} implies that each fibre contains a subcurve isomorphic to  $C_0^{-}$, and locally $f : C \rightarrow S$ is  a deformation of each of its fibres preserving the nodes given by the marked points on  $C_0^{-}$. Hence the first implication of the claim follows from the above discussion.

Conversely, suppose that such commutative diagrams exist locally. In particular, this implies that for any closed point $s\in S$, its fibre $C_s$ contains a subcurve isomorphic to $C_0^-$, in such a way that the marked points of $C_0^-$ are nodes on $C_s$. By lemma \ref{314a}, the closure $C_s^{+}$ of the complement of $C_0^{-}$ in the  fibre $C_s$  is either irreducible,  or $C_s$ has the maximal number of nodes. 

If $C_s^{+}$ is irreducible, then it follows from our discussion in remark \ref{37} that its topological type is uniquely determined by the number $m$ of its marked points, up to a reordering of the marked points via the action of the group $\Gamma(C_0^-;\CQ)$. Therefore the point $\theta_f(s)= [C_s]$ representing the fibre lies in the same irreducible component of $\msbar_g^{(3g-4)}$ as the fibre of a deformation of $C_0$ which fixes the subcurve $C_0^{-}$, and this component is by definition $\dppbar$. 

If $C_s^{+}$ is not irreducible, then $C_s$  has $3g-3$ nodes. Let $d : \tilde{C}^{+} \rightarrow B$ be a deformation of $C_s^{+}$ such that all fibres are irreducible $m$-pointed stable curves of genus $g$, except the one fibre $\tilde{C}^{+}_{{b_0}}$, for some $b_0 \in B$, which is  isomorphic to  
$C_s^{+}$. To see that this is possible consult remark \ref{37}, where all types of curves were listed, which can occur as subcurves $C_s^{+}$.   Restricting $B$ we can achieve furthermore that after glueing with the trivial family $C_0^{-}\times B \rightarrow B$ along the $m$ marked points we obtain a stable curve of genus $g$, where each fibre has at least $3g-4$ nodes. By the previous paragraph the induced morphism 
\[ \theta_d : \:\: B \rightarrow \msbar_g, \]
restricted to $B \setminus \{b_0\}$, factors through $\dppbar$, and hence $\theta_f(s) = \theta_d(b_0) \in \dppbar$ as well. 
\ebew

\begin{cor}\label{K324}
Let $f: C \rightarrow S$ be a stable curve of genus $g$, with reduced base $S$, such that for
all closed points $s\in S$ the fibre $C_s$ has exactly $3g-4$ nodes. Suppose that
the induced morphism $\theta_f : S \rightarrow \msbar_g$ factors through
$\dpp$. Then there exists a distinguished subcurve
$f^{-} : C^{-} \rightarrow S$ of $f: C \rightarrow S$, which in the \'etale topology is
locally isomorphic to the trivial curve $C_0^{-} \times S \rightarrow
S$ as an $m/\Gamma(C_0^-;\CQ)$-pointed stable curve. 
\end{cor}

\proof
By the previous lemma \ref{315} such a curve exists locally in the \'etale topology. Since each
fibre $C_s$ of $f: C \rightarrow S$ has exactly $3g-4$ nodes, there
is a unique irreducible component of $C_s$ which has not the maximal
number of nodes when considered as an $m$-pointed curve. The closure
of the complement of this component in $C_s$ is the unique subcurve
isomorphic to $C_0^{-}$. Hence the images of the trivial curves $C_0^{-}\times
S_\alpha  \rightarrow S_\alpha$ in $f: C\rightarrow S$, as in lemma \ref{315}, glue together to up
to automorphisms of $C_0^{-}$, which means up to a permutation of
their marked points by the action of the group $\Gamma(C_0^-;\CQ)$. 

This defines a subcurve $f^-:C^-\rightarrow S$ of $f:C\rightarrow S$ over $S$. Note that lemma \ref{sectnod} implies that there exists an \'etale covering of $f^-:C^-\rightarrow S$ by an $m$-pointed stable curve, so $f^-:C^-\rightarrow S$ is indeed an $\mgamma$-pointed stable curve, and it is \'etale locally isomorphic to a trivial curve. 
\ebew

The above corollary \ref{K324} deals only with stable curves $f:C\rightarrow S$, where the base $S$ is reduced. In our applications, we will need an analogous statement for arbitrary schemes $S$, and for this we need some extra assumptions.

\begin{cor}\label{K324a}
Let $f: C \rightarrow S$ be a stable curve of genus $g$, such that for
all closed points $s\in S$ the fibre $C_s$ has exactly $3g-4$ nodes. Suppose that
the induced morphism $\theta_f : S \rightarrow \msbar_g$ factors through
$\dpp$. Let $(E,p,\phi)$ be the triple associated to $f: C \rightarrow S$ as in proposition \ref{28}, consisting of a principal $\PGL(N+1)$-bundle $p:E\rightarrow S$, and a morphism $\phi: E\rightarrow  \overline{H}_{g,n,0}$. Suppose that $\phi$ factors through $K(P_1)$. 
Then there exists a distinguished subcurve
$f^{-} : C^{-} \rightarrow S$ of $f: C \rightarrow S$, which in the \'etale topology is
locally isomorphic to the trivial curve $C_0^{-} \times S \rightarrow
S$ as an $m/\Gamma(C_0^-;\CQ)$-pointed stable curve. 
\end{cor}

\proof
Let ${\cal C}_{g,n,0}\rightarrow \overline{H}_{g,n,0}$ denote the universal curve, and ${\cal C}_{\overline{K}(P_1)}\rightarrow \overline{K}(P_1)$ its pullback to $
\overline{K}(P_1)$. Since $\overline{K}(P_1)$ is reduced, corollary \ref{K324} 
shows the existence of a unique subcurve ${\cal C}_{\overline{K}(P_1)}^-\rightarrow \overline{K}(P_1)$, which is locally isomorphic to the trivial curve $C_0^-\times \overline{K}(P_1)\rightarrow \overline{K}(P_1)$ as an $m/\Gamma(C_0^-;{\cal Q})$-pointed curve. Let ${\cal C}_E$ and ${\cal C}_E^-$ denote the  respective pullbacks of ${\cal C}_{\overline{K}(P_1)}$ and ${\cal C}_{\overline{K}(P_1)}^-$ to $E$. By the definition of $E$, one has $f: C={\cal C}_E/\PGL(N+1) \rightarrow E/\PGL(N+1)=S$, and $C^- := {\cal C}_E^-/\PGL(N+1) \rightarrow S$  is the unique subcurve, which is locally isomorphic to $C_0^-\times S\rightarrow S$ as an $m/\Gamma(C_0^-;{\cal Q})$-pointed curve.
\ebew

\begin{rk}\em
In general  this subcurve $f^{-} : C^{-} \rightarrow S$ is not
an $m$-pointed curve, as the markings of the
fibres might not be given by global sections. It is however an $\mgamma$-pointed stable curve of genus $g^+$, with $\Gamma= \Gamma(C_0^-;\calq)$. 
\end{rk}

\begin{rk}\em\label{323a}
If $C_1^{-},\ldots,C_r^{-}$ are the connected components of $C_0^{-}$, then for the groups of automorphisms of pointed curves holds
\[ \Aut(C_0^{-})  = \Aut(C_1^{-})\times \ldots \times \Aut(C_r^{-}) .\]
Note that  $\Aut(C_0^{-})$ is in a natural way a subgroup of $\Aut(C_0)$ which acts trivially on $C_0^{+}$. Recall that we identified $\Aut(C_0)$ with the subgroup in $\PGL(N+1)$ stabilizing the fixed embedding of  $C_0$ into $\PP^N$, so we can also view $\Aut(C_0^{-})$ as a subgroup of $\PGL(N+1)$.

Note that this embedding is compatible with the embedding of the group $\Aut(C_0^-;\CQ)$, which also contains $\Aut(C_0^-)$ as a subgroup, into $\PGL(N+1)$, as in remark \ref{R311b}. 
\end{rk}

\begin{rk}\em\label{323b}
The group $\Aut(C_0^-;\CQ)$  of automorphisms of $C_0^{-}$, stabilizing the subset $\calq$ as a set, as introduced in definition  \ref{D318a}, has an induced
action on the projective bundle $ \PP \left( \bigoplus_{i=1}^r 
(\pr_2)_\ast(\omega_{C_i^{-}\times S / S}(S_1+\ldots+S_m))^{\otimes
n}\right)$, where $S_i := \{Q_i\} \times S$ for $i=1,\ldots,m$ are the
sections of marked points of the trivial $m$-pointed curve $\pr_2: C_0^{-}\times S \rightarrow S$. 

Let $f : C\rightarrow S$ be a stable curve of genus $g$, such that the
induced morphism $f: S \rightarrow \msbar_g$ factors through
$\dpp$.  
Doing a similar construction as in \ref{const_part_2}, we obtain for the associated projective bundles  an inclusion
\[ \PP \left( \bigoplus_{i=1}^r 
(\pr_2)_\ast(\omega_{C_i^{-}\times S / S}(S_1+\ldots+S_m))^{\otimes n} \right) \, / \,    \Aut(C_0^-;\CQ) \hookrightarrow  \]
\[ \hspace*{2cm} \hookrightarrow
\PP
f_\ast(\omega_{C/ S})^{\otimes n} \, / \, \Aut(C_0^-) . \]
Note that by lemma \ref{315} there is at first only an inclusion of projective
bundles on an \'etale cover of $S$. Only after dividing out by the action of the group of automorphisms of $C_0^{-}$, which stabilize the set $\calq$,  these inclusions glue globally. 
Recall that the group $\Aut(C_0^{-})$ is a normal subgroup of $\Aut(C_0^-;\CQ)$,
with quotient group isomorphic to $\Gamma(C_0^-;\CQ)$, see remark \ref{R319}. Locally the inclusion of  $\PP ( \bigoplus_{i=1}^r 
(\pr_2)_\ast(\omega_{C_i^{-}\times S / S}(S_1+\ldots+S_m))^{\otimes n} )$ into $\PP f_\ast(\omega_{C/S})^{\otimes n}$ is invariant with respect to reorderings of the sections of the marked points by elements of $\Gamma(C_0^-;\CQ)$.  
\end{rk}

After all of these  preliminaries, and having introduced most of our notation, let us now define the moduli stacks we are interested in. 

\begin{defi}\em
$(i)$ The fibred category  $\FDPtilde$\label{m7} is defined as the following full subcategory of $\cmgbar$. For any scheme $S$, the set of objects of the fibre category $\FDPtilde(S)$ is equal to the set of those stable curves $f:C\rightarrow S$ of genus $g$, such that the induced morphism $\theta_f: S\rightarrow \msbar_g$ factors through the subscheme $\dpobar$. 

$(ii)$ The stack $\FDP$ is defined as the reduced stack underlying  $\FDPtilde$.\label{n25}
\end{defi}

\begin{rk}\em\label{R323}
The fibred category $\FDP$ can be considered as a  $2$-functor. It is  a substack of $\mgstack{g}$, and a Deligne-Mumford stack itself, as we will see below. It follows from lemma \ref{315} that  $\FDP$ is a moduli stack for curves $f:C\rightarrow \Spec(k)$, which contain a subcurve isomorphic to $C_0^-$, where the marked points on $C_0^-$ correspond to nodes on $C$, with $\dpobar$ as its moduli space. 
\end{rk}

\begin{prop}\label{pr341}
The stack $\FDPtilde$ is the preimage substack of the scheme $\dpobar$ in $\mgstack{g}$ under the canonical map from $\mgstack{g}$ to $\msbar_g$. In other words, the following diagram is Cartesian.
\[ \diagram
\FDPtilde \dto \rto & \mgstack{g}\dto \\
\dpobar \rto & \msbar_g .
\enddiagram \]
\end{prop}

\proof
Note that we only have to check that $\FDPtilde$ satisfies the universal property of fibre products for all schemes $S$.  Consider a commutative diagram
\[ \diagram
S^\bullet \xto[rrd]^\beta \xto[rdd]_\alpha \xdotted[dr]|>\tip^{?}\\
& \FDPtilde \rto \dto & \mgstack{g} \dto\\
& \dpobar \rto& \msbar_g.
\enddiagram \]
The morphism $\alpha$ is represented by a morphism $a : S \rightarrow \dpobar$ of schemes, and by Yoneda's lemma $\beta$ defines an element $b \in \mgstack{g}(S)$, hence a stable curve $b: B \rightarrow S$ of genus $g$. Therefore there is a morphism $S \rightarrow \msbar_g$ into the moduli space, which factors through $\dpobar$ via $\alpha$. This implies that $b: B \rightarrow S $ is an element of $\FDPtilde(S)$, which again by Yoneda's lemma defines a morphism $S^\bullet \rightarrow \FDPtilde$, and this morphism  is necessarily unique. 

Since $\FDPtilde$ is by definition a full subcategory of $\mgstack{g}$, condition $(ii)$ of proposition \ref{A0118} is automatically satisfied. 
\ebew

\begin{notation}\em
We denote by $\FDPopen$\label{n26} the reduction of the preimage substack of $\dpo$ in $\mgstack{g}$. Analogously, we define $\FDPc$\label{m8} and $\FDPcopen$ as the reductions of the preimages of $\dxbar$ and $\dx$.
\end{notation}

\begin{lemma}\label{140}
There are  isomorphisms of stacks
\[ \FDP \: \cong \: \left[ \, \overline{K}(C_0;P_1) / \PGL(N+1) \right],\]
and
\[ \FDPc \: \cong \: \left[ \, \kxbar / \PGL(N+1) \right],\]
and analogous quotient representations of the open substacks $\FDPopen$ and $\FDPcopen$. 
\end{lemma}

\proof
Using proposition \ref{pr341}, 
this follows immediately from the definition of $\overline{K}(C_0;P_1)$ and $\kxbar$ by standard arguments on Cartesian products and quotient stacks. 
\ebew

\begin{cor}
The stack $\FDPopen$ is smooth and irreducible, and of dimension one.
\end{cor}

\proof
The scheme ${K}(C_0;P_1)$ is an atlas of $\FDPopen $, considered as Artin stack. Hence the claim follows immediately from the fact that ${K}(C_0;P_1)$ is smooth and irreducible, see remark \ref{R336}.
\ebew

\begin{rk}\em
From lemma \ref{140} the following characterization of $\cdpobar$ can be derived. For any scheme $S$, the fibre category $\cdpobar(S)$ has as objects such stable curves $f:C\rightarrow S$ of genus $g$, where  the induced morphism $\theta_f : S\rightarrow \mgbar$ factors through $\dpobar$, and the morphism $\phi: E \rightarrow \hilb{g,n,0}$ factors through $\overline{K}(C_0;P_1)$. Here, $(E,p,\phi)$ is the triple associated to the curve $f:C\rightarrow S$ as in proposition \ref{28}.
\end{rk}

\chapter{The moduli substack}\label{par4}

This is the central part of our work. Recall that there is a natural morphism $\overline{\Theta}: \msbar_{g^+,m} \rightarrow \mg$ between moduli spaces of stable curves, as in \ref{D322}. By proposition \ref{312}, this morphism induces an isomorphism between $\dpobar$ and the subscheme $\mgprimebarquot$ of the moduli space of $\mgamma$-pointed stable curves of genus $g$, for $\Gamma = \Gamma(C_0^-;\calq)$. We now want to consider the analogous situation for the corresponding moduli stacks.

Let $\cmgprimebarquot$ denote the substack of $\cmgbarquot$, which is the reduction of the preimage of $\mgprimebarquot$. It turns out that the natural morphism of stacks $
\cmsbar_{g^+,m} \rightarrow \cmgbar$ does not induce an isomorphism between the stacks $\cmgprimebarquot$ and $\cdpobar$, even though their moduli spaces are isomorphic. 

Instead, it turns out that in order to relate the two stacks, we need to consider a finite quotient of the stack $\cmgprimebarquot$, see remark \ref{Rk46}. There is no global morphism from this quotient to $\cdpobar$, so we must restrict ourselves to an open and dense substack $\cmgammaxbar$, as defined in \ref{not004}, which maps  to an open and dense substack $\cdxbar$ of $\cdpobar$, compare proposition \ref{MainProp}. 

This morphism is even  an isomorphism between the two stacks.  This is the content of our main theorem \ref{MainThm}, which says that the stack $\cdpobar$ is generically, i.e. on an open and dense substack,  isomorphic to a quotient stack, which is at the same time a finite quotient of a natural substack of the moduli stack of $\mgamma$-pointed stable curves of genus $g^+$.

\begin{notation}\label{not004}\em
We denote by  $\cmgprimebarquot$\label{k1} the substack of $\cmgbarquot$, which is the reduction of the preimage of $\mgprimebarquot$ under the canonical morphism. The open substack corresponding to the open subscheme $\mgprimequot$ is denoted by $\cmgprimequot$.

Analogously, we  denote by  $\cmgammaxbar$\label{k2} the substack of $\cmgbarquot$, which is the reduction of the preimage of $\mgammaxbar$, and the  open substack corresponding to the open subscheme $\mgammax$ by $\cmgammax$.
\end{notation}

\begin{abschnitt}\em
Let us briefly recall the setup from above. Throughout,  $f_0 : C_0  \rightarrow \Spec(k)$ is a fixed stable curve of genus $g$ with $3g-3$ labeled nodes. The node $P_1$ distinguishes a decomposition of $C_0$ into two $m$-pointed prestable curves $C_0^+$ and $C_0^-$ of genus $g^+$ and $g^-$, respectively, which meet in their marked points. Recall that $C_0^+$ was defined as the union of all those irreducible components of $C_0$ which contain the node $P_1$, and $C_0^-$ as the closure of its complement. Deformations of $C_0$, preserving the $3g-4$ nodes different from $P_1$,  correspond to certain deformations of the curve $C_0^+$ as an $m$-pointed stable curve. Hence the irreducible component $\overline{D}(C_0;P_1)$ of $\msbar_g^{(3g-4)}$, which is distinguished by such  deformations, corresponds to a subscheme $\msbar_{g^+,m}'$ of the moduli space $\msbar_{g^+,m}$. 

The subscheme $C_0^-$ remains invariant throughout all deformations. By glueing deformations of $C_0^+$ to the fixed curve $C_0^-$ the ordering of the marked points disappears. To describe which differently marked $m$-pointed curves of genus $g^+$ result in the same curve of genus $g$ after glueing, one needs to look at automorphisms of $C_0^-$, considered as a prestable curve with $m$ distinguished points, but without ordering of the points. To do this, we introduced the notion of $\mgamma$-pointed stable curves. We obtained an isomorphism of $\overline{D}(C_0;P_1)$ with $\msbar_{g^+,m}'/\Gamma(C_0^-;\CQ)$. 

The aim of this section is to understand this isomorphism on the open part $\dpo$ on the level of the corresponding moduli stacks $\cdpo$, or at least on an open and dense substack of it. 
\end{abschnitt}
 
\begin{rk}\em\label{21}
Recall that we can view $\PGL(N^{+}+1)$ as a subgroup of $\PGL(N+1)$, such that the morphism  
$\thx : \hxbar \rightarrow \hilb{g,n,0}$ of proposition \ref{const_prop} is $\PGL(N^{+}+1)$-equivariant. By trivial extension, the group $\Aut(C_0^{-})$ of automorphisms of the $m$-pointed stable curve $C_0^{-}$ can be considered as   a subgroup of $\Aut(C_0)$. Hence it is via 
\[ \Aut(C_0) \cong \Stab_{\PGL(N+1)}([C_0]) \]
also a subgroup of $\PGL(N+1)$. Under this identification, the action of $\Aut(C_0^{-})$ restricted to  $C_0^{+}$ is trivial. 
So there is an injective group homomorphism
\[ \PGL(N^{+}+1) \times \Aut(C_0^{-}) \quad \longrightarrow \quad \PGL(N+1) .\]
Note that by remark \ref{R311} the embedding of $\PGL(N^+ +1)$ into $\PGL(N+1)$ is induced by the distinguished embeding of $\PP^{N^+}$ into $\PP^N$ as in remark \ref{const_part_1}. Analogously, the group $\Aut(C_0^-)$ is a subgroup of $\PGL(N^- +1)$, which lies also inside $\PGL(N+1)$, and the inclusion is determined by the embedding of a certain $\PP^{N^-}$ into $\PP^N$. By definition, elements of $\Aut(C_0^-)$ fix the $m$ marked points of $C_0^-$, which span the intersection $\PP^{N^-} \cap \PP^{N^+}$ in $\PP^N$ by remark \ref{const_part_3}. Thus such an element acts as identity on the intersection, and therefore elements of $\Aut(C_0^-)$ and elements of $\PGL(N^+ +1)$ commute inside $\PGL(N+1)$. 
\end{rk}

\begin{lemma}\label{L42}
The morphism $\thx: \hxbar \rightarrow \hilb{g,n,0}$ is
equivariant with respect to the action of $\, \PGL(N^{+}+1) \times \Aut(C_0^{-})$, where $\Aut(C_0^{-})$ acts trivially on $\hxbar$.
\end{lemma}

\proof
The equivariance of $\thx$ with respect to the action of $\PGL(N^{+}+1)$ was alredy stated in proposition \ref{const_prop}. It remains to show that the image of $\hxbar$ under $\thx$ is pointwise fixed under the action  of $\Aut(C_0^{-})$. Recall that  $\thx$ is constructed in \ref{const_part_4} via  compatible trivializations as in the following diagram:
\[ \diagram
\PP u^{+}_\ast(\omega_{{{\cal
C}}_{g^{+},n,m}/\hilb{g^{+},n,m}}(\sum_{i=1}^m{\cal S}_i ))^{\otimes n}  \rto|<\hole|<<\ahook \dto^{\cong}& \PP(u_0)_\ast(\omega_{{\cal
C}_0/\hilb{g^{+},n,m}})^{\otimes n} \dto_{\cong}\\
\PP(f_0^{+})_\ast (\omega_{C_0^{+}/k}(\sum_{i=1}^m Q_i))^{\otimes
n}  \times \hilb{g^{+},n,m} \rto|<\hole|<<\ahook &
\PP(f_0)_\ast(\omega_{C_0/k})^{\otimes n} \times
\hilb{g^{+},n,m},
\enddiagram \]
restricted to the subscheme $\hxbar$. 
By definition of the inclusion of $\Aut(C_0)$ into $\PGL(N+1)$,  the embedding of $C_0$ into $\PP^N$ is fixed under the action of $\Aut(C_0)$, and hence in particular fixed under the action of $\Aut(C_0^{-})$. Therefore the embedding of the projective bundles is invariant under the action of $\Aut(C_0^{-})$, and thus $\thx$ is also invariant with respect to the action of $\Aut(C_0^{-})$. 
\ebew

\begin{rk}\em
The action of $\PGL(N^+ +1)$  on $\hilb{g^+,n,m}$ induces an action of   $\PGL(N^+ +1)$  on $\hilb{g^+,n,m/\Gamma(C_0^-;\CQ)}=\hilb{g^+,n,m} / \Gamma(C_0^-;\CQ) $. Indeed, the action of  $\PGL(N^+ +1)$ on
$\PP u_\ast^+(\omega_{{\cal C}_{g^+,n,m}/\hilb{g^+,n,m}}({\cal S}_1+\ldots +{\cal S}_m))^{\otimes n} $ does not depend on the ordering of the marked points. Note that the open subscheme $\hxbar$ is fixed under the actions of both $\PGL(N^+ +1)$ and $\Gamma(C_0^-;\CQ)$. Hence lemma \ref{L42} implies that the induced morphism
\[ \Theta: \: \: \hgammadetailxbar  \rightarrow \hilb{g,n,0} \]
is equivariant with respect to the action of  $\PGL(N^{+}+1) \times \Aut(C_0^{-})$, where $\Aut(C_0^{-})$ acts trivially on $\hgammadetailxbar= \hxbar/\Gamma(C_0^-;\CQ) $.
\end{rk}

The following proposition provides a stack version of the natural morphism $\mgprimebarquot \rightarrow \dpobar$. More precisely, it generalizes the restricted morphism  $\mgammaxbar \rightarrow \dxbar$, because it turns out that to make things work, we have to restrict ourselves to the locus of curves  where the group of automorphisms splits as in  \ref{def_split}.

\begin{prop}\label{MainProp}
There is a morphism of Deligne-Mumford stacks
\[ \Lambda : \quad  \left[\hgammadetailxbar \, / \, \PGL(N^{+}+1) \times \Aut(C_0^{-}) \right] \: \: \longrightarrow  \: \:  \FDPc \]
where $\Aut(C_0^{-})$ acts trivially on $\hgammadetailxbar$. 
\end{prop}

\begin{notation}\em
As a shorthand notation we will from now on use $A:= \Aut(C_0^{-})$, $\Gamma := \Gamma(C_0^-;\CQ)$  and $P := \PGL(N^{+}+1)$. We will also abbreviate the subschemes $\kxbar$ and $\dxbar$ by $\kpxbar$ and $\dpxbar$, respectively. 
\end{notation}

\begin{rk}\em\label{Rk46}
Recall that by proposition \ref{P217} there is an isomorphism of stacks $\cmgbarquot \cong [ \hgbarquot/\PGL(N^+ +1)]$, and therefore an isomorphism $\cmgammaxbar  \cong [ \hgammaxbar/\PGL(N^+ +1)]$. By composing the morphism $\Lambda$ of proposition \ref{MainProp} with the canonical quotient morphism $\cmgammaxbar \rightarrow \cmgammaxbar / \Aut(C_0^-) = [ \hgammaxbar/\PGL(N^+ +1) \times\Aut(C_0^-)]$, we obtain a morphism of stacks 
\[ \cmgammaxbar \:\: \longrightarrow \:\: \cmgbar ,\]
which factors through $\cdxbar$, and which induces the natural morphism $\mxbar \rightarrow \dxbar$ on the corresponding  moduli spaces. 
\end{rk}

\proofof{proposition \ref{MainProp}}
 Viewing stacks as categories fibred in groupoids over the category of
 schemes, the morphism $\Lambda$  will be given as a functor respecting the fibrations.

$(i)$ Let a scheme $S$ be fixed. Consider a tuple $(E',p',\phi') \in $ \linebreak
$[ \, \hgammaxbar \, / \, P\times A ](S)$, i.e. a principal $P\times A$-bundle $p':E'\rightarrow S$, together with a $P\times\ A$-equivariant morphism $\phi':E' \rightarrow \hgammaxbar$. Via the injective group homomorphism $P\times A \rightarrow \PGL(N+1)$ from remark \ref{21}
we can define the extension $E$ of $E'$ by
\[ E := (E' \times \PGL(N+1) ) / (P\times A) \]
as a principal $\PGL(N+1)$-bundle $p:E\rightarrow S$ over $S$. By corollary \ref{314} we have $\kpxbar = \thx(\hxbar)\cdot\PGL(N+1) = \Theta(\hgammaxbar)\cdot\PGL(N+1)$. Therefore we can define a morphism
\[ \begin{array}{cccc}
\phi : & E& \rightarrow & \kpxbar\\
& e & \mapsto& \Theta(\phi'(e'))\cdot \gamma, 
\end{array}\]
where  $e\in E$ is represented by a pair $(e',\gamma)\in E'\times \PGL(N+1)$.
This map is indeed well-defined, as $\Theta \circ\phi'$ is equivariant with respect to $P\times  A$.  For all $g \in \PGL(N+1)$ we have 
\[ \phi(eg) = \phi(e) g\]
for all $e\in E$, so $\phi$ is $\PGL(N+1)$-equivariant. Hence 
\[ (E,p,\phi) \in \left[ \, \kpxbar / \PGL(N+1) \right](S) \cong \FDPc(S),\]
where the last isomorphism holds by lemma \ref{140}.

$(ii)$ Let $\varphi': (E_2',p_2',\phi_2') \rightarrow (E_1',p_1',\phi_1')$ be a morphism in the fibre category $
[ \, \hgammaxbar \, / \, P\times A ](S)$. In other words, $\varphi'$ is given by a bundle homomorphism $f': E_1' \rightarrow E_2'$, such that $p_1'=p_2'\circ f'$ and $\phi_1' = \phi_2'\circ f'$ hold. Then there exists a natural extension $f : E_1 \rightarrow E_2$ of $f'$, where $E_1$ and $E_2$ are constructed from $E_1'$ and $E_2'$ as above, such that $p_1=p_2\circ f$. One easily verifies that this extension satisfies $\phi_1 = \phi_2\circ f$ as well. 

Hence we have constructed a functor
\[ \Lambda_S : \: \:  \left[\hgammaxbar \, / \, P\times A \right](S) \: \rightarrow \: \FDPc(S), \]
for any scheme $S$, and this in turn defines a morphism of stacks
\[ \Lambda : \: \: \left[\hgammaxbar \, / \, P\times A \right] \: \rightarrow \: \FDPc \]
in the usual way. 
\ebew

We now arrive at our central theorem on one-dimensional substacks of the moduli stack of Deligne-Mumford stable curves. Note that here we are restricting our attention to the open substack $\cdx$ only. Before we formulate the result, it may be helpful to recall some of the notation. 

\begin{rk}\em
Let $f: C\rightarrow \Spec(k)$ be a stable curve of genus $g$, which is represented by a point $[C]\in \dpo$. By lemma \ref{newsplit},  there is an exact sequence of groups
\[ \id \rightarrow \Aut(C_0^-) \rightarrow \Aut(C) \rightarrow \Aut_\Gamma(C^+;\calq) \rightarrow \id ,\]
where $C^+$ denotes the unique $\mgamma$-pointed stable subcurve of $C$, which is the closure of the complement of the unique subcurve of $C$ isomorphic to $C_0^-$, and $\Gamma= \Gamma(C_0^-;\calq)$. As we noted in remark \ref{newsplitrk}, this sequence splits for almost all $[C]\in \dpo$. We denoted by
\[ \dx \:\: \subset \:\: \dpo  \label{n26a}\]
the dense open subscheme parametrizing those stable curves $f:C\rightarrow \Spec(k)$, for which there exists a splitting $\Aut(C) = \Aut_\Gamma(C^+;\calq)\times \Aut(C_0^-)$, and by 
\[ \cdx \:\: \subset \:\: \cdpo \label{n26b}\]
the dense open substack, which is defined as the reduction of the preimage of $\dx$ in $\cmgbar$. 

In notation \ref{D322}, we defined $\mx$ and $\hx$ as the corresponding subschemes in $\mgprime$ and $\hgprime$. They are both invariant under the action of $\Gamma$, and we denote their respective quotients by $\mgammax$ and $\hgammax$.\label{n26c} Clearly, the subscheme $\hgammax$ is invariant under the action of the group  $\PGL(N^+ +1)$ on $\hgprimebar$ as well. 

By  $\kx$\label{n26d} we denote the corresponding subscheme of $K(C_0;P_1)$ as in notation  \ref{n24}. 
\end{rk}

\begin{thm}\label{MainThm}
There is an isomorphism of Deligne-Mumford stacks
\[ \left[ \hgammadetailx  \, / \, \PGL(N^{+}+1)\times \Aut(C_0^{-}) \right] \: \: \cong \: \:  \cdx \]
where $\Aut(C_0^{-})$ acts trivially on $\hgammadetailx$. 
\end{thm}

\begin{rk}\em\label{rk00411}
By the isomorphism of proposition \ref{P217}, we can view the quotient stack $[\overline{H}_{g^{+},n,m/\Gamma(C_0^-;\CQ)}'  / \PGL(N^{+}+1)\times \Aut(C_0^{-})]$ also as the stack quotient $\cmgprimebarquot / \Aut(C_0^-)$, or equivalently, as a product of the moduli stack $\cmgprimebarquot$ with the classifying stack $\mbox{B}\Aut(C_0^-)$ of the group $\Aut(C_0^-)$. Thus theorem \ref{MainThm} states that generically the stack $\cdpobar$ is isomorphic to the product of stacks $\cmgprimebarquot \times \mbox{B} \Aut(C_0^-)$. 

With respect to the open substack $\cdpo$, theorem \ref{MainThm} is in fact the strongest result possible. Indeed, consider the case of a curve $C_0$, for which the stacks $\cdx$ and $\cdpo$ are not the same. For a curve $f: C\rightarrow \Spec(k)$, which is an object of $\cdpo(\Spec(k)$, but not of $\cdx(\Spec(k)$, the isomorphism must fail, since the splitting of the stack into a product would imply a corresponding splitting of the automorphism group $\Aut(C)$, contradicting the definition of $\cdx$. 

Later on, in the following chapters, we will discuss when and how the result of theorem \ref{MainThm} extends to the boundary $\cdpobar \smallsetminus \cdpo$. 
\end{rk}

\proofof{theorem \ref{MainThm}}
The proof is by explicit construction of an inverse morphism $\Xi$ to the
restriction of the morphism $\Lambda$ from  proposition
\ref{MainProp} to the open substack $[ \, \hgammax \, / \,
P\times A  \, ]$ in such a way, that the composition $\Lambda \circ \Xi$ is isomorphic to, and the composition $\Xi\circ \Lambda$ is equal to the respective identity functor. 

$(i)$  Consider $(E,p,\phi) \in \cdx(S)$ for some scheme $S$. By lemma \ref{140} we may assume that $p:E \rightarrow S$ is a principal $\PGL(N+1)$-bundle, together with a $\PGL(N+1)$-equivariant morphism $\phi: E \rightarrow \kx$. Consider the universal stable curve over $\hilb{g,n,0}$ restricted to $\kx$
\[ u : \: \:  {\cal C}_{g,n,0} \rightarrow \kx, \]
 and denote its pullback to $E$ via $\phi$ by 
\[ u_E : \: \: {\cal C}_E \rightarrow E .\]
Taking quotients by the action of $\PGL(N+1)$  induces a morphism
\[ f_S : \quad {\cal C}_S := {\cal C}_E / \PGL(N+1) \rightarrow E / \PGL(N+1) = S ,\]
which is a stable curve of genus $g$ over $S$. The induced morphism  $\theta_f : S \rightarrow \msbar_g$ factors through $\dpx$. Therefore,  by corollary \ref{K324a},  there exists a unique subscheme $f_S^-: {\cal C}_S^- \rightarrow S$ of $f_S: {\cal C}_S\rightarrow S$, and an \'etale covering $\{s_\alpha: S_\alpha\rightarrow S\}_{\alpha\in A}$ of $S$, such that for the pullback ${\cal C}_{S_\alpha}^-$ of ${\cal C}_S^-$ to $S_\alpha$, for each $\alpha \in A$,  there is a commutative diagram
\[ \diagram
{{\cal C}}^{-}_{S,\alpha} 
\dto_{f^{-}_{S,\alpha}}\rrto&& {{\cal C}}_S \dto^{f_S}\\
S_\alpha \rrto_{s_\alpha}&& S ,
\enddiagram\]
where ${{\cal C}}^{-}_{S,\alpha} = C_0^{-} \times S_\alpha$ is the trivial $m$-pointed stable curve over $S_\alpha$. 

Recall that $C_0^{-}$ decomposes into connected components
$C_1^{-},\ldots,C_r^{-}$, where $C_i^{-}$ is an $m_i$-pointed
stable curve of genus $g_i^{-}$ for $i=1,\ldots,r$, using the
notation of \ref{38}.

For  $i=1,\ldots,r$ let $f_{S,\alpha,i}^{-} : {\cal
C}_{S,\alpha,i}^{-} \rightarrow S_\alpha$ denote the corresponding
subcurve of $f_{S,\alpha}^{-} : {\cal C}_{S,\alpha}^{-} \rightarrow
S_\alpha$. Let $f_{{S_\alpha}} : {\cal C}_S | S_\alpha \rightarrow S_\alpha$ denote the pullback of $f_S$ to $S_\alpha$.
As in construction \ref{const_part_1} there is a natural inclusion of projective bundles
\[ \PP( f_{S,\alpha,i}^{-})_\ast(\omega_{{\cal C}_{S,\alpha,i}^{-}/S_\alpha}(S_{i,1}+\ldots+S_{i,m_i}))^{\otimes n} \: \hookrightarrow \:
\PP(f_{{S_\alpha}})_\ast(\omega_{{\cal C}_S|S_\alpha/S_\alpha})^{\otimes n}, \]
where $S_{i,1},\ldots,S_{i,m_i}$ denote the sections of the $m_i$ marked points on ${\cal C}_{S,\alpha,i}^{-}$. 

Put  $N_i^{-} := (2g_i^{-}-2+m_i)n-g_i^{-}+1$. There is a principal
$\PGL(N^{-}_i+1)$-bundle $p^{-}_{\alpha,i} : E_{\alpha,i}^{-}
\rightarrow S_\alpha$ on $S_\alpha$ associated to the projective
bundle $\PP( f_{S,\alpha,i}^{-})_\ast(\omega_{{\cal
C}_{S,\alpha,i}^{-}/S_\alpha}(S_{i,1}+\ldots+S_{i,m_i}))^{\otimes
n}$. It is a subbundle of the principal $\PGL(N+1)$-bundle $p_\alpha:
E_\alpha \rightarrow S_\alpha$ associated to the projective bundle
$\PP(f_{{S_\alpha}})_\ast (\omega_{{\cal C}_S|S_\alpha/S_\alpha})^{\otimes n}$
over $S_\alpha$, which is nothing else but the pullback of $p:
E\rightarrow S$ to $S_\alpha$ by proposition
\ref{28}.  

If one defines $f^{-}_{E,\alpha,i}: {\cal C}^{-}_{E,\alpha,i} \rightarrow E_{\alpha,i}^{-}$ by the Cartesian diagram
\[ \diagram 
{\cal C}^{-}_{E,\alpha,i} \rto \dto_{f_{E,\alpha,i}^{-}} & {\cal C}^{-}_{S,\alpha,i} \dto^{f_{S,\alpha,i}^{-}} \\
E_{\alpha,i}^{-} \rto_{p^{-}_{\alpha,i}} & S_\alpha,
\enddiagram \]
then it is an $m_i$-pointed stable curve of genus $g^{-}_i$ over
$E_{\alpha,i}^{-}$. Let $\tilde{S}_{i,1},\ldots,\tilde{S}_{i,m_i}$
denote the divisors of marked points on ${\cal C}^{-}_{E,\alpha,i}$. 
The projective bundle
\[ \PP(f_{E,\alpha,i}^{-})_\ast (\omega_{{\cal C}^{-}_{E,\alpha,i}/E_{\alpha,i}^{-}}(\tilde{S}_{i,1}+\ldots+\tilde{S}_{i,m_i}))^{\otimes n}  \]
is canonically isomorphic to the pullback bundle 
\[ (p_{\alpha,i}^{-})^\ast \PP(f_{S,\alpha,i}^{-})_\ast( \omega_{{\cal
C}_{S,\alpha,i}^{-}/S_\alpha}(S_{i,1}+\ldots+S_{i,m_i}))^{\otimes
n}, \]
 and  thus has a natural trivialization. This is true because any projective bundle pulled back to its associated principal bundle becomes trivial. 
Hence by the universal property of the Hilbert scheme there is an induced morphism
\[ \theta_{\alpha,i}^{-}: \:\: E_{\alpha,i}^{-} \longrightarrow  \hilb{g^{-}_i,n,m_i}, \]
which is $\PGL(N^{-}_i+1)$-equivariant.

Put $P^{-} := \PGL(N^{-}_1+1) \times \ldots \times \PGL(N^{-}_r +1)$. Taking everything together, we obtain a principal $P^{-}$-bundle
\[ p^{-}_\alpha : \quad E^{-}_\alpha := E_{\alpha,1}^{-}\times_{S_\alpha} \ldots\times_{S_\alpha} E_{\alpha,r}^{-} \: \longrightarrow \: S_\alpha \]
over $S_\alpha$, together with a $P^{-}$-equivariant morphism
\[ \theta_\alpha^{-}:= \theta_{\alpha,1}^{-}\times
\ldots\times\theta_{\alpha,r}^{-}  : \quad E^{-}_\alpha \: \longrightarrow \: \hilb{g^{-},n,m}^{-} ,\]
where $\hilb{g^{-},n,m}^{-} := \hilb{g^{-}_1,n,m_1}\times \ldots
\times  \hilb{g^{-}_r,n,m_r}$. 

Put ${\cal C}_{E,\alpha}^- := {\cal C}^{-}_{E,\alpha,1}\times \ldots\times{\cal C}^{-}_{E,\alpha,r}$. 
By construction the induced  morphism on the quotients
\[ {\cal C}_{S,\alpha}^{-}={\cal C}^{-}_{E,\alpha} / P^{-}
\: \longrightarrow \: E_\alpha^{-} / P^{-}  \cong S_\alpha \]
is an $m$-pointed prestable curve of genus $g^{-}$ over $S_\alpha$. Hence there is an induced  morphism
\[ \overline{\theta}_\alpha^-: \:\:\: S_\alpha \: \longrightarrow \: \msbar_{{g^{-},m}}^{-} := \msbar_{{g^{-}_1,m_1}}\times\ldots\times \msbar_{{g^{-}_r,m_r}}. \]
Putting the above  morphisms together we obtain a commutative diagram
\[ \diagram 
E^{-}_\alpha \rto^{\theta_\alpha^{-}} \dto_{p_\alpha^{-}} & \hilb{g^{-},n,m}^{-} \dto\\
S_\alpha \rto_{\overline{\theta}_\alpha^-} & \msbar_{{g^{-},m}}^{-}.
\enddiagram \]
Note that by construction all fibres of ${\cal C}^{-}_{S,\alpha} \rightarrow S_\alpha$ are isomorphic to $C_0^{-}$, hence $\overline{\theta}_\alpha^-: S_\alpha \rightarrow  \msbar_{{g^{-},m}}^{-}$ is constant. So there is in fact a commutative diagram
\[ \diagram 
E^{-}_\alpha \rto^{\theta_\alpha^{-}} \dto_{p_\alpha^{-}} & F \dto\\
S_\alpha \rto & \{ [C_0^{-}]\},
\enddiagram \]
where $F\cong P^{-} / A$ is the fibre of the canonical
morphism $\hilb{g^{-},n,m}^{-} \rightarrow \msbar_{{g^{-},m}}^{-}$
over the point $[C_0^{-}]$ representing the isomorphism class of $C_0^-$.

Fix an embedding of $C_0^{-}$ into $\PP^{N_1^{-}}\times \ldots\times
\PP^{N_r^{-}}$, corresponding to a point in $\hilb{g^{-},n,m}^{-}$,
which shall also be denoted by $[C_0^{-}]$. Then the decomposition of $A =
\Aut(C_0^{-})$ into a product as in remark \ref{323a} can be
considered as contained in the respective decomposition of $P^{-}$. 
A standard computation now shows that 
\[ E^{-}_\alpha/A \cong F \times S_\alpha \]
is a trivial fibration. Indeed, we can define a trivializing section
$\sigma_\alpha: S_\alpha \rightarrow E^{-}_\alpha/ A$  as follows.
For $s \in S_\alpha$ choose some $e\in (p_\alpha^{-})^{-1}(s) \subset
E_\alpha^{-}$, such that $\theta_\alpha^{-}(e) = [C_0^{-}]\in F
\subset \hilb{g^{-},n,m}^{-} $. Under the identification $F \cong
P^{-}/A$ the point  $[C_0^{-}]$ corresponds to the class $[\id]$ of
the identiy $\id \in P^{-}$. Hence by the $P^{-}$-equivariance of $
\theta_\alpha^{-}$ the point $e$ is uniquely defined up to the action of
$A$, and thus $\sigma_\alpha(s) := [e]\in E_\alpha^{-}/A$ is well defined. 

$(ii)$ So far we have been working on an \'etale cover of $S$. Passing to the global picture we have to take into account the actions of certain groups of automorphisms. 
 
Note that the local bundles $p_\alpha^{-}: E^{-}_\alpha \rightarrow
S_\alpha$ do not necessarily glue to give a global bundle on
$S$. However, it follows from remark \ref{323b} that we can glue
$E^{-}_\alpha / \Aut(C_0^-;\CQ)$ to obtain a bundle over $S$, which we denote by abuse of notation by
$p^{-}: E^{-}/\Aut(C_0^-;\CQ) \rightarrow S $, and which admits a section
\[ \sigma: \:\: S \rightarrow  E^{-}/\Aut(C_0^-;\CQ)\]
by glueing the trivializing sections $\sigma_\alpha$. 
By construction, there is an inclusion of the quotient $E^{-}/\Aut(C_0^-;\CQ)$ into
 $E/A$, compare  remark \ref{323b}.  By composition with the natural quotient morphism, we obtain a section 
\[ \tilde{\sigma} : \:\: S \rightarrow E / P\times  A .\]

However, this section defines a reduction of the principal
$\PGL(N+1)$-bundle $p : E \rightarrow S$ to a principal $P \times
 A$-bundle 
\[ p^{+} : \:\: E^{+} \rightarrow S,\]
 which is a subbundle of $E$, such that the extension
$(E^{+}\times \PGL(N+1))/( P\times  A)$ is isomorphic to $E$. Compare also remark \ref{rkB141} of the appendix. 
Note that this property makes this construction the inverse of the
construction of proposition \ref{MainProp}. In fact, we define $E^{+}$
as the preimage of the image of the section $\tilde{\sigma}$ under the
natural quotient map $E \rightarrow E/P\times  A$.  

$(iii)$ 
Recall that our aim is to construct an object
$(E^{+},p^{+},\phi^{+})\in $ \linebreak $[\hgammax/P\times
A](S)$. To do this, we stil need to define an $P\times
A$-equivariant morphism $\phi^{+}:E^{+} \rightarrow
\hgammax$. 

Locally in the \'etale topology, the universal curve $u: {\cal C}_{g,n,0}\rightarrow \kx $
over $\kx$ is covered by stable curves $u_\beta: {\cal C}_\beta
\rightarrow K_\beta$, where $\{k_\beta:K_\beta\rightarrow \kx\}_{\beta \in B}$ is an \'etale 
cover of $\kx$ as in lemma \ref{315}, such that $u_\beta$ contains
a trivial subcurve $u_\beta^{-} : C_0^{-}\times K_\beta \rightarrow
K_\beta$. Let $u_\beta^{+} : {\cal C}_\beta^{+} \rightarrow K_\beta$
denote the closure of the complement of this trivial subcurve in
$u_\beta: {\cal C}_\beta \rightarrow K_\beta$. 
This is an $m$-pointed stable curve of genus $g^{+}$ over $K_\beta$. 
By construction, $E^{+}
\subset E$ is a subbundle, which is stabilized under the action of the subgroup $P$ in $\PGL(N+1)$. Thus the $\PGL(N+1)$-equivariant morphism $\phi: E\rightarrow \kx$ induces a
$P$-equivariant morphism
\[\tilde{\phi} : \: \:   E^{+} \rightarrow \kx .\]
Consider the pullback 
\[ \tilde{\phi}_\beta : \quad E^{+}_\beta := k_\beta^\ast E^+   \rightarrow K_\beta .\]
Define  ${\cal C}_{E,\beta}^{+}$ by the Cartesian diagram
\[ \diagram
{\cal C}_{E,\beta}^{+} \dto_{f_{E,\beta}^{+}} \rrto && {\cal C}^{+}_\beta\dto^{u_\beta^{+}}\\
E^{+}_\beta  \rrto_{\tilde{\phi}_\beta} &&              K_\beta.
\enddiagram  \]
Then $f_{E,\beta}^{+}: {\cal C}_{E,\beta}^{+} \rightarrow E^{+}_\beta$ is an $m$-pointed stable curve of genus $g^{+}$ over $E^{+}_\beta$, and the projective bundle 
\[ \PP(f_{E,\beta}^{+})_\ast(\omega_{{\cal C}_{E,\beta}^{+}/E^{+}_\beta}(S_1+\ldots+S_m))^{\otimes n} \]
on $E_\beta^{+}$ is isomorphic to the pullback of the projective bundle
\[ \PP(u_\beta^{+})_\ast (\omega_{{\cal C}^{+}_\beta/K_\beta}({\cal S}_1+\ldots+{\cal S}_m))^{\otimes n} \]
on $K_\beta$ along the morphism $\tilde{\phi}_\beta: E_\beta^{+}\rightarrow K_\beta$. Here again $S_1,\ldots,S_m$ and ${\cal S}_1,\ldots,{\cal S}_m$ denote the divisors of the $m$ marked points on ${\cal C}_{E,\beta}^{+}$ and ${\cal C}_{\beta}^{+}$, respectively. 

$(iv)$ 
Note that there is a distinguished trivialization of the projective
bundle  $\PP(u_\beta^{+})_\ast (\omega_{{\cal
C}^{+}_\beta/K_\beta}({\cal S}_1+\ldots+{\cal S}_m))^{\otimes n}$
over $K_\beta$.\footnote{In fact, the triviality of the pullback $\tilde{\phi}_\beta^\ast \PP(u_\beta^{+})_\ast (\omega_{{\cal
C}^{+}_\beta/K_\beta}({\cal S}_1+\ldots+{\cal S}_m))^{\otimes n}$ shows that the principal $P$-bundle $E_\beta^+/A$ is isomorphic to the principal $P$-bundle associated to the projective bundle $\PP(u_\beta^{+})_\ast (\omega_{{\cal
C}^{+}_\beta/K_\beta}({\cal S}_1+\ldots+{\cal S}_m))^{\otimes n}$ over $K_\beta$.} 

This can be seen as follows. The curve $u_\beta: {\cal C}_\beta\rightarrow K_\beta$ is constructed as a pullback of the universal curve on $\overline{H}_{g,n,0}$, restricted to $\kx$. Therefore there exists a natural embedding of the subcurve $u_\beta^+: {\cal C}_\beta^+\rightarrow K_\beta$ into $\pr_2: \PP^N \times K_\beta\rightarrow K_\beta$. To distinguish a trivialization of $\PP(u_\beta^{+})_\ast (\omega_{{\cal
C}^{+}_\beta/K_\beta}({\cal S}_1+\ldots+{\cal S}_m))^{\otimes n}$ over $K_\beta$, it suffices to distinguish a corresponding embedding of the subcurve $u_\beta^+: {\cal C}_\beta^+\rightarrow K_\beta$ into $\pr_2: \PP^{N^+} \times K_\beta\rightarrow K_\beta$. 

Consider first the curve $u^+: {\cal C}_{g^+,n,m}^\times \rightarrow \hx$, which is the restricition of the universal curve on the full Hilbert scheme. By definition, this curve is embedded into $ \PP^{N^+}\times \hx = \PP u^+_\ast(\omega_{{\cal C}_{g^+,n,m}'/\hx}({\cal S}_1+\ldots+{\cal S}_m))^{\otimes n}$, which in turn is embedded into $\PP(u_0)_\ast(\omega_{{\cal C}_0/\hx})^{\otimes n} = \PP^N \times \hx$ by the natural inclusion 
\[ \PP u^+_\ast(\omega_{{\cal C}_{g^+,n,m}^\times/\hx}({\cal S}_1+\ldots+{\cal S}_m))^{\otimes n} \:\: \subset \:\: 
\PP(u_0)_\ast(\omega_{{\cal C}_0/\hx})^{\otimes n} . \]
Compare construction \ref{const_part_1} for details and for the notation. Note that the embedding is fibrewise independent of the ordering of the labels of the marked points. We used this in proposition \ref{const_prop} for the construction of the  morphism $\thx: \hxbar \rightarrow \overline{H}_{g,n,0}$. Recall that $\thx$ induces a morphism $\Theta: \hgammaxbar \rightarrow \overline{H}_{g,n,0}$, see proposition \ref{312}. 

For each closed point $[C]\in \Theta(\hgammax) = \thx(\hx)$ representing a stable curve $C$ embedded into $\PP^N$, there is a unique subcurve isomorphic to $C_0^-$, and thus a unique subcurve $C^+$, which is an $\mgamma$-pointed stable curve of genus $g^+$, satisfying $\Theta([C^+]) = [C]$. Thus, the restricition of $u^+_\beta:  {\cal C}_\beta^+\rightarrow K_\beta$ to the subscheme $K_\beta' \subset K_\beta$ lying over $\Theta(\hgammax)$ has a distinguished embedding into $\pr_2 : \PP^{N^+}\times K_\beta' \rightarrow K_\beta'$, namely the one given by the embedding of the universal curve over $\hx$. We want to extend this embedding to all of $K_\beta$, using corollary \ref{314}, which implies that $\kx = \Theta(\hgammax) \cdot \PGL(N+1)$. 

Let $[C]\in \Theta(\hgammax)$ be a closed point representing a stable curve $C$ embedded into $\PP^N$, with a decomposition $C= C^+ \cup C_0^-$ as above, together with the embedding of $C^+$ into $\PP^{N^+}$. For an element $\gamma\in \PGL(N+1)$, the distinguished $\mgamma$-pointed stable subcurve of genus $g^+$ of the translated embedded curve $\gamma C$ is the curve $\gamma C^+$, by the uniqueness statement of corollary \ref{K324}. The obvious embedding of $\gamma C^+$ into $\gamma \PP^{N^+} \subset \PP^N$ is in fact independent of the choice of $\gamma$. 

Indeed, let $\gamma'\in \PGL(N+1)$ be another element such that $\gamma'([C])=\gamma([C])$. Then $\gamma'\circ\gamma^{-1}\in\Aut(C)$. By the definition of $\cdx$, this group splits as
$\Aut(C)=\Aut_\Gamma(C^+;\calq)\times \Aut(C_0^-)$. Since the embedding of $C^+$ into $\PP^{N^+}$ is independent of the ordering of the labels of the marked points, and since $\Aut(C^+) = \Stab_{\PGL(N^+ +1)}([C^+])$, any elements of $\Aut(C)$ lying in $\Aut_\Gamma(C^+;\calq)$ leave the embedding invariant. By construction \ref{const_part_2}, there is a fixed embedding $C_0^- \subset \PP^{N^-} \subset \PP^N$, and thus  inclusions $\Aut(C_0^-)\subset \PGL(N^- +1) \subset \PGL(N+1)$. Note that $\PGL(N^- +1)$ acts trivially on $\PP^{N^+}$ outside the intersection $\PP^{N^+}\cap \PP^{N^-}$. However, elements of $\Aut(C_0^-)$ fix by definition the marked points of $C_0^-$, which span the intersection $\PP^{N^+}\cap \PP^{N^-}$ as a subspace of $\PP^N$, and thus fix the subspace $\PP^{N^+}$. Hence elements of $\Aut(C)$ lying in $\Aut(C_0^-)$ fix the embedding of $C^+$ as well. 

Therefore, for all closed points $[C]\in \kx$, representing an embedded stable curve $C\subset \PP^N$, there is a well-defined embedding of the distinguished $\mgamma$-pointed stable subcurve $C^+ \subset C$ into $\PP^{N^+}$, and so there exists an embedding of  $u_\beta^+: {\cal C}_\beta^+\rightarrow K_\beta$ into $\pr_2: \PP^{N^+} \times K_\beta\rightarrow K_\beta$, as claimed.  

$(v)$ 
So finally we obtain for each index $\beta\in B$  an $m$-pointed stable curve $f_{E,\beta}^{+}: {\cal C}_{E,\beta}^{+} \rightarrow E^{+}_\beta$, together with a trivialization of the projective bundle $\PP(f_{E,\beta}^{+})_\ast(\omega_{{\cal C}_{E,\beta}^{+}/E^{+}_\beta}(S_1+\ldots+S_m))^{\otimes n}$. By the universal property of the Hilbert scheme, this defines a $P$-equivariant morphism
\[ \phi^{+}_\beta : \: \: E_\beta^{+} \rightarrow \hx.\]

Up to a permutation of the $m$ marked points under the action of $\Gamma$ these morphisms glue together, so we obtain a morphism
\[ {\phi}^{+} : \quad {E}^{+} \rightarrow
\hx/\Gamma = \hgammax,\]
which is $P\times  A$-equivariant, with $ A$
acting trivially on $\hgammax$.

$(vi)$
Taking everything together, we obtain an object 
\[ (E^{+}, p^{+},\phi^{+}) \in \left[  \hgammax/P\times A \right](S). \]
Note that locally and up to the action of $A$ this triple just represents the $\mgamma$-pointed stable curve over $S_\alpha$, which is the closure of the complement of the trivial curve $C_0^{-}\times S_\alpha \rightarrow S_\alpha$ in $C_\alpha\rightarrow S_\alpha$, as in lemma \ref{315}. 

Up to isomorphism, this construction of the reduction $E^{+}$ of $E$
is inverse to the extension used in the proof of proposition
\ref{MainProp}. In fact one has the equality $\Xi(\Lambda(E^{+},p^+,\phi^+)) =
(E^{+},p^+,\phi^+)$, and an isomorphism $\Lambda(\Xi(E,p,\phi)) \cong (E,p,\phi)$. 

$(vii)$ It finally remains to define the inverse functor on morphisms
\[ \varphi : \quad (E_2,p_2,\phi_2) \rightarrow (E_1,p_1,\phi_1) \]
in $\cdx(S)$. Such a morphism is given by a morphism of principal $\PGL(N+1)$-bundles $f : E_1 \rightarrow E_2$, such that the diagram
\[ \diagram
E_1 \xto[rrd]^{\phi_1} \drto^f \xto[rdd]_{p_1}\\
& E_2 \rto_{\phi_2} \dto^{p_2}& \kx\\
&S
\enddiagram \]
commutes. We need to show that $f$ restricts to a morphism of principal $P\times A$-bundles
\[ f^{+} : \: \: E_1^{+} \rightarrow E_2^{+} ,\]
where $E_1^{+}$ and $E_2^{+}$ are the $P\times  A$-subbundles of $E_1$ and $E_2$ constructed as above. Note that this restriction, if it exists, necessarily satisfies 
\[ \phi_1 = \phi_2\circ f^+ .\]
For $i=1,2$, the morphism $\phi_i$ factors as $\phi_i=\Theta\circ\phi_i^+$, where $\Xi(E_i,p_i,\phi_i)=(E_i^+,p_i^+,\phi_i^+)$. The injectivity of $\Theta : \hgammax \rightarrow \kx$, which follows from proposition \ref{312},  implies the identity 
\[ \phi_1^+ = \phi_2^+\circ f^+,\]
so $f^+:E_1^+\rightarrow E_2^+$ defines indeed a morphism between the triples $(E_2^+,p_2^+,\phi_2^+)$ and $(E_1^+,p_1^+,\phi_1^+)$. 

It is clear from the construction that this inverts the definition of $\Lambda$ on morphisms. Thus $\Xi$ is shown to be an inverse functor to $\Lambda$, up to isomorphism, provided that we can prove the existence of $f^+$. 

To show that such a restriction $f^+:E_1^+\rightarrow E_2^+$ exists, we  need to verify 
\[ f(E_1^{+}) \subseteq E_2^{+} .\]
Recall that  $E_1^{+}$ and $E_2^{+}$ are the preimages of 
the images of sections $\tilde{\sigma}_1 : S \rightarrow
E_1/P\times A$ and $\tilde{\sigma}_2: S \rightarrow
E_2/P\times A$ under the canonical quotient maps $E_1 \rightarrow E_1/P\times A$ and $E_2 \rightarrow E_2/P\times\ A$.

Let $i=1,2$. Recall that we used the morphism $\phi_i : E_i \rightarrow \kx$ and the pullback diagram of the universal curve
\[ \diagram
{\cal C}_{E_i}^{-} \rto \dto_{u_{E_i}} & {\cal C}_{g,n,0}\dto^u\\
E_i \rto& \kx
\enddiagram \]
to define a morphism 
\[ \theta_i^{-} : \: \:  E_i^{-} \rightarrow \hilb{g^{-},n,m}^{-}.\]
Strictly speaking, this is true only locally in the \'etale topology, or up to the action of
$\Aut(C_0^-;\CQ)$. Since the claim is local on the base $S$, we omit
for simplicity the indices ``$\alpha$'' we were working with in part $(i)$ of the
proof. 

By the universal property of the Hilbert space, we have the identity
\[ (\ast) \qquad \qquad \theta_1^- = \theta_2^-\circ f.    \]
We saw that in fact
\[ \theta_i^- : \:\: E_i^{-} \rightarrow F ,\]
for $i=1,2$, where $F\cong P^{-}/A$ is the fibre of the quotient morphism $\hilb{g^{-},n,m}^{-} \rightarrow \msbar_{g^{-},m}^{-}$ over the point representing the isomorphism class of $C_0^{-}$. 

Again, we identify here the image  $[\id]$ of the unit $\id \in
P^{-}/A$ with the point $[C_0^{-}]$ representing the fixed embedding of $C_0^{-}$. 

There is a section $\sigma_i: S \rightarrow E_i^{-}/A$, sending a point $s\in S$ to the equivalence class of points $e\in E_i^{-}$ in the fibre over $s$  with $\theta_i^-(e) = [\id]\in F$. This is well-defined because $\theta_i^-$ is $P^{-}$-equivariant. 
Let $\tilde{\sigma}_i$ denote the morphism from $S$ to $E_i/P\times A$ induced by $\sigma_i$. 
Now the bundle $E_i^{+}$ is defined as the unique subbundle of $E_i$,
which fits into the Cartesian diagram
\[ \diagram
E_i^{+}\rrto \dto_{p_i^{+}} &&E_i\dto^{q_i}\\
S\rrto_{\tilde{\sigma}_i}&&E_i/P\times  A ,
\enddiagram\]
where $q_i : E_i \rightarrow E_i/P\times  A $ denotes the natural quotient map. 
In other words,
for an element $e\in E_i$ holds $e\in E_i^{+}$ if and only if
\[ q_i (e) = \tilde{\sigma}_i(p_i^{+}(e)),\]
which in turn is equivalent to
\[ \theta_i^-(e) = [\id]\]
by the definition of $\tilde{\sigma}_i$. 
So the relation $f(E_1^{+}) \subseteq E_2^{+}$ translates into the condition
\[ \theta_2^-(f(e)) = [\id] \]
for all $e\in E_1^{+}$. However, above we found the identity $(\ast) \:\: \theta_1^- = \theta_2^-\circ f$, which concludes the proof.
\ebew

\begin{rk}\em
Recall that $\cmgammax$ denotes the substack of $\cmgbarquot$, which is the reduction of the preimage of
$\mgammax$. Analogously as in remark \ref{R323}, the scheme
$\dx$ is a moduli space of the stack
\[ \cdx \: \cong \: \left[ \hgammax/P\times A\right] ,\]
where the isomorphism holds by theorem \ref{MainThm}, 
and the scheme $\mgammax$ is a moduli space of the stack
\[ \cmgammax  \: \cong \:
\left[ \hgammax /P\right] .\]
There is an isomorphism $\dx \cong \mgammax $, but in general the stack $\cdx$ is not isomorphic to the stack
$\cmgammax$. 

This may seem surprising at a first glance. However, obviously there are in general more automorphisms of $C_0$ that there are automorphisms of $C_0^{+}$. Moduli stacks, in contrast to moduli spaces, are designed not to parametrize isomorphism classes of stable curves only, but also to ``remember'' their automorphisms. Therefore the moduli stack $\cdpobar$ should be ``richer'' than  the moduli stack  $\cmgbarquot$. 
\end{rk}

\begin{cor}
There is an isomorphism of stacks
\[ \cdx \cong \cmgammax \]
if and only of $\Aut(C_0^{-})$ is trivial.
\end{cor}

\proof
This follows immediately from theorem \ref{MainThm} and the theory of quotient stacks.
\ebew

\begin{rk}\em
To getter a better grasp of the idea that in general the stacks $\cdx$ and $\cmgammax$ cannot be isomorphic without taking the additional quotient with respect to the action of $\Aut(C_0^{-})$, consider the natural morphism of stacks
\[\tilde{\Xi} : \:\:\:  \cdx \rightarrow \cmgammax.\]
As a functor it is constructed by sending a stable curve $f: C\rightarrow S$ of genus $g$, which is an object of $\cdx(S)$, to the unique $\mgamma$-pointed stable subcurve $f' : C' \rightarrow S$, which is the closure of the complement of the distinguished subcurve $f^-: C^-\rightarrow S$ with fibres isomorphic to $C_0^-$ as in corollary \ref{K324a}. Suppose that $\Aut(C_0^-)$ is not trivial. Then there exists a nontrivial automorphism of $C$ which is the identity when restricted to $C'$. The functor $ \tilde{\Xi}$ maps this automorphism to the identity on $f' : C' \rightarrow S$ in $\cmgammax(S)$. Hence the morphism   $\tilde{\Xi}$ is not faithful as a functor, and thus in particular not an isomorphism of stacks.
\end{rk}

\begin{rk}\em
In general, there is a commutative diagram
\[ \begin{array}{ccccc}
\hgammax& =   & \hgammax & \hookrightarrow &H_{g,n,0}\\
 \downarrow &&\downarrow &&\downarrow \\
\cmgammax &\longleftarrow & \cdx &\hookrightarrow & \cmg\\
\downarrow &&\downarrow &&\downarrow \\
\mgammax & \cong &\dx &\hookrightarrow &\ms_{g}.
\end{array} \]
The canonical morphism
\[ \cdx \cong \left[ \hgammax / P\times A\right] \:  \longrightarrow \: \left[\hgammax / P \right] \cong \cmgammax \]
is of degree equal to $\frac{1}{\# \Aut(C_0^{-})}$.
\end{rk}

\begin{rk}\em
$(i)$ It is important to note that unlike the situation of the moduli spaces, where we have an isomorphism $\mgprimebarquot \cong \dpobar$, the result of theorem \ref{MainThm} does not hold true when extended to the closed stack $\cdpobar$, even in the cases where the two stacks $\cdx \subset \cdpo$ are the same. It is the topic of the following chapters to analyze the boundary behavior of $\cdpobar$ in detail. \\
$(ii)$ If the substack $\cdx$ is strictly smaller that the stack $\cdpo$, then the isomorphism of theorem \ref{MainThm} does not even extend to an isomorphism between  $\cmgprimequot /\Aut(C_0^-)$ and $\cdpo$. This can easily be seen. Consider a curve $f: C\rightarrow \Spec(k)$, which is an object of $\cdpo(\Spec(k))$. It decomposes naturally into two subvurves. One of them is isomorphic to $C_0^-$, the other subcurve  $f^+: C^+\rightarrow \Spec(k)$ is given by an  object of the stack product $\cmgprimequot\times \mbox{B}\Aut(C_0^-)(\Spec(k))$, see \ref{rk00411}. If an isomorphism exists, the automorphism groups of both objects must be isomorphic, i.e. there must be a splitting $\Aut(C)\cong  \Aut_\Gamma(C^+;\CQ) \times \Aut(C_0^-)$. 
\end{rk}


%
%

\chapter{Partial compactification}\label{par5}

\begin{rk}\em
The key part in the proof of theorem \ref{MainThm} above was   to construct to a given triple $(E,p,\phi) \in \cdx(S)$, i.e. to  a given principal $\PGL(N+1)$-bundle $p : E\rightarrow S$, a principal $\PGL(N^+ +1)\times\Aut(C_0^-)$-bundle 
$p^+: E^+ \rightarrow S$, which is compatible with the inclusion $\PGL(N^+ +1)\times\Aut(C_0^-) \hookrightarrow \PGL(N+1)$. Geometrically, this means to find for a given stable curve $f:C\rightarrow S$ of some genus $g$ a subcurve $f^+:C^+\rightarrow S$, which is an $\mgamma$-pointed stable curve of genus $g^+$, or equivalently, a subcurve $f^-:C^-\rightarrow S$, which is an $\mgamma$-pointed stable curve of genus $g^-$, where $\Gamma = \Gamma(C_0^-;\calq)$. 

If the induced morphism $\theta_f : S\rightarrow \msbar_g$ factors through $D(P_1)$, then each fibre of $f$ has exactly $3g-4$ nodes. Therefore in each fibre $C_s$, with $s\in S$,  there is an unique irreducible component of $C_s$ which has not the maximal number of nodes when considered as a pointed stable curve itself. In this case the required decomposition of $f:C \rightarrow S$  into two subcurves exists. 

However, for curves with fibres having $3g-3$ nodes, this need no longer be true, even if their groups of automorphisms split. In general there is no way to distinguish a subcurve as above as there may be no unique choice fibrewise. For example this happens  for our fixed curve $f_0: C_0\rightarrow \Spec(k)$ if and only if the subcurve $C_0^+$, distinguished by the choice of $P_1$, is not stabilized under the action of $\Aut(C_0)$. To deal with this difficulty  we will replace the subcurve $C_0^+$ by an equivariant subcurve $C_0^{++}$. 
\end{rk}

\begin{defi}\em
Two nodes $P_1$ and $P_2\in C_0$ are called {\em equivalent} if there exists an automorphism of $C_0$ interchanging them.
\end{defi}

\begin{rk}\em
We saw above that $\overline{D}(P_1) = \overline{D}(P_2)$ holds if and only if $P_1$ and $P_2$ are equivalent nodes on $C_0$.
\end{rk}

\begin{defi}\em
Let $[P_1]$ denote the equivalence class of a node on a stable curve $C_0\rightarrow\Spec(k)$ of genus $g$ with $3g-3$ nodes. Let $P_1\in C_0$ be a node representing $[P_1]$, and let $C_0^{+}$ be the subscheme of $C_0$, which is the union of all irreducible components of $C_0$ containing $P_1$ as before. We define $C_0^{++}$\label{n28}  as the (reduced) subscheme of $C_0$, which is the orbit of $C_0^{+}$ under the action of $\Aut(C_0)$. Define $C_0^{--}$ as the closure of the complement of $C_0^{++}$ in $C_0$. If we want to emphasize the defining node we will also use the notations $C_0^{++}(P_1)$ and $C_0^{--}(P_1)$ instead of $C_0^{++}$ and $C_0^{--}$.  
\end{defi}

\begin{rk}\em
Equivalently, $C_0^{++}$ is the union of all irreducible components of $C_0$ containing at least one node of the equivalence class $[P_1]$. 

The intersection of the subcurves  $C_0^{--}$ and $C_0^{++}$ consists of finitely many points ${\cal R} = \{R_1,\ldots,R_{\mu}\}$, all of which are nodes of $C_0$. Note that the chosen enumeration of the nodes of $C_0$ distinguishes an ordering of these marked points on $C_0^{--}$ and $C_0^{++}$. Thus the curves $f_0^{++}: C_0^{++} \rightarrow \Spec(k)$ and $f_0^{--}: C_0^{--} \rightarrow \Spec(k)$ are in a natural way $\mu$-pointed prestable curves of genus $g^{++}$ and $g^{--}$, respectively. Both of them decompose into families of pointed stable curves  $C_i^{++} \rightarrow \Spec(k)$\label{n29} and $C_j^{--} \rightarrow \Spec(k)$\label{n30}, which  are  $\mu_i^{++}$-pointed stable curves of genus $g_i^{++}$, for $i=1,\ldots,r^{++}$, and $\mu_j^{--}$-pointed stable curves of genus $g_j^{--}$, for $j=1,\ldots,r^{--}$, each with the maximal number of nodes possible. 
\end{rk}

\begin{lemma}
All connected components $C_1^{++},\ldots,C_{r^{++}}^{++}$ of $\: C_0^{++}$ are isomorphic to each other. In particular $g_1^{++}=\ldots=g_{r^{++}}^{++}$ and $\mu_1^{++}=\ldots=\mu_{r^{++}}^{++}$.
\end{lemma}

\proof
Compare the connected component containing the subcurve $C_0^+$ with a second connected component. Without loss of generality assume that these components are $C_1^{++}$ and $C_2^{++}$. By definition, there exists an automorphism $\phi  \in \Aut(C_0)$ such that $C_2^{++}$ contains $\phi(C_0^{+})$ as a subscheme. Since $\phi$ maps connected components to connected components, $\phi$ restricts to an isomorphism of $C_1^{++}$ with $C_2^{++}$.
\ebew 

\begin{rk}\em
By construction, $C_0^{--}$ and $C_0^{++}$ are invariant under the action of $\Aut(C_0)$ as subschemes of $C_0$. 
\end{rk}

\begin{notation}\em
Consider the group $\Aut(C_0^{--};{\cal R})$\label{n31} of those automorphisms of the scheme $C_0^{--}$ over $\Spec(k)$, which stabilize the subset ${\cal R} = 
\{R_1,\ldots,R_\mu\}$ in $C_0^{--}$. The natural homomorphism from this group into the permutation group $\Sigma_\mu$ acting on this subset defines a subgroup $\Gamma(C_0^{--};\calr)$\label{n32} of $\Sigma_\mu$. As in remark \ref{R319}, there is an exact sequence of groups
\[ \id \rightarrow \Aut(C_0^{--}) \rightarrow \Aut(C_0^{--};\calr) \rightarrow \Gamma(C_0^{--};\calr) \rightarrow \id ,\]
where $ \Aut(C_0^{--})$ denotes the group of automorphisms of $C_0^{--}$  considered as a $\mu$-pointed prestable curve. As in definition \ref{D318a}, for any subgroup $\Gamma \subset \Sigma_\mu$, we denote by $\Aut_\Gamma(C_0^{--})$ the subgroup of those elements $\phi\in  \Aut(C_0^{--};\calr)$ for which  there exists a permutation $\gamma\in \Gamma$ such that
\[ \phi(R_i) = R_{\gamma(i)} \]
for all $i=1,\ldots,\mu$. In particular the following groups are the same.
\[ \Aut_{\Gamma(C_0^{--};\calr)}(C_0^{--}) = \Aut_{\Sigma_m}(C_0^{--}) = 
\Aut(C_0^{--};\calr) .\]
\end{notation}

\begin{rk}\em
It is not possible to give a comprehensive list of all possible types of subcurves $C_0^{++}$ as we did in remark \ref{37}. Note that unlike $C_0^+$, the curve $C_0^{++}$ needs not to be connected. 

Obviously, $C_0^{+} \subset C_0^{++}$, and  $C_0^{--} \subset C_0^{-}$. It may  happen that $C_0^{++} = C_0$, and thus $C_0^{--}=\emptyset$, as in example \ref{E57b} below. Any of the cases $m< \mu,m= \mu$ and  $m> \mu$ can occur. Throughout this section we will always assume $C_0^{--}\neq \emptyset$.
\end{rk}

\begin{ex}\em\label{E57a}
Consider the stable curve $C_0$ of genus $g=3$ as pictured below.

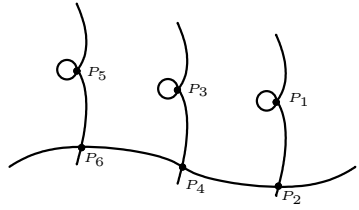
\begin{center}
\setlength{\unitlength}{0.00016667in}
\begingroup\makeatletter\ifx\SetFigFont\undefined%
\gdef\SetFigFont#1#2#3#4#5{%
  \reset@font\fontsize{#1}{#2pt}%
  \fontfamily{#3}\fontseries{#4}\fontshape{#5}%
  \selectfont}%
\fi\endgroup%
{\renewcommand{\dashlinestretch}{30}

}
\end{center}

The subcurve $C_0^+$ distinguished by the node $P_1$ is given by

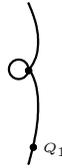
\begin{center}
\setlength{\unitlength}{0.00016667in}
\begingroup\makeatletter\ifx\SetFigFont\undefined%
\gdef\SetFigFont#1#2#3#4#5{%
  \reset@font\fontsize{#1}{#2pt}%
  \fontfamily{#3}\fontseries{#4}\fontshape{#5}%
  \selectfont}%
\fi\endgroup%
{\renewcommand{\dashlinestretch}{30}
%
}
\end{center}

considered as an $1$-pointed stable curve of genus $g^+=1$. Clearly this subcurve is not invariant under the action of $\Aut(C_0)$. The corresponding subcurve $C_0^{++}$ is given by

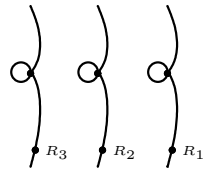
\begin{center}
\setlength{\unitlength}{0.00016667in}
\begingroup\makeatletter\ifx\SetFigFont\undefined%
\gdef\SetFigFont#1#2#3#4#5{%
  \reset@font\fontsize{#1}{#2pt}%
  \fontfamily{#3}\fontseries{#4}\fontshape{#5}%
  \selectfont}%
\fi\endgroup%
{\renewcommand{\dashlinestretch}{30}
%
}
\end{center} 

as a prestable $3$-pointed curve, represented by a point $[C_0^{++}]\in \msbar_{1,1}\times  \msbar_{1,1}\times \msbar_{1,1}$. 
\end{ex}

\begin{ex}\em\label{E57b}
Consider the same curve $C_0$ as in example \ref{E57a}, but now with focus on the node $P_2$. For the subcurve $C_0^{+}(P_2)$ one obtains

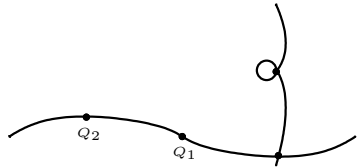
\begin{center}
\setlength{\unitlength}{0.00016667in}
\begingroup\makeatletter\ifx\SetFigFont\undefined%
\gdef\SetFigFont#1#2#3#4#5{%
  \reset@font\fontsize{#1}{#2pt}%
  \fontfamily{#3}\fontseries{#4}\fontshape{#5}%
  \selectfont}%
\fi\endgroup%
{\renewcommand{\dashlinestretch}{30}

}
\end{center} 

as a $2$-pointed stable curve of genus $g^+=1$. 
Clearly the invariant curve in this case is $C_0^{++}(P_2) = C_0$.
\end{ex}

\begin{notation}\em
Let $g' \ge 0$ and $\mu' \ge 0$ be given such that $2g'-2+\mu' > 0$. As in remark \ref{M23}, for sufficiently large $n'$, there is a subscheme $\overline{H}_{g',n',\mu'}$\label{n33} of the Hilbert scheme Hilb$_{\PP^{N'}}^{P_{g',n',\mu'}}\times (\PP^{N'})^{\mu'}$, parametrizing $\mu'$-pointed prestable curves $f: C \rightarrow S$ of genus $g'$, embedded into $\PP^{N'}$, where $N' = (2g'-2+\mu')n'-g'$, such that 
\[ (\omega_ {C/S}(S_1+\ldots+S_{\mu'}))^{\otimes n'} = \oka_{\PP^{N'}}(1) | C .\]
Here $S_1,\ldots,S_{\mu'}$ denote the divisors defined by the  $\mu'$ sections of marked points. Recall that we only consider such prestable curves, where for any fibre all connected components are pointed stable curves. 

As in \ref{def_split}, we define a reduced subscheme 
\[ \hxxbar \: \subset \: \hilb{g^{++},n,\mu} \label{k3} \]
of the Hilbert scheme $\hilb{g^{++},n,\mu}$ as  the locus of points $[C^{++}]$ parametrizing prestable curves  $f^{++}:C^{++}\rightarrow \Spec(k)$ with the following property: if $f: C\rightarrow \Spec(k)$ denotes the stable curve of genus $g$, which is obtained by glueing $C^{++}$ with $C_0^{--}$ in the $\mu$ marked points, then  $C$ has exactly $3g-4$ nodes, and its group of automorphisms splits naturally as
\[ \Aut(C) = \Aut_{\hat{\Gamma}}(C^{++};\calr)\times \Aut(C_0^{--}) ,\]
with $\hat{\Gamma} := \Gamma(C_0^{--};\calr)$. 

Note that the subscheme $\hxxbar$ does not only depend on the parameters $g^{++}$, $n$ and $\mu$, but also on the fixed curve $C_0$ and the chosen node $P_1$. 
\end{notation}

\begin{rk}\em\label{nonemptysub}
For certain choices of $C_0$ and $P_1$, the subscheme $\hxxbar$ may be the empty scheme. For example, consider the stable curve $C_0$ of genus $g=5$, consisting of a string of three rational lines, intersecting in a total of two points $P_1$ and $P_2$, with five nodal rational curves appropriately attached. 

For our considerations below we will always assume that $\hxxbar$ is not empty.
\end{rk}

\begin{rk}\em\label{R513}
Analogously to what we did in construction \ref{const_part_1} and \ref{const_part_2}, we can relate the schemes  $\overline{H}_{g^{++},n,\mu}$ and  $\overline{H}_{g^{--},n,\mu}$ to the scheme $\overline{H}_{g,n,0}$, where $n$ is assumed to be sufficiently large. By glueing  the curve $C_0^{--}$ along the $\mu$ marked points to each  curve represented by a point in $\hxxbar$ one can construct a $\PGL(N^{++}+1)$-equivariant morphism 
\[ \thxx :\quad  \hxxbar \: \longrightarrow \: \overline{H}_{g,n,0}. \label{n34} \]
By the same arguments as in the proof of proposition \ref{P318}, we may assume that the morphism $\thxx$ factors through the quotient 
\[ \hgammadetailxxbar \: := \: \hxxbar/ \Gamma(C_0^{--};\calr) . \label{n35}\] 
Note that the action of $\Gamma(C_0^{--};\calr)$ commutes with the action of the group $\PGL(N^{++}+1)$ on $ \hxxbar$. 

During this section we will from now on use the abbreviation $\hat{\Gamma} :=  \Gamma(C_0^{--};\calr)$. 
\end{rk}

\begin{notation}\em
Let $\hat{D}(C_0;P_1)$\label{n36} denote the partial compactification of the subscheme $D(C_0;P_1)$ at the point $[C_0]$ inside $\msbar_g$, i.e. as a set 
\[ \hat{D}(C_0;P_1) := D(C_0;P_1) \cup \{[C_0]\} .\]
We denote by $\dxxbar$\label{k4}  the reduced subscheme of $\hat{D}(C_0;P_1)$ which is the locus of those stable curves, for which the group of automorphisms splits as in notation \ref{n33}.

Below, we will always assume that $\dxxbar$ is non-empty, and thus an open and dense subscheme of $\dpobar$. There are, however, cases where $\dxxbar$  is the empty scheme, compare remark \ref{nonemptysub}.

As before we will often omit $C_0$ from the notation.
\end{notation}

\begin{rk}\em
The case we are interested in below is where $\dxxbar$ contains the boundary point $[C_0]$. Recall that we were not able to say something about this point in chapter \ref{par4}. 

However, there are some cases where one has $\dxxbar \subset \dpo$. By definition, this happens, if and only if the group of automorphisms $\Aut(C_0)$ of the curve $C_0$ does not split as a product $\Aut_{\hat{\Gamma}}(C_0^{++};\calr)\times \Aut(C_0^{--})$. 

For example, consider a decomposition of a stable curve $C_0$ of genus $g=3$ into two curves $C_0^{--}$ and $C_0^{++}$, where $C_0^{--}$  consists of two disjoint $1$-pointed nodal rational curves, and $C_0^{++}$ is the union of two rational lines meeting in two points, and with one marked point on each component. Then, as in remark \ref{rk00328}, the group $\Aut(C_0)$  does not split.

As another example, consider the stable curve $C_0$ of genus $g=3$, where $C_0^{++}$  consists of three disjoint $1$-pointed nodal rational curves, and $C_0^{--}$ is a $3$-pointed rational line. Then clearly  $\Aut(C_0^{--})$ is trivial, so  $\Aut(C_0)$  splits. Hence $[C_0]\in\dxxbar$. 
\end{rk}

\begin{lemma}\label{47}
Let $f : C\rightarrow S$ be a stable curve of genus $g$, with reduced base $S$, such that for all $s\in S$ the fibre $C_s$ has at least $3g-4$ nodes. Suppose that the induced morphism $\theta_f : S \rightarrow \msbar_g$ factors through $\overline{D}(P_1)$. Then there exists an \'etale covering $\{S_\alpha\}_{\alpha\in A}$ of $S$, and for each $\alpha \in A$  a commutative diagram
\[\diagram
C_0^{--} \times S_\alpha \drto_{\pr_2}\rrto|<\hole|<<\ahook && C_\alpha \rrto\dlto&& C \dto^f\\
&S_\alpha \xto[rrr]&&& S 
\enddiagram\]
with 
\[ \rho_i(s) = (R_i,s) \]
for $i=1,\ldots, \mu$ and for all $s \in S_\alpha$, where $\rho_1,\ldots ,\rho_{\mu} : S_\alpha \rightarrow
C_\alpha$ denote the sections of nodes as in lemma \ref{sectnod}, and $R_1,\ldots,R_\mu$ are the marked points on $C_0^{--}$.
\end{lemma}

\proof
The lemma is an immediate consequence of lemma \ref{315}, since $C_0^{--}$ is a subcurve of $C_0^-$, and marked points on $C_0^{--}$ correspond to marked points or nodes of $C_0^-$. 
\ebew

\begin{rk}\em
Note that unlike the situation in lemma \ref{315},  the inverse of this lemma is no longer true. This becomes obvious from example \ref{E57b}, where $C_0^{--} = \emptyset$. However, there is still an equivalent of corollary \ref{K324}. 
\end{rk}

\begin{lemma}\label{L416}
Let $f : C\rightarrow S$ be a stable curve of genus $g$, with reduced base $S$, such that for all $s\in S$ the fibre $C_s$ has at least $3g-4$ nodes. Suppose that the induced morphism $\theta_f : S\rightarrow \msbar_g$ factors through $\hat{D}(C_0;P_1)$. Then there exists a naturally distinguished subcurve $f^{--}:C^{--}\rightarrow S$ of $f:C\rightarrow S$, which in the \'etale topology is locally isomorphic to the trivial curve $\pr_2:C_0^{--}\times S \rightarrow S$ as a $\mu/\Gamma(C_0^{--};\calr)$-pointed curve.
\end{lemma}

\proof
It is clear that locally in the \'etale topology  there are subcurves isomorphic to the trivial curve with fibre $C_0^{--}$. We need to show that up to a reordering of the labels of their sections of marked points these local subcurves fit together to give a global subcurve of $f: C\rightarrow S$. In other words, for any $s\in S$ we need to distinguish in a natural way a subcurve $C_s^{--}$ of the fibre $C_s$, which is isomorphic to $C_0^{--}$. 

If $C_s$ has $3g-3$ nodes, then by assumption $\theta_f(s) = [C_0]$, as this is the only point of $\hat{D}(C_0;P_1)$ representing a curve with the maximal number of  nodes. Hence the inclusion $C_0^{--} \subset C_0$ determines a subcurve $C_s^{--}$ up to an automorphism of $C_0$. Since $C_0^{--}$ is invariant under the action of $\Aut(C_0)$ by construction,   $C_s^{--}$ is uniquely determined by $C_0^{--}$. 

Now suppose that $C_s$ has $3g-4$ nodes. Then there is a unique irreducible component $C_s^{+}$ of $C_s$, which has not the maximal number of nodes when considered as a pointed stable curve. Necessarily, any subcurve $C_s^{--}$ isomorphic to $C_0^{--}$ must be contained in the closure $C_s^-$ of the complement of $C_s^+$ in $C_s$. Note that $C_s^-$ is isomorphic to $C_0^-$ as an $m$-pointed curve. Thus there is an embedding of $C_s^-$ into $C_0$. It follows from proposition \ref{312} that this embedding in uniquely determined up to an automorphism of $C_0$. The subcurve $C_0^{--}$ in $C_0$ is invariant under the action of $\Aut(C_0)$, so its preimage $C_s^{--}$ in $C_s^-$ is independent of the chosen embedding of $C_s^-$ into $C_0$. Thus there is a distinguished subcurve $C_s^{--}$ in $C_s$, isomorphic to $C_0^{--}$. Clearly the order of the marked points on the subcurve $C_s^{--}$ is only determined up to a permutation in ${\Gamma}(C_0^{--};\calr)$. 
\ebew

\begin{notation}\em
Let $\hat{K}(C_0;P_1)$\label{n37} denote the reduction of the  preimage of $\hat{D}(C_0;P_1)$ under the canonical quotient morphism $\pi : \overline{H}_{g,n,0} \rightarrow \msbar_g$. We denote by $\kxxbar$\label{k5} the reduced subscheme corresponding to $\dxxbar$.  

Thus we have 
\[ \kxxbar \subset \hat{K}(C_0;P_1) \subset \overline{K}(C_0;P_1) . \]
By its definition, the morphism $\thxx$ from remark \ref{R513} is a surjection
\[ \thxx : \quad \hxxbar \: \rightarrow  \: \kxxbar .\]
Analogously to corollary \ref{314} we have 
\[  \kxxbar  = \thxx(\hxxbar ) \cdot \PGL(N+1).\]
\end{notation}

\begin{cor}\label{cor_005}
Let $f: C \rightarrow S$ be a stable curve of genus $g$, such that for
all closed points $s\in S$ the fibre $C_s$ has exactly $3g-4$ nodes. Suppose that
the induced morphism $\theta_f : S \rightarrow \msbar_g$ factors through
$\dxxbar$. Let $(E,p,\phi)$ be the triple associated to $f: C \rightarrow S$ as in proposition \ref{28}, consisting of a principal $\PGL(N+1)$-bundle $p:E\rightarrow S$, and a morphism $\phi: E\rightarrow  \overline{H}_{g,n,0}$. Suppose that $\phi$ factors through $\kxxbar$. 
Then there exists a distinguished subcurve
$f^{--} : C^{--} \rightarrow S$ of $f: C \rightarrow S$, which in the \'etale topology is
locally isomorphic to the trivial curve $C_0^{--} \times S \rightarrow
S$ as an $\mu/\Gamma(C_0^{--};\CR)$-pointed stable curve. 
\end{cor}

\proof The proof of this corollary is analogous to that of corollary \ref{K324a}. \ebew

\begin{rk}\em\label{R516}
Let $\msbar^{\,ps}_{g^{++},\mu}$\label{n40} denote the moduli space of $\mu$-pointed prestable curves of genus $g^{++}$. Note that $\msbar^{\,ps}_{g^{++},\mu}$ can be considered as a quotient of $\overline{H}_{g^{++},n,\mu}$ by the action of the group $\PGL(N^{++}+1)$. Let $\mxxbar$\label{n41} denote the reduced image of $\hxxbar$ in this moduli space, and put 
\[ \hgammaxxbar := \hxxbar /  \Gamma(C_0^{--};\calr) .\label{n42}\]
\end{rk}

The following analogue of proposition \ref{312} holds true.

\begin{lemma}\label{L0516a}
The morphism $\thxx : \hxxbar \: \longrightarrow \: \overline{H}_{g,n,0}$ from remark \ref{R513} induces an injective morphism
\[ \hat{\Theta} : \:\: \hgammaxxbar \rightarrow \hat{K}(C_0;P_1) ,\label{n43}\]
and an isomorphism
\[ \mgammaxxbar \:  \cong \: \dxxbar .\]
\end{lemma}

\proof
To see the second part of the claim, note that for each stable curve $C\rightarrow \Spec(k)$ parametrized by a point in $\dxxbar$, there is a distinguished subcurve isomorphic to $C_0^{--}$ as a scheme, see lemma \ref{L416} above. If $[C_0]\in \dxxbar$, then, unlike the situation in proposition \ref{312}, this is still true for the boundary point $[C_0]$. The order of the labels of the marked points of this subcurve is only determined up to a permutation in $\Gamma(C_0^{--};\calr)$. The closure of the complement of this subcurve determines thus a $\mu/\Gamma(C_0^{--};\calr)$-pointed stable curve represented by a unique point in $\mgammaxxbar$. 

The arguments for the first part are analogous.
\ebew

\begin{rk}\em\label{R517}
Glueing with the curve $C_0^{--}$ along the marked points defines a morphism 
\[ \msbar^{\,ps}_{g^{++},\mu} \longrightarrow \msbar_g .\]
By definition, $\mxxbar$ is the reduced preimage of $\dxxbar$ under this morphism. In particular $\mxxbar$ is contained in the one-dimensional part of the boundary stratification of $\msbar^{\,ps}_{g^{++},\mu}$. Note that the closure of the curve $\mxxbar$ is in general  neither smooth nor irreducible.
\end{rk}

\begin{lemma}\label{L518}
The irreducible components of the closure of the subscheme $\mxxbar$ inside $\msbar_{g^{++},\mu}^{\,ps}$ are smooth curves. Intersections of  two irreducible components are transversal.
\end{lemma}

\proof
It follows from the standard description  of the moduli space $\msbar_{g^{++},\mu}^{\,ps}$ that the closure of  $\mxxbar$ is given locally as the union of lines, meeting pairwise transversally.   To prove the lemma it therefore suffices to show that there are no nodes on a fixed  irreducible component of the closure of $\mxxbar$. 

Let $C\rightarrow\Spec(k)$ be a stable curve of genus $g^{++}$, which has the maximal number of nodes. Let $\calr$ be the set of nodes on  $C$, and let $P_1,P_2\in\calr$ be two nodes. Suppose that the two irreducible components of the one-dimensional part of the boundary stratum of $\msbar_{g^{++},\mu}$, which are distinguished by such deformations of $C$, which preserve all nodes except $P_1$ or $P_2$, respectively, are the same. 
To prove that this irreducible component of the closure of $\mxxbar$ is smooth, we need to show that there exists an automophism $\rho\in\Aut(C)$, such that $\rho(P_1) = P_2$. 

If $m=1$, then a curve representing a general point of $\mxxbar$ is given by the following picture.

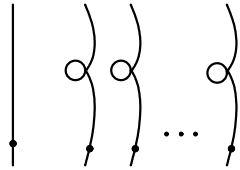
\begin{center}
\setlength{\unitlength}{0.00016667in}
\begingroup\makeatletter\ifx\SetFigFont\undefined%
\gdef\SetFigFont#1#2#3#4#5{%
  \reset@font\fontsize{#1}{#2pt}%
  \fontfamily{#3}\fontseries{#4}\fontshape{#5}%
  \selectfont}%
\fi\endgroup%
{\renewcommand{\dashlinestretch}{30}

}
\end{center}

In words, a general curve is a disjoint union of one $1$-pointed smooth curve of genus one, and finitely many $1$-pointed nodal rational curves. Recall that we assumed throughout that $g \ge 3$. In this case the claim is trivial, since the label of the marked point of the smooth component is fixed. 

If $m=2$, then there are two possibilities. Suppose first that $C_0^+$ is a union of two smooth rational curves, meeting in two points. Then a curve representing a general point of $\mxxbar$ is given by the disjoint union of exactly one special ``bracelet'', and a finite (maybe zero) number of identical general bracelets. Here a general bracelet is constructed by glueing a finite number of copies of $C_0^+$ together in their marked points, as in the picture below. 

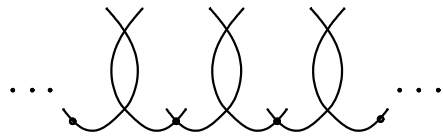
\begin{center}
\setlength{\unitlength}{0.00016667in}
\begingroup\makeatletter\ifx\SetFigFont\undefined%
\gdef\SetFigFont#1#2#3#4#5{%
  \reset@font\fontsize{#1}{#2pt}%
  \fontfamily{#3}\fontseries{#4}\fontshape{#5}%
  \selectfont}%
\fi\endgroup%
{\renewcommand{\dashlinestretch}{30}

}
\end{center}

Note that a bracelet may be open or closed. To obtain the  one special bracelet one has to replace one copy of $C_0^+$ with a bracelet by a $2$-pointed nodal curve of genus $1$, as in the next picture.   

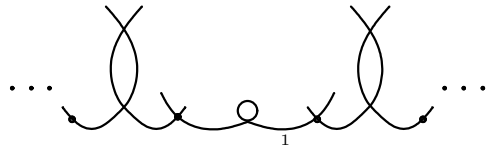
\begin{center}
\setlength{\unitlength}{0.00016667in}
\begingroup\makeatletter\ifx\SetFigFont\undefined%
\gdef\SetFigFont#1#2#3#4#5{%
  \reset@font\fontsize{#1}{#2pt}%
  \fontfamily{#3}\fontseries{#4}\fontshape{#5}%
  \selectfont}%
\fi\endgroup%
{\renewcommand{\dashlinestretch}{30}
%
}
\end{center}

The number ``1'' next to an irreducible component indicates that this component has geometric genus $1$, whereas all other components are rational. 

Curves of this type can only degenerate in two ways to curves with the maximal number of nodes, and in both cases the claim holds. 

Suppose now that $C_0^+$ consists of one smooth rational curve with three marked points on it, where in one of them a nodal rational curve is attached. If a curve representing a general point of $\mxxbar$ is given by bracelets as before, then we are done by the same reasoning as above. The schematic picture of a general bracelet in this case is as follows.  

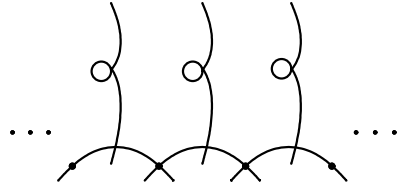
\begin{center}
\setlength{\unitlength}{0.00016667in}
\begingroup\makeatletter\ifx\SetFigFont\undefined%
\gdef\SetFigFont#1#2#3#4#5{%
  \reset@font\fontsize{#1}{#2pt}%
  \fontfamily{#3}\fontseries{#4}\fontshape{#5}%
  \selectfont}%
\fi\endgroup%
{\renewcommand{\dashlinestretch}{30}

}
\end{center}

Again, the one special bracelet is obtained by replacing one copy of $C_0^+$ by a $2$-pointed nodal curve of genus $1$, as shown below.

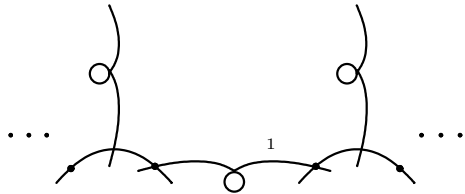
\begin{center}
\setlength{\unitlength}{0.00016667in}
\begingroup\makeatletter\ifx\SetFigFont\undefined%
\gdef\SetFigFont#1#2#3#4#5{%
  \reset@font\fontsize{#1}{#2pt}%
  \fontfamily{#3}\fontseries{#4}\fontshape{#5}%
  \selectfont}%
\fi\endgroup%
{\renewcommand{\dashlinestretch}{30}
%
}
\end{center}

If the connected components of a general curve are not bracelets, and $g\ge 4 $, then a schematic picture of the curve is the following.

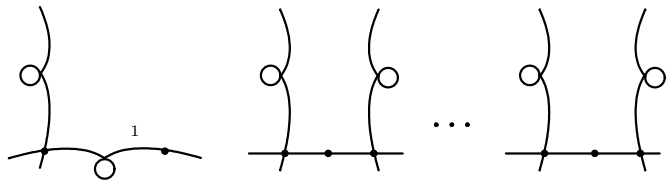
\begin{center}
\setlength{\unitlength}{0.00016667in}
\begingroup\makeatletter\ifx\SetFigFont\undefined%
\gdef\SetFigFont#1#2#3#4#5{%
  \reset@font\fontsize{#1}{#2pt}%
  \fontfamily{#3}\fontseries{#4}\fontshape{#5}%
  \selectfont}%
\fi\endgroup%
{\renewcommand{\dashlinestretch}{30}
%
}
\end{center}

Again, the claim follows by looking at the only two possible degenerations. 

Finally, if $g=3$, then the general curve in $\mxxbar$  may be of the type pictured below.

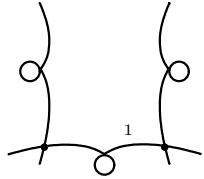
\begin{center}
\setlength{\unitlength}{0.00016667in}
\begingroup\makeatletter\ifx\SetFigFont\undefined%
\gdef\SetFigFont#1#2#3#4#5{%
  \reset@font\fontsize{#1}{#2pt}%
  \fontfamily{#3}\fontseries{#4}\fontshape{#5}%
  \selectfont}%
\fi\endgroup%
{\renewcommand{\dashlinestretch}{30}
%
}
\end{center}

In this case the claim is clearly true.

Now assume $m=4$. Let $C_1^+$ denote the subcurve of $C$ which is the union of the two lines meeting in $P_1$, and let $C_1^-$ denote the closure of its complement in $C$. Analogously, we define $C_2^+$ and $C_2^-$ for the node $P_2$. Thus there are two  decompositions 
\[ C = C_1^+ + C_1^- = C_2^+ + C_2^- ,\]
where $C_1^- \cong C_2^-$ holds by the assumption that the deformations of $P_1$ and $P_2$ distinguish the same curve in the boundary stratum of $\msbar_{g^{++},\mu}^{\,ps}$. Call this isomorphism $\lambda : C_1^- \rightarrow C_2^-$. By lemma \ref{L326a} there exists an isomorphism $\gamma:C_1^+\rightarrow C_2^+$, and an automorphism $\tau$ of $C_2^-$, which glue together to give an isomorphism 
\[ \rho=(\gamma,\tau\circ\lambda) : \:\:\:\: C_1^+ + C_1^- \longrightarrow C_2^+ + C_2^- ,\]
i.e. an automorphism $\rho \in \Aut(C)$, which clearly satisfies $\rho(P_1) = P_2$.
\ebew

\begin{ex}\em
Consider once again the curve $C_0$ of genus $g=3$, with node $P_1$, as pictured below.

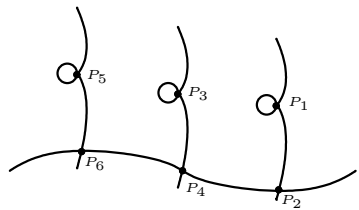
\begin{center}
\setlength{\unitlength}{0.00016667in}
\begingroup\makeatletter\ifx\SetFigFont\undefined%
\gdef\SetFigFont#1#2#3#4#5{%
  \reset@font\fontsize{#1}{#2pt}%
  \fontfamily{#3}\fontseries{#4}\fontshape{#5}%
  \selectfont}%
\fi\endgroup%
{\renewcommand{\dashlinestretch}{30}

}
\end{center}

Then $\mxxbar = \hat{M}_{1,3}^{\times}$ consists of three lines, each of them isomorphic to $\msbar_{1,1}$, meeting in the point $[C_0^{++}]$. Note that the quotient $\mgammaxxbar=  \hat{M}_{1,3/\Sigma_3}^\times$  is isomorphic to $\msbar_{1,1}$ and smooth.
\end{ex}

\begin{rk}\em\label{R0520}
By construction, the finite morphism $\mxxbar\rightarrow \dxxbar$ factors through an isomorphism 
$ \mgammaxxbar \rightarrow \dxxbar$, see lemma \ref{L0516a}. Since $ \overline{D}(C_0;P_1)$ is smooth and irreducible, the inverse morphism on the open and dense  subscheme extends to all of $ \overline{D}(C_0;P_1)$. Hence  there is in fact an isomorphism of schemes between $\overline{D}(C_0;P_1)$ and the closure of $\mgammaxxbar$ inside $  \msbar_{g^{++},\mu/\hat{\Gamma}}$.
 In particular, the closure of  $\mgammaxxbar$ is a smooth irreducible curve. 
\end{rk}

\begin{defi}\em
We denote by \[\FDPxx \label{n45}\] the reduction of the preimage substack of the subscheme $\dxxbar$ under the canonical morphism $\cmgbar \rightarrow \msbar_g$.
\end{defi}

\begin{rk}\em\label{R0526}
As in proposition \ref{MainProp}, one can construct a morphism of stacks
\[ \left[\hgammaxxbar / \PGL(N^{++}+1)\times\Aut(C_0^{--}) \right] \: \rightarrow  \: \FDPxx ,
\]
where $\Aut(C_0^{--})$ acts trivially on $\hgammaxxbar$. 

Let $\cmgammaxxbar$\label{n46} denote the stack associated to $\mu/\Gamma(C_0^{--};\CR)$-pointed prestable curves of genus $g^{++}$, which are parametrized by points in $\mgammaxxbar$. Recall that we assume throughout that all connected components of a prestable curve are stable. Proposition \ref{P217} generalizes to the case of prestable curves, so there is an isomorphism
\[ \cmgammaxxbar  \cong \left[\hgammaxxbar / \PGL(N^{++}+1) \right] .\]
Thus there is in particular a morphism of stacks 
\[ \hat{\Omega}: \:\: \cmgammaxxbar  \rightarrow \FDPxx ,\]
which will be shown to be finite and surjective in lemma \ref{L6107} below.
\end{rk}

\begin{thm}\label{49}
There is an isomorphism of Deligne-Mumford stacks
\[ \left[\hgammaxxbar / \PGL(N^{++}+1) \times\Aut(C_0^{--}) \right] \: \cong \: \FDPxx ,
\]
where $\Aut(C_0^{--})$ acts trivially on $\hgammaxxbar$. 
\end{thm}

\proof
The proof is analogous to the proof of theorem \ref{MainThm}.
Therefore we will only indicate the main steps here. 
In the proof we use the abbreviations $A^{--}:= \Aut(C_0^{--})$, $P^{++} :=  \PGL(N^{++}+1)$ and $\hat{\Gamma}:= \Gamma(C_0^{--};\calr)$. 

$(i)$ A morphism $\hat{\Lambda}$ from $
[\hgammaxxbar / P^{++} \times A^{--}]$ to $\FDPxx$ is constructed as in proposition \ref{MainProp}. Note that there is again an injective group homomorphism
\[ P^{++}\times A^{--} \rightarrow \PGL(N+1) . \]
By remark \ref{R513}, 
the morphism induced by $\thxx$ on the scheme $\hgammaxxbar$ 
\[ \hat{\Theta} : \:\: \hgammaxxbar \rightarrow \overline{H}_{g,n,0} \]
is equivariant with respect to the action of $ P^{++}\times A^{--}$, where $A^{--}$ acts trivially on $\hgammaxxbar$. Using this, one can define for any scheme $S$ and any object $(E',p',\phi') \in [\hgammaxxbar / P^{++} \times A^{--}](S)$ an object $(E,p,\phi)\in [\kxxbar/\PGL(N+1)] \cong \FDPxx(S)$. For this latter isomorphism compare lemma \ref{140}. 
Here $p: E\rightarrow S$ is the extension of the principal $P^{++}\times A^{--}$-bundle $p':E'\rightarrow S$ by
\[ E := (E' \times \PGL(N+1)) / (P^{++}\times A^{--}) .\]
The morphism $\phi: E \rightarrow \kxxbar$ is constructed from $\phi': E' \rightarrow  \hgammaxxbar$ and $\hat{\Theta}$, as in the proof of proposition \ref{MainProp}.  

$(ii)$ The main difficulty in constructing an inverse morphism in the proof of theorem \ref{MainThm} was, geometrically speaking, to distinguish an appropriate subcurve of a given stable curve. Here we are in a situation which is modified in such a way that we are looking for an invariant subcurve, which will resolve this difficulty over the additional boundary point.

Consider an object $(E,p,\phi)\in \FDPxx(S)$, for some scheme $S$. So $p:E\rightarrow S$ is a principal $\PGL(N+1)$-bundle, and $\phi:E\rightarrow \kxxbar$ is a $\PGL(N+1)$-equivariant morphism. This triple represents a stable curve $f_S: {\cal C}_S\rightarrow S$, such that the induced morphism $\overline{\theta}_{f_S}: S \rightarrow \msbar_g$ factors through $\dxxbar$. By corollary \ref{cor_005}, there exists an \'etale cover $\{S_\alpha\}_{\alpha\in A}$ of $S$, such that the restriction of $f_S: {\cal C}_S\rightarrow S$ to each $S_\alpha$ contains a subcurve $f_{S,\alpha}^{--} : {\cal C}_{S,\alpha}^{--}\rightarrow S_\alpha$ which is isomorphic to the trivial curve $\pr_2: C_0^{--}\times S_\alpha\rightarrow S_\alpha$. Consider the restrictions $p_\alpha: E_\alpha\rightarrow S_\alpha$ of 
the principal $\PGL(N+1)$-bundle $p:E\rightarrow S$. As in the proof of theorem \ref{MainThm},  the subcurves $f_{S,\alpha}^{--} : {\cal C}_{S,\alpha}^{--}\rightarrow S_\alpha$ define subbundles $p_\alpha^{--}: E_\alpha^{--}\rightarrow S_\alpha$ of $p_\alpha: E_\alpha\rightarrow S_\alpha$, which are principal $ P^{++}\times A^{--}$-bundles. By the choice of the cover $\{S_\alpha\}_{\alpha\in A}$  according to corollary \ref{cor_005},  the local $\mu$-pointed prestable subcurves $f_{S,\alpha}^{--} : {\cal C}_{S,\alpha}^{--}\rightarrow S_\alpha$ are restrictions of a global $\mu/\hat{\Gamma}$-pointed prestable subcurve $f_S^{--}: {\cal C}_S^{--} \rightarrow S$ of $f_S: {\cal C}_S\rightarrow S$. Hence the local 
principal $ P^{++}\times A^{--}$-bundles $p_\alpha^{--}: E_\alpha^{--}\rightarrow S_\alpha$ glue together, up to the action of $\Aut_{\Sigma_\mu}(C_0^{--})$, due to the non-uniqueness of the labels of the $\mu$ sections of marked points of the local curves. Compare also remark \ref{323a}.  As in the proof of theorem \ref{MainThm} we denote this global object by
\[ p^{--} : \:\: E^{--}/\Aut_{\Sigma_\mu}(C_0^{--})\rightarrow S .\]
For the same reasons as in the proof of theorem \ref{MainThm} there exists a section from $S$ to $E^{--}/\Aut_{\Sigma_\mu}(C_0^{--})$, which by composition defines a section
\[ \hat{\sigma} : S \rightarrow E/(  P^{++}\times A^{--}) .\]
This section now defines a reduction of the principal $\PGL(N+1)$-bundle $p:E\rightarrow S$ to a principal $ P^{++}\times A^{--}$-subbundle $p^{++}: E^{++}\rightarrow S$. 

Using corollary \ref{cor_005} again, and the universal property of the Hilbert scheme, one can construct local morphisms 
\[ \phi_\beta^{++} : E^{++}_\beta \rightarrow  \hxxbar \]
on a covering $\{p_\beta^{++} : E_\beta^{++} \rightarrow S_\beta\}_{\beta\in B}$ of $p:E\rightarrow S$. This construction is analogous to the one in the proof of theorem \ref{MainThm}. Again, the morphisms glue together, up to the action of $\Gamma^{--}$. So one finally obtains a morphism
\[ \phi^{++} : E^{++} \rightarrow \hgammaxxbar, \]
which is $P^{++}\times A^{--}$-equivariant. Thus an object $(E^{++},p^{++},\phi^{++})$ of the quotient stack $   [\hgammaxxbar /P^{++}\times A^{--}](S)$ has been constructed. 

The construction of the functor from $\FDPxx$ to $[\hgammaxxbar/P^{++}\times A^{--}]$ on morphisms is straightforward. Note that this uses the injectivity of the morphism $ \hat{\Theta} : \hgammaxxbar \rightarrow \kxxbar$ provided by lemma \ref{L0516a}. 

The remainder of the proof is completely analogous to the proof of theorem \ref{MainThm}.
\ebew

\begin{notation}\em
Note that the subschemes $\dxbar$ and $\dxxbar$ are both defined as loci of curves in $\dpobar$ with splitting groups of automorphisms, but where the splitting is with respect to different subgroups. Thus their intersections with  $\dpo$, which are  $\dx$ and $\dxxbar\cap \dpo$,  need not be the same. 

We define their mutual intersection by
\[ \doo \: \: := \: \: \dx \: \cap \: \dxxbar \label{cc1}
\]
as a reduced subscheme of $\dpo$. If $\dxxbar$ is non-empty, as we assume throughout, then this is an open and dense subscheme of $\dpobar$.  

The reduced substack of $\FDPopen$ corresponding to $\doo$ shall be denoted by $\cdoo$\label{cc2}. 

We define by $\koo$\label{cc3} the reduction of the preimage of $\doo$ inside $H_{g,n,0}$. Via the morphisms $\thx$ and $\thxx$, we define reduced subschemes  
\[\hoo \subset H_{g^+,n,m}   \text{ and } \hoohat \subset {H}_{g^{++},n,\mu} \label{cc4b}\]
respectively, together with their quotients
\[ \hoogamma \subset H_{g^+,n,\mgamma} \label{cc5a} \text{ and } \hoogammahat \subset {H}_{g^{++},n,\mu/\hat{\Gamma}} \, .\]
\end{notation}

\begin{cor}\label{C49}
There is an isomorphism of Deligne-Mumford stacks
\[ \left[\hoogammahat / \PGL(N^{++}+1)\times\Aut(C_0^{--}) \right] \: \cong \: \cdoo ,
\] 
where $\Aut(C_0^{--})$ acts trivially on $\hoogammahat$.
\end{cor}

\proof
This is an immediate consequence of theorem \ref{49}.
\ebew

\begin{rk}\em
In view if corollary \ref{C49} one may think of our considerations in this section as an invariant reformulation of theorem \ref{MainThm}. The statement of the result  no longer depends on a particular choice of a node $P_1\in C_0$, but only on the equivalence class of this node.  This may seem conceptually more appealing, since after all there is no a priory ordering of the nodes on $C_0$. However, replacing the subcurve $C_0^+$ by the larger invariant subcurve $C_0^{++}$ reduces the usefulness of the theorem in applications. In worst-case situations, such as described in example \ref{E57b}, one finds $C_0^{++} = C_0$, so that $\Aut(C_0^{--})$ is the trivial group, just as ${\Gamma}(C_0^{--};\calr)$ is trivial. The morphism $\thxx$ in this case is just the identity on $\overline{H}_{g,n,0}$, and hence $\hoogammahat = \kx $. Thus corollary \ref{C49} reduces to the statement
\[ [ \kx / \PGL(N+1) ] \cong \cdx ,\]
which is a tautology in view of lemma \ref{140}.
\end{rk}

%
%

\chapter{Glueing}

\begin{rk}\em\label{R51}
As in the previous chapters, let $\overline{D}\subset \msbar_g^{3g-4)}$ be an irreducible component of the one-dimensional boundary stratum of $\msbar_g$. Then $ \overline{D} = \overline{D}(C_0;P_1)$ for some stable curve $f_0:C_0\rightarrow \Spec(k)$ of genus $g$ with $3g-3$ nodes. Let $P_1,\ldots,P_{3g-3}$ be an enumeration of the nodes of $C_0$. 

In chapter \ref{par4} we gave a description of the open substack $\FDPopen$ of the stack $\overline{\cal D}= \FDP$ by refering to  $m/\Gamma(C_0^-;\calq)$-pointed stable curves of genus $g^+$. In chapter \ref{par5} we described the substack of $\FDP$, corresponding to the partial compactification $\hat{D}(C_0;P_1)$ of the scheme $D(C_0;P_1)$, which in turn was obtained by adding the point $[C_0]$. Here we used $\mu/\Gamma(C_0^{--};\CR)$-pointed prestable curves of genus $g^{++}$. More precisely, in both cases we had to restrict ourselves  to substacks $\cdx$ and $\FDPxx$, respectively, corresponding to curves with splitting automorphism groups. 

The presentation of $\overline{D}$ as $\overline{D}(C_0;P_1)$ is not unique. Indeed, if $D$ denotes the locus of curves with less than the maximal number of nodes in $\overline{D}$, then the complement $\overline{D} \smallsetminus D$ consists of finitely many points $[C_1],\ldots,[C_k]$, each representing a stable curve over $\Spec(k)$ with $3g-3$ nodes.  For $i=1,\ldots,k$ there is a node $P_{i,1} \in C_i$ such that
\[ \overline{D}=\overline{D}(C_i;P_{i,1}) \,  .\]
The finitely many partial compactifications $\hat{D}(C_i;P_{i,1})$ cover all of $\overline{D}$. Thus the stack $\overline{\cal D}$ is covered by finitely many open substacks 
\[ {\hat{\cal D}(C_i;P_{i,1})} \: \: \subset \: \: \FDP \label{lb1}\]
which are the reductions of the preimages of $\hat{D}(C_i;P_{i,1})$  for $i=1,\ldots,k$.

However, while $g^+$ and $m$ depend only on the component $\overline{D}$, this is no longer true for $g^{++}$ and $\mu$. Indeed, for $i=1,\ldots,k$, let
\[ C_i = C_i^+ + C_i^- \]
denote the decomposition of $C_i$ defined by the node $P_{i,1}$, such that $C_0^+$ consists of the union of all irreducible components containing $P_{i,1}$. 
Then one has isomorphisms 
\[ C_1^- \cong C_2^- \cong \ldots \cong C_k^- \]
of $m/\Gamma(C_0^-;\CQ)$-pointed  curves. Therefore the number $m$ of marked points and the genus of the complementary subcurve are fixed. Consider now  for $i=1,\ldots,k$ the $\Aut(C_i)$-invariant decomposition
\[ C_i = C_i^{++ }+ C_i^{-- },\]
where $C_i^{++}$ is given as the reduced orbit of $C_i^+$ under the action of $\Aut(C_i)$. In general, neither the genus nor the number of points of $C_i^{--}$ will be the same for different $1\le i\le k$. 

Still, the different partial compactifications $\hat{\cal D}(C_i;P_{i,1})$ for $i=1,\ldots,k$ glue together, since their restrictions to ${{\cal D}(C_i;P_{i,1})}= {{\cal D}} $ are the same. To understand these restrictions better, we need to compare the stacks $\cmgammadetailxbar$ and $\cmgammadetailxxbar$. Recall that we have isomorphisms 
\[\cmgammadetailxbar \: \cong \: [\hgammadetailxbar /\PGL(N^+ +1)]
,\] 
and 
\[ \cmgammadetailxxbar \: \cong \: [\hgammadetailxxbar /\PGL(N^{++} +1)]\]
as a consequence of proposition \ref{P217}. 
In particular, there is a morphism between the open substacks 
\[ 
\begin{array}{c}\left[\hoogammadetail / \PGL(N^++1)\times\Aut(C_0^-)\right] \\
\downarrow \\
 \left[\hoogammahatdetail/ \PGL(N^{++}+1)\times\Aut(C_0^{--}))\right]. 
\end{array}\]
Recall that we assume throughout that we are in the case where the scheme $\dxxbar$ is not the empty scheme, compare notation \ref{n36}, and thus both $\hoogammadetail$ and $\hoogammahatdetail$ are non-empty. 

From theorem \ref{MainThm} and corollary \ref{C49} it follows that both of the above quotient stacks are isomorphic to the dense open substack $\cdoo$ of  $\FDPopen$. Without this knowledge, it is not a priory obvious that such a morphism exist, as  that the inclusions
\[\PGL(N^+ +1) \subset \PGL(N^{++}+1) \]
and
\[ \Aut(C_0^-) \supset \Aut(C_0^{--}) \]
seem to be acting in opposite directions.  To understand the relation of these quotient stacks directly, without refering to $\cdoo$, we need again more notation. 
\end{rk}

\begin{notation}\em
For our fixed stable curve $f_0:C_0\rightarrow \Spec(k)$ we defined above two decompositions
\[ C_0 = C_0^+ + C_0^- = C_0^{++} +
 C_0^{--} \]
into prestable subcurves with marked points. By definition we have $C_0^+ \subset C_0^{++}$ and $C_0^{--} \subset C_0^-$. We define a subcurve $C_0^{+-}$\label{n47} of $C_0$ as the closure of the complement of $C_0^+$ in $C_0^{++}$, considered as a prestable curve with $\mu^{+-}$ marked points. 

Let ${\cal Q}=\{Q_1,\ldots,Q_m\}$ denote the set of marked points of $C_0^+$. Note that $\cal Q$ is equal to the set of marked points on $C_0^-$. There is a disjoint decomposition
\[ {\cal Q} = \CQ^{+-} \cup \CQ^{--}, \label{n48}\]
where $ \CQ^{+-}$ denotes the subset of those marked points of $C_0^+$, which are glueing points with $C_0^{+-}$, or equivalently, those marked points which are nodes of $C_0^{++}$. The set $\CQ^{--}$ consists of those marked points of $C_0^+$, which are also marked points of $C_0^{++}$, or equivalently, $\CQ^{--} = C_0^+ \cap C_0^{--}$. Similarly, we define a decomposition of the set $\CR=\{R_1,\ldots,R_\mu\}$ of marked points of $C_0^{--}$ into subsets $\CR = \CR^{+-} \cup \CR^{+}$. Here $\CR^{+-}$ is the subset of those marked points, which are glueing points with $C_0^{+-}$, and $\CR^+$ is the set of those points, which remain marked points when considered as marked points of  $C_0^-$. Obviously, the latter are just the glueing points of $C_0^{--}$ with $C_0^+$, so that $\CR^+ = \CQ^{--}$. Hence we have a disjoint decomposition 
\[ \CR = \CR^{+-} \cup  \CQ^{--} . \label{n49}\]
In particular, the set $\CR^{+-}$ is the set of those marked points of $C_0^{--}$ which are nodes of $C_0^-$. Clearly, $\CQ^{+-} \cup \CQ^{--}\cup\CR^{+-}$ is a distribution of the set of all nodes of $C_0$.
Schematically, we have the following picture of $C_0$. 

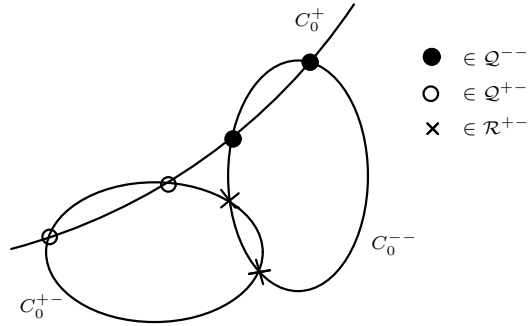
\begin{center}
\setlength{\unitlength}{0.00025000in}
\begingroup\makeatletter\ifx\SetFigFont\undefined%
\gdef\SetFigFont#1#2#3#4#5{%
  \reset@font\fontsize{#1}{#2pt}%
  \fontfamily{#3}\fontseries{#4}\fontshape{#5}%
  \selectfont}%
\fi\endgroup%
{\renewcommand{\dashlinestretch}{30}
\begin{picture}(10651,6872)(0,-10)
\thicklines
\put(-3256.354,13945.096){\arc{25577.401}{0.5896}{1.3115}}
\put(3096,1521){\ellipse{4650}{3000}}
\put(6173,3159){\ellipse{3000}{4950}}
\put(848,1847){\ellipse{310}{310}}
\put(3398,2972){\ellipse{310}{310}}
\texture{55888888 88555555 5522a222 a2555555 55888888 88555555 552a2a2a 2a555555 
	55888888 88555555 55a222a2 22555555 55888888 88555555 552a2a2a 2a555555 
	55888888 88555555 5522a222 a2555555 55888888 88555555 552a2a2a 2a555555 
	55888888 88555555 55a222a2 22555555 55888888 88555555 552a2a2a 2a555555 }
\put(4785,3947){\blacken\ellipse{320}{320}}
\put(4785,3947){\ellipse{320}{320}}
\put(6435,5597){\blacken\ellipse{320}{320}}
\put(6435,5597){\ellipse{320}{320}}
\put(9041,5709){\blacken\ellipse{336}{336}}
\put(9041,5709){\ellipse{336}{336}}
\put(9023,4922){\ellipse{310}{310}}
\path(9173,4284)(8948,4059)
\path(9173,4284)(8948,4059)
\path(8948,4284)(9173,4059)
\path(8948,4284)(9173,4059)
\path(4523,2859)(4823,2484)(4898,2484)
\path(4523,2484)(4898,2784)
\path(5123,1059)(5573,1134)
\path(5123,1059)(5573,1134)
\path(5273,1359)(5423,834)
\put(198,234){\makebox(0,0)[lb]{\smash{{{\SetFigFont{7}{8.4}{\rmdefault}{\mddefault}{\updefault}$C_0^{+-}$}}}}}
\put(6098,6384){\makebox(0,0)[lb]{\smash{{{\SetFigFont{7}{8.4}{\rmdefault}{\mddefault}{\updefault}$C_0^+$}}}}}
\put(7748,1584){\makebox(0,0)[lb]{\smash{{{\SetFigFont{7}{8.4}{\rmdefault}{\mddefault}{\updefault}$C_0^{--}$}}}}}
\put(9698,5559){\makebox(0,0)[lb]{\smash{{{\SetFigFont{7}{8.4}{\rmdefault}{\mddefault}{\updefault}$\in{\cal Q}^{--}$}}}}}
\put(9698,4059){\makebox(0,0)[lb]{\smash{{{\SetFigFont{7}{8.4}{\rmdefault}{\mddefault}{\updefault}$\in{\cal R}^{+-}$}}}}}
\put(9698,4809){\makebox(0,0)[lb]{\smash{{{\SetFigFont{7}{8.4}{\rmdefault}{\mddefault}{\updefault}$\in{\cal Q}^{+-}$}}}}}
\end{picture}
}
\end{center}

We furthermore denote by 
\[  \CN := \CQ^{+-} \cup \CR^{+-} \]
the set of all of the $\mu^{+-} := \#{\cal N}$  marked points of $C_0^{+-}$.
\end{notation}

\begin{notation}\em
Let $f: C \rightarrow \Spec(k)$ be a stable curve with marked points. Let $Y_1,\ldots,Y_r$ and $Z_1,\ldots,Z_s$ be finite subsets of $C$. Then we denote by
\[ \Aut(C;Y_1,\ldots,Y_r,\underline{Z}_1,\ldots,\underline{Z}_s) \label{n50}\]
the group of those automorphisms of $C$ as a scheme, which stabilize the sets $Y_1,\ldots,Y_r$ as subsets, and the sets  $Z_1,\ldots,Z_s$ pointwise. 

If $\lambda_i := \# Y_i$ denotes the number of elements of $Y_i$, $i=1,\ldots,r$,  then there is a natural homomorphism to the group of permutations
\[  \Aut(C;Y_1,\ldots,Y_r,\underline{Z}_1,\ldots,\underline{Z}_s) \: \rightarrow \: \Sigma_{\lambda_1} \times \ldots \Sigma_{\lambda_r}. \]
We denote by 
 $\Gamma(C;Y_1,\ldots,Y_r,\underline{Z}_1,\ldots,\underline{Z}_s)$\label{n51} the image of this homomorphism, or write  simply $ \Gamma(C;Y_1,\ldots,Y_r)$.  Note that this notation generalizes  our earlier conventions. 
\end{notation}

\begin{lemma}\label{Gbeta}
The action of $\Aut(C_0^-;\CQ)$ stabilizes the subcurve $C_0^{--}$ as a subscheme with distinguished, but not marked points. In particular, there is a natural homomorphism of groups
\[  \Aut(C_0^-;\CQ) \rightarrow \Aut(C_0^{--};\CR) .\]
\end{lemma}

\proof
The claim of the lemma is equivalent to the claim that $C_0^{+-}$ is stabilized by the action of $\Aut(C_0^-;\CQ)$. Recall that 
\[ C_0^{+-} = C_0^- \cap C_0^{++} ,\]
and $C_0^{++}$ is the union of all translates $\rho(C_0^+)$ for $\rho \in \Aut(C_0)$. Hence it suffices to show that for all $\phi\in \Aut(C_0^-;\CQ)$ and all  $\rho\in \Aut(C_0)$ holds that $\phi(\rho(C_0^+)\cap C_0^-) \subset C_0^{++}$.  

If $m=1$ or $m=2$ then any automorphism $\phi\in\Aut(C_0^-;\CQ)$ extends to an automorphism of $C_0$. In these cases we may assume that $\phi\circ \rho \in \Aut(C_0)$, and thus clearly $\phi(\rho(C_0^+))\subset C_0^{++}$. So from now on let $m=4$. 

Suppose first that $\rho(C_0^+) \subset C_0^-$. Denote the closure of the complement of $\rho(C_0^+)$ in $C_0^-$ by $C_0^\vee$. We have a decomposition
\[ C_0^- = \phi(C_0^- )= \phi(\rho(C_0^+))+ \phi(C_0^\vee),\]
and thus, since $C_0 = \rho(C_0)$,
\[ \begin{array}{rcl}
C_0 &=& \phi(\rho(C_0^+)) + \phi(C_0^\vee) + C_0^+ \\ 
&=& \rho^{-1}(C_0^+) + \rho^{-1}(C_0^-) .
\end{array}\]
We claim that 
\[ \rho^{-1}(C_0^-) \cong \phi(C_0^-) + C_0^+ \]
as $m$-pointed prestable curves. Once proven, this claim implies by corollary \ref{C333} that there is an automorphism $\gamma\in\Aut(C_0)$, such that $\gamma(\rho^{-1}(C_0^+)) = \phi(\rho(C_0^+))$. Hence $\phi(\rho(C_0^+))\subset C_0^{++}$, and we are done. 

Let us prove the claim. Since $\rho\in\Aut(C_0)$ and $\phi\in\Aut(C_0^-;\CQ)$, one has
\[ \rho^{-1}(C_0^-) \cong C_0^- =\rho(C_0^+) + C_0^\vee \cong \phi(\rho(C_0^+)) + \phi(C_0^\vee) \]
as $m$-pointed curves. Furthermore,
\[ \phi(\rho(C_0^-)) \cong \rho(C_0^+)\cong C_0^+ \]
as  $4$-pointed curves. Thus
\[ \phi(\rho(C_0^+)) + \phi(C_0^\vee) \cong C_0^+ + \phi(C_0^\vee)\]
as $m$-pointed curves, and  the claim follows. 

Suppose now that $\rho(C_0^+) \not\subset C_0^-$. Without loss of generality we may assume that  $\rho(C_0^+) \neq C_0^+$, since otherwise there is nothing to prove. In other words, $\rho(C_0^+)\cap C_0^-$ consist of one line with three distinguished points.  There are at least one and at most three marked points of $C_0^-$ on $\rho(C_0^+)\cap C_0^-$. If there are 3 marked points, then $\rho(C_0^+)\cap C_0^-$ is disjoint from the remainig part of $C_0^-$, and any automorphism of $C_0^-$ stabilizing the set $\CQ$ of marked points  stabilizes  $\rho(C_0^+)\cap C_0^-$. In particular, $\phi(\rho(C_0^+)\cap C_0^-) \subset C_0^{++}$.

If  $\rho(C_0^+)\cap C_0^-$ contains two marked points, then these two points lie necessarily on different irreducible components of $C_0^+$. There are two cases. If  $\rho(C_0^+)\cap C_0^-$ is the only irreducible component of $C_0^-$ containig two marked points, then any $\phi\in\Aut(C_0^-;\CQ)$ must stabilize this component, and we are done. Otherwise there are exactly two irreducible components containing two marked points. If $\phi$ does not stabilize $\rho(C_0^+)\cap C_0^-$, then it interchanges these two components.  
Since $C_0$ is a stable curve with the maximal number of nodes, there exists an automorphism $\gamma\in\Aut(C_0^-;\CQ)$, stabilizing  $\rho(C_0^+)\cap C_0^-$ and interchanging the two marked points on it, and fixing the complement of  $\rho(C_0^+)\cap C_0^-$ in $C_0^-$ pointwise. Up to such an automorphism, the automorphism $\phi\in \Aut(C_0^-;\CQ)$ extends to all of $C_0$. Thus
$\phi(\rho(C_0^+)\cap C_0^-) = \phi(\gamma(\phi(\rho(C_0^+)\cap C_0^-))$ is contained in $C_0^{++}$. 

It remains to consider the case where  $\rho(C_0^+)\cap C_0^-$ contains exactly one marked point of $C_0^-$. Let $C_{0,1}^+$ and $C_{0,2}^+$ denote the irreducible components of $C_0^+$, such that $C_{0,1}^+$ contains the unique marked point on $\phi(\rho(C_0^+))$. We have decompositions
\[ \begin{array}{rcl}
C_0 &=& \phi(\rho(C_0^+)\cap C_0^-) + C_{0,1}^+ + \phi(C_0^\vee) + C_{0,2}^+ \\
&=& \rho^{-1}(C_0^+) + \rho^{-1} (C_0^-) ,
\end{array}\]
where $C_0^\vee$ denotes the closure of the complement of $\rho(C_0^+)\cap C_0^-$ in $C_0^-$. Note that both $ \phi(\rho(C_0^+)\cap C_0^-) + C_{0,1}^+$ and $
 \rho^{-1}(C_0^+)$ are connected $4$-pointed rational curves. As before, it suffices to show that 
\[ \phi(C_0^\vee) + C_{0,2}^+ \cong \rho^{-1}(C_0^-) \]
as $m$-pointed curves, because we then may apply corollary \ref{C333} again, and are done.  One has
\[\begin{array}{rcl}
\rho^{-1}(C_0^-)\cong C_0^-&=& (\rho(C_0^+)\cap C_0^-) + C_0^\vee \\
&\cong&\phi(\rho(C_0^+)\cap C_0^-) + \phi(C_0^\vee) 
\end{array}\]
as $m$-pointed curves. Since both $C_{0,2}^+$ and $\phi(\rho(C_0^+)\cap C_0^-)$ are $3$-pointed lines, they are obviously isomorphic, and this concludes the proof. 
\ebew

\begin{lemma} 
There is a natural commutative diagram of groups 
\[ \diagram
\id \dto &\id\dto\\
\Aut(C_0^{+-};\CQ^{+-},\underline{\CR}^{+-})\rdouble\dto^{\alpha_1} &\Aut(C_0^{+-};\CQ^{+-},\underline{\CR}^{+-})\dto^{\alpha_2}\\
\Aut(C_0^-;\CQ^{+-},\underline{\CR}^{+-},\CQ^{--})\rto|<\hole|<<\ahook\dto^{\beta_1} &\Aut(C_0^-;\CQ)\dto^{\beta_2}\\
\Aut(C_0^{--},\underline{\CR}^{+-},\CQ^{--})\rto|<\hole|<<\ahook\dto&\Aut(C_0^{--};\CR).\\
\id 
\enddiagram \]
Both colums of the diagram are exact sequences. 
\end{lemma}

\proof
The morphisms $\alpha_1$ and $\alpha_2$ are given by trivial extension of automorphisms from $C_0^{+-}$ to all of $C_0^-$. To see the  exactness of the first column note that the sequence splits. By extending automorphisms in  $\Aut(C_0^{--};\underline{\CR}^{+-},\CQ^{--})$ trivially to all of $C_0^-$ we obtain a section of  $\beta_1$.

The existence of the homomorphism $\beta_2$  follows from lemma \ref{Gbeta}. Let $\varphi\in\ker(\beta_2)$. Since $\beta_2(\varphi)$ is the identity map on $C_0^{--}$, we must in particular have that $\varphi$ fixes the set $\CR^{+-}$ pointwise, so $\varphi \in \Aut(C_0^-;\CQ^{+-},\underline{\CR}^{+-},\CQ^{--})$. Since $\beta_1(\varphi)=\id$, the automorphism $\varphi$ must in fact be contained in the group  $\Aut(C_0^{+-};\CQ^{+-},\underline{\CR}^{+-})$.
\ebew

\begin{rk}\em
$(i)$ The group $\Aut(C_0^{+-};\CQ^{+-},\underline{\CR}^{+-})$ is in a natural way a normal subgroup of $\Aut(C_0^-;\CQ)$.\\
$(ii)$ The natural homomorphism $\Aut(C_0^-;\CQ) \rightarrow \Aut(C_0^{--};\CR)$ is in general not surjective. As an example consider the case of the stable curve $C_0$ of genus $g=3$ with distinguished node $P_1$ as pictured below.

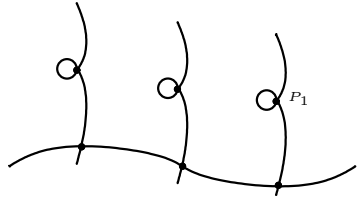
\begin{center}
\setlength{\unitlength}{0.00016667in}
\begingroup\makeatletter\ifx\SetFigFont\undefined%
\gdef\SetFigFont#1#2#3#4#5{%
  \reset@font\fontsize{#1}{#2pt}%
  \fontfamily{#3}\fontseries{#4}\fontshape{#5}%
  \selectfont}%
\fi\endgroup%
{\renewcommand{\dashlinestretch}{30}

}
\end{center}

$(iii)$ The natural homomorphism $\Aut(C_0) \rightarrow \Aut(C_0^{--};\CR)$ is also not surjective, as one can see in the following example of a stable curve $C_0$ of genus $g=4$ with distinguished node $P_1$ as below. 

\begin{center}
\setlength{\unitlength}{0.00016667in}
\begingroup\makeatletter\ifx\SetFigFont\undefined%
\gdef\SetFigFont#1#2#3#4#5{%
  \reset@font\fontsize{#1}{#2pt}%
  \fontfamily{#3}\fontseries{#4}\fontshape{#5}%
  \selectfont}%
\fi\endgroup%
{\renewcommand{\dashlinestretch}{30}
%
}
\end{center}
\end{rk}

\begin{rk}\label{diamond}\em
Let $f : C\rightarrow \Spec(k)$ be a stable curve represented by a point in $\dpobar$. By lemma \ref{L416}, there exists a distinguished decomposition of $C$ into two subcurves, one of them isomorphic to $C_0^{--}$, the other one a prestable curve $C^{++}$ of genus $g^{++}$ with $\mu$ marked points. There is a unique irreducible component $C^+$ of $C^{++}$ which has not the maximal number of nodes, when considered as a pointed stable curve itself. The closure of its complement in $C^{++}$ is isomorphic to $C_0^{+-}$, while $C^+$ is represented by a point in $\mgprimebar$. The marked points on the curve $C^+$ correspond naturally to the points of the sets $\CQ^{--}$ and $\CQ^{+-}$. 

It follows from remark \ref{37} that for almost all points of $\dpobar$ the group of automorphisms of the corresponding curve $C^{++}$ splits as
\[ \Aut(C^{++};\underline{\CQ}^{--},\underline{\CR}^{+-}) = \Aut_{\Gamma^\diamond}(C^{++};\underline{\CQ}^{--},{\CQ}^{+-}) \times \Aut(C_0^{+-};\underline{\CQ}^{+-},\underline{\CR}^{+-}), \]
where $\Gamma^\diamond = \Gamma(C_0^{+-};{\CQ}^{+-},\underline{\CR}^{+-})$. 
\end{rk}

\begin{notation}\em
We define the subscheme 
\[ \dii \: \: \subset \: \: \dpo  \label{cd1} \]
as  reduced subscheme of $\doo$ which is the locus of all curves with splitting automorphism group of $C^{++}$ as in remark \ref{diamond}. Since we consider only those cases where $\doo$ is non-empty, this is an open and dense subscheme of $\dpobar$.  

The reduced substack of $\FDPopen$ corresponding to $\dii$ shall be denoted by $\cdii$\label{cd2}. 

We define by $\kii$\label{cd3} the reduction of the preimage of $\dii$ inside $H_{g,n,0}$. Via the morphisms $\thx$ and $\thxx$, we define reduced subschemes  
\[\hii \subset H_{g^+,n,m}   \text{ and } \hiihat \subset {H}_{g^{++},n,\mu} \label{cd4b}\]
respectively, together with their quotients
\[ \hiigamma \subset H_{g^+,n,\mgamma}  \text{ and } \hiigammahat \subset {H}_{g^{++},n,\mu/\hat{\Gamma}} \, .\label{cd5b}\]
As before, we have $\Gamma = \Gamma(C_0^-;\CQ)$ and $\hat{\Gamma}(C_0^{--};\CQ)$. 

The corresponding subschemes of the moduli spaces are denoted by $\mii$, $\miihat$, $\miigamma$ and $\miigammahat$, respectively. 
\end{notation}

\begin{rk}\em
Using an analogous construction as in \ref{const_part_1}, we can define a morphism
\[ \Theta^{+-} : \quad \hii \longrightarrow \hiihat . \label{n52}\]
To do this we first glue the universal embedded stable curve $u_{g^+,n,m}: {\cal C}_{g^+,n,m}\rightarrow \hii$ to the trivial curve $\pr_2: C_0^{+-}\times \hii \rightarrow \hii$ along the sections of those marked points which are indexed by points in $\CQ^{+-}$. We denote the curve thus obtained by $u^{++}: {\cal C}_0^{++} \rightarrow \hii$. Now choose once and for all an embedding of this curve into $\mathbb P^N\times \hii\rightarrow  \hii$. The universal property of the Hilbert scheme induces a morphism $\Theta^{+-} : \hii \rightarrow  \overline{H}_{g^{++},n,\mu}$, which actually maps into $\hiihat$. 

Similar to what we did earlier, we can view $\PGL(N^{++}+1)$ as a subgroup of $\PGL(N+1)$, so that with respect to this inclusion the morphism $\Theta^{+-}$ is $\PGL(N^{++}+1)$-equivariant. Without loss of generality we may assume that the embeddings of $C_0^{+-}$ and of  $\mathbb P^{N^{+}}$ into $\mathbb P^N$ are invariant of the fibre of $\pr_2: \PP^N\times\hii\rightarrow \hii$.
\end{rk}

\begin{prop}\label{P614}
The morphism  $\Theta^{+-} :\hii \longrightarrow \hiihat$ factors through an injective morphism
\[ \Theta^\diamond : \quad \hii/\Gamma(C_0^{+-};\CQ^{+-},\underline{\CR}^{+-})  \longrightarrow \hiihat . \label{n53}\]
\end{prop}

\proof
Let $[C_1^+],[C_2^+]\in  \hii$ be represented by embedded curves $C_1^+$ and $C_2^+$ in $\mathbb P^{N^{+}}$, and such that $\Theta^{+-}([C_1^+])= \Theta^{+-}([C_2^+]) $. By construction of $\Theta^{+-}$, the embedding of $C_0^{+-}$ into $\mathbb P^N$ is independent of the fibre. Thus $C_1^+$ and $C_2^+$ are identically embedded into $\PP^N$, as long as we disregard the labels of their marked points. Since the embedding of $\PP^{N^{+}}$ into $\PP^N$ is also independent of the fibre, the embedded curves $C_1^+$ and $C_2^+$ are identical in $\PP^{N^{+}}$, at least up to the order of the labels of their marked points. 

Marked points on $C_1^+$ and $C_2^+$, which are not glueing points with $C_0^{+-}$, are marked points of the curve representing $\Theta^{+-}([C_1^+])= \Theta^{+-}([C_2^+]) $, and hence they agree on $C_1^+$ and $C_2^+$. Therefore  $C_1^+$ and $C_2^+$ can differ as marked curves at most by a permutation of the labels of their glueing points, i.e. those points which are identified with the points of $\CQ^{+-}$ on $C_0^{+-}$. Such a permutation must be induced by an automorphism $\gamma$ of $C_0^{+-}$, which fixes those  marked points of $C_0^{+-}$ which are not glueing points. Hence $\gamma$ can only be an element of $\Gamma(C_0^{+-};\CQ^{+-},\underline{\CR}^{+-})$. 
\ebew

\begin{lemma}\label{L515}
Let $\pi: \overline{H}_{g^{++},n,\mu} \rightarrow \overline{H}_{g^{++},n,\mu/\Gamma(C_0^{--};\CR)}$ denote the canonical quotient morphism. The composed morphism  $\pi\circ {\Theta}^\diamond$ from the quotient scheme $  H^\diamond_{g^+,n,m/\Gamma(C_0^{+-};\CQ^{+-},\underline{\CR}^{+-})}$ to  $  \hat{H}{}^\diamond_{g^{++},n,\mu/\Gamma(C_0^{--};\CR)}$ factors through a morphism
\[ \overline{\Theta}{}^\diamond: \quad {H}^\diamond_{g^+,n,m/\Gamma(C_0^-;\CQ)} \longrightarrow \hat{ H}{}^\diamond_{g^{++},n,\mu/\Gamma(C_0^{--};\CR)}. \label{n54}\]
\end{lemma}

\proof
Let $C_1\rightarrow \Spec(k)$ and $C_2\rightarrow \Spec(k)$ be two embedded stable curves of genus $g^+$ with $m$ marked points, which differ by a reordering of the labels of their marked points by a permutation $\gamma\in \Gamma(C_0^-;\CQ)$. Recall that $\gamma$ is induced by an automorphism $\rho\in \Aut(C_0^-;\CQ)$. Let ${\cal P}_i$ denote the set of marked points of $C_i$, for $i=1,2$. They are disjoint unions of subsets ${\cal P}_i^{+-}$ and ${\cal P}_i^{--}$, where ${\cal P}_i^{+-}$ denotes the subset of those points, in which $C_i$ is glued to $C_0^{+-}$ when the morphism $\Theta^{+-}$ is applied, and  ${\cal P}_i^{--}$ its complement. In other words, if the points of ${\cal P}_i$ are enumerated as $P_{i,1},\ldots,P_{i,m}$, then $P_{i,j}\in {\cal P}_i^{+-}$ if and only if $Q_{j}\in \CQ^{+-}$. 
 By lemma \ref{Gbeta} the subcurve $C_0^{+-}$of $C_0^-$ is stabilized under the action of $\Aut(C_0^-;\CQ)$. In particular, the subset $\CQ^{+-}$, and hence the sets  $ {\cal P}_i^{+-}$, are stabilized. The automorphism $\rho\in\Aut(C_0^-;\CQ)$ restricts to an automorphism in $\Aut(C_0^{+-};\CN)$, where $\CN$ was defined as the set of marked points on $C_0^{+-}$. By construction, the embedding of $C_0^{+-}$ into $\PP^N$ is invariant under the action of $\Aut(C_0^{+-};\CN)$, and thus the embedded curves representing $\Theta^{+-}([C_1])$ and $\Theta^{+-}([C_2])$ are equal as embedded curves, up to a permutation of the labels of their marked points. These marked points are in one-to-one correspondence with the points of $\CR$. The permutation of the labels is given by a permutation $\sigma\in\Gamma(C_0^-;\CQ)$. Let $\sigma$ be induced by an automorphism $\gamma\in\Aut(C_0^-;\CQ)$. Then $\gamma$ restricts to an automorphism in $\Aut(C_0^{--};\CR)$, which in turn defines a permutation in $\Gamma(C_0^{--};\CR)$. Hence $\Theta^{+-}([C_1])$ and $\Theta^{+-}([C_2])$ are equal as $m/\Gamma(C_0^{--};\CR)$-pointed curves. Therefore $\overline{\Theta}{}^\diamond$ is well-defined. 
\ebew

\begin{lemma}\label{L6010}
The morphism 
\[ \overline{\Theta}{}^\diamond \: : \quad {H}^\diamond_{g^+,n,m/\Gamma(C_0^-;\CQ)} \longrightarrow \hat{H}{}^\diamond_{g^{++},n,\mu/\Gamma(C_0^{--};\CR)}\]
is injective.
\end{lemma}

\proof
Consider two embedded $m/\Gamma(C_0^-;\CQ)$-pointed stable curves $C_1^+\rightarrow \Spec(k)$ and $C_2^+\rightarrow \Spec(k)$ of genus $g^+$. Suppose that $C_1^+ + C_0^{+-} = C_2^+ + C_0^{+-}$ as embedded $\mu/\Gamma(C_0^{--};\CR)$-pointed stable curves in $\PP^{++}$. Here ``+'' denotes again glueing in the marked points $\CQ^{+-}$ of $C_0^{+-}$. 

By the construction of $\overline{\Theta}{}^\diamond$, the embedding of $C_0^{+-}$ into $\PP^{N^{++}}$ is fixed. Hence the two $\mu$-pointed curves  $C_1^+ + C_0^{+-}$ and $  C_2^+ + C_0^{+-}$ may differ at most by a permutation of the labels of those marked points which lie on $C_1^+$ and $C_2^+$, respectively. In other words, they differ by a permutation in $\Gamma(C_0^{--};\underline{\CR}^{+-},\CQ^{--})$. However, an automorphism of $C_0^{--}$, inducing such a permutation, can be trivially extended to an automorphism of $C_0^-$, which in turn defines a permutation in $\Gamma(C_0^-;\underline{\CQ}^{+-},\CQ^{--})$. In particular, the $m$-pointed curves $C_1^+$ and $C_2^+$ differ at most by a permutation in $\Gamma(C_0^-;\CQ)$, and hence are the same as $m/\Gamma(C_0^-;\CQ)$-pointed curves.
\ebew

\begin{rk}\em
The morphisms $\Theta^{+-}$, ${\Theta}^\diamond$ and $\overline{\Theta}{}^\diamond$ from above induce a commutative diagram of morphisms between moduli spaces as below.
\[ \diagram
\mii \rrto \dto&&\miihat\dto^=\\
\ms_{g^+,m/\Gamma(C_0^{+-};\CQ^{+-},\underline{\CR}^{+-})}^\diamond \rrto \dto&&\miihat\dto\\
\ms_{g^+,m/\Gamma(C_0^-;\CQ)}^\diamond \rrto^\cong \drto_\cong&&\hat{M}_{g^{++},\mu/\Gamma(C_0^{--};\CR)}^\diamond\dlto^\cong\\
& \dii . 
\enddiagram\]
The vertical arrows are all given by the action of finite groups, where the action is not free in general. The corresponding morphisms will therefore be finite but ramified.  The isomorphisms with $\dii$ were already dealt with in proposition \ref{312} and remark \ref{R0520}. By composition, this implies that the  last horizontal morphism of the diagram is an isomorphism as well. 
\end{rk}

\begin{lemma}\label{L617}
The natural morphism
\[ \ms_{g^+,m/\Gamma(C_0^{+-};\CQ^{+-},\underline{\CR}^{+-})}^\diamond \longrightarrow \miihat \]
is injective. Thus $ \msbar_{g^+,m/\Gamma(C_0^{+-};\CQ^{+-},\underline{\CR}^{+-})}'$ is isomorphic to an irreducible component $\msbar_{g^{++},\mu}^\#$ of the closure of $\miihat$ in  $\msbar_{g^{++},\mu}^{ps}$.
\end{lemma}

\proof
Let $C_1\rightarrow \Spec(k)$ and $C_2\rightarrow \Spec(k)$ be two $m$-pointed stable curves of genus $g^+$, such that $C_1+C_0^{+-}$ and $C_2+C_0^{+-}$ are isomorphic as $\mu$-pointed stable curves of genus $g^{++}$. 
Here ``$+$'' denotes glueing along the marked points indexed by points of $\CQ^{+-}$. 
If $C_1$ and $C_2$ have both less than the maximal number of nodes possible, then any isomorphism from $C_1+C_0^{+-}$ to $C_2+C_0^{+-}$ restricts to an isomorphism from $C_1$ to $C_2$. Such an isomorphism may permute the labels of the marked points by an element in $\Gamma(C_0^{+-};\CQ^{+-},\underline{\CR}^{+-})$. This shows that the restricted morphism from $\ms_{g^+,m/\Gamma(C_0^{+-};\CQ^{+-},\underline{\CR}^{+-})}^\diamond$  to $\hat{M}_{g^{++},\mu}^\diamond$ is injective. 

Because of remark \ref{R0323}, the moduli space $\msbar_{g^{+},m}'$ is smooth and irreducible, and so is its quotient $\msbar_{g^+,m/\Gamma(C_0^{+-};\CQ^{+-},\underline{\CR}^{+-})}'$. Likewise, the irreducible components of the closure of $\miihat$ are smooth by lemma \ref{L518}.  Therefore the injectivity of the restricted finite morphism between the open subschemes  extends to the closure of the dense open subscheme $ \ms_{g^+,m/\Gamma(C_0^{+-};\CQ^{+-},\underline{\CR}^{+-})}^\diamond$, which is the scheme $ \msbar_{g^+,m/\Gamma(C_0^{+-};\CQ^{+-},\underline{\CR}^{+-})}'$. 
\ebew

\begin{rk}\em \label{R618}
The above $\PGL(N^{++}+1)$-equivariant morphism of schemes $\Theta^{+-}: \hii\rightarrow \hiihat$ induces a morphism of stacks
\[ \left[ \hii/\PGL(N^++1) \right] \longrightarrow \left[ 
\hiihat / \PGL(N^{++}+1) \right]. \]
The morphism $\Theta^{+-}$ factors through $ {H}{}^\diamond_{g^+,n,m/\Gamma(C_0^{+-};\CQ^{+-},\underline{\CR}^{+-})}$, which implies that the induced morphism of stacks factors through the quotient stack $[ {H}{}^\diamond_{g^+,n,m/\Gamma(C_0^{+-};\CQ^{+-},\underline{\CR}^{+-})}/\PGL(N^++1)]$. In general, the factorized morphism is not an ismorphism of stacks. 
 \end{rk}

\begin{prop}\label{P6014}
There is an injective morphism of stacks
\[ \begin{array}{c}
\left[{H}{}^\diamond_{g^+,n,m/\Gamma(C_0^{+-};\CQ^{+-},\underline{\CR}^{+-})} \, / \, \PGL(N^++1)\times\Aut(C_0^{+-})\right] \\
\downarrow  \\
 \left[ 
\hiihat \,  / \, \PGL(N^{++}+1) \right] .
\end{array}\]
\end{prop}

\proof 
The construction of such a morphism $\Psi$  is straightforward, compare proposition \ref{MainProp}. To show injectivity, we have show for all schemes $S$ that the restriction $\Psi_S$ is full and faithful as a functor between the respective fibre categories, see \cite[2.2]{LM}. Fix a scheme $S$, and two objects $(E_1,p_1,\theta_1)$ and $(E_2,p_2,\theta_2)$ in $[{H}{}^\diamond_{g^+,n,m/\Gamma(C_0^{+-};\CQ^{+-},\underline{\CR}^{+-})}/\PGL(N^++1)\times\Aut(C_0^{+-})](S)$. A morphism $\phi$ from $(E_1,p_1,\theta_1)$ to  $(E_2,p_2,\theta_2)$ is given by a morphism  $\lambda: E_2\rightarrow E_1$ of principal $\PGL(N^+ +1)\times\Aut(C_0^{+-})$-bundles over $S$, satisfying $\theta_2=\theta_1\circ \lambda$. By definition, the morphism $\Psi_S(\phi)$ is given by the $\PGL(N^{++}+1)$-extension $\overline{\lambda}$ of $\lambda$
\[ \overline{\lambda} : \:\: \:  E_2(\PGL(N^{++}+1)) \rightarrow  E_1(\PGL(N^{++}+1)) \]
between the $\PGL(N^{++}+1)$-extensions of $E_2$ and $E_1$. Clearly, $\lambda$ is uniquely determined by the restriction of $\overline{\lambda}$, so the functor $\Psi_S$ is faithful.

To show that the functor $\Psi_S$ is also full, we need to verify that any morphism $\phi$ from $\Psi_S(E_1,p_1,\theta_1)$ to $\Psi_S(E_2,p_2,\theta_2)$ arises in this way. Consider a morphism $\tilde{\lambda}: E_2(\PGL(N^{++}+1)) \rightarrow  E_1(\PGL(N^{++}+1))$ of principal $\PGL(N^{++}+1)$-bundles over $S$, which is compatible with the respective morphisms of the principal $\PGL(N^{++}+1)$-bundles into $\hiihat$. Once we demonstrate that the image of the restriction of $\tilde{\lambda}$ of the subbundle $E_2$ is contained in $E_1$, we are done. It is enough to do this locally, or even fibrewise over $S$. Indeed, if such a restriction $\lambda$ exists, then the identity $ \tilde{\Theta}^+\circ \theta_2 =\tilde{\Theta}^+\circ\theta_1\circ\tilde{\lambda}$ implies the equality $ \theta_2 =\theta_1\circ{\lambda}$. Recall that by proposition \ref{P614} the  morphism $ \tilde{\Theta}^+$ is injective. Therefore $\lambda$ defines a morphism $\phi'$ between the triples $(E_1,p_1,\theta_1)$ and $(E_2,p_2,\theta_2)$, such that $\Psi_S(\phi')=\phi$.

Let $f_1^+: C_1^+\rightarrow S $ and $f_2^+: C_2^+\rightarrow S$ denote the $m/\Gamma(C_0^{+-};\CQ^{--},\underline{\CR}^{+-})$-pointed stable curves of genus $g^+$ corresponding to the two triples $(E_1,p_1,\theta_1)$ and $(E_2,p_2,\theta_2)$. Then the triples $\Psi_S(E_1,p_1,\theta_1)$ and $\Psi_S(E_2,p_2,\theta_2)$ correspond to the two curves $f_1: C_1\rightarrow S$ and $f_2: C_2\rightarrow S$, where $C_1\rightarrow S$ is obtained from $C_1^+\rightarrow S$ by glueing with the trivial curve $C_0^{+-}\times S\rightarrow S$ along the sections of those marked points, which lie in $\CQ^{+-}$, and analogously for $C_2$. The bundle $E_1(\PGL(N^{++}+1))$ is the principal $\PGL(N^{++}+1)$-bundle associated to the projective bundle $\PP((f_1)_\ast(\omega_{C_1/S}({\cal S}_{1,1}+ \ldots + {\cal S}_{1,\mu}))^{\otimes n}))$, where ${\cal S}_{1,1}, \ldots , {\cal S}_{1,\mu}$ denote the divisors of marked points of $C_1$. 

Locally, there exists a principal $\PGL(N^+ +1)$-subbundle $E_1^+$, which is the principal bundle associated to the projective bundle $\PP((f_1^+)_\ast(\omega_{C_1^+/S}({\cal S}^+_{1,1}+ \ldots + {\cal S}^+_{1,m}))^{\otimes n}))$, where ${\cal S}^+_{1,1}, \ldots , {\cal S}^+_{1,m}$ denote the divisors of marked points of $C_1^+$. Recall that we fixed the embedding of $C_0^{+-}$ into $\PP^{N^{++}}$. Thus, via the inclusion
\[ \omega_{C_1^+/S}({\cal S}^+_{1,1}+ \ldots + {\cal S}^+_{1,m}) \hookrightarrow \omega_{C_1/S}({\cal S}_{1,1}+ \ldots + {\cal S}_{1,\mu}) \]
the embedding of $E_1^+$ into $E_1(\PGL(N^{++}+1))$ is uniquely determined up to the action of $\Aut(C_0^{+-})$.  By construction, for the extension of $E_1^+$ holds
\[ E_1^+(\PGL(N^+ +1)\times\Aut(C_0^{+-}) ) = E_1 ,\]
as subbundles of $E_1(\PGL(N^{++}+1))$. 

Note that such a subbundle may not exist globally. Indeed, for a fibre $C_{1,s}^+$ of $f_1^+: C_1^+\rightarrow S$, the assignment of its $m$ marked points to nodes and marked points of the fibre $C_{1,s} = C_{1,s}^+ + C_0^{+-}$ of $f_1:C_1 \rightarrow  S$ is only determined up to the action of the group $\Gamma(C_0^{+-};\CQ^{--},\underline{\CR}^{+-})$. 

Analogously, there is locally a principal $\PGL(N^+ +1)$-bundle $E_2^+$, which is associated to the projective bundle $\PP((f_2^+)_\ast(\omega_{C_2^+/S}({\cal S}^+_{2,1}+ \ldots + {\cal S}^+_{2,m}))^{\otimes n}))$,
and which is locally a subbundle of $E_2$. 

The morphism $\phi$ corresponds to a morphism $\tilde{\phi}: C_2\rightarrow C_1$ of prestable curves over $S$. Since $C_0^{+-}$ does have the maximal number of nodes possible, and the fibres of $f_1^+: C_1^+\rightarrow S$ and $f_2^+:C_2^+\rightarrow S$ don't, the morphism $\tilde{\phi}$ restricts to a morphism $\phi' : C_2^+ \rightarrow C_1^+$ over $S$. The induced morphism on the projective bundles thus restricts to a morphism 
\[  \begin{array}{c}
\PP((f_2^+)_\ast(\omega_{C_2^+/S}({\cal S}^+_{2,1}+ \ldots + {\cal S}^+_{2,m}))^{\otimes n})) \\
\downarrow\\
\PP((f_1^+)_\ast(\omega_{C_1^+/S}({\cal S}^+_{1,1}+ \ldots + {\cal S}^+_{1,m}))^{\otimes n})).
\end{array}\]
This induces a morphism $\lambda': E_2^+ \rightarrow E_1^+$ of the locally existing principal $\PGL(N^+ +1)$-bundles. By construction, its $\PGL(N^{++}+1)$-extension is $\tilde{\lambda}$. Thus in particular $\tilde{\lambda}(E_2)\subset E_1$, and we are done. 
\ebew

\begin{notation}\em
Let $\overline{H}{}^\#_{g^{++},n,\mu}$ denote the subscheme of $\overline{H}{}_{g^{++},n,\mu}$, which is the reduction of the preimage of the irreducible component $\msbar_{g^{++},\mu}^\#$ of the closure of $\miihat$ as in lemma \ref{L617}. Its intersection with $\hiihat$, i.e. the open subscheme, which  parametrizes curves with less than the maximal number of nodes and splitting groups of automorphisms, shall be denoted by ${H}^\#_{g^{++},n,\mu}$. 
\end{notation}

\begin{prop}\label{P53a}
There is an isomorphism of stacks
\[ \begin{array}{c}
\left[{H}{}^\diamond_{g^+,n,m/\Gamma(C_0^{+-};\CQ^{+-},\underline{\CR}^{+-})}/\PGL(N^++1)\times\Aut(C_0^{+-})\right] \\
\downarrow \cong \\
 \left[ 
 {H}^\#_{g^{++},n,\mu} / \PGL(N^{++}+1) \right] .
\end{array}\]
\end{prop}

\proof
Note that the injective morphism of proposition \ref{P6014} clearly factors through $[ 
{{H}}^\#_{g^{++},n,\mu} / \PGL(N^{++}+1)]$. The proof of the existence of an inverse morphism works
 analogously as in the  proofs of theorem \ref{MainThm} and  theorem \ref{49}. Consider an object of $[ 
{{H}}^\#_{g^{++},n,\mu} / \PGL(N^{++}+1)](S)$ for some scheme $S$. Such an object is given by a  $\mu$-pointed stable curve $f:C\rightarrow S$ of genus $g^{++}$, which is represented by a point of $\msbar_{g^{++},\mu}^\#$ and has at most $3g-4$ nodes. Essentially, what we have to do is to find  a natural subcurve $f^+: C^+ \rightarrow S$, which is an $m/\Gamma(C_0^{+-};\CQ^{+-},\underline{\CR}^{+-})$-pointed stable curve of genus $g^+$, and such that glueing with the trivial curve $C_0^{+-}\times S\rightarrow S $ gives $f:C\rightarrow S$ again. Here glueing means glueing along the sections of those marked points which lie in $\CQ^{+-}$. 
Because we are working on the open moduli substacks of curves with less than the maximal number of nodes, such a subcurve exists, up to the labelling of the marked points. Indeed, for any $\mu$-pointed stable curve $f: C \rightarrow \Spec(k)$ of genus $g^{++}$, which has one node less than the maximal possible number, there is a unique irreducible component $C^+$ of $C$ which has not the maximal number of nodes when considered as a pointed stable curve itself. Note that  because  a  fibre $C_s$  of $f:C\rightarrow S$ is  represented by a point of $\msbar_{g^{++},\mu}^\#$, the distribution of the labels of the marked points of $C_s$ is compatible with the decomposition of the fibre as $C_s= C_s^+ + C_0^{+-}$. 

Suppose that the stable curve $f:C\rightarrow S$ is represented by an object $(E,p,\theta)\in [ 
{{H}}^\#_{g^{++},n,\mu} / \PGL(N^{++}+1)](S)$, i.e. by a $\PGL(N^{++}+1)$-principal bundle $p:E\rightarrow S$, together with a  $\PGL(N^{++}+1)$-equivariant morphism $\theta: E\rightarrow {{H}}^\#_{g^{++},n,\mu}$. Using the subcurve $f^+:C^+\rightarrow S$, we obtain locally $\PGL(N^{+}+1)$-subbundles $p^+:E^+
\rightarrow S_\alpha$, where $\{S_\alpha\}_{\alpha\in A}$ is a covering of $S$. As before, these local subbundles glue together, up to the action of $\Aut(C_0^{+-})$. Like in the proof of theorem \ref{MainThm}, from this  we can  construct an object $(E^+, p^+,\theta^+)$ of $[{H}^\diamond_{g^+,n,m/\Gamma(C_0^{+-};\CQ^{+-},\underline{\CR}^{+-})}/\PGL(N^++1)\times\Aut(C_0^{+-})](S)$, which maps to the original triple  $(E,p,\theta)$.

Similarly, the construction of the functor on morphisms is analogous to the construction in the proof of theorem \ref{MainThm}. The injectivity of the morphism ${\Theta}^{\diamond}$ between the corresponding atlases is provided by proposition \ref{P614}. 
\ebew

\begin{rk}\em\label{R52}
Composition of the morphism of proposition   \ref{P6014} with the natural morphism from $[\hiihat / \PGL(N^{++}+1) ]$ to the quotient stack $[\hat{{ H}}{}^\diamond_{g^{++},n,\mu/\Gamma(C_0^{--};\CR)} / \PGL(N^{++}+1) ]$
we obtain a morphism of stacks
\[ 
\begin{array}{c}
\left[ {H}^\diamond_{g^+,n,m/\Gamma(C_0^{+-};\CQ^{+-},\underline{\CR}^{+-})} /\PGL(N^++1)\times\Aut(C_0^{+-})\right] \\
\downarrow \\
 \left[ 
\hat{{H}}{}^\diamond_{g^{++},n,\mu/\Gamma(C_0^{--};\CR)} / \PGL(N^{++}+1) \right]. \end{array}\]
By lemma \ref{L515} the morphism $ {H}^\diamond_{g^+,n,m/\Gamma(C_0^{+-};\CQ^{+-},\underline{\CR}^{+-})} \rightarrow \hat{{H}}{}^\diamond_{g^{++},n,\mu/\Gamma(C_0^{--};\CR)}$ facors through $ {H}{}^\diamond_{g^+,n,m/\Gamma(C_0^-;\CQ)}$.

 Therefore the above morphism of stacks factors through the quotient stack
$
[ \,  {H}^\diamond_{g^+,n,m/\Gamma(C_0^-;\CQ)} /\PGL(N^++1)\times\Aut(C_0^{+-})]$.
\end{rk}

\begin{prop}\label{P54}
There is an injective morphism  of stacks
\[ \begin{array}{c}
 \left[ {H}^\diamond_{g^+,n,m/\Gamma(C_0^-;\CQ)} \, / \, \PGL(N^++1)\times\Aut(C_0^{+-})\right]\\
\downarrow \\
  \left[ 
\hat{{H}}{}^\diamond_{g^{++},n,\mu/\Gamma(C_0^{--};\CR)} \,  / \, \PGL(N^{++}+1) \right] .
\end{array}\]
\end{prop}
  
\proof
The proof is completely analogous to that of  proposition \ref{P6014}. Note that on order to show that the relevant functor is full, one needs the injectivity of $\overline{\Theta}^\diamond$, which is provided by lemma \ref{L6010}.
\ebew

\begin{rk}\em
Taking everything together we obtain the follwing commutative diagram of quotient stacks. For typographical reasons we use the abbreviations  $P^+ :=  \PGL(N^++1)$ and $P^{++} :=  \PGL(N^++1)$, as well as $A^- :=  \Aut(C_0^{-})$, $A^{+-} :=  \Aut(C_0^{+-})$ and $A^{--} :=  \Aut(C_0^{--})$. 
\pagebreak
\[ 
\left[ {H}^\diamond_{g^+,n,m/\Gamma(C_0^{+-};\CQ^{+-},\underline{\CR}^{+-})} / P^+\times A^{+-} \right] \:\: \longrightarrow \:\: 
\left[ \hat{{H}}{}^\diamond_{g^{++},n,\mu} / P^{++} \right]\]
\[ \qquad \downarrow  \hspace{7cm} \downarrow \qquad \]
\[ \left[  {H}^\diamond_{g^+,n,m/\Gamma(C_0^{-};\CQ)} / P^+\times A^{+-} \right] \:\:  \longrightarrow   \:\: 
\left[ \hat{{H}}{}^\diamond_{g^{++},n,\mu/\Gamma(C_0^{--};\CR)} / P^{++} \right]  \]
\[ \qquad \downarrow  \hspace{7cm} \downarrow \qquad \]
\[ \left[ {H}^\diamond_{g^+,n,m/\Gamma(C_0^{-};\CQ)} /P^+\times A^{-} \right] \:\:  \longrightarrow \:\: 
\left[ \hat{{H}}^\diamond_{g^{++},n,\mu/\Gamma(C_0^{--};\CR)} / P^{++}\times A^{--} \right]  \]
All vertical arrows of the diagram are induced by group actions, all horizontal
 arrows are injective. The last horizontal arrow is an isomorphism. 
\end{rk}

\begin{rk}\em
Suppose that there are two boundary points $[C_1]$ and $[C_2]\in \overline{D}(P_1)$, represented by  stable curves $f_1: C_1\rightarrow \Spec(k)$ and $f_2: C_2\rightarrow \Spec(k)$, both of genus $g$ and with $3g-3$ nodes. Then  
\[  \overline{D}(P_1) =  \overline{D}(C_1;P_1') = \overline{D}(C_2;P_1'') \]
for suitable nodes $P_1' \in C_1$ and $P_1'' \in C_2$. Let $C_1= C_1^{++}+ C_1^{--}$ and $C_2= C_2^{++}+ C_2^{--}$ denote the  decompositions into invariant subcurves determined by the equivalence classes of $P_1'$ and $P_1''$. Suppose that $C_i^{++}$ is a $\mu_i$-pointed prestable curve of genus $g_i^{++}$, with $\CR_i$ as its set of marked points, for $i=1,2$.  Then the respective partial compactifications belonging to $[C_1]$ and $[C_2]$ glue together over the open substacks via an isomorphism
\[ \begin{array}{c}
 \left[  \hat{{H}}{}^\diamond_{g_1^{++},n,\mu_1/\Gamma(C_1^{--};\CR_1)} \: / \: \PGL(N_1^{++}+1)\times\Aut(C_1^{--})\right]\\
\downarrow\cong\\
  \left[  \hat{{H}}{}^\diamond_{g_2^{++},n,\mu_2/\Gamma(C_2^{--};\CR_2)} \: / \: \PGL(N_2^{++}+1)\times\Aut(C_2^{--})\right] .
\end{array} \]
\end{rk}

%
%

\chapter{The global picture}

\begin{rk}\em
In the previous chapters we approached the stack $\FDP$ by studying substacks
\[ \cdx \subset \FDPxx \subset \FDP .\]
The substack $\cdx$, corresponding to stable curves of genus $g$  with exactly $3g-4$ nodes and splitting groups of automorphims, which are represented by points of $\dx$, can be described as a quotient stack as in theorem \ref{MainThm}. This construction depends on a decomposition of $C_0 = C_0^+ + C_0^-$, which is determined by the fixed node $P_1 \in C_0$. Replacing this decomposition by a decomposition $C_0 = C_0^{++} + C_0^{--}$, which is equivariant with respect to the action of $\Aut(C_0)$, we obtain a description of $\FDPxx$, the preimage substack belonging to the partial compactification $\dxxbar$ of $\dpo$ at the point $[C_0] \in \msbar_g$, as a quotient stack, see theorem \ref{49}. 

It follows from proposition \ref{MainProp} and proposition \ref{P217} that there is a morphism of stacks 
\[ {\Omega} : \quad \cmgammadetailx  \longrightarrow \cdx, \label{n55}\]
which is the composition of the morphism $\Lambda $ from poposition \ref{MainProp} with the canonical morphism
\[  \begin{array}{c}
\left[ \overline{H}_{g^+,n,m/\Gamma(C_0^-;\CQ)} \: / \:\PGL(N^++1)\right] \\
\downarrow \\
\left[ \overline{H}_{g^+,n,m/\Gamma(C_0^-;\CQ)} \: / \: \PGL(N^++1)\times\Aut(C_0^-)\right],
\end{array} \]
using the isomorphism $ \cmgammadetailx \cong [ \hgammadetailx / \PGL(N^++1)]$. 
\end{rk}

\begin{lemma}
Considered as a functor between fibred categories, the morphism $\Lambda$ is a faithful  functor.
\end{lemma}

\proof
Recall that $\cdx \cong [\kx / \PGL(N+1)]$ by lemma \ref{140}.
Consider two objects $(E_1^+,p_1^+,\theta_1^+) \in [\hgammadetailx / \PGL(N^++1)\times \Aut(C_0^-)](S_1)$ and $(E_2^+,p_2^+,\theta_2^+) \in [ \hgammadetailx / \PGL(N^++1)\times \Aut(C_0^-)](S_2)$, together with a pair of morphisms $(f_i,\varphi_i)$, $i=1,2$,  between them. Suppose that $\Lambda((\varphi_1,f_1))=\Lambda((\varphi_2,f_2))$. In other words, for $i=1,2$ there is a  commutative diagram
\[\diagram
&&&& \hgammadetailx \\
E_1^+ \rrto_{\varphi_i} \xto[rrrru]^{\theta_1^+}\drto_{p_1^+}&& f_i^\ast(E_2^+) \rto\dlto & E_2^+ \urto_{\theta_2^+}\dto^{p_2^+}\\
&S_1\rrto^{f_i}&&S_2.
\enddiagram\]
If $E_i := (E_i^+\times\PGL(N+1))/(\PGL(N^++1)\times\Aut(C_0^-))$ denotes the extension of $E_i^+$ by $\PGL(N+1)$, then $(\overline{\varphi}_i,f_i)= \Lambda((\varphi_i,f_i))$ is given by the commutative diagram
\[\diagram
&&&& \kx \\
E_1 \rrto_{\overline{\varphi}_i} \xto[urrrr]^{\theta_1}\drto_{p_1}&& f_i^\ast(E_2) \rto\dlto & E_2 \urto_{\theta_2}\dto^{p_2}\\
&S_1\rrto^{f_i}&&S_2.
\enddiagram\]
By construction, $E_i^+$ is a subbundle of $E_i$, and $\overline{\varphi}_i$ is the $\PGL(N+1)$-equivariant extension of $\varphi_i$. Thus $\varphi_1=\varphi_2$ if $\overline{\varphi}_1=\overline{\varphi}_2$. Recall that  $\theta_i = \Theta\circ \theta_i^+$, where $\Theta : \hgammadetailx \rightarrow \kx$ is the injective morphism from proposition \ref{312}. This shows that the morphism of stacks $\Lambda$ is faithful when considered as a functor of fibred categories.
\ebew

\begin{rk}\em\label{R63}
In general, $\Lambda$ is not a full functor. Indeed, consider a stable curve $f: C\rightarrow \Spec(k)$ of genus $g$ with $3g-4$ nodes, for which there exists an automorphism $\gamma\in \Aut(C)$ which acts non-trivially on the distinguished subcurve isomorphic to $C_0^-$. Let $f^+:C^+\rightarrow \Spec(k)$ denote the complementary subcurve. Then there is certainly no automorphism of $f^+: C^+\rightarrow \Spec(k)$ in $\cmgammadetailx(\Spec(k))$ which gets sent to $\gamma$ by the morphism  $\Lambda$.
\end{rk}  

\begin{rk}\em
Note that the morphism 
\[ \Omega : \: \cmgammadetailx = [ \hgammadetailx/\PGL(N^+ +1)]  \longrightarrow \cdx \]
is the composition of  a representable,  finite and surjective morphism with an isomorphism. This alone is enough to prove the following lemma. However the proof given below helps to understand the morphism $\Omega$ geometrically.
\end{rk}

\begin{lemma}\label{L61}
The morphism $\Omega :\cmgammadetailx \rightarrow \cdx$ is representable, surjective, finite and \'etale,  and its  degree is equal to the order of the group $\Aut(C_0^-)$. 
\end{lemma}

\proof
Let $S$ be a scheme and let $s: S\rightarrow \cdx$ be a morphism. Then $s\in\cdx(S)$ is given by a triple $(E,p,\theta)\in[\kx/\PGL(N+1)](S)$. In the proof of theorem \ref{MainThm} we constructed a principal $\PGL(N^+ +1)\times\Aut(C_0^-)$-bundle $p^+: E^+\rightarrow S$ as a subbundle of $p: E\rightarrow S$,  such that for the $\PGL(N+1)$-extension of $E^+$ holds $E^+(\PGL(N+1))= E$. 

Define $\tilde{E} := E^+/ \PGL(N^+ +1)$, so that $\tilde{p}: \tilde{E} \rightarrow S$ is a principal $\Aut(C_0^-)$-bundle over $S$. By the universal property of the Hilbert scheme, the morphism $\theta: E\rightarrow \kx$ defines a stable curve $f_E: C_E\rightarrow E$ of genus $g$. By corollary \ref{K324} there exists a distinguished subcurve $f_E^+: C_E^+\rightarrow E$, which is locally isomorphic to the trivial curve over $E$ with fibre $C_0^-$, considered as an $m/\Gamma(C_0^-;\calq)$-pointed curve. Thus there is an induced morphism
\[ \theta_{f_E^+} : \:\: E \rightarrow \hgammadetailx .\]
Composition with the inclusion $E^+\subset E$ gives thus a $\PGL(N^+ +1)$-equivariant morphism
\[ \theta^+ : \:\: E^+ \rightarrow \hgammadetailx .\]
The quotient morphism  $p^+: E^+\rightarrow \tilde{E}$ is a principal $\PGL(N^+ +1)$-bundle over $\tilde{E}$. In particular, we have constructed an object 
\[ (E^+,p^+,\theta^+) \in [\hgammadetailx/\PGL(N^+ +1) ] (\tilde{E}) ,\]
or equivalently, a morphism $\phi: \tilde{E}\rightarrow  \cmgammadetailx$. In fact, the diagram
\[ \diagram
\tilde{E}\rrto^{\tilde{p}} \dto_\phi&& S\dto^s\\
\cmgammadetailx  \rrto_\Omega && \cdx
\enddiagram \]
is Cartesian, as we will explain in a minute. Thus, by definition, the morphism $\Omega$ is representable. Since $\tilde{p}$ is surjective, unramified and finite of degree equal to the order of $\Aut(C_0^-)$, so is $\Omega$. Since both stacks ars smooth of the same dimension, unramified implies \'etale. 

The main point in proving that the above diagram is Cartesian is to find for each given object 
\[ \left( (E_B,p_B,\theta_B),b,\lambda\right) \in \cmgammadetailx \times_{\cdx} S^\bullet (B)\]
an object $h \in \tilde{E}^\bullet(B)$, making the usual diagrams commute. Here $p_B: E_B\rightarrow B$ is a principal $\PGL(N^+ +1)$-bundle over a scheme $B$, $\theta_B: E_B \rightarrow \hgammadetailx$ is a $\PGL(N^+ +1)$-equivariant morphism, $b: B\rightarrow S$ is a morphism of schemes, and $\lambda: E_B(\PGL(N+1)) \rightarrow b^\ast E$ is an isomorphism of principal $\PGL(N+1)$-bundles. 

Note that by the uniqueness property of corollary \ref{K324} the subbundle $E_B$ maps under $\lambda$ to the subbundle $b^\ast E^+$, compare also the proof of proposition \ref{P6014}.  Since $E_B /\PGL(N^+ +1) = B$, there is an induced composed morphism 
\[ h : \:\: B \rightarrow b^\ast E^+/\PGL(N^+ +1) \rightarrow \tilde{E}, \]
as desired. It is straightforward to verify that this assignment defines an isomorphism 
\[ \tilde{E}^\bullet  \:\: \cong \:\: \cmgammadetailx \times_{\cdx} S^\bullet \]
between stacks. 
\ebew

\begin{rk}\em\label{R61}
Let $\phi: \tilde{E}\rightarrow  \cmgammadetailx$ be the morphism constructed in the proof of lemma \ref{L61}. The $m/\Gamma(C_0^-;\calq)$-pointed stable curve of genus $g^+$ corresponding to $\phi$ can be constructed as follows. Consider the morphism $\theta^+: E^+\rightarrow \hgammadetailx$. By the universal property of the Hilbert scheme this determines an embedded stable curve $C_{E^+} \rightarrow E^+$. By the $\PGL(N^+ +1)$-equivariance of $\theta^+$ the quotient 
\[  C_{\tilde{E}} := C_{E^+}/\PGL(N^+ +1) \rightarrow E^+/\PGL(N^+ +1) = \tilde{E} \]
exists as an $m/\Gamma(C_0^-;\CQ)$-pointed stable curve of genus $g^+$ over $\tilde{E}$, and this curve is represented by the morphism $\phi$. 

In fact, since the action of $\Aut(C_0^-)$ on $\hgammadetailx$ is trivial, the quotient curve 
\[ C_S^+ := C_{\tilde{E}} / \Aut(C_0^-) \rightarrow \tilde{E}/\Aut(C_0^-) = S \] 
exists as an $m/\Gamma(C_0^-;\CQ)$-pointed stable curve of genus $g^+$ over $S$. Glueing this curve with the trivial curve $C_0^-\times S \rightarrow S$ along the sections of marked points gives a stable curve $C_S\rightarrow S$ of genus $g$ over $S$, which is the curve corresponding to the  morphism $s: S \rightarrow \cdx$. Formally we have thus
\[ \Omega\circ\phi = s \circ\tilde{p} ,\]
confirming the commutativity of the diagram. 
\end{rk}

The above lemma \ref{L61} extends to the partial compactifications. Let $\hat{\Omega}$ be the morphism of remark \ref{R0526}. 

\begin{lemma}\label{L6107}
The morphism $\hat{\Omega} : \cmgammadetailxxbar  \rightarrow \FDPxx$ is representable, surjective, finite and unramified,  and its  degree is equal to the order of the group $\Aut(C_0^{--})$. 
\end{lemma}

\proof
The proof is completely analogous to that of lemma \ref{L61}. As we did there, one can construct for a given triple $(E,p,\theta)\in\FDPxx(B)$ a Cartesian diagram
\[ \diagram
\hat{E}\rrto^{\hat{p}} \dto_\phi&& S\dto^s\\
\cmgammadetailxxbar \rrto_{\hat{\Omega}} && \FDPxx.
\enddiagram \]
Here $\hat{E}:= E^{++}/\PGL(N^{++}+1)$ as a principal $\Aut(C_0^{--})$-bundle over $S$, where $E^{++}$ is the principal $\PGL(N^{++}+1)\times\Aut(C_0^{--}) $-subbundle of $E$ constructed in the proof of theorem \ref{49}. 
\ebew

\begin{rk}\em
The $\mu/\Gamma(C_0^{--};\CR)$-pointed prestable curve of genus $g^{++}$ over $\hat{E}$ corresponding to the morphism $\phi \in  \cmgammadetailxxbar(\hat{E})$ can be described analogously as in   remark \ref{R61}. 
\end{rk}

\begin{notation}\em
Let $\hat{M}_{g^+,m/\Gamma(C_0^-;\CQ)}^{\: \dagger}$\label{n57} denote the preimage of the partial compactification $\dxxbar$ of $\dx$ at the point $[C_0]$ under the isomorphism $\overline{M}_{g^+,m/\Gamma(C_0^-;\CQ)}\rightarrow \dpobar$. 

Let $\hat{H}^{\: \dagger}_{g^+,m(\Gamma(C_0^-;\CQ))}$\label{n56} denote the reduction of its preimage in $\overline{H}_{g^+,m/\Gamma(C_0^-;\CQ)}$.
\end{notation}

The following proposition is very useful. It allows us to decide whether the description of $\cdx$ as a quotient stack as in theorem \ref{MainThm} extends to the boundary of $\FDP$ or not.

\begin{prop}\label{P62}
The isomorphism 
\[ \Lambda: \quad \left[\hgammadetailx / \PGL(N^++1)\times\Aut(C_0^-)\right] \longrightarrow \cdx \]
 extends to an isomorphism from  $[\hat{H}^{\: \dagger}_{g^+,m/\Gamma(C_0^-;\CQ)}/ \PGL(N^++1)\times\Aut(C_0^-)]$ to $  \FDPxx$ if and only if $C_0^+$ is invariant with respect to the action of $\Aut(C_0)$ as an $m/\Gamma(C_0^-;\calq)$-pointed curve. 
\end{prop}

\proof
If $C_0^+$ is invariant, then $C_0^+ = C_0^{++}$ as  $m(\Gamma(C_0^-;\calq))$-pointed curve, and hence   $[\hat{H}^{\: \dagger}_{g^+,m/\Gamma(C_0^-;\CQ)}/ \PGL(N^++1)\times\Aut(C_0^-)]$ is equal to the quotient stack  $[\hat{{H}}^\times_{g^{++},\mu/\Gamma(C_0^{--};\CR)}/ \PGL(N^{++}+1)\times\Aut(C_0^{--})]$. Thus the claim follows from theorem \ref{49}. 

Suppose now that $C_0^+$ is not invariant. For each $\gamma\in \Aut(C_0)$ choose a different point $x_\gamma\in \PP^1$. Define $U_\gamma$ as the open complement of the point $x_\gamma$ in $\PP^1$, and put $U_0:= \bigcap\limits_{\gamma\in\Aut(C_0)} U_\gamma$. We now glue the trivial families $f_\gamma:= \pr_2: \:  C_0\times U_\gamma\rightarrow U_\gamma$ over $U_0$ using the automorphisms $\gamma$ fibrewise. This defines a stable curve $f: C\rightarrow \PP^1$, for which there exists no subcurve $f^+: C^+\rightarrow \PP^1$ of genus $g^+$. Thus there is no curve in $
[\hat{H}^{\: \dagger}_{g^+,m/\Gamma(C_0^-;\CQ)}/ \PGL(N^++1)\times\Aut(C_0^-)](\PP^1)$ which maps to $f:C\rightarrow \PP^1$ under the morphism $\Lambda$. Therefore $\Lambda$ cannot be an equivalence of categories. 
\ebew

\begin{rk}\em 
$(i)$ Let the number of glueing points between $C_0^+$ and $C_0^-$ be $m=1$. Then by remark \ref{37} and proposition \ref{312} the curve $\overline{D}(P_1)$ is isomorphic to $\msbar_{1,1}$, so there is exactly  one point of $\overline{D}(P_1)$ representing a curve with exactly $3g-3$ nodes. Hence
\[ \hat{D}(P_1) =  \overline{D}(P_1) \]
holds in this case, and therefore we have in particular $\FDPhat = \FDP$. \\
$(ii)$ If $m=2$, then there are always two such points on $\overline{D}(P_1)$, corresponding to the two points of $\Delta_0 \subset \msbar_{1,2}$ representing the two singular, reducible curves of genus $1$ with two nodes. These two points are fixed under any permutation in $\Sigma_2$. 
In particular, there are always two partial compactifications needed to cover $\overline{D}(P_1)$. \\
$(iii)$ If $m=3$, then $\overline{D}(P_1)$ is a finite quotient of $\msbar_{0,4}$, which has 3 boundary points. Depending on the symmetries of $C_0$, there may be $1$, $2$ or $3$ points on $\overline{D}(P_1)$ representing stable curves with $3g-3$ nodes.
\end{rk}

\begin{rk}\em
For all $m=1,2$ or $4$, the subscheme $\msbar_{g^+,m}'$ of the moduli space $\msbar_{g^+,m}$ is smooth and irreducible. Its preimage $\overline{H}_{g^+,n,m}'$ in the Hilbert scheme is also  smooth and irreducible. Since the canonical morphism $\overline{H}_{g^+,n,m}'\rightarrow \cmgprimebar = [\overline{H}_{g^+,n,m}'/\PGL(N^+ +1)] $ is an atlas in the sense of an Artinian stack, one can conclude that the stack $\cmgprimebar$ is smooth and irreducible as well, and of dimension one. 

Note that in general neither the scheme $\mxxbar$, nor $\hxxbar$, nor the stack $\cmxxbar$ are smooth or irreducible. 
\end{rk}

By the definition of a Deligne-Mumford stack, for any given morphism $t: T\rightarrow \cmgammadetailx$, where $T$ is a scheme, the composition ${\Omega}\circ t : T\rightarrow \cdx$ is representable, even if ${\Omega}$ is not representable. The following proposition shows explicitely how to construct a representing morphism of schemes.

\begin{prop}\label{P6112}
Let $S$ and $T$ be schemes. Let  $s: S\rightarrow \cdx$ and $t: T\rightarrow  \cmgammadetailx$  be  morphisms, which are represented by the triples $(E_s,p_S,\theta_S)\in  \cdx(S)$ and $(E_T,p_T,\theta_T)\in\cmgammadetailx(T)$, respectively. Let $\Theta: {H}^\times_{g^+,n,m/\Gamma(C_0^-;\CQ)}\rightarrow {H}_{g,n,0}$ be the morphism constructed earlier. Put 
\[ R := E_T \times_{\Theta\circ\theta_T, {H}_{g,n,0}, \theta_S} E_S ,\]
and let $\overline{p}_S : R/\PGL(N^+ +1)\rightarrow S$ and $\overline{p}_T : R/\PGL(N^+ +1)\rightarrow T$ denote the projections induced by $p_S$ and $p_T$. 
Then the diagram
\[ \diagram
R/\PGL(N^+ +1) \xto[rrr]^{\overline{p}_S} \dto_{
\overline{p}_T} &&& S\dto^s\\
T \rto_t & \cmgammadetailx  \rrto_{{\Omega}} && \cdx
\enddiagram \]
is Cartesian.
\end{prop}

\proof
By definition, the scheme $R$ can be described as a closed subscheme of $E_T\times E_S$ by 
\[ R = \{ (e_T,e_S) \in E_T\times E_S: \:\: \Theta\circ\theta_T(e_T) = \theta_S(e_S)\} .\]
The defining equation is invariant under the diagonal action of $\PGL(N^+ +1)$ on $E_T\times E_S$, so there is an induced action of $\PGL(N^+ +1)$ on $R$, and the quotient $R':=R/\PGL(N^+ +1)$  exists as a scheme. 

To determine the triple representing the composed morphism $s\circ\overline{p}_S$, consider the diagram
\[\diagram
\overline{p}_S^\ast E_S \rto^{\tilde{p}_S} \dto & E_S \rto^{\theta_S} \dto^{p_S} &  {H}_{g,n,0}\\
R' \rto_{\overline{p}_S} & S .
\enddiagram\]
Note that the pullback bundle $ \overline{p}_S^\ast E_S$ is trivial. 
Thus the triple associated to $s\circ\overline{p}_S$ is 
\[ (R'\times\PGL(N+1),\pr_1,\sigma(\theta_S\circ\tilde{p}_S,\id_{\PGL(N+1)})) \in\cdx(R'), \]
where $\sigma$ denotes the action of $\PGL(N+1)$ on $ {H}_{g,n,0}$. 

Similarly, consider the diagram
\[\diagram
\overline{p}_T^\ast E_T(\PGL(N+1)) \rto^{\tilde{p}_T}\dto& E_T(\PGL(N+1)) \dto\rto& \hgammadetailx\\
R'\rto_{\overline{p}_T} & T, 
\enddiagram\]
where $E_T(\PGL(N+1))$ denotes the $\PGL(N+1)$-extension of $E_T$ with respect to the fixed embedding of $\PGL(N^+ +1)$ into $\PGL(N+1)$. From this we obtain the triple
\[ ( R'\times\PGL(N+1), \pr_1,\sigma(\Theta'\circ \overline{\theta}_T\circ\tilde{p}_T, \id_{\PGL(N+1)})) \in \cdx(R'),\]
representing the composed morphism ${\Omega}\circ t\circ\overline{p}_T$. Here $\overline{\theta}_T$ denotes the $\PGL(N+1)$-extension of $\theta_T$. 
By comparing the two triples we find that the diagram of the proposition commutes. Indeed, the defining equation of $R$ implies the identity $\Theta\circ \overline{\theta}_T\circ\tilde{p}_T= \theta_S\circ\tilde{p}_S$. 

In particular, by the universal property of the fibre product, there is a morphism of stacks
\[ R' \rightarrow T \times_{\cdx} S .\] 
By the representability of ${\Omega}\circ t$ this is even a morphism of schemes. We have to construct an inverse  morphism 
\[ T \times_{\cdx} S \rightarrow R' .\] 
Consider an object $(f,g,\eta)\in  (T \times_{\cdx} S)^\bullet (B)$ for some scheme $B$. Here, $f: B\rightarrow T$ and $g: B\rightarrow S$ are morphisms of schemes, and $\eta$ is an isomorphism between the triple representing $s\circ g$ and ${\Omega}\circ t\circ f$ in $\cdx(B)$. Thus $\eta$ is given by an isomorphism of principal $\PGL(N+1)$-bundles
\[ \lambda : \:\: f^\ast E_T(\PGL(N+1)) \rightarrow g^\ast E_S \]
over $B$, which is compatible with the morphisms into ${H}_{g,n,0}$. 

Let us at first assume that $f^\ast E_T$ is trivial over $B$. Then there is a section $\nu: B \rightarrow f^\ast E_T$, and by composition with the natural morphism  $\overline{f}: f^\ast E_T\rightarrow E_T$, we obtain a morphism
\[ \tilde{f}:= \overline{f}\circ \nu : \:\:\:  B \rightarrow E_T .\]
Since there is an inclusion $E_T \subset E_T(\PGL(N+1))$, we have furthermore a morphism 
\[ \tilde{g}:= \overline{g}\circ \lambda\circ\nu: \:\: \:  B \rightarrow E_S ,\]
where  $\overline{g}: g^\ast E_S\rightarrow E_S$ denotes the natural morphism. 

By assumption, $\lambda$ satisfies the compatibility condition 
\[ \theta_S\circ \overline{g}\circ \lambda \: = \: \Theta\circ \theta_T \circ\overline{f} \]
on $f^\ast E_T$. 
From this we derive the identity 
\[ \theta_S \circ \tilde{g} \: = \: \Theta\circ\theta_T\circ \tilde{f} , \]
which shows that the pair $(\tilde{f},\tilde{g})$ is indeed an object of $(R')^\bullet(B)$. 

Consider now the case where $f^\ast E_T$ is not globally trivial. We may choose an \'etale covering $\{B_\alpha\}_{\alpha\in A}$ of $B$, such that the restricted subbundles $f^\ast E_T|B_\alpha$ are trivializable. The local trivializations are only determined up to the action of $\PGL(N^+ +1)$. Therefore the morphism $\tilde{f}_\alpha: B_\alpha\rightarrow E_T$ constructed above, and thus the pair  $(\tilde{f}_\alpha,\tilde{g}_\alpha)$, is only determined up to the action of $\PGL(N^+ +1)$. However, the pair $(\tilde{f}_\alpha,\tilde{g}_\alpha)$ determines a well-defined object of $(R/\PGL(N^+ +1))^\bullet(B_\alpha)$. Clearly, these objects glue to give a global object of $(R/\PGL(N^+ +1))^\bullet(B)$. It is straightforeward to verify that this construction gives the desired inverse morphism.
\ebew 

\begin{rk}\em
An alternative, and more geometric argument for the commutativity of the diagram of proposition \ref{P6112} can be given as follows. Let $C_T\rightarrow T$ denote the $m$-pointed stable curve of genus $g^+$ corresponding to $t\in\ \cmgammadetailx(T)$, and let $C_S\rightarrow S$  be the stable curve of genus $g$ corresponding to $s\in\cdx(S)$. Then the composed morphism $s\circ\overline{p}_S$ corresponds to the curve $\overline{p}_S^\ast C_S \rightarrow R'$, and ${\Omega}\circ t\circ \overline{p}_T$ corresponds to $\overline{p}_T^\ast( C_T + (C_0^-\times T)) \rightarrow R'$, where ``+''  denotes glueing of $C_T\rightarrow T$ with the trivial curve $C_0^-\times T\rightarrow T$ along the sections of marked points. The defining equation of $R$ implies that there is an isomorphism $\overline{p}_T^\ast( C_T + (C_0^-\times T)) \cong \overline{p}_S^\ast C_S$ of stable curves over $R'$. 
\end{rk}

As applications of proposition \ref{P6112}, we want to prove the following corollaries, which are of course both implied by lemma \ref{L61}. 

\begin{cor}\label{C6114}
The morphism ${\Omega}: \cmgammadetailx \rightarrow \cdx$ is surjective.
\end{cor}

\proof
Consider the commutative diagram of proposition \ref{P6112} with $T:= \hgammadetailx$ and $S:=  \kx$. By corollary \ref{314} we have $\kx = \Theta( \hgammadetailx) \cdot \PGL(N+1)$. Therefore, for each $x\in \kx$ there exists a point $y\in  \hgammadetailx$ and an element $\gamma\in\PGL(N+1)$, such that $x= \Theta(y)\cdot \gamma$. Thus for each $e_x\in E_S$ with $p_S(e_x)=x$, there exists an $e_y\in E_T$ with $p_T(e_y)=y$, and an element $\gamma\in\PGL(N+1)$,  such that $\Theta\circ\theta_T(e_y) =\theta_S(e_x\cdot \gamma)$. In other words, the morphism $\overline{p}_S: R/\PGL(N^+ +1) \rightarrow S$ is surjective, and hence the composed morphism from $ R/\PGL(N^+ +1)$ to $\cdx$ is surjective. Since ${\Omega}\circ t\circ \overline{p}_T = s\circ\overline{p}_S$,  the morphism ${\Omega}$ must be surjective as well.
\ebew

\begin{cor}\label{C6115}
The morphism ${\Omega}: \cmgammadetailx \rightarrow \cdx$ is \linebreak quasi-finite.
\end{cor}

\proof
Both stacks $\cmgammadetailx$ and $\cdx$ are irreducible and one-dimensional. Hence there are \'etale atlases $t: T\rightarrow \cmgammadetailx$ and $s: S\rightarrow \cdx$, where $S$ and $T$ are irreducible curves, and they can even be assumed to be smooth. Thus the induced morphism $\overline{p}_S: R/\PGL(N^+ +1) \rightarrow S$ of proposition \ref{P6112} is a surjective morphisms between irreducible curves, and hence finite. Since $R/\PGL(N^+ +1) \rightarrow \cmgammadetailx $ is an \'etale atlas, we are done. 
\ebew

\begin{rk}\em\label{R6116}
The finiteness of the morphism $\overline{p}_S: R/\PGL(N^+ +1) \rightarrow S$ in corollary \ref{C6115} can be understood geometrically as follows. Let the morphisms $s: S\rightarrow \cdx$ and $t:T\rightarrow \cmgammadetailx$ correspond to stable curves $C_S\rightarrow S$ and $C_T\rightarrow T$, respectively. A point $x\in S$ determines a fibre $C_{S,x}\rightarrow \Spec(k)$ as a stable curve of genus $g$. By corollary \ref{C6114} there exists an element $e_x\in E_{S,x}$ in the fibre of $p_S: E_S\rightarrow S$ over $x$, such that $\theta_S(e_X) \in \overline{K}(P_1)$. 

The composed morphism $T\rightarrow\cmgammadetailx\rightarrow \mgammadetailx$ is a surjective morphism of irreducible curves, and thus finite. Hence there are only finitely many points $y_1,\ldots,y_k\in T$, such that for the corresponding fibres $C_{T,y_i}$ of $C_T\rightarrow T$ holds $C_{T,y_i} + C_0^- \cong C_{S,x}$. Here again, ``+''  denotes glueing in the marked points. Fix one $1\le i\le k$, and let $E_{T,y_i}$ denote the fibre of $E_T\rightarrow T$ over $y_i$. 
For all $e\in E_{T,y_i}$
there exists an element $\gamma_e \in\PGL(N +1)$, such that 
\[ \Theta\circ\theta_T(e) = \theta(e_x)\cdot \gamma_e \]
holds. 
If the curve $C_{S,x}$ has exactly $3g-4$ nodes, then up to the action of $\PGL(N^+ +1)$, the element $\gamma_e$ is uniquely determined up to multiplication with  an element of $\Aut(C_0^-)$. In this case, the fibre of $\overline{p}_S: R/\PGL(N^+ +1) \rightarrow S$ over $x$ is isomorphic to $\Aut(C_0^-)^k$ as a set. Note that the number $k$ need not be constant, but may vary with $x$. 

If the curve $C_{S,x}$ has $3g-3$ nodes, then there may be more than one embedding of $C_{T,y_i}$ into $C_{S,x}$. In this case, as a set the fibre of $\overline{p}_S$ over $x$ is still finite, and the numer if its elements is bounded by the order of  $\Aut(C_{S,x})$.
\end{rk}

\begin{rk}\em
For any subgroup $\Gamma\subset \Sigma_m$, an $m$-pointed stable curve of genus $g^+$ defines obviously also an $\mgamma$-pointed stable curve of the same genus. Thus there is a natural morphism of stacks
\[\cmgplus \rightarrow \cmgbarquot .\]
Recall that both stacks have a presentation as quotient stacks by proposition \ref{28} and proposition \ref{P217}. Since $\Gamma$ acts freely on $\overline{H}_{g^+,n,m}$, and since its action commutes with the action of $\PGL(N^+ +1)$, the presentations as quotient stacks immediately give an isomorphism of stacks
\[ \cmgplus  / \Gamma \:\: \cong \:\: \cmgbarquot .\]
In particular, the morphism $\cmgplus \rightarrow \cmgbarquot$ is  surjective, unramified and finite of degree equal to the order of $\Gamma$.
\end{rk}

\begin{rk}\em\label{R6117a}
By the above observations, there is in particular a finite and surjective morphism
\[   \omega: \:\:\: \cmxbar \rightarrow \cdx .\]
Note that there is an exact sequence of groups
\[ \id \rightarrow \Aut(C_0^-)\rightarrow \Aut(C_0^-;\CQ)\rightarrow \Gamma(C_0^-;\CQ)\rightarrow \id.\]
Hence the degree of $\omega$ equals the order of the group $\Aut(C_0^-;\CQ)$. 

If $C_0^+$ is invariant with respect to the action of $\Aut(C_0)$ as an $m/\Gamma(C_0^-;\CQ)$-pointed subcurve of $C_0$, then $\Aut(C_0^+)$ is a normal subgroup of $\Aut(C_0)$. Indeed, there is an exact sequence of groups
\[ \id \rightarrow \Aut(C_0^+)\rightarrow \Aut(C_0)\rightarrow\Aut(C_0^-;\CQ)\rightarrow \id\]
in this case. 
\end{rk}

\begin{rk}\em
The morphism $\omega: \cmxbar \rightarrow \cdx$ is in general not representable. Thus the morphism ${\Omega}$ is also in general not representable. 

This can be seen as follows. 
Let $s: S\rightarrow \cdx$ be an \'etale atlas of $\cdx$, where $S$ is a smooth and irreducible curve. An object of $\cmxbar\times_{\cdx} S^\bullet(\Spec(k))$ is a triple $(C^+,x,\lambda)$, consisting of an $m$-pointed stable curve $C^+\rightarrow \Spec(k)$ of genus $g^+$, a geometric point $x:\Spec(k)\rightarrow S$ of $S$, and an isomorphism $\lambda: C^+ + C_0^- \rightarrow C_{S,x}$, if $C_S\rightarrow S$ is the curve corresponding to the  atlas $s: S\rightarrow \cdx$, and $C_{S,x}$ its fibre over $x$. Note that as abstract curves, $C^+$ and $C_{S,x}$ are only defined up to an isomorphism. Hence the isomorphisms $\lambda$ are in one-to-one correspondence with the elements of the quotient $\Aut(C_{S,x})/ \Aut(C^+)$, for some fixed embedding $C^+ \subset C_{S,x}$. 

Suppose that $\omega$ is representable. Then $T:= \cmxbar\times_{\cdx} S$ is a scheme, and there is a canonical morphism $\tau: T\rightarrow S$ representing $\omega$. By remark \ref{R6117a}, this is a finite morphism of degree equal to the order of $\Aut(C_0^-;\CQ)$ between two smooth curves. We saw above that the fibre of $\tau$ over a point $x: \Spec(k)\rightarrow S$ is bijective to the set $\Aut(C_{S,x})/ \Aut(C^+_x)$, where $C_x^+$ denotes the closure of the complement of some embedding of $C_0^-$ into $C_{S,x}$. However, in general  there exist stable curves of genus $g$, represented by points of $\dx$, where the number of elements of the quotient $\Aut(C_{S,x})/ \Aut(C^+_x)$ exceeds the number of elements of $\Aut(C_0^-;\CQ)$. This is a contradiction. 
\end{rk}

\begin{rk}\em\label{R722}
Let us conclude with a short summary. Consider an irreducible component $\overline{D}$ of the one-dimensional boudary stratum $\msbar_g^{(3g-4)}$ of the moduli space $\msbar_g$. There exists a stable curve $C_0$ of genus $g$ with $3g-3$ nodes, together with an enumeration of its nodes, such that 
\[ \overline{D}=\overline{D}(C_0;P_1) ,\]
or shorter, $ \overline{D}=\overline{D}(P_1)$. Let $\tilde{\cal D}$ denote the substack  of the moduli stack $\cmgbar$ defined by the  Cartesian diagram
\[ \diagram 
\tilde{\cal D} \rto \dto & \cmgbar\dto\\
\overline{D}\rto& \msbar_g.
\enddiagram\]
The underlying reduced substack $\overline{\cal D}$ 
is  irreducible and of dimension one. 

The node $P_1$ distinguishes a subcurve $C_0^+$ of $C_0$ as an $m/\Gamma(C_0^-;\calq)$-pointed stable curve of some genus $g^+$. There is a closed subscheme $\msbar_{g^+,m}'$ of the moduli space $\msbar_{g^+,m}$, which is a smooth irreducible curve, so that there is an isomorphism
\[ \overline{D} \:\: \cong \:\: \msbar_{g^+,m/\Gamma(C_0;\calq)}' .\]
In general, the corresponding stacks $\overline{\cal D}$ and $\cmgprimebarquot$ are not isomorphic. However, on dense open substacks, there exists a finite and surjective morphism of smooth Deligne-Mumford stacks
\[ \cmgammax \: \longrightarrow \:\: \cdx,\]
which is of degree equal to the order of the group $\Aut(C_0^-)$, and even 
representable and unramified. 

For the open substack $\cdx$ one has an isomorphism of quotient stacks
\[ \cdx \:\: \cong \:\: \left[ \hgammadetailx /\PGL(N^+ +1)\times\Aut(C_0^-)\right] .\]
Using the isomorphism $[ \hgammadetailx /\PGL(N^+ +1)]\cong \cmgammadetailx$, and the interpretation of a quotient stack as a stack quotient, as explained in the appendix, we can write this as an isomorphism of stacks
\[  \cdx \:\:\: \cong \:\:\: \cmgammadetailx / \Aut(C_0^-) .\]

Consider the partial compactification $\hat{D}$ of $D$ at the point representing the curve $C_0$, i.e. $\hat{D}= \hat{D}(C_0;P_1) = D(C_0;P_1)\cup\{[C_0]\}$. 
Let $\hat{\cal D}$ denote the reduced substack of $\overline{\cal D}$ defined by $\hat{D}$. 
Then the above isomorphism can be suitably extended  to $\hat{\cal D}$
if and only if the subcurve $C_0^+$ is invariant with respect to the action of $\Aut(C_0)$, up to a reordering of its marked points by a permutation in $\Gamma(C_0^-;\calq)$.  

In general, over  the partial compactification $\hat{D}$, there is an isomorphism of open substacks
\[  \FDPxx \:\:\: \cong \:\:\:  \cmgammadetailxxbar / \Aut(C_0^{--}) .\]
Note that for certain boundary components $\overline{D}$ this stack may be empty.  

The curve $\overline{D}\subset\msbar_g$ has finitely many boundary points, i.e. points representing stable  curves with the maximal number of nodes. Hence there are finitely many stable curves $C_1,\ldots,C_k$ of genus $g$ with $3g-3$ nodes, together with nodes $P_{i,1}\in C_i$ for $i=1,\ldots ,k$, such that
\[ \overline{D}=\overline{D}(C_0;P_1) =\overline{D}(C_1;P_{1,1}) = \ldots = 
\overline{D}(C_k;P_{k,1}),\]
and such that the partial compactifications $\hat{D}(C_i;P_{i,1})$ for $i=1,\ldots ,k$, cover  $\overline{D}$. In particular, the corresponding open substacks ${\hat{\cal D}(C_i;P_{i,1})}$ cover all of the stack $\FDP$. 
\end{rk}

%

%

\chapter{Examples and Outlook}

\section{The case of genus three}\label{Sect81}

To illustrate the  description of one-dimensional boundary substacks of the moduli space of Deligne-Mumford stable curves as given above, we want to consider the case of genus $g=3$ in detail. 

\begin{rk}\em\label{R8301}
There are exactly five different stable curves of genus $g=3$ with $3g-3=6$ nodes. We denote them by $C_1,\ldots,C_5$, in the same order as they are listed in Faber's paper \cite{Fa}. Schematically they can be represented as in the following pictures.

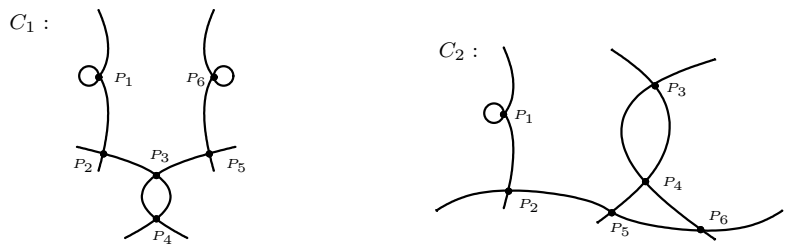
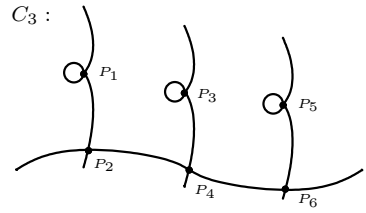
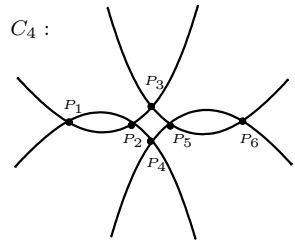
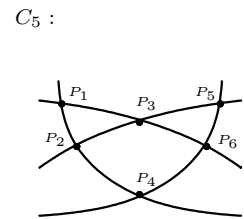
\begin{center}
\setlength{\unitlength}{0.00016667in}
\begingroup\makeatletter\ifx\SetFigFont\undefined%
\gdef\SetFigFont#1#2#3#4#5{%
  \reset@font\fontsize{#1}{#2pt}%
  \fontfamily{#3}\fontseries{#4}\fontshape{#5}%
  \selectfont}%
\fi\endgroup%
{\renewcommand{\dashlinestretch}{30}

}
\end{center}

\end{rk}

\begin{rk}\em
There are eight different types of stable curves of genus $3$ with exactly $5$ nodes. Following Faber's notation, we denote their types by $(a),\ldots,(h)$. They are represented by the pictures below.

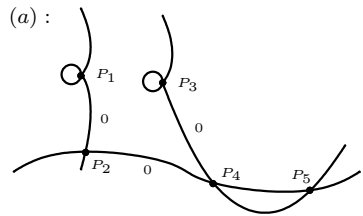
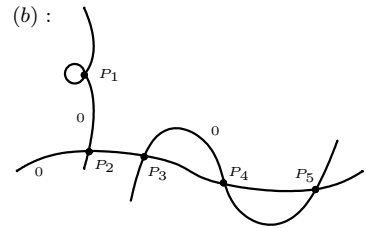
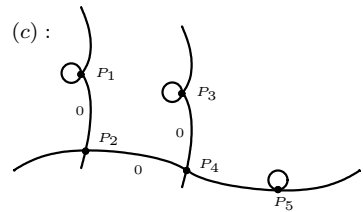
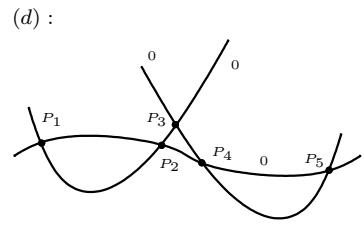
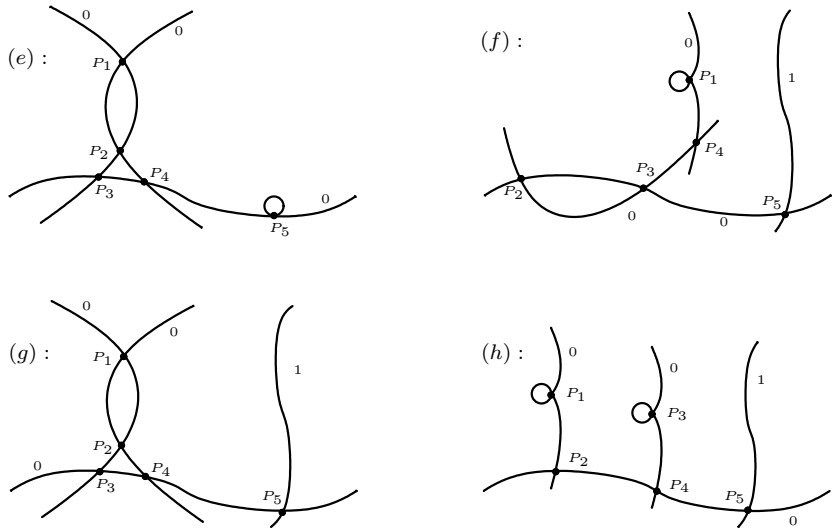
\begin{center}
\setlength{\unitlength}{0.00016667in}
\begingroup\makeatletter\ifx\SetFigFont\undefined%
\gdef\SetFigFont#1#2#3#4#5{%
  \reset@font\fontsize{#1}{#2pt}%
  \fontfamily{#3}\fontseries{#4}\fontshape{#5}%
  \selectfont}%
\fi\endgroup%
{\renewcommand{\dashlinestretch}{30}
%
}
\end{center}

The number attached to an irreducible component gives its geometric genus.
\end{rk}

\begin{rk}\em
The type of a curve containing five nodes does not vary in the open part $D$  of an  irreducible component $\overline{D}$  of the boundary stratum $\msbar_3^{(5)}$. Each such component $\overline{D}$ can be written as $\overline{D}= \overline{D}(C_i \, ;P_j)$ for some $1\le i\le 5$ and  $1\le j\le 6$, but this presentation of $\overline{D}$ is not unique. The following table lists for each curve $C_1,\ldots,C_5$ all boundary components distinguished by such deformations of this curve, which preserve all but one node.
\[\begin{array}{lccccc}
C_i \qquad & P_j&m & D(C_i\, ;P_j) & \Aut(C_i^-) & \Gamma(C_i^-;\CQ) \\
\hline\\
C_1&P_1,P_6& 1& (f)&4& 1\\
   &P_2,P_5& 2& (a)&2& 2\\
   &P_3,P_4& 2& (c)&4& 2\\
C_2&P_1    & 1& (g)&4& 1\\
   &P_2    & 2& (e)&2& 2\\
   &P_3,P_4& 2& (a)&2& 2\\
   &P_5,P_6& 4& (b)&2& 6\\
C_3&P_1,P_3,P_5&1&(h)&8&1\\
   &P_2,P_4,P_6&2&(c)&4&2\\
C_4&P_1,P_2,P_5,P_6&2&(e)&2&2\\
   &P_4,P_5&        4&(d)&1&8\\
C_5&P_1,\ldots,P_6&4&(d)&1&8 
\end{array}\]
The first column states  the curve $C_i$ we are considering. The second column lists all equivalent nodes on the curve $C_i$, where the equivalence is taken with respect to the action of $\Aut(C_i)$. The third column states the number of marked points of $C_i^-$. The fourth column tells the type of a generic curve in $\overline{D}(C_i;P_j)$. The fifth and sixth column give the orders of the associated groups of automorphisms. 
\end{rk}

\begin{rk}\em
In the next table we will give a complete list of all possible presentations of all irreducible components $\overline{D}$ of $\msbar_3^{(5)}$.
\[\begin{array}{lcccc}
\overline{D} \qquad & \msbar_{g^+,m}'& \Text{cover}& \Text{boundary}&\Text{neighbours}\\
\hline\\ 
(a)&\msbar_{1,2}'& 4& C_1&(c),(f)\\
 &&&C_2&(b),(e),(g)\\
(b)&\msbar_{0,4}&12&C_2&(a),(e),(g)\\
(c)&\msbar_{1,2}'&8&C_1&(a),(f)\\
 &&&C_3&(h)\\
(d)&\msbar_{0,4}&8&C_4&(e)\\
 &&&C_5&-\\
(e)&\msbar_{1,2}'&4&C_2&(a),(b),(g)\\
 &&&C_4&(d)\\
(f)&\msbar_{1,1}&4&C_1&(a),(c)\\
(g)&\msbar_{1,1}&4&C_2&(a),(b),(e)\\
(h)&\msbar_{1,1}&8&C_3&(c)
\end{array}\]
The first column specifies the irreducible component $\overline{D}$  by indicating the type of a general point. The second column gives the moduli scheme $\msbar_{g^+,m}'$ covering $\overline{D}$ as in proposition \ref{312}, and the third column states the degree of the finite morphism $\cmxbar \rightarrow \cdx$ as in remark \ref{R6117a}. In column four, the  boundary points of $\overline{D}$ are listed, i.e. points corresponding to stable curves $C_i$ with $6$ nodes. The irreducible components meeting $\overline{D}$ in $[C_i]$ are given in column five, again by the types of their general curves.  
\end{rk}

\begin{ex}\rm
Consider the irreducible component $\overline{D}= \overline{D}(C_2;P_1)$. A general point of this component represents a curve of type $(g)$. Since $m=1$, all relevant groups of automorphisms split. In particular we have the equality $D^\times(C_2;P_1) = D(C_2;P_1)$. Note that the subcurve $C_2^+$ determined by $P_1$ is invariant with respect to the action of $\Aut(C_2)$. Thus the isomorphism  ${{\cal M}_{g^+,m/\Gamma(C_2^-;\CQ)}'} /\Aut(C_2^-) \cong {\cal D}(C_2;P_1)$   extends to an isomorphism over the partial compactification at the point $[C_2]$, see proposition \ref{P62}. Since there is only one boundary point on $\overline{D}(C_2;P_1)$, the partial compactification $\hat{D}(C_2;P_1)$ is in fact equal to  
$\overline{D}(C_2;P_1)$. Thus we have an isomorphism of stacks
\[ {\overline{\cal D}(C_2;P_1)}  \:\: \cong \:\: {\cal M}_{1,1} / {\mathbb Z}_2 \times {\mathbb Z}_2 .\]
\end{ex}

\begin{ex}\rm
A general point of the irreducible component $\overline{D}=\overline{D}(C_3; P_2)$ represents a curve of type $(c)$. There is an isomorphism
\[ {\cal D}^\times(C_3;P_2)  \: \: \cong \:\: {{\cal M}^\times_{1,2/\Sigma_2}} / {\mathbb Z}_2 \times {\mathbb Z}_2 .\]
Since the curve $C_3^+$ determined by the node $P_2$ is not invariant with respect to the action of $\Aut(C_3)$, this isomorphism does not extend to the partial compactification
$ {\hat{D}(C_3;P_2)}$.

The component $\overline{D}$ contains two boundary points, which are $[C_3]$ and $[C_1]$. Therefore one has an alternative presentation $\overline{D}=\overline{D}(C_1; P_3)$. Here, the subcurve $C_1^+$ determined by $P_3$ is invariant with respect to the action of $\Aut(C_1)$. 

However, the group of automorphisms $\Aut(C_1)$ does not split as a product  $\Aut_{\hat{\Gamma}}(C_1^+;\{P_2,P_5\}) \times \Aut(C_1^-)$. Therefore $\hat{D}{}^\times(C_1;P_3) = {D}{}^\times(C_1;P_3)$, and no additional information about the boundary point $[C_1]$ is gained. 
\end{ex}

\begin{ex}\em
Consider the irreducible component $\overline{D}=\overline{D}(C_2;P_5)$, with general point of type $(b)$. Again we have an equality $D^\times(C_2;P_5) = D(C_2;P_5)$. There is an isomorphism of stacks
\[ {\cal D}(C_2;P_5)  \: \: \cong \:\: {{\cal M}_{0,4/\Sigma_3}} / {\mathbb Z}_2 .\]
Note that there is only one boundary point on $\overline{D}$, so the partial compactification $\hat{D}$ is equal to $\overline{D}$. Because the subscheme $C_2^+$, determined by the node $P_5$, is not invariant with respect to the action of the group $\Aut(C_2)$, the above isomorphism does not extend to a global isomorphism, see proposition \ref{P62}.

The invariant subcurve $C_2^{++}$ determined by $P_5$ is an $1$-pointed stable curve of genus $2$, with four nodes. In particular, one has  $\hat{D}{}^\times(C_2;P_5)  =\hat{D}{}(C_2;P_5)  =\overline{D}(C_2;P_5)  $. One finds $\Aut(C_2^{--}) = {\mathbb Z}_2$ and $\Gamma(C_2^{--};\calr) = \{\id\}$. By theorem \ref{49}, there is a global isomorphism of stacks
\[ \overline{\cal D}(C_2;P_5)  \: \: \cong \:\: \hat{\cal M}{}^\times_{2,1} / {\mathbb Z}_2 .\]
\end{ex}

\begin{ex}\em
Finally we want to look at the case  of the irreducible component $\overline{D}$, where general points represent curves of type $(a)$. There are two boundary points on $\overline{D}$, which are $[C_1]$ and $[C_2]$. Thus there exist two presentations
\[ \overline{D}=  \overline{D}(C_2;P_3) =  \overline{D}(C_1;P_2). \]
Note that $D^\times(C_2;P_3) = D(C_2;P_3)$, as well as $\hat{D}{}^\times(C_2;P_3)  =\hat{D}{}(C_2;P_3) $.  
The subcurve $C_2^+$ defined by $P_3$ is invariant with respect to the action of $\Aut(C_2)$, so for the partial compactification at the point $[C_2]$ there is an isomorphism
\[ {\hat{D}(C_2;P_3)} \: \: \cong \:\: {\hat{\cal M}{}^\times_{1,2/\Sigma_2}} / {\mathbb Z}_2 .\]
The subcurve $C_1^+$ defined by $P_2$ is not invariant with respect to the action of $\Aut(C_1)$. Therefore this isomorphism cannot be extended to the global stack. 

Using invariant subschemes does not improve the situation either. Clearly one hat $C_2^{++} = C_2^+$, but one finds also $C_1^{++} = C_1$. Hence no new insight into ${\hat{\cal D}(C_1;P_2)}$ is gained. 
\end{ex}

\vfill

\pagebreak

\section{Higher-dimensional boundary strata}

There is a natural generalization of the results described so far. Let us briefly outline how the techniques developed above can be extended to describe boundary strata of the stack $\cmgbar$ of dimensions other  than one. There will be a subsequent paper to describe the situation in more detail.

\begin{abschnitt}\em
Let $f_0: C_0\rightarrow \Spec(k)$ be a fixed stable curve of genus $g$, with $3g-3$ nodes. Choose an enumeration $P_1,\ldots,P_{3g-3}$ of its nodes, and choose a number $0\le d \le 3g-3$. 

Let $f: C\rightarrow S$ be a deformation of $f_0: C_0\rightarrow \Spec(k)$ with irreducible and reduced base $S$, and with central fibre $C_{s_0}\cong C_0$ for some $s_0\in S$. Suppose that there are sections $\rho_{d+1},\ldots,\rho_{3g-3}: S\rightarrow C$, such that for all $i=d+1,\ldots,3g-3$ holds that\\
\hspace*{5mm} $(i)$ \quad for all $s\in S$, the point $\rho_i(s)$ is a node of the fibre $C_s$,  and\\
\hspace*{5mm} $(ii)$ \quad $\rho_i(s_0) = P_i$ under the above isomorphism.\\
Assume furthermore that a generic fibre $C_s$ of $f$ has no other nodes than $\rho_{d+1}(s),\ldots,\rho_{3g-3}(s)$. Then the induced morphism $\theta_f : S\rightarrow \msbar_g$ maps into the boundary stratum $\msbar_g^{(3g-3-d)}$. Since $S$ is irreducible, the morphism $\theta_f$ distinguishes one $d$-dimensional irreducible component of $\msbar_g^{(3g-3-d)}$, and this component shall be denoted by $\overline{D}(C_0;P_1,\ldots,P_d)$. 
Compare also lemma \ref{L31} for the one-dimensional case.
\end{abschnitt}

\begin{abschnitt}\em
Define a substack ${\overline{\cal D}(C_0;P_1,\ldots,P_d)}$ of $\cmgbar$ as the reduction of the  preimage substack of $\overline{D}(C_0;P_1,\ldots,P_d)$ under the canonical morphism $\cmgbar \rightarrow \msbar_g$. Analogously, let $D(C_0;P_1,\ldots,P_d)$ denote the open subscheme of $\overline{D}(C_0;P_1,\ldots,P_d)$, which  parametrizes  curves with exactly $3g-3-d$ nodes, and let ${{\cal D}(C_0;P_1,\ldots,P_d)}$ denote the corresponding open substack of the stack ${\overline{\cal D}(C_0;P_1,\ldots,P_d)}$.

Generalizing notation \ref{N35}, we define $C_0^+ = C_0^+(P_1,\ldots,P_d)$ as the subcurve of $C_0$ which is the union of all irreducible components of $C_0$ containing at least one of the nodes $P_1,\ldots,P_d$. We define $C_0^-$ as the closure of the complement of $C_0^+$ in $C_0$. The points of intersection of $C_0^+$ and $C_0^-$ are nodes of $C_0$, so they posess a natural ordering determined by the choosen enumeration of the nodes of $C_0$. Let $\calq = \{ Q_1,\ldots,Q_m\}$ be the ordered set of intersection points. In this way $C_0^+$ and $C_0^-$ can be considered as $m$-pointed prestable curves of some genus $g^+$ and $g^-$, respectively. In general neither of the two subcurves will be stable. However, in both cases the connected components are pointed stable curves, each of them with the maximal number of nodes possible. If $d=0$, then  we define $C_0^+ := \emptyset$ and $C_0^- := C_0$. 
\end{abschnitt}

\begin{abschnitt}\em
Clearly there is no finite list of all existing types of curves $C_0^+$ as in remark \ref{37}. The combinatorics of possible configurations complicates the description of $\overline{D}(C_0;P_1,\ldots,P_d)$ and its associated stack ${\overline{\cal D}(C_0;P_1,\ldots,P_d)}$. In particular the description of the boundary of $\overline{D}(C_0;P_1,\ldots,P_d)$, i.e. the locus of curves with more nodes than the $3g-3-d$ nodes of a generic curve, becomes quite involved. 

However, the ``boundary'' ${\overline{\cal D}(C_0;P_1,\ldots,P_d)} \smallsetminus  {{\cal D}(C_0;P_1,\ldots,P_d)}$ can be covered by closed substacks ${\overline{\cal D}(C_0;P_1,\ldots,\hat{P}_i,\ldots ,P_d)}$ of dimension $d-1$, with $i= 1,\ldots,d$. Here $\hat{P}_i$ indicates the omission of the $i$-th node $P_i$. Note that there is a natural chain of inclusions
\[ {\overline{\cal D}(C_0;\emptyset)} \subset 
{\overline{\cal D}(C_0;P_1)} \subset
{\overline{\cal D}(C_0;P_1,P_2)} \subset
\ldots \subset
{\overline{\cal D}(C_0;P_1,\ldots,P_{3g-3})} = \cmgbar.\]
The stack ${\overline{\cal D}(C_0;\emptyset)}$ can be thought of as the fibre of the canonical morphism $\cmgbar\rightarrow \msbar_g$ over the point $[C_0]$ representing the fixed curve $C_0$, compare also lemma \ref{LDNULL} below.  
\end{abschnitt}

\begin{abschnitt}\em
Let $\overline{H}_{g^+,n,m}^{ps}$ denote the subscheme of the Hilbert scheme, which para\-metrizes $m$-pointed prestable curves of genus $g^+$, whose connected components are stable. As in construction \ref{const_prop} one can define a morphism 
\[ \Theta^\times: \:\: \overline{ H}{}^\times_{g^+,n,m}\rightarrow \overline{H}_{g,n,0},\]
on a reduced open subscheme $\overline{ H}{}^\times_{g^+,n,m}$ of $\overline{H}_{g^+,n,m}^{ps}$, which is defined as the locus of curves for which the group of automorphisms splits appropriately. To construct the morphism, one first glues the universal curve $u^+: {\cal C}_{g^+,n,m}\rightarrow \overline{H}{}^\times_{g^+,n,m}$ over $\overline{ H}{}^\times_{g^+,n,m}$ along the $m$ sections of marked points to the trivial curve $\pr_2: C_0^-\times \overline{H}{}^\times_{g^+,n,m}\rightarrow 
\overline{H}{}^\times_{g^+,n,m}$ to obtain a stable curve $u_0: {\cal C}_0 \rightarrow 
\overline{H}{}^\times_{g^+,n,m}$ of genus $g$. Then one chooses a suitable embedding of this curve into $\PP^N\times\overline{H}{}^\times_{g^+,n,m}\rightarrow \overline{H}{}^\times_{g^+,n,m}$, which in turn induces  the desired morphism $\Theta^\times: \overline{H}{}^\times_{g^+,n,m}\rightarrow \overline{H}_{g,n,0}$ via the universal property of the Hilbert scheme. Again, we assume that $n$ is chosen large enough to work simultaneously for both  $\overline{H}_{g,n,0}$ and $ \overline{H}{}^\times_{g^+,n,m}$. 

The embedding can be chosen in such a way that the restriction to the subcurve $\pr_2: C_0^-\times \overline{ H}{}^\times_{g^+,n,m}\rightarrow 
\overline{H}{}^\times_{g^+,n,m}$ is independent of the fibre. We have then furthermore a natural inclusion of groups
\[ \PGL(N^+ +1)\times \Aut(C_0^-) \:\: \subset  \:\: \PGL(N+1) ,\]
with $N^+ := (2g^+ -2+m)n-g^+ +1$, compare remark \ref{21}. 
By construction, the morphism $\Theta^\times$ is $\PGL(N^+ +1)$-equivariant.
\end{abschnitt}
 
\begin{abschnitt}\em
The permutation group $\Gamma(C_0^-;\calq)$ acts freely on $\overline{H}^{ps}_{g^+,n,m}$, and the the morphism $\Theta^\times$ factors through a morphism
\[ \Theta : \:\: \overline{H}{}^\times_{g^+,n,m} /\Gamma(C_0^-;\calq) \rightarrow \overline{H}_{g,n,0} .\]
\end{abschnitt}

The following fact is a straightforward generalization of corollary \ref{K324}. 

\begin{lemma}\label{L801}
Let $f: C\rightarrow S$ be a stable curve of genus $g$, with reduced base $S$, such that for all $s\in S$ the fibre $C_s$ has exactly $3g-3-d$ nodes. Suppose that the induced morphism $\theta_f : S\rightarrow \msbar_g$ factors through ${D}(C_0;P_1,\ldots,P_d)$. Then there exists a unique subcurve $f^-: C^-\rightarrow S$ of $f:C\rightarrow S$, which is locally isomorphic to the trivial curve $\pr_2: C_0^-\times S \rightarrow S$ as an $m/\Gamma(C_0^-;\calq)$-pointed prestable curve. 
\end{lemma}

\proof 
To prove this lemma, pick for each fibre $C_s$ of $f$ all of its irreducible components, which have the maximal number of nodes when considered as pointed stable curves themselves. We obtain a subcurve $C_s^-$, which is a prestable curve of some genus $g_1$, with $m_1$ unordered distinguished points, which are the points of intersection of $C_s^-$ with the closure of its complement in $C_s$.
 
Since $f: C\rightarrow S$ is flat, the subcurves of the individual fibres  glue together to a global subcurve $f^-: C^-\rightarrow S$ over $S$. Locally, there exist sections of nodes of $f:C \rightarrow S$, inducing locally an ordering of the marked points in the fibres, and thus locally $f^-: C^-\rightarrow S$ is an $m_1$-pointed prestable curve. 

Subcurves of neighbouring fibres are deformations of each other, preserving the number of their nodes. Since the number of nodes in each fibre of $f^-: C^-\rightarrow S$ is maximal, neighbouring fibres are in fact isomorphic as $m_1$-pointed prestable curves. Thus $f^-: C^-\rightarrow S$ is locally isomorphic to the trivial curve $\pr_2: C_1^-\times S \rightarrow $, where $C_1^-$ is some $m_1$-pointed prestable curve of genus $g_1$.  By the definition of $\overline{D}(C_0;P_1,\ldots,P_d)$, we must have $m_1=m$, and the curve $C_1^-$  must be isomorphic to $C_0^-$ as an $m/\Gamma(C_0^-;\calq)$-pointed curve.
\ebew

\begin{abschnitt}\em
Let ${D}^\times(C_0;P_1,\ldots,P_d)$ denote the reduced subscheme of the scheme $\overline{D}(C_0;P_1,\ldots,P_d)$, which is the locus of curves for which the groups of automorphisms split appropriately. The corresponding reduced substack of $\overline{\cal D}(C_0;P_1,\ldots,P_d)$ shall be denoted by ${\cal D}^\times(C_0;P_1,\ldots,P_d)$

Define the scheme $\overline{K}(C_0;P_1,\ldots,P_d)$ as the reduction of the preimage of $\overline{D}(C_0;P_1,\ldots,P_d)$ under the canonical morphism $\overline{H}_{g,n,0}\rightarrow \msbar_g$, and let ${K}(C_0;P_1,\ldots,P_d)$ and $K^\times (C_0;P_1,\ldots,P_d)$ denote the subschemes defined by ${D}(C_0;P_1,\ldots,P_d)$ and ${D}^\times(C_0;P_1,\ldots,P_d)$, respectively. 

We define $ {H}^\times_{g^+,n,m}$ as the reduction of the preimage of $\overline{K}(C_0;P_1,\ldots,P_d)$ under the morphism $\Theta^\times$.  Its image under the canonical morphism $\pi: \overline{H}_{g^+,n,m} \rightarrow \msbar_{g^+,m}$ is denoted by $\ms^\times_{g^+,m}$. 

The corresponding reduced  substack of $\cmgplus$ shall be denoted by $\cmx$.
\end{abschnitt}

Using the above lemma it is easy to show the following proposition.

\begin{prop}\label{P802}
The morphism $\Theta^\times :\overline{H}{}^\times_{g^+,n,m}\rightarrow \overline{H}_{g,n,0} $ induces a morphism 
\[ \Theta : \:\:  {H}^\times_{g^+,n,m}/\Gamma(C_0^-;\calq) \rightarrow 
{K}^\times(C_0;P_1,\ldots,P_d) \]  
which is injective, and an isomorphism
\[ {M}^\times_{g^+,m} / \Gamma(C_0^-;\calq) \: \cong \: {D}^\times(C_0;P_1,\ldots,P_d).\]
\end{prop}

For the proof compare also the proof of proposition \ref{312}. Note that the hard bit in that proof was to show injectivity on the boundary points. Since we restrict ourselves here to the open subscheme we can avoid this difficulty. 

\begin{abschnitt}\em
From proposition \ref{28} it follows that there is a representation of the stack ${\overline{\cal D}(C_0;P_1,\ldots,P_{d})}$ as a quotient stack
\[ {\overline{\cal D}(C_0;P_1,\ldots,P_{d})} \cong 
\left[ \, \overline{K}(C_0;P_1,\ldots,P_{d}) / \PGL(N^+ +1) \right] ,\]
compare also lemma \ref{140}. There is 
 an analogous representation of the open substacks   ${{\cal D}(C_0;P_1,\ldots,P_{d})}$ and ${{\cal D}^\times(C_0;P_1,\ldots,P_{d})}$. Note that the substack ${{\cal D}(C_0;P_1,\ldots,P_{d})}$ is smooth, since the scheme $ {K}(C_0;P_1,\ldots,P_{d})$ is smooth, and  ${D}(C_0;P_1,\ldots,P_{d})$ is smooth as well. 

However, in general neither the scheme $ \overline{K}(C_0;P_1,\ldots,P_{d})$ nor  its quotient $ \overline{D}(C_0;P_1,\ldots,P_{d})$ are  smooth, and ${\overline{\cal D}(C_0;P_1,\ldots,P_{d})}$ is not a smooth stack. 
\end{abschnitt}

Using the morphism $\Theta : { H}^\times_{g^+,n,m}/\Gamma(C_0^-;\calq) \rightarrow \overline{H}_{g,n,0}$ and the inclusion $\PGL(N^+ +1)\subset \PGL(N+1)$ it is straightforward to construct a morphism of stacks from $[\, { H}^\times_{g^+,n,m/\Gamma(C_0^-;\calq)} /\PGL(N+1)]$ to ${{\cal D}^\times(C_0;P_1,\ldots,P_{d})}$. From   this morphism one can derive the following generalization of our main theorem \ref{MainThm}.

\begin{thm}\label{T809}
There is an isomorphism of Deligne-Mumford stacks
\[ \left[{H}^\times_{g^+,n,m/\Gamma(C_0^-;\calq)} /\PGL(N+1)\times\Aut(C_0^-)\right] \:\:\: \cong \:\:\: {{\cal D}^\times(C_0;P_1,\ldots,P_{d})},\]
where $\Aut(C_0^-)$ acts trivially on ${H}^\times_{g^+,n,m/\Gamma(C_0^-;\calq)}$.
\end{thm}

The proof is analogous to that of theorem \ref{MainThm}. As before, there are two  main ingredients to the construction of the inverse  functor from ${{\cal D}^\times(C_0;P_1,\ldots,P_{d})}$ to $[{H}^\times_{g^+,n,m/\Gamma(C_0^-;\calq)} /\PGL(N+1)\times\Aut(C_0^-)]$.  The first is the existence of locally trivial subschemes as provided by lemma \ref{L801}, which allows to define the  functor on objects. The second is the injectivity of the morphism $\Theta:   {H}^\times_{g^+,n,m}/\Gamma(C_0^-;\calq) \rightarrow 
{K}^\times(C_0;P_1,\ldots,P_d)$, which follows from proposition \ref{P802}, and which is necessary to define the functor on morphisms. 

\begin{rk}\em
The theorem shows in particular that there is a surjective, finite and unramified morphism 
\[ \cmgammadetailx  \rightarrow {{\cal D}^\times(C_0;P_1,\ldots,P_{d})}, \] 
which is of degree equal to the order of $\Aut(C_0^-)$.
\end{rk} 

\begin{abschnitt}\em
For the sake of completeness, let us take a closer look at the  boundary stratum $\msbar_g^{(3g-3)}$, and the corresponding zero-dimensional substacks of $\cmgbar$. Irreducible components $\overline{D}$ of $\msbar_g^{(3g-3)}$ are isolated points $[C_0]$, where $C_0\rightarrow \Spec(k)$ is a stable curve with the maximal number of nodes. By our earlier convention we have thus
\[ \overline{D} = \overline{D}(C_0;\, \emptyset) .\]
An object of ${\overline{\cal D}(C_0;\, \emptyset)}(S)$, for some scheme $S$, is a stable curve $f: C\rightarrow S$ of genus $g$, where each fibre is isomorphic to $C_0$. As in lemma \ref{140}, there is an isomorphism
\[ {\overline{\cal D}(C_0;\, \emptyset)} \: \cong \: [\, \overline{K}(C_0;\, \emptyset) / \PGL(N+1) ] ,\]
where $\overline{K}(C_0;\, \emptyset)$ denotes the fibre of the canonical morphism $\overline{H}_{g,n,0}\rightarrow \msbar_g$ over the point $[C_0]$. Recall that $\overline{K}(C_0;\, \emptyset) \cong \PGL(N+1) / \Aut(C_0)$. 
\end{abschnitt}

\begin{lemma}\label{LDNULL}
There is an isomorphism of stacks
\[ {\overline{\cal D}(C_0;\, \emptyset)} \:\: \cong  \:\: \left[ \, \Spec(k) / \Aut(C_0) \, \right] ,\]
where $\Aut(C_0)$ acts trivially on $\Spec(k)$. In other words, $ {\overline{\cal D}(C_0;\, \emptyset)}$ is isomorphic to the classifying stack of the group $\Aut(C_0)$.
\end{lemma}

\proof
Fix a geometric point $x: \Spec(k)\rightarrow \overline{K}(C_0;\, \emptyset)$. This distinguishes an embedding of $C_0$ into $\PP^N$. At the same time, this defines an embedding of the group $\Aut(C_0)$ into $\PGL(N+1)$, such that 
\[ \Aut(C_0) = \Stab_{\PGL(N+1)}(x) .\]
Let $S$ be a scheme. Consider an object $(E',p',\theta')\in [ \, \Spec(k) / \Aut(C_0)](S)$. We define a principal $\PGL(N+1)$-bundle $p: \: E:=E'(\PGL(N+1)) \rightarrow S$ by extension. There is a well-defined $\PGL(N+1)$-equivariant morphism $\theta: E\rightarrow \overline{K}(C_0;\, \emptyset)$, which is given by
\[ \theta([e',\gamma]) := x\cdot \gamma \]
for a point $[e',\gamma] \in E= (E'\times \PGL(N+1))/\Aut(C_0)$, with $e'\in E'$ and $\gamma\in \Aut(C_0)$. Together with the obvious construction on morphisms this induces a functor from $[ \, \Spec(k) / \Aut(C_0)]$ to ${\overline{\cal D}(C_0;\, \emptyset)}$. 

Conversely, take an object $(E,p,\theta)\in {\overline{\cal D}(C_0;\, \emptyset)}(S)$, for some scheme $S$. Let $E_s$ denote a fibre of $p:E\rightarrow S$ over some point $s\in S$. Because of its $\PGL(N+1)$-equivariance, the restricted morphism $\theta_s : E_s\rightarrow \overline{K}(C_0;\, \emptyset)$ is surjective. By the choice of the embedding of $\Aut(C_0)$ into $\PGL(N+1)$, the preimage $\theta_s^{-1}(x)$ defines a point $x_s$ in $E_s/\Aut(C_0)$. Therefore we can construct a section 
\[ \sigma: \:\:  S \rightarrow E/\Aut(C_0)\]
by $\sigma(s) := x_s$. This implies that there exists a principal $\Aut(C_0)$-subbundle $p': E'\rightarrow S$ of $E$, which extends to $p: E\rightarrow S$. By construction, the restriction $\theta':E'\rightarrow \overline{K}(C_0;\, \emptyset)$ of $\theta$ maps constantly to the point $x$. Thus $(E',p',\theta')$ is an object of $[\Spec(k) / \Aut(C_0)](S)$. 

It is easy to see that under this construction morphisms between triples in ${\overline{\cal D}(C_0;\, \emptyset)}$ restrict to morphisms of triples in $[\Spec(k) / \Aut(C_0)]$. Thus we obtain  a functor from  ${\overline{\cal D}(C_0;\, \emptyset)}$ to $[\Spec(k) / \Aut(C_0)]$, which is inverse to the functor constructed above, up to isomorphism, and hence an isomorphism of stacks.
\ebew

\begin{ex}\em
We want to conclude this section with an example which takes this generalization even a bit further, to the case of moduli stacks of pointed prestable curves.

Consider the curve $C_0$ of genus $g=3$, with nodes $p_i$ for $i=1,\ldots,3g-3=6$ as pictured below. 

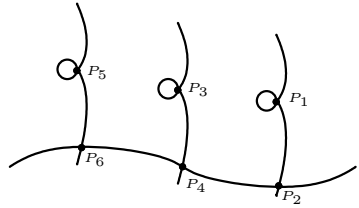
\begin{center}
\setlength{\unitlength}{0.00016667in}
\begingroup\makeatletter\ifx\SetFigFont\undefined%
\gdef\SetFigFont#1#2#3#4#5{%
  \reset@font\fontsize{#1}{#2pt}%
  \fontfamily{#3}\fontseries{#4}\fontshape{#5}%
  \selectfont}%
\fi\endgroup%
{\renewcommand{\dashlinestretch}{30}

}
\end{center}

On the level of moduli spaces we obtain isomorphisms
\[ \begin{array}{lcl}
\overline{D}(C_0;P_1) &\cong & \msbar_{1,1} ,\\
\overline{D}(C_0;P_2) &\cong & \msbar_{1,2/\Sigma_2}' ,\\
\overline{D}(C_0;P_1,P_2) &\cong & \msbar_{1,2/\Sigma_2} .
\end{array}\] 
For the open substacks we obtain as in section \ref{Sect81} isomorphisms
\[ \begin{array}{lcl}
{{\cal D}(C_0;P_1)} &\cong & {{\cal M}_{1,1}} /{{\mathbb Z}_2}^3 ,\\[2mm]
{{\cal D}^\times(C_0;P_2)} &\cong & {{\cal M}^\times_{1,2/\Sigma_2}}  /{{\mathbb Z}_2}^2,\\[2mm]
{{\cal D}^\times(C_0;P_1,P_2)} &\cong & {{\cal M}^\times_{1,2/\Sigma_2}}  /{{\mathbb Z}_2}^2.
\end{array}\] 
Note that the representations by quotients do not extend to the closed stacks. 

We will try to make this description more geometric.  
Let $f: \msbar_{1,2} \rightarrow \msbar_{1,1}$ denote the canonical morphism, which is more or less defined by forgetting the second marked point on each curve. Recall that components of curves, which become unstable by the omission of one distinguished  point, are contracted to a (marked) point. 
 The fibre over the unique boundary point of $\msbar_{1,1}$ is isomorphic to $\msbar_{1,2}'$. Recall that by definition $\msbar_{1,2}' =\Delta_0$ is the closure of the locus of irreducible singular  stable curves in $\msbar_{1,2}$.  

The morphism $f$ admits a section $j: \msbar_{1,1} \rightarrow \msbar_{1,2}$. The embedding is constructed by glueing to a given $1$-pointed stable curve of genus $1$ a $3$-pointed line in one of its marked points. The image of $j$ is a divisor $\Delta_1$ in $ \msbar_{1,2}$, and for the boundary of $ \msbar_{1,2}$ holds
\[  \msbar_{1,2}^{(1)} = \Delta_0 \cup \Delta_1 .\] 
In fact, $\Delta_1$ is the closed locus of reducible singular $2$-pointed stable curves of genus $1$. 

Both divisors $\Delta_0$ and $\Delta_1$ are invariant with respect to the action of $\Sigma_2$ on $\msbar_{1,2}$. Let $\Delta_0^\diamond$ and $\Delta_1^\diamond$  denote their images in $\msbar_{1,2/\Sigma_2}$, and note that they meet transversally in exactly two points. 

Consider the curve $C_1 := C_0^+(P_2)$, considered as a $2/\Sigma_2$-pointed stable curve of genus $1$, and note that it has the maximal number of nodes possible for such a curve. Choose an enumeration of its nodes as  shown in  the picture below. 

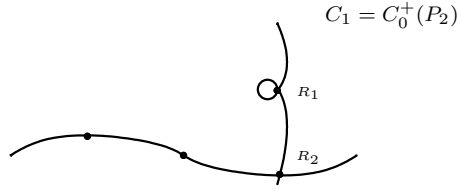
\begin{center}
\setlength{\unitlength}{0.00016667in}
\begingroup\makeatletter\ifx\SetFigFont\undefined%
\gdef\SetFigFont#1#2#3#4#5{%
  \reset@font\fontsize{#1}{#2pt}%
  \fontfamily{#3}\fontseries{#4}\fontshape{#5}%
  \selectfont}%
\fi\endgroup%
{\renewcommand{\dashlinestretch}{30}

}
\end{center}

Note also that the two distinguished points are not labelled. 

Extending our previous notation in an obvious way, we obtain identities
\[ \begin{array}{lcl}
\Delta_0 &=& \overline{D}(C_1;R_2), \:\:\: \text{ and} \\
\Delta_1 &=& \overline{D}(C_1;R_1).
\end{array}\]
Consider the decomposition of $C_1$ defined by the node $R_1$. We obtain for $C_1^+$ a $1$-pointed nodal stable curve of genus $0$, and for $C_1^-$ a smooth stable curve of genus $0$ with the one marked point $R_2$ and two more unlabelled distinguished points. Formally, $C_1^-$ can be described as a $3/\Sigma_2$-pointed stable curve. Let ${\cal X}$   denote the set of the two unmarked distinguished points. Then we have
\[ \begin{array}{lcl}
\Aut(C_2^-;{\cal X}) &=& {\mathbb Z}_2, \:\:\: \text{ and} \\[2mm]
\Gamma(C_0^-; \{R_2\}) &=& \{\id\} .
\end{array}\]
By a generalization of theorem  \ref{MainThm} one obtains isomorphisms of stacks
\[ \begin{array}{lcl}
{{\cal D}(C_1;R_1)} &\cong & {{\cal M}_{1,1}} /{{\mathbb Z}_2} , \:\: \text{ and}\\[2mm]
{{\cal D}(C_1;R_2)} &\cong & {{\cal M}_{1,2/\Sigma_2}'} .
\end{array}\] 
Note that in the second case the decomposition of $C_1$ defined by $R_2$ is the trivial one with $C_1^+=C_1$. 

We can now replace our initial diagram
\[ \diagram
{{\cal M}_{1,1}} \dto\rrto^{finite} && {{\cal M}_{1,1}/{{\mathbb Z}_2}^3} \dto \rto|<\hole|<<\ahook & {{\cal M}_3}\dto\\
{ M}_{1,1} \rrto^\cong  &&  D(C_0;P_1) \rto|<\hole|<<\ahook & \msbar_3
\enddiagram \]
by the two  much more refined diagrams below. The first diagram shows  explicit descriptions of  the open substacks. 
\[\begin{array}{cccccccccccc}
{{\cal M}_{1,1}} &\rightarrow&
{{\cal M}_{1,1}}/{{\mathbb Z}_2}\subset
{{\cal M}_{1,2/\Sigma_2}} \: \supset \: {{\cal M}_{1,2/\Sigma_2}}\rightarrow 
{{\cal M}_{1,2/\Sigma_2}}/{{\mathbb Z}_2}^2& \subset&
{{\cal M}_3}\\[3mm]
\downarrow&&\downarrow\hspace{2cm}\downarrow\hspace{2cm}\downarrow\hspace{2cm}\downarrow&&\downarrow\\[3mm]
{M}_{1,1}& \cong&
D(C_1;R_1) \subset
\msbar_{1,2/\Sigma_2} \: \supset \: {{M}_{1,2/\Sigma_2}}\cong
D(C_0;P_1,P_2)& \subset&
\msbar_3.
\end{array}\]
The second diagram depicts finite morphisms between closed stacks.
\[\diagram
{{\cal M}_{1,1}} \dto\rto^{fin.}&
{{\cal D}(C_1;R_1)}\dto\rto|<\hole|<<\ahook&
{{\cal M}^\times_{1,2/\Sigma_2}} \dto\rto^{fin.}&
{{\cal D}^\times(C_0;P_1,P_2)}\dto\rto|<\hole|<<\ahook&
{{\cal M}_3}\dto\\
{\ms}_{1,1} \rto^\cong &
\overline{D}(C_1;R_1) \rto|<\hole|<<\ahook&
\ms_{1,2/\Sigma_2}^\times\rto^\cong&
{D}^\times(C_0;P_1,P_2)\rto|<\hole|<<\ahook&
\ms_3.
\enddiagram\]
Of course, there are analogous diagrams for ${D}(C_0;P_2)$. 
In particular, one has isomorphisms
\[ \begin{array}{lcr}
{{\cal D}(C_0;P_1)} &\cong & {{\cal D}(C_1;R_1)} \: / \: {{\mathbb Z}_2}^2 ,\\[2mm]
{{\cal D}^\times(C_0;P_2)} &\cong & {{\cal D}^\times(C_1;R_2)}  \: / \: {{\mathbb Z}_2}^2,\\[2mm]
{{\cal D}^\times(C_0;P_1,P_2)} &\cong & {{\cal M}^\times_{1,2/\Sigma_2}}  \: / \: {{\mathbb Z}_2}^2.
\end{array}\] 
In other words, the description of the boundary component ${\overline{\cal D}(C_0;P_1,P_2)}$ in $\overline{\cal M}_3$ can be reduced to the description of $\overline{\cal M}_{1,2}$ using a 4-to-1 covering. 
\end{ex}

\begin{abschnitt}\em
This example may serve as an illustration of the way in which  components of the boundary strata of the moduli stack $\stackf_{\msbar_g}$ might  be analyzed by an inductive reduction to moduli stacks of pointed curves of some smaller  genus $g' \le g$.
\end{abschnitt}

\vfill\pagebreak

\section{A remark on intersection products}

The underlying theme of this paper was to describe in what a non-intuitive way a substack ${\overline{\cal D}}$ of the moduli stack $\cmgbar$ behaves, when $\overline{D}$ is an irreducible component of some boundary stratum $\msbar_g^{(3g-3-d)}$. In fact, if $\overline{M}{}^+$ denotes the corresponding subscheme of the moduli space of $m/\Gamma(C_0^-;\CQ)$-pointed prestable curves of genus $g^+$, isomorphic to $\overline{D}$, then generically the corresponding moduli stack $\overline{\cal M}{}^+$ differs from the stack ${\overline{\cal D}}$ by a finite covering of degree equal to the order of the group $\Aut(C_0^-)$. Recall that the induced action of $\Aut(C_0^-)$ on the moduli scheme $\overline{ M}{}^+$ is trivial. 

We want to give a first indication, why this unexpected ``deviation'' from the geometric situation is indeed not a weakness, but a strenght of the stack description. In fact, automorphisms of curves represented by points of $\msbar_g$ profoundly influence the intersection theory of the moduli space. In the scheme description the information on the automorphism is lost, and can only be recovered with some effort. This can be done for example by applying the theory of $\mathbb Q$-varieties, as in Faber's paper \cite{Fa}. Dealing with stacks however, this information is readily availiable. We refer to the papers  \cite{Vi},\cite{EG} and \cite{Kr} for three excellent accounts on this subject. 

Our description of substacks contained in the boundary of $\cmgbar$ is not designed to be applied in intersection theory. In fact, it is a very rigid construction, as it depends heavily on features of stable curves, which are only availiable in the individual components of the boundary stratification of $\msbar_g$. Still, in cases where the boundary components behave in a nice way, something can be said.  
Let us illustrate this in one small example. 

\begin{ex}\em
Consider the moduli space $\msbar_3$ of stable curves of genus $g=3$, which is of dimension $6$. In this case we can rely on  a very good description of the bondary stratification, which is provided by  Faber's work \cite{Fa}. 
Let $f: C_3 \rightarrow \Spec(k)$ be the stable curve of genus $g=3$, with $6$ nodes as in the list of  remark \ref{R8301}, with an enumeration $P_1,\ldots,P_6$ of its nodes as indicated there.  Consider the $3$-dimensional irreducible components $\overline{D}(C_3;P_2,P_4,P_6)$ and  $\overline{D}(C_3;P_1,P_3,P_5)$ contained in the boundary stratum $\msbar_3^{(3)}$. A general fibre of $\overline{D}(C_3;P_2,P_4,P_6)$ is given schematically by the picture below.

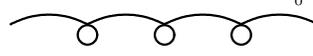
\begin{center}
\setlength{\unitlength}{0.00016667in}
\begingroup\makeatletter\ifx\SetFigFont\undefined%
\gdef\SetFigFont#1#2#3#4#5{%
  \reset@font\fontsize{#1}{#2pt}%
  \fontfamily{#3}\fontseries{#4}\fontshape{#5}%
  \selectfont}%
\fi\endgroup%
{\renewcommand{\dashlinestretch}{30}

}
\end{center}

Recall that $P_2,P_4$ and $P_6$ denote the nodes, which are not preserved. The   number in the picture denotes the geometric genus of the irreducible component next to it. 
 Similarly, $\overline{D}(C_3;P_1,P_3,P_5)$ can be drawn as follows. 

\begin{center}
\setlength{\unitlength}{0.00016667in}
\begingroup\makeatletter\ifx\SetFigFont\undefined%
\gdef\SetFigFont#1#2#3#4#5{%
  \reset@font\fontsize{#1}{#2pt}%
  \fontfamily{#3}\fontseries{#4}\fontshape{#5}%
  \selectfont}%
\fi\endgroup%
{\renewcommand{\dashlinestretch}{30}
%
}
\end{center}

In Faber's list \cite[p.343, table 6]{Fa}, these components  are denoted by $(a)$ and $(i)$, respectively. There are obvious inclusions
\[ \overline{D}(C_3;P_1,P_3,P_5) \supset \overline{D}(C_3;P_1,P_3)\supset \overline{D}(C_3;P_1)\supset \overline{D}(C_3;\emptyset) .\]
One sees easily that there are no other boundary components contained in $\overline{D}(C_3;P_1,P_3,P_5)$. Thus, set theoretically we have for the intersection
\[ \overline{D}(C_3;P_1,P_3,P_5) \cap \overline{D}(C_3;P_2,P_4,P_6) = \overline{D}(C_3;\emptyset) = \{[C_3]\} .\]
The two $3$-dimensional subschemes meet in fact transversally in the point $[C_3]$. 
From \cite[p. 412, table 4]{Fa} one can read off a presentation of the $\QQ$-classes $(a)_\QQ$ and $(i)_\QQ$ associated to $\overline{D}(C_3;P_2,P_4,P_6)$ and  $\overline{D}(C_3;P_1,P_3,P_5)$ in terms of generating elements of the Chow ring of $\msbar_g$. From the multiplication table \cite[p. 418, table 10]{Fa} one computes for the intersection product 
\[  (a)_\QQ \cdot(i)_\QQ = \frac{1}{48}. \]
There are isomorphisms between the Chow groups of the stack $\cmgbar$ and the Chow groups of its moduli space $\msbar_g$, see \cite[prop. 6.1]{Vi}. The $\QQ$-classes of subschemes of $\msbar_g$ correspond to the classes of the respective substacks of $\cmgbar$. In our case, the substacks ${\overline{\cal D}(C_3;P_1,P_3,P_5)}$ and ${\overline{\cal D}(C_3;P_2,P_4,P_6)}$ are both smooth and meet transversally in the substack ${\overline{\cal D}(C_3;\emptyset)}$. Therefore, for the classes in the respective Chow groups holds
\[ [{\overline{\cal D}(C_3;P_1,P_3,P_5)}].[{\overline{\cal D}(C_3;P_2,P_4,P_6)}] = [{\overline{\cal D}(C_3;\emptyset)} ] .\]
By lemma \ref{LDNULL} we have an isomorphism
\[ {\overline{\cal D}(C_3;\emptyset)} \cong [ \Spec(k) / \Aut(C_3)], \]
and the canonical morphism ${\overline{\cal D}(C_3;\emptyset)} \rightarrow \{[C_3]\} $ is of degree $\frac{1}{\#\Aut(C_3)}$. Note that $\Aut(C_3) \cong \Sigma_3\times {\ZZ_2}^3$. Thus, the $\QQ$-class $[C_3]_\QQ$ in $A^6(\msbar_3) \cong \QQ$ corresponding to the class of ${\overline{\cal D}(C_3;\emptyset)}$ is just $\frac{1}{48}$, which is in accordance with Faber's result. 
\end{ex}

\begin{abschnitt}\em
It should not be concealed that this example has been chosen with care, and that it is not representative for the behaviour of bondary strata. In general it will not be the case  that individual components meet transversally in the expected dimension, or that they are smooth in a neighbourhood of their intersection locus. Our construction may be helpful in certain cases. For a complete analysis of intersection products on $\cmgbar$ it seems not suitable, but, as we said before, it is not designed to be. 
\end{abschnitt}

\appendix
\setcounter{chapter}{0}

%
%

\chapter{Proof of lemma \ref{L326A}}\label{AppA}

We will prove the lemma by induction on the number $\nu = \nu(C)$ of such nodes of $C$ which are intersections of two smooth rational irreducible components. Incidentially, this is the number of subcurves of $C$, which are represented by points of $\msbar_{0,4}$. 

Throughout, let ``$+$'' denote glueing of two $k$-pointed curves along their marked points according to the labels of the marked points. To simplify the notation, we will denote the group $\Gamma(C_i^-;P_1^{(i)},\ldots,P_{\mu+k-4}^{(i)},\{Q_1^{(i)},\ldots,Q_k^{(i)}\}) $ by $\Gamma(C_i^-)$, for $i=1,2$.

Let $\nu =1$. Then necessarily $C_1^+ = C_2^+$ as sets. Since $C_1^+ + C_1^- = C = C_2^+ + C_2^-$, the the order of the marked points can differ at most by a permutation induced by an automorphism of the closure of the complement, which is $C_1^-= C_2^-$, considered as a set. So in this case the claim is clear. 

Assume now $\nu > 1$, and assume that the claim has been  proved for all numbers less than $\nu$. Note that we may also assume that $C_1^+$ and $C_2^+$ are different as sets, since otherwise we could argue as we did above for the case $\nu =1$  to show the claim.  

$(i)$ We consider first the case where  $C_1^+$ and $C_2^+$ have no common irreducible component. Since the intersection $C_1^+ \cap C_2^+$ is at most finite,  we must have
\[ C_1^+ \subset C_2^- .\]
Define
\[ \tilde{C}_2^- := C_2^- - C_1^+ .\]
Here ``$-$'' denotes taking the closure if the complement, considered as an appropriately pointed curve. Thus
\[ C = C_1^- + C_1^+ = \tilde{C}_2^- + C_1^+ + C_2^+ .\]
From this we obtain
\[ C_1^- = \tilde{C}_2^- + C_2^+ .\]
Using the isomorphism $\phi: C_2^- \rightarrow C_1^-$ we get also 
$C_1^- = \phi(\tilde{C}_2^- ) + \phi(C_1^+)$, 
and therefore the equality
\[ \phi(\tilde{C}_2^- ) + \phi(C_1^+) =  \tilde{C}_2^- + C_2^+ .\]
We now put
\[ \tilde{C} := C_1^- . \]
Note that $\nu(\tilde{C}) < \nu(C)$. Clearly,  there is an isomorphism $\tilde{\phi}$ from  $\phi(\tilde{C}_2^-)$ to $\tilde{C}_2^-$, considered as pointed curves, which is just given by the inverse of the restriction of $\phi$. Both $\phi(C_1^+)$ and $C_2^+$ are represented by points in $\msbar_{0,4}$. From the induction hypothesis it now follows that $\phi(C_1^+)$ is isomorphic to $C_2^+$, considered as $4$-pointed curve, up to a permutation of the labels of the marked points given by an element $\tilde{\sigma} \in \Gamma(\tilde{C}_2^-)$. 

Recall that by definition $\tilde{\sigma}$ is induced by some automorphism $\rho$ of $\tilde{C}_2^-$, which fixes all the marked points of $\tilde{C}_2^-$, except possibly those which form the intersection $\tilde{C}_2^+ \cap C_2^+$. 

If $C_1^+ \cap C_2^+ = \emptyset$, then in particular all points of the intersection $\tilde{C}_2^+ \cap C_1^+$ are fixed under $\rho$. Indeed, if $P\in \tilde{C}_2^+ \cap C_1^+$, then $P\not\in C_2^+$, and hence $P\not\in \tilde{C}_2^+ \cap C_2^+$. 
Thus $\rho$ can be extended trivially onto $\tilde{C}_2^+ + C_1^+ = C_2^-$. Therefore $\rho$ can be considered an an automorphism of $C_2^-$, and thus $\tilde{\sigma}$ can be considered as an element of $\Gamma(C_2^-)$. Hence $\phi(C_1^+)$ is isomorphic to $C_2^+$, considered as $4$-pointed curve, up to a permutation of the labels of the marked points given by an element $\tilde{\sigma} \in \Gamma({C}_2^-)$. Clearly $\phi(C_1^+)$ and $C_1^+$ are isomorphic as $4$-pointed curves, so finally we find that $C_1^+$ is isomorphic to $C_2^+$, considered as $4$-pointed curve, up to a permutation in $ \Gamma({C}_2^-)$.

If the intersection $C_1^+ \cap C_2^+ $ is not empty, then in order to be able to extend $\rho$ to an automorphism of $C_2^-$, we need to verify that all points $P\in \tilde{C}_2^+ \cap C_1^+ $ are fixed under $\rho$. This follows because $\tilde{C}_2^+ \cap C_1^+ \cap C_2^+$ is empty. Indeed, by assumption such a point $P$ is a marked point of $\tilde{C_2}^- + C_2^+$, and so  in particular not a point of $\tilde{C}_2^- \cap C_2^+$. So as before we can conclude that $C_1^+$ is isomorphic to $C_2^+$, considered as $4$-pointed curve, up to a permutation in $ \Gamma({C}_2^-)$.

$(ii)$ We still have to consider the second case, where $C_1^+$ and $C_2^+$ have exactly one common component. Suppose that
\[ C_1^+ = C_{1,a}^+ + C_{1,b}^+ \quad \text{and} \quad  C_2^+ = C_{2,a}^+ + C_{2,b}^+ ,\]
where $C_{1,a}^+,C_{1,b}^+,C_{2,a}^+,C_{2,b}^+$, are $3$-pointed smooth rational curves, and ``+'' denotes glueing in one marked point. Suppose that
\[ C_{1,a}^+=C_{2,a}^+ \]
as sets. Since $C= C_1^+ + C_1^- = C_2^+ + C_2^-$, and $C_{1,b}^+ \neq C_{2,b}^+$, we must have
\[ C_{1,b}^+ \subset C_2^- \quad \text{and} \quad C_{2,b}^+ \subset C_1^- .\]

Note that in this case for  the number $k$ of  glueing points of $C_1^+$ and $C_1^-$ holds  $1 \le k \le 4$. The case $k=0$ cannot occur: there is always one glueing point of $C_{1,a}^+$, which by assumption is equal to $C_{2,a}^+$, with $C_{2,b}^+$, which is contained in $C_1^-$. This point is thus a point in $C_1^+ \cap 
C_1^-$. 

If $k$ = 1, then any isomorphism $\phi : C_2^- \rightarrow C_1^-$ extends to an isomorphism from $C_2^- + C_2^+ $ to $C_1^- + C_1^+ $. In particular there is an isomorphism between $C_1^+$ and $C_2^+$ as claimed. 

The same holds true if $k=2$. Indeed, there are only two configurations possible, in which the two glueing points can be distributed on $C_1^+ = C_{1,a}^+ + C_{1,b}^+$. Either there is exactly one glueing point on each irreducible component, or both glueing points lie on $C_{1,a}^+$. Recall that there is always at least one glueing point on $C_{1,a}^+$, which is the point given by the intersection of $C_{1,a}^+=C_{2,a}^+$ with $C_{2,b}^+$. Since by assuption $C_{1,a}^+=C_{2,a}^+$ as subsets of $C$, the configuration on $C_2^+$ must be the same as on $C_1^+$. Thus $C_1^+$ and $C_2^+$ are isomorphic as $2$-pointed curves. Note that in this case, too, the isomorphism $\phi : C_2^- \rightarrow C_1^-$ extends to an automorphism of $C$. 

We are finally left with the cases $k = 3,4$. 

We may assume that $\tilde{C}_2^- \neq \emptyset$. Indeed, if $k=4$, then there is always the one  point of $C_{1,a}^+$ lying on $\tilde{C}_2^-$, which is the unique node of $C_{1,a}^+$ which is neither contained in $C_{1,b}^+$ nor in $C_{2,b}^+$. If $k=3$, and $\tilde{C}_2^- = \emptyset$, then $C_2^- = C_{1,b}^+$. Hence $C_2^-$ is a line with three marked points $A_1,A_2,A_3$ on it. Any permutation of these three points is induced by an automorphism of the line. In other words, we have $\Gamma(C_2^-) = \Sigma_3$. So in particular $C_1^+$ is isomorphic to $C_2^+$ as a $4$-pointed curve up to a permutation in $\Gamma(C_2^-)$. 

Note that the same arguments work if we assume that $\tilde{C}_2^-$ is nonempty, but has empty intersection with both $C_1^+$ and $C_2^+$.
 
To deal with the remaining cases  we need the following claim, which will be proven  by a lenghty case-by-case analysis in \ref{PCLAIM} below.

{\bf Claim.}\label{PG213} For the isomorphism   $\phi : C_2^- \rightarrow C_1^-$ holds at least one of the following statements.\\
$(i)$ The isomorphism   $\phi : C_2^- \rightarrow C_1^-$ extends to an automorphism of $C$.\\
$(ii)$ The intersection $\phi(C_{1,b}^+) \cap C_{1,a}^+$  consists of exactly one point.\\
$(iii)$ The intersection $\phi^{-1}(C_{2,b}^+) \cap C_{1,a}^+$  consists of exactly one point.\\
$(iv)$ There exists a connected component $C_A^-$ of $C_2^-$, which is glued to $C_2^+$ in exactly one point $P$, such that $P$ is a fixed point of $\phi$.\\
$(v)$  The group $\Gamma(C_2^-)\subset \Sigma_k$  contains a subgroup isomorphic to the  permutation group $\Sigma_3$.

Let us assume that we know the statement of the claim to be true. Clearly, if the isomorphism $\phi : C_2^- \rightarrow C_1^-$ extends to an automorphism of $C$ then in particular $C_1^+$ and $C_2^+$ are isomorphic, and the lemma holds in this case. We don't even have to resort to induction. Also, if $\Sigma_3 \subset \Gamma(C_2^-)$, then clearly $C_1^+$ is isomorphic to $C_2^+$ as a $4$-pointed curve, up to a permutation in $\Gamma(C_2^-)$.

We will finish the proof of the remaining cases by induction on the number of connected components of $C_2^-$. Note that for our purposes connected components of $C$, which intersect neither $C_1^+$ nor $C_2^+$, can safely be ignored.

$(\alpha)$ If there is only one irreducible component of $C_2^-$, then statement $(iv)$ of the claim cannot apply. Without loss of generality we may assume that  $\phi(C_{1,b}^+) \cap C_{1,a}^+ \neq \emptyset$.  Otherwise replace $\phi$ by $\phi^{-1}$, and interchange all indices $1$ and $2$. Remember that $C_{1,a}^+ = C_{2,a}^+$ as sets.

Recall that $C_{1,b}^+ \subset C_2^-$. We define
\[ \tilde{C_2}^- := C_2^- - C_{1,b}^+ .\]
From this we obtain for the curve $C-C_{1,b}^+$ two decompositions
\[ C_1^- + C_{1,a}^+ = \tilde{C_2}^- + C_{2,a}^+ + C_{2,b}^+ .\]
We denote this curve by $\tilde{C} :=  C_1^- + C_{1,a}^+$. Note that $\nu(\tilde{C}) < \nu(C)$, since the node, where $C_{1,a}^+$ and $C_{1,b}^+$ intersect, is removed. It also holds $\nu(\tilde{C}) \ge 1$, since $\tilde{C}$ contains $C_2^+$. Using the isomorphism $\phi : C_2^- \rightarrow C_1^-$ we obtain
\[ C_1^ - = \phi(\tilde{C}_2^-) + \phi(C_{1,b}^+) ,\]
and hence
\[ \phi(\tilde{C}_2^-) + \phi(C_{1,b}^+) + C_{1,a}^+ = \tilde{C}_2^- + C_{2,a}^+ + C_{2,b}^+ .\]
By definition $ C_{2,a}^+ + C_{2,b}^+ = C_2^+$. We  also assumed that $\phi(C_{1,b}^+)$ meets $C_{1,a}^+$ in exacly one point, so $\tilde{C}_1^+ := \phi(C_{1,b}^+) + C_{1,a}^+$ is represented by a point in $\msbar_{0,4}$. Note that the curve $\tilde{C}_1^+$ is isomorphic to $C_1^+ = C_{1,b}^+ + C_{1,a}^+$ as a $4$-pointed curve. Indeed, to know the isomorphism class of a stable singular rational $4$-pointed curve it is enough to know the labels of the marked points on one of the two irreducible comonents. Thus the last equation above becomes
\[ \phi(\tilde{C}_2^-) +  \tilde{C}_1^+ = \tilde{C}_2^- + C_2^+ .\]
By induction we may conclude that $\tilde{C}_1^+$ is isomorphic to $C_2^+$, up to a permutation induced by an automorphism $\tilde{\rho}$ on $\tilde{C}_2^-$, which fixes the marked points of  $\tilde{C}_2^-$, except possibly those in the intersection $\tilde{C}_2^- \cap C_2^+$. 

Note that $\tilde{C}_2^-$ differs from $C_2^-$ by the  line $C_{1,b}^+$ glued to $\tilde{C}_2^-$ in at most two marked points. These points are in particular no points of the intersection $\tilde{C}_2^- \cap C_2^+$, and are thus fixed under $\tilde{\rho}$. Therefore the morphism $\tilde{\rho}$ extends to an automorphism of $C_2^-$. So finally we obtain that $\tilde{C}_1^+$, and hence $C_1^+$, is isomorphic to $C_2^+$ as a $4$-pointed curve, up to a permutation in $\Gamma(C_2^-)$.

$(\beta)$ Suppose now that there is more than one connected component of $C_2^-$. If one of the statements $(ii)$ or $(iii)$ applies, then we are done for the same reasons as above. Otherwise consider the curve
\[ \overline{C} := C - C_A^- .\]
Clearly the image of the subcurve $C_A^-$ under the isomorphism $\phi:C_2^-\rightarrow C_1^-$ is a connected component of $C_1^-$. Since there is a point $P\in C_A^-$ with $\phi(P)=P$, we must have $\phi(C_A^-) = C_A^-$. 
Put
\[ \overline{C}_1^- := C_1^- - C_A^- \Text{ \: \: and \: \: }  \overline{C}_2^- := C_2^- - C_A^- .\]
Note that the decomposition 
\[ \overline{C} = \overline{C}_1^- + C_1^+ = \overline{C}_2^- + C_2^+ \]
satisfies all assumptions of the lemma, and the number of irreducible components of $\overline{C}$ is one less than the number of irreducible components of ${C}$. Hence induction applies, and we can conclude that $C_1^+$ is isomorphic to $C_2^+$ as a $4$-pointed curve, up to a permutation in $\Gamma(\overline{C}_2^-)$. By construction of $\overline{C}_2^-$, any automorphism on it can be trivially extended to an automorphism of $C_2^-$. Therefore $C_1^+$ is isomorphic to $C_2^+$ as a $4$-pointed curve, up to a permutation in $\Gamma({C}_2^-)$, as stated in the lemma. 
\ebew 

The proof of the lemma is completed only up to the claim stated above. Before we can verify this claim we need some preparations.

\begin{defi}\em
Let $f:C\rightarrow \Spec(k)$ be an $m'$-pointed prestable curve of some genus $g'$. Let $P_1,P_2 \in C$. We define the {\em distance} between $P_1$ and $P_2$ as the minimal number of irreducible components of a connected subcurve of $C$ which contains both $P_1$ and $P_2$. This number shall be denoted by $d_C(P_1,P_2)$\label{n58}. If $P_1$ and $P_2$ lie on different connected components of $C$ we put $d_C(P_1,P_2):= \infty$. 
\end{defi}

\begin{rk}\em\label{R329}
One should think of  $d_C(P_1,P_2)$ as the minimal length of a path connecting the points $P_1$ and $P_2$. Note that  $d_C(P_1,P_2)=1 $ if and only if $P_1$ and $P_2$ lie on the same irreducible component of $C'$. In particular $d_C(P_1,P_1)=1$. 

If $C \subset C'$, then $d_{C'}(P_1,P_2) \le d_C(P_1,P_2)$. Intuitively, there may be a shortcut in $C'$, which is not availiable in $C$. More generally, if $\phi : C \rightarrow C'$ is any morphism, then  
\[ d_{C'}(\phi(P_1),\phi(P_2)) \le d_C(P_1,P_2) .\]
If  $\phi : C \rightarrow C'$ is an an isomorphism, then $d_{C'}(\phi(P_1),\phi(P_2)) =  d_C(P_1,P_2)$. 

Suppose that $C \subset C'$, and that the intersection of $C$ with the closure of the complement of $C$ in $C'$ consists of at most one point. Then $d_{C'}(P_1,P_2) = d_C(P_1,P_2)$.
\end{rk}

\begin{abschnitt}\em\label{PCLAIM}
\proofof{the claim on page \pageref{PG213}}
Recall that we are considering decompositions of the curve $C$ into subcurves $C= C_i^- + C_i^+$, for $i=1,2$. The number of glueing points of $C_i^-$ and $C_i^+$ is given by $k$. The curves $C_i^+$ decompose into lines $C_{i,a}^+ + C_{i,b}^+$, where $C_{1,a}^+ = C_{2,a}^+$ as sets, and $C_{1,b}^+ \neq C_{2,b}^+$. By definition, $C_2^- = \tilde{C}_2^- + C_{1,b}^+$. 

$(i)$ Consider the case $k=3$ first. The following pictures show schematically  all possible configurations for the curve $C$. The closed curve always represents the pointed curve $\tilde{C}_2^-$. Marked points, which are not glueing points of $\tilde{C}_2^-$ with either $C_{1,b}^+$ or $C_2^+$ are omitted from the pictures.

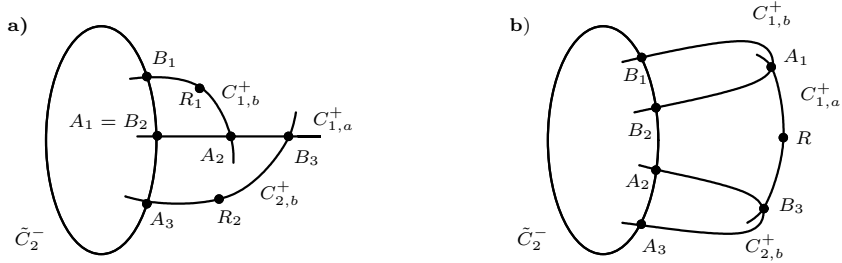
\begin{center}
\setlength{\unitlength}{0.00016667in}
\begingroup\makeatletter\ifx\SetFigFont\undefined%
\gdef\SetFigFont#1#2#3#4#5{%
  \reset@font\fontsize{#1}{#2pt}%
  \fontfamily{#3}\fontseries{#4}\fontshape{#5}%
  \selectfont}%
\fi\endgroup%
{\renewcommand{\dashlinestretch}{30}

}
\end{center} 

Denote the glueing points of $C_2^-$ with $C_2^+$ by $A_1,A_2,A_3$, and 
the glueing points of $C_1^-$ with $C_1^+$ by $B_1,B_2,B_3$. 

In case $a)$ there is one additional distinguished point $R_1$ on $C_{1,b}^+$, which is not a glueing point, and also one point $R_2$ on $C_{2,b}^+$. It may happen that $B_1 = A_3$.  In case $b)$ there is  one additional distinguished point $R$ on $C_{1,a}^+$. Here the intersection $\{A_2,A_3\} \cap \{ B_1,B_2\}$ may be nonempty. 

In case $a)$ we obtain schematically the following picture of the two decompositions of $C$. 
 
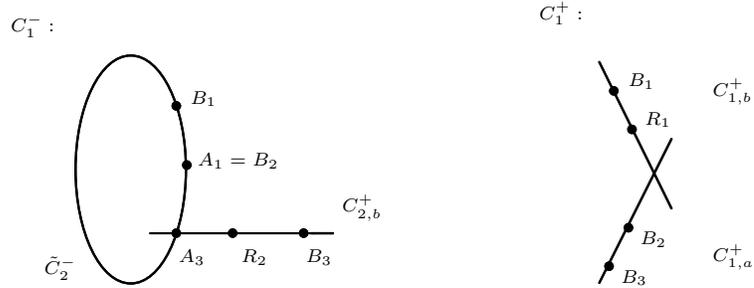
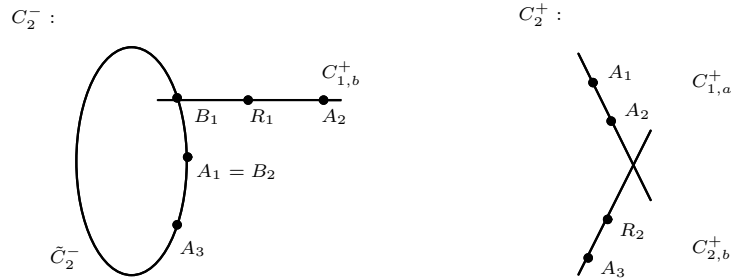
\begin{center}
\setlength{\unitlength}{0.00016667in}
\begingroup\makeatletter\ifx\SetFigFont\undefined%
\gdef\SetFigFont#1#2#3#4#5{%
  \reset@font\fontsize{#1}{#2pt}%
  \fontfamily{#3}\fontseries{#4}\fontshape{#5}%
  \selectfont}%
\fi\endgroup%
{\renewcommand{\dashlinestretch}{30}

}
\end{center} 

Note that $A_1 = B_2$.
A priory there are six permutations $\gamma \in \Sigma_3$, which could be induced by morphisms $\phi: C_2^- \rightarrow C_1^-$ such that $\phi(A_i) = B_{\gamma(i)}$ for $i=1,2,3$. 

If $\phi(A_2)= B_2$ or if $\phi(A_2)=B_3$, then $\phi(A_2)$  is an element of $\phi(C_{1,b}^+)\cap C_{1,a}^+$. There can be no other point in the intersection. Thus the claim holds in this case. So assume now that $\phi(A_2)=B_1$. If $\phi(A_1)=B_3$ then  $\phi^{-1}(B_3)=A_1$ is the unique point in the intersection $\phi^{-1}(C_{2,b}^+)\cap C_{1,a}^+$, and again the claim is true. 

It remains to consider isomorphisms  $\phi: C_2^- \rightarrow C_1^-$ with  $\phi(A_2) = B_1$, $\phi(A_1) = B_2$, and $\phi(A_3) = B_3$. 

In the special case, where $A_3=B_1$, the schematic picture of $C$ is as follows.

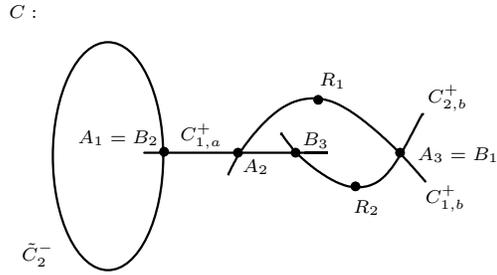
\begin{center}
\setlength{\unitlength}{0.00016667in}
\begingroup\makeatletter\ifx\SetFigFont\undefined%
\gdef\SetFigFont#1#2#3#4#5{%
  \reset@font\fontsize{#1}{#2pt}%
  \fontfamily{#3}\fontseries{#4}\fontshape{#5}%
  \selectfont}%
\fi\endgroup%
{\renewcommand{\dashlinestretch}{30}

}
\end{center} 

Hence the isomorphism $\phi$ can be visualized by the next picture.

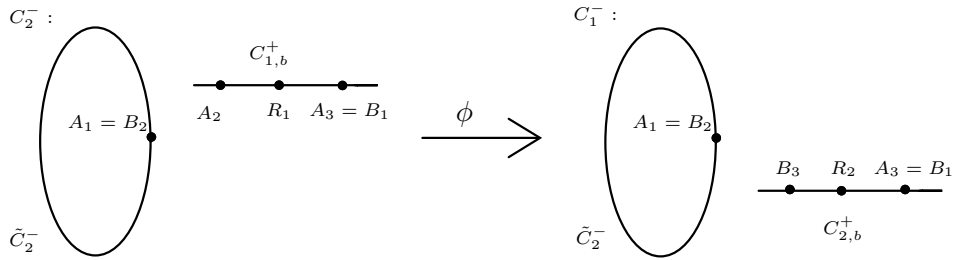
\begin{center}
\setlength{\unitlength}{0.00016667in}
\begingroup\makeatletter\ifx\SetFigFont\undefined%
\gdef\SetFigFont#1#2#3#4#5{%
  \reset@font\fontsize{#1}{#2pt}%
  \fontfamily{#3}\fontseries{#4}\fontshape{#5}%
  \selectfont}%
\fi\endgroup%
{\renewcommand{\dashlinestretch}{30}
%
}
\end{center} 

One sees that the isomorphism $\phi$ can be extended first to an isomorphism from $C_2^- + C_{1,a}^+$ to $C_1^- + C_{1,a}^+$, and then 
to an isomorphism $\overline{\phi}: C_2^-+ C_2^+ \rightarrow C_1^- + C_1^+$, i.e. to an automorphism of $C$. 

We claim that in the cases where $A_3 \neq B_1$ no isomorphism $\phi: C_2^- \rightarrow C_1^-$ with  $\phi(A_2) = B_1$, $\phi(A_1) = B_2$, and $\phi(A_3) = B_3$ exists. Note that here $A_1$ is a fixed point of $\phi$. Consider the connected component $C_{A_1}^-$ of $C_2^-$ which contains $A_1$. Suppose that  $C_{A_1}^-$ contains also the point $A_2$. Then we can  compute the distance of the points $A_1$ and $A_2$ on $C_2^-$ in two different way. First consider 
\[ \begin{array}{ccl}
d_{C_2^-}(A_1,A_2) &=&  d_{C_1^-}(\phi(A_1),\phi(A_2))\\
&=& d_{C_1^-}(B_2,B_1)\\
&=& d_{C_1^-}(A_1,B_1).
\end{array} \]
However, looking at the geometry of the curve $C_2^-$, one reads off the equality
\[ d_{C_2^-}(A_1,A_2) = d_{C_2^-}(A_1,B_1) + 1,\]
a contradiction. Note that because of the special structure of the curves the distances of $A_1$ and $B_1$ are the same on $C_1^-$ and on $C_2^-$, compare remark \ref{R329}. 

Now suppose that $A_3$ is contained in  $C_{A_1}^-$. Then one computes
\[ \begin{array}{ccl}
d_{C_2^-}(A_1,A_3) &=&  d_{C_1^-}(\phi(A_1),\phi(A_3)) \\
&=& d_{C_1^-}(A_1,B_3)\\
 &=& d_{C_1^-}(A_1,A_3)+ 1\\ 
 &=& d_{C_2^-}(A_1,A_3)+ 1,
\end{array}\]
which is a contradiction. The last two equalities follow from the special geometry of the curves $C_1^-$ and $C_2^-$. 

If neither $A_2$ nor $A_3$ lie on $C_{A_1}^-$, then the subcurve $C_{A_1}^-$ satisfies part $(iv)$ of the claim. This settles the case $a)$. 

In case $b)$ the schematic picture of the decompositions of $C$ is

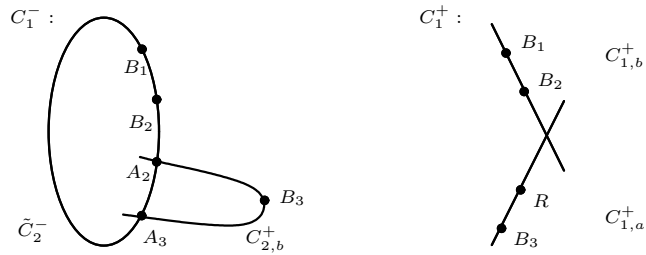
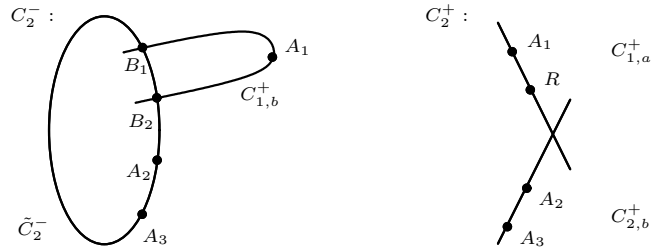
\begin{center}
\setlength{\unitlength}{0.00016667in}
\begingroup\makeatletter\ifx\SetFigFont\undefined%
\gdef\SetFigFont#1#2#3#4#5{%
  \reset@font\fontsize{#1}{#2pt}%
  \fontfamily{#3}\fontseries{#4}\fontshape{#5}%
  \selectfont}%
\fi\endgroup%
{\renewcommand{\dashlinestretch}{30}

}
\end{center} 

Again there are a priory six cases to consider. If $\phi(A_1)=B_3$, then
necessarily $\phi(C_{1,b}^+) = C_{2,b}^+$. The isomorphism $\phi:C_2^-\rightarrow C_1^-$ can then be extended to an automorphism of $C$. 

Let us now look at the cases where  $\phi(A_1)=B_1$ or $\phi(A_1)=B_2$. 

Consider first the special case where $\{A_2,A_3\}\cap \{B_1,B_2\} \neq \emptyset$. If the two sets are equal, then $\tilde{C}_2^-$ is disjoint from $C_1^+\cup C_2^+$, and we already dealt with this case.  Recall that we found in this exceptional situation that $\Gamma(C_2^-) = \Sigma_3$. 

So let us now assume that the intersection of the two sets consists of exactly one point. Suppose without loss of generality that this point is $A_2=B_2$, the second case being analogous.  Look at the following schematic picture of the curve $C$.

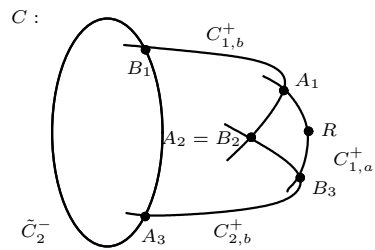
\begin{center}
\setlength{\unitlength}{0.00016667in}
\begingroup\makeatletter\ifx\SetFigFont\undefined%
\gdef\SetFigFont#1#2#3#4#5{%
  \reset@font\fontsize{#1}{#2pt}%
  \fontfamily{#3}\fontseries{#4}\fontshape{#5}%
  \selectfont}%
\fi\endgroup%
{\renewcommand{\dashlinestretch}{30}

}
\end{center} 

Assume first that for the isomorphism $\phi:C_2^-\rightarrow C_1^-$ holds  $\phi(A_1)=B_1$ and $\phi(A_2)=B_2$. Then one computes
\[ 1 = d_{C_2^-}(A_1,A_2) = d_{C_1^-}(B_1,B_2) \ge 2 ,\]
which is of course a contradiction. We obtain the same contradiction under the assumption  $\phi(A_1)=B_1$ and $\phi(A_2)=B_3$, and for  $\phi(A_1)=B_2$ and $\phi(A_2)=B_1$. 

The remaining case is  $\phi(A_1)=B_2$ and $\phi(A_2)=B_3$, so that $\phi(A_3)=B_1$. In this case the isomorphism $\phi: C_2^- \rightarrow C_1^-$ is extendable to an isomorphism from $C_2^- + C_{2,b}^+ $ to $C_1^-+ C_{1,b}^+$, and then further to an isomorphism $\overline{\phi}: C_2^- + C_{2,b}^+ + C_{1,a}^+ \rightarrow C_1^-+ C_{1,b}^+ + C_{1,a}^+$, i.e. to an automorphism of $C$. Note that we have in particular $\overline{\phi}(C_{1,a}^+)= C_{1,a}^+$, and also that we can construct $\overline{\phi}$ in such a way that $\overline{\phi}(R) = R$ holds. 

For further use we will keep in mind that in the situation where $\{A_2,A_3\}\cap\{B_1,B_2\}$ consists of exactly one point, an isomorphism $\phi: C_2^- \rightarrow C_1^-$ exists only under the very special circumstances described above. 

We still need to deal with the case where for the isomorphism $\phi:C_2^-\rightarrow C_1^-$ holds  $\phi(A_1)\in \{B_1,B_2\}$, while  $\{A_2,A_3\}\cap \{B_1,B_2\} =  \emptyset$. We assume without loss of generality that $\phi(A_1)=B_1$ and $\phi(A_2)=B_2$ , the other  cases being analogous.  For this situation we need some preparations.

First, note that if such an isomorphism $\phi$ exists, then there exist an integer $\ell$ such that $\phi^{\ell+1}(A_1)=B_3$. Indeed, assume otherwise. For all $i\in\NN$, for which the point $\phi^i(A_1)$ exists, it is a marked point or a node of $C_1^-$. If  $\phi^i(A_1)\neq B_3$, we must have $\phi^i(A_1) \in \tilde{C}_2^-$ for all $i\in \NN$. Hence there is an infinite sequence $\{\phi^i(A_1)\}_{i\in\NN}$ of nodes and marked points in $C$. If $j\le i$, then $\phi^i(A_1)=\phi^j(A_1)$ implies $\phi^{i-j}(A_1)=A_1$ by the injectivity of $\phi$. As the target set of $\phi$ is  $C_1^-$, and $A_1$ is not an element of $C_1^-$, we must have $i=j$. But then all of the infinitely many  elements of the sequence $\{\phi^i(A_1)\}_{i\in\NN}$ are different, which is absurd. 

Note that for this number $\ell$ also holds $\phi^\ell(C_{1,b}^+) = C_{2,b}^+$, as well as $\phi^\ell(A_1)=A_3$, $\phi^\ell(B_1)=B_3$ and $\phi^\ell(B_2)=A_2$. 

We claim that $\ell > 1$. Indeed, if we had  $\ell =1$, then in particular $\phi(B_1)= B_3$. However, by assumption $\{A_2,A_3\}\cap  
 \{B_1,B_2\} = \emptyset$, so $B_1$ is a node of $C_2^-$. Since $\phi$ is an isomorphism, $\phi(B_1)$ must be a node of $C_1^-$, whereas $B_3$ is a marked point.

The subcurve $C_{1,b}^+$ of $C_2^-$ can be thought of as a ``handle'' attached to the curve $\tilde{C}_2^-$. It is a line, glued to $\tilde{C}_2^-$ in exacly two nodes $B_1$ and $B_2$, with one additional marked point $A_1$ ``on top''. Since $\ell > 1$, we have $\phi(C_{1,b}^+ ) \subset C_1^- - C_{2,b}^+$. We now define
\[ C_{1,b}^{(1)} := \phi(C_{1,b}^+) \Text{ \: \: and \: \: } \tilde{C}^{(1)}_2 := \tilde{C}_2^- -  C_{1,b}^{(1)} .\]
The isomorphism $\phi: C_2^-=\tilde{C}_2^-+C_{1,b}^+ \rightarrow \tilde{C}_2^- + C_{2,b}^+ = C_1^- $ clearly restricts to an isomorphism
\[ \phi^{(1)} : \; \tilde{C}^{(1)}_2   +    C_{1,b}^{(1)} \rightarrow  
\tilde{C}^{(1)}_2+ C_{2,b}^+.\]
Note that $\phi(B_1)\neq B_1$, since $B_1=\phi(A_1)$, and $\phi$ is  injective. We have also $\phi(B_1)\neq B_2$. Otherwise $C_2^-$ would schematically look like follows.

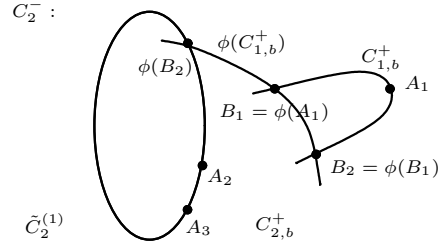
\begin{center}
\setlength{\unitlength}{0.00016667in}
\begingroup\makeatletter\ifx\SetFigFont\undefined%
\gdef\SetFigFont#1#2#3#4#5{%
  \reset@font\fontsize{#1}{#2pt}%
  \fontfamily{#3}\fontseries{#4}\fontshape{#5}%
  \selectfont}%
\fi\endgroup%
{\renewcommand{\dashlinestretch}{30}

}
\end{center} 

From $B_2=\phi(B_1)=\phi^2(A_1)$ follows $\phi(B_2)= \phi^2(B_1)$. Thus $\phi^2(C_{1,b}^+) =  \phi(C_{1,b}^+)$, and hence $\phi^2(B_2)= B_1$. This implies $\phi^2(B_2)=\phi(A_1)$ , and by the injectivity of $\phi$  the equality $\phi(B_2)=A_1\not\in \tilde{C}_2^-$, a contradiction. 

If the intersection of $\{\phi(B_1),\phi(B_2)\}$ and $\{A_2,A_3\}$ contains exactly one point, then we also have  a contradiction. Indeed, we saw above that such an isomorphism $\phi^{(1)}$ can only exist under very special conditions, which are not satisfied here. To see this, let $T$ denote the intersection point of $\{\phi(B_1),\phi(B_2)\}$ and $\{A_2,A_3\}$. From our earlier discussion it follows that a necessary condition for the existence of $\phi^{(1)}$ is the equality $\phi^{(1)}(T)=B_3$. So we must have 
$T = A_3 =\phi(B_1)$. A schematic picture of $C$ is thus as follows. 

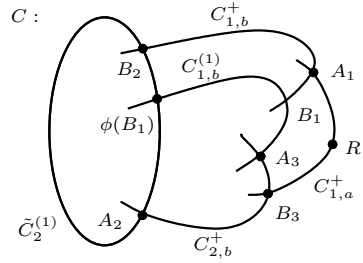
\begin{center}
\setlength{\unitlength}{0.00016667in}
\begingroup\makeatletter\ifx\SetFigFont\undefined%
\gdef\SetFigFont#1#2#3#4#5{%
  \reset@font\fontsize{#1}{#2pt}%
  \fontfamily{#3}\fontseries{#4}\fontshape{#5}%
  \selectfont}%
\fi\endgroup%
{\renewcommand{\dashlinestretch}{30}

}
\end{center} 

The three marked points of $\tilde{C}_2^+ = \tilde{C}_2^{(1)}+ C_{1,b}^{(1)}$ are $B_1$, $A_3$ and $A_2$, as in the picture below.  

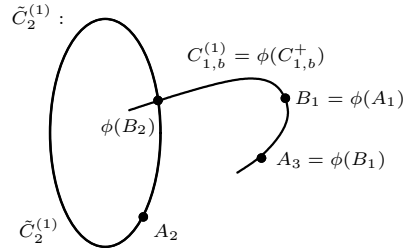
\begin{center}
\setlength{\unitlength}{0.00016667in}
\begingroup\makeatletter\ifx\SetFigFont\undefined%
\gdef\SetFigFont#1#2#3#4#5{%
  \reset@font\fontsize{#1}{#2pt}%
  \fontfamily{#3}\fontseries{#4}\fontshape{#5}%
  \selectfont}%
\fi\endgroup%
{\renewcommand{\dashlinestretch}{30}
%
}
\end{center} 

Our earlier discussion also implies that on the marked points the isomorphism 
 $\phi^{(1)}$ must satisfy $\phi^{(1)}(A_2) = \phi(B_2)$ and $ \phi^{(1)}(B_1)= A_3$. But the first of the two  equations contradicts the injectivity of 
$\phi$ and our initial assumption $\{A_2,A_3\}\cap\{B_1,B_2\}=\emptyset$.

If the two sets are equal, then necessarily $\phi(B_1)= A_3$ and $\phi(B_2)= A_2$. This implies that $\phi^{(1)}( C_{1,b}^{(1)}) = C_{2,b}^+$. In particular, in this case $\ell = 2$ holds. Schematically we have the following picture of $\tilde{C}_2^-$.
There may be other connected components of $\tilde{C}_2^-$, disjoint from $C_1^+\cup C_2^+$. However, they have no effect on the claim, so we may ignore them. 

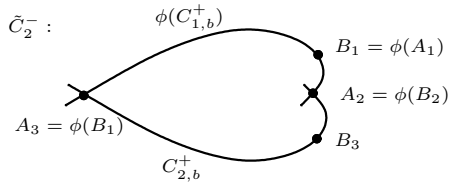
\begin{center}
\setlength{\unitlength}{0.00016667in}
\begingroup\makeatletter\ifx\SetFigFont\undefined%
\gdef\SetFigFont#1#2#3#4#5{%
  \reset@font\fontsize{#1}{#2pt}%
  \fontfamily{#3}\fontseries{#4}\fontshape{#5}%
  \selectfont}%
\fi\endgroup%
{\renewcommand{\dashlinestretch}{30}

}
\end{center} 

Note that this situation cannot occur. Indeed, there is no way in which this curve can be derived from a stable curve be removing a handle, which by assumption is not attched in the point $B_3$. The only point availiable for glueing is $B_1=\phi(A_1)$, which is not enough. 

If the intersection of the two sets $\{\phi(B_1),\phi(B_2)\}$ and $\{A_2,A_3\}$ is empty, then $\phi(C_{1,b}^+)$ is again a handle glued to $ \tilde{C}^{(1)}_2$. Note that the restricted morphism $\phi^{(1)}$ maps the  point on top of the handle to one of its glueing points. We obtain the following picture.

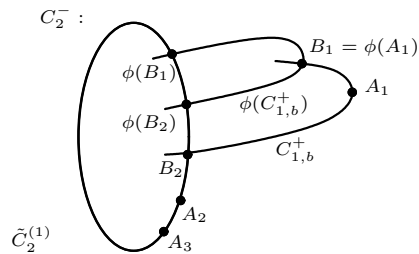
\begin{center}
\setlength{\unitlength}{0.00016667in}
\begingroup\makeatletter\ifx\SetFigFont\undefined%
\gdef\SetFigFont#1#2#3#4#5{%
  \reset@font\fontsize{#1}{#2pt}%
  \fontfamily{#3}\fontseries{#4}\fontshape{#5}%
  \selectfont}%
\fi\endgroup%
{\renewcommand{\dashlinestretch}{30}

}
\end{center} 

Note that $(\phi^{(1)})^{\ell-1}( C_{1,b}^{(1)}) = C_{2,b}^+$. We now proceed inductively with respect to the number $\ell$, by removing further handles, until we obtain a contradiction. 

Summing things up, we conclude that no isomorphism $\phi: C_2^- \rightarrow C_1^-$ exists with $\phi(A_1)=B_1$, except in the special cases where $\{A_2,A_3\} = \{B_1,B_2\}$. In the latter cases we found $\Gamma(C_2^-) = \Sigma_3$. The analogous result holds for isomorphisms with $\phi(A_1)=B_2$. This settles case b).

$(ii)$ We finally need to consider the case $k=4$. The schematic picture in this case is as follows.

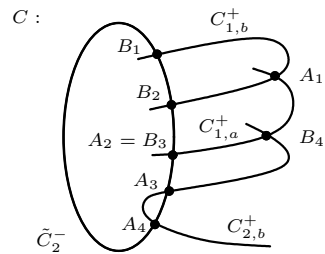
\begin{center}
\setlength{\unitlength}{0.00016667in}
\begingroup\makeatletter\ifx\SetFigFont\undefined%
\gdef\SetFigFont#1#2#3#4#5{%
  \reset@font\fontsize{#1}{#2pt}%
  \fontfamily{#3}\fontseries{#4}\fontshape{#5}%
  \selectfont}%
\fi\endgroup%
{\renewcommand{\dashlinestretch}{30}

}
\end{center}

We denote the glueing points of $C_2^-$ with $C_2^+$ by $A_1,A_2,A_3,A_4$, and 
the glueing points of $C_1^-$ with $C_1^+$ by $B_1,B_2,B_3,B_4$. Note that $A_2 = B_3$. 
The sets $\{A_2,A_3\}$ and $\{B_1,B_2\}$ may have nonempty intersection. 
The two decompositions of $C$ are represented by the next pictures.

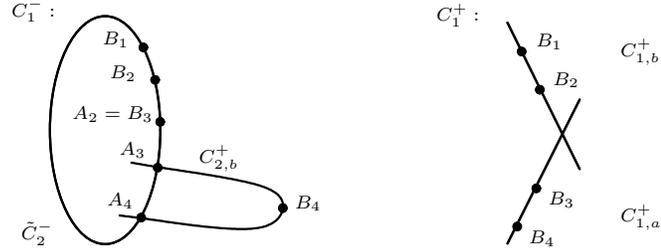
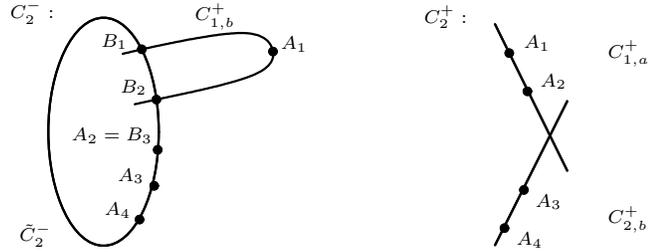
\begin{center}
\setlength{\unitlength}{0.00016667in}
\begingroup\makeatletter\ifx\SetFigFont\undefined%
\gdef\SetFigFont#1#2#3#4#5{%
  \reset@font\fontsize{#1}{#2pt}%
  \fontfamily{#3}\fontseries{#4}\fontshape{#5}%
  \selectfont}%
\fi\endgroup%
{\renewcommand{\dashlinestretch}{30}

}
\end{center} 

In this situation there are 24 possible permutations. If $\phi(A_1)$ equals either $B_3$ or $B_4$, then $\phi(C_{1,b}^+)\cap C_{1,a}^+ \neq \emptyset$. Since there can be no other points of intersection we are done. Similarly, $\phi^{-1}(C_{2,b}^+)\cap C_{1,a}^+ \neq \emptyset$ holds  if $\phi^{-1}(B_4) $ equals either $A_1$ or $A_2$, or equivalently, if either $\phi(A_1)$ equals $B_4$ or  $\phi(A_2)$ equals $B_4$.  
 
So  are left with those 8 permutations, where $\phi(A_1)\in\{B_1,B_2\}$ and $\phi(A_2) \neq B_4$. 

Let us begin with the case  where  $\phi(\{A_1,A_2\})\subset \{B_1,B_2\}$. We can extend the isomorphism $\phi: C_2^- \rightarrow C_1^-$ to an isomorphism $\phi':  C_2^- + C_{1,a}^+  \rightarrow C_1^- + C_{1,b}^+$ in such a way that ${\phi}'(B_4) = A_1$. This is illustrated in the following  picture.

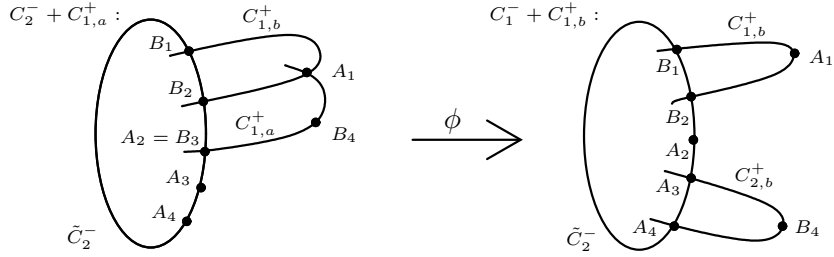
\begin{center}
\setlength{\unitlength}{0.00016667in}
\begingroup\makeatletter\ifx\SetFigFont\undefined%
\gdef\SetFigFont#1#2#3#4#5{%
  \reset@font\fontsize{#1}{#2pt}%
  \fontfamily{#3}\fontseries{#4}\fontshape{#5}%
  \selectfont}%
\fi\endgroup%
{\renewcommand{\dashlinestretch}{30}

}
\end{center} 

 By assuption on $\phi$ we also have 
 $\phi(\{A_3,A_4\})\subset \{B_3,B_4\}$. Thus we can extend the isomorphism ${\phi}'$ to an isomorphism ${\phi}''$ from $ C_2^- + C_{1,a}^+ + C_{2,b}^+$ to $ C_1^- + C_{1,b}^+ + C_{1,a}^+$, with $\phi''(C_{2,b}^+) = C_{1,a}^+$. In other words, $\phi$ extends to an automorphism of $C$, as claimed.

It finally remains to consider the cases where $\phi(A_1)\in\{B_1,B_2\}$ and $\phi(A_2)=B_3$. Note that since $A_2=B_3$, the point $A_2$ is a fixed point of $\phi$.  Assume without loss of generality $\phi(A_1)=B_1$, the second case being completely analogous. We apply some sort of surgery to the curve $C$. Let $\hat{n} : \hat{C}\rightarrow C$ denote the local normalization of $C$ at the point $A_2$. In other words, $\hat{n}$ is an isomorphism except over the point $A_2$, where $\hat{n}^{-1}(A_2)=\{A_{2,a},A_{2,b}\}$. Denote the preimage of $C_1^+$ by $\hat{C}_1^+$, the preimage of 
$C_2^+$ by $\hat{C}_2^+$, and so on. There is a decomposition of $\hat{C}$ by
\[ \hat{C}_1^- + \hat{C}_1^+ = 
\hat{C}_2^- + \hat{C}_2^+ ,\] 
satisfying the assumptions of the lemma, but where the number of glueing points is $\hat{k}=3$. The schematic picture is as follows.

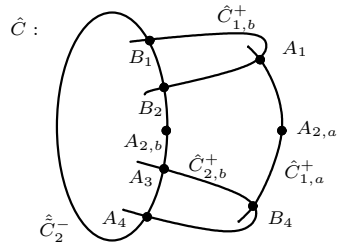
\begin{center}
\setlength{\unitlength}{0.00016667in}
\begingroup\makeatletter\ifx\SetFigFont\undefined%
\gdef\SetFigFont#1#2#3#4#5{%
  \reset@font\fontsize{#1}{#2pt}%
  \fontfamily{#3}\fontseries{#4}\fontshape{#5}%
  \selectfont}%
\fi\endgroup%
{\renewcommand{\dashlinestretch}{30}

}
\end{center} 

We have studied this situation already in part $b)$ of our proof for the case $k=3$ above. There is an isomorphism $\hat{\phi}:\hat{C}_2^- \rightarrow \hat{C_1}^-$ with $\hat{\phi}(A_1)=B_1$ if and only if there is  an isomorphism ${\phi}:{C}_2^- \rightarrow {C_1}^-$.  We saw above that there are  only two cases in which such an isomorphism $\hat{\phi}$ can exist. In the first case  $\hat{\tilde{C}}_2^-$ has empty intersection with $\hat{C}_1^+\cup\hat{C}_2^+$. In this situation  $\hat{C}_2^- = \tilde{C}_2^- \stackrel{\cdot}{\cup} C_{1,b}^+$, and this disjoint union is isomorphic to $C_2^-$. Therefore $\Aut(C_2^-;P_1^{(2)},\ldots,P_\mu^{(2)},\{Q_1^{(2)},\ldots,Q_4^{(2)}\})$ contains the automorphism group of a line with three unlabelled distinguished points, so $\Sigma_3 \subset \Gamma(C_2^-)$.  

In the second case, the isomorphism ${C}_2^- \rightarrow C_1^-$ extends to an automorphism $\overline{\phi}$ of $\hat{C}$, which stabilizes the component $\hat{C}_{1,a}$, with fixed point $\overline{\phi}(A_{2,a})=A_{2,a}$. As an automorphism of pointed curves, the morphism $\overline{\phi}$ permutes only glueing points, whereas all other marked points of $\hat{C}$ are fixed. So in particular $\overline{\phi}(A_{2,b})=A_{2,b}$ holds. Therfore  $\overline{\phi}$ descends to an automorphism of $C$, extending $\phi$. 

This finally concludes the proof. 
\ebew
\end{abschnitt}

%
%

%
%
%
%
%

\chapter{Stacks and Quotient Stacks}

In this appendix we are collecting some basic facts about algebraic stacks, with special focus on quotient stacks. Most of the material here is well known and widely used. Our reason to present it again is mainly to fix our notation, but also the experience that explicit references are sometimes hard to locate. 

This is particularly true for quotient stacks. As this work relies heavily on their formal properties, it seems necessary to discuss them in detail. A number of the statements made in the section on quotient stacks  has not, or not in a quoteable way, appeared in the literature before. 

Many of the definitions and results below can be formulated  in wider generality. However, to faciliate access, we will concentrate on those features of the theory,  which seem most relevant to our studies of moduli problems. 

A highly recommendable introduction to stacks can be found in the appendix to a paper of Vistoli \cite{Vi}. Further  sources, from which I borrowed freely, apart from the book of Laumon and Moret-Bailly \cite{LM},  are the papers \cite{DM},\cite{Ed}, \cite{Gi}, the book  \cite{Hakim}, and many others as listed in the bibliography. Only recently another paper by Vistoli \cite{Vi3} has  become availiable.




%
%

\section{Abstract Preliminaries}

Throughout let $\catc$ be a category in which fibre products exist.

\begin{abschnitt} \em
By definition, a {\em $2$-category}  is a category $\catc$, where in addition for all pairs of objects $X,Y \in \ob(\catc)$ the set of morphisms $\mor_\catc(X,Y)$ between them  equals the set of objects of a category $\catm_{X,Y}$. Morphisms in the category $\catm_{X,Y}$ are called {\em $2$-morphisms}, and one usally writes $\eta: f \Rightarrow g$ for a $2$-morphism $\eta$ from $f: X\rightarrow Y$ to $g: X\rightarrow Y$. Part of the definition is also the existence of a {\em composition functor}
\[ \mu_{X,Y,Z}: \:\:\:  \catm_{X,Y}\times\catm_{Y,Z} \rightarrow \catm_{X,Z} ,\]
for $X,Y,Z\in\ob(\catc)$, which on objects $f\in\catm_{X,Y}$ and $g\in\catm_{Y,Z}$ is just the usual composition of morphisms 
\[ \mu_{X,Y,Z}(f,g) = g\circ f.\]
Furthermore, this composition functor is associative, and compositions with  identity morpisms are trivial. 
\end{abschnitt}

\begin{rk}\em
The composition functor defines a pairing ``$\ast$'' of $2$-morphisms. If $f,g:X\rightarrow Y$ and $f',g':Y\rightarrow Z$ are given in $\catc$, together with two $2$-morphisms $\eta:f\Rightarrow g$ and $\eta':f'\Rightarrow g'$, then there exist a well-defined $2$-morphism 
\[ \eta'\ast\eta := \mu_{X,Y,Z}(\eta,\eta')  : \:\:\: f'\circ f \Longrightarrow g'\circ g, \]
satisfying certain natural compatibility conditions with respect to the composition of morphisms in $\catc$. These conditions can be found in the book of Hakim \cite{Hakim}, from which our notation is taken.  
\end{rk}

\begin{ex}\em
$(i)$ The prototype of all $2$-categories is the category $(\Cat)_2$ of categories. For a pair of categories $\cata$ and $\catb$, the category of morphisms $\catm_{\cata,\catb}$ between them consists of the set of all functors from $\cata$ to $\catb$ as the set of objects, and the set of $2$-morphisms between two functors is  given by the set of all natural transformations between them. \\
$(ii)$ Similarly, the full subcategory of groupoids, together with natural transformations, has the structure of a $2$-category $\Gp_2$. Recall that a {\em groupoid}  is a category, where all morphisms are isomorphisms. Therefore, all of the  $2$-morphisms in $\Gp_2$ are isomorphisms, too. 
\end{ex}

\begin{defi}\em \label{appA3}
A category $\catg$ is called a {\em category fibred in groupoids over
$\catc$} if there is a functor
\[ \mathbf p : \:\: \catg \rightarrow \catc \]
with the following properties.\\
$(i)$  For all morphisms $f : X' \rightarrow X$ in $\catc$ and all
objects $G \in \ob(\catg)$ with $\mathbf p(G) = X$, there exists an object  $G'
\in \ob(\catg)$ together with a morphism $\phi \in
\mor_\catg(G',G)$, such that
\[ \mathbf p(G') = X', \quad \Text{and} \quad \mathbf p (\phi) = f .\]
$(ii)$  For all commutative diagrams
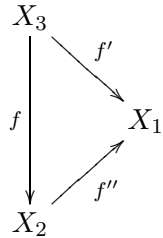
\[ \diagram
X_3 \ddto_f\drto^{f'}\\
&X_1\\
X_2\urto_{f''}
\enddiagram \]
in $\catc$, and all $G_i\in \ob(\catg)$ with $\mathbf p(G_i)=X_i$ for
$i=1,2,3$, and all $ \phi' \in \mor_\catg(G_3,G_1)$, and
$\phi''\in\mor_\catg(G_2,G_1)$ with $\mathbf p(\phi')=f'$ and
$\mathbf p(\phi'')=f''$, there exists a unique morphism $\phi \in
\mor_\catg(G_3,G_2)$ such that 
\[  \phi'' \circ \phi = \phi' \quad \Text{and} \quad \mathbf p (\phi) = f .\]
\end{defi}

\begin{rk}\em
Condition $(ii)$ of definition \ref{appA3} can be reformulated as follows. Consider two diagrams
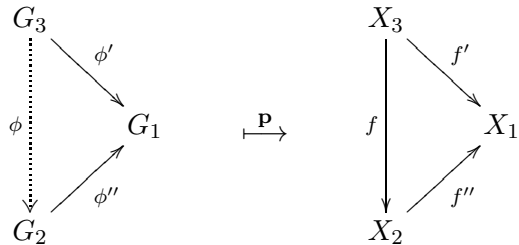
\[ \diagram
G_3 \drto^{\phi'} \dddotted|>\Tip_{\phi}&&&X_3 \ddto_f\drto^{f'}\\
&G_1&\stackrel{\mathbf p}{\longmapsto}&&X_1\\
G_2\urto_{\phi''} &&&X_2\urto_{f''}
\enddiagram \]
where the projection functor $\mathbf p$ sends the left diagram to the right commutative diagram. Then there exists a unique morphism $\phi: G_3\rightarrow G_2$, making the left diagram commuative, and satisfying $\mathbf p(\phi) = f$.  
\end{rk}

\begin{rk}\em
Let $\mathbf p: \catg \rightarrow \catc$ be a category fibred in
groupoids over $\catc$. For any $X \in \ob(\catc)$ we define the {\em
fibre category} $\catg(X)$ by
\[ \ob(\catg(X)) := \{ G \in \ob(\catg): \mathbf p(G) = X \}, \]
and for $G,G' \in \ob(\catc(X))$
\[ \mor_{\catg(X)}(G,G') := \{ \phi \in \mor_\catg(G,G'): \mathbf
p(\phi) = \id_X \} .\]
Often we write shorter $G \in \catg(X) $ instead of $G \in
\ob(\catg(X))$. The second axiom of definition \ref{appA3} implies that a morphism $\phi : G'\rightarrow G$ in $\catg$ is an isomorphism if $\mathbf p(\phi)$ is an isomorphism in $\catc$. Hence  the fibre
category $\catg(X)$ is in fact a groupoid for any $X\in\ob(\catc)$. 
\end{rk}

\begin{defi}\em\label{AD15}
$(i)$ Let $\mathbf p: \catg \rightarrow \catc$ be a category fibred in
groupoids over $\catc$. A {\em normalized lifting} consists of the
choice of an element $G' \in \ob(\catg)$ and a morphism $\phi\in
\mor_\catg(G',G)$ for any tuple $(G,X',X,f)$ as in axiom $(i)$ of definition \ref{appA3}. 
The element $G'$ is denoted by
$f^\ast G  \in \ob(\catg(X'))$, 
and hence $ \phi:
f^\ast G\rightarrow G$.
We require furthermore that for all $G \in \ob(\catg)$ holds 
$ \id_X^\ast G = G$,
with  $\phi = \id_G$. \\
$(ii)$ Let $\mathbf p: \catg \rightarrow \catc$ and $\mathbf p': \catg' \rightarrow \catc$
be categories fibred in groupoids over $\catc$. A {\em morphism of categories fibred in groupoids over $\catc$ between $\catg$ and $\catg'$} is a  functor 
\[ \mathbf F : \:\: \catg \rightarrow \catg' \]
such that
$ \mathbf p'\circ \mathbf F = \mathbf p $. The induced functor between  the fibre categories is denoted by
\[ \mathbf F_X : \:\: \catg(X) \rightarrow \catg'(X), \]
if $X\in\ob(\catc)$. Two such  categories fibred in groupoids over $\catc$   are called {\em isomorphic} if there exists an equivalence of categories $ \mathbf F :  \catg \rightarrow \catg'$ which is a morphism of fibred categories, such that for all $X\in\ob(\catc)$ the functor $\mathbf F_X$ is an equivalence of categories. We denote an isomorphism of fibred categories by ``$\equiv$'', indicating that this notion is weaker than that of an isomorphism of the underlying categories. \\
$(iii)$ The $2$-category $\Gpc_2$ of categories fibred in groupoids over $\catc$ is the subcategory of $(\Cat)_2$, which has as $2$-morphisms only those natural transformations, which are compatible with the projection functors to $\catc$. 
\end{defi}

\begin{abschnitt}\em
From now on we will assume for any fibred category that a normalized lifting has been chosen. 
\end{abschnitt}

\begin{abschnitt}\em
By definition, a {\em $2$-functor} between two $2$-categories $\cata$ and $\catb$ consists of a map
\[ \mathbf P: \:\: \ob(\cata) \rightarrow \ob(\catb) ,\]
and for all $X,Y\in\ob(\cata)$ a functor 
\[  \mathbf P_{X,Y} : \:\: \mor_\cata(X,Y)\rightarrow \mor_\catb(\mathbf P(X),\mathbf P(Y)) ,\]
together with natural transformations
\[ \zeta_{X,Y,Z} : \:\:\: \mu_{\mathbf P(X),\mathbf P(Y),\mathbf P(Z)} \circ (\mathbf P_{X,Y}\times \mathbf P_{Y,Z} ) \Rightarrow \mathbf P_{X,Z}\circ \mu_{X,Y,Z} \]
for all $X,Y,Z\in\ob(\cata)$. The natural transformations $\zeta_{X,Y,Z}$ are required to be isomorphisms of functors, and to satisfy certain natural compatibility conditions with respect to the composition of morphisms. For full details on these conditions see for example \cite[I.1.5]{Hakim}. Often $2$-functors are also called {\em lax functors} or {\em pseudo functors}. 

The set of all $2$-functors from $\cata$ to $\catb$ will be denoted by $\Lax(\cata,\catb)$. 
\end{abschnitt}

\begin{abschnitt}\em
Any category can be viewed as a $2$-category with trivial $2$-morphisms, i.e. identities, and therefore any functor between two categories can be  viewed as a $2$-functor between $2$-categories. 

Let $\catc$ be a category, and let an object $X\in\ob(\catc)$ be given. The {\em functor of points} 
\[ X^\bullet : \:\:  \catc^{op} \rightarrow \Set \]
is given by $X^\bullet(Y):= \mor_\catc(Y,X)$ on objects $Y\in\catc$, and by $X^\bullet(f) := f^\ast:\mor_\catc(Y',X)\rightarrow\mor_\catc(Y,X)$ on morphisms $f\in\mor_\catc(Y,Y')$. The notation $ \catc^{op}$ stands for the {\em opposite category} of $\catc$, reminding us that the functor $X^\bullet $ is contravariant.  

A morphism $g: X_1\rightarrow X_2$ in $\catc$ induces a natural transformation of functors  $\mathbf g: X_1^\bullet \rightarrow  X_2^\bullet$ by composition. 
\end{abschnitt}

\begin{abschnitt}\em
By definition, a {\em sieve} $\mathbf R$ on $X$ is a subfunctor $\mathbf R:  \catc^{op} \rightarrow \Set$ of $X^\bullet$. A {\em Grothendieck topology} $J$ on $\catc$ is defined by specifying for each $X\in\ob(\catc)$ a set $J(X)$ of sieves on $X$, with certain natural additional properties. We do not want to give the general definition here, but concentrate only on those two cases that are needed  below. 

$(i)$ Consider for example the category $\catc = (\Af)$ of affine schemes. The {\em Zariski topology} $J_{Zar}$  on $(\Af)$ is defined as follows. 
For a given scheme $X\in \ob(\catc)$, a sieve $\mathbf R$ is contained in $J_{Zar}(X)$ if and only if there exists a so-called {\em covering family} $\{X_\alpha\}_{\alpha\in A}$ of Zariski open subsets $i_\alpha: X_\alpha\hookrightarrow X$ of $X$, such that $\bigcup_{\alpha\in A} X_\alpha = X$, with the property that for all schemes $Y\in\ob(\catc)$ holds
\[ 
\mathbf R(Y) = \left\{
f: Y\rightarrow X:   \begin{array}{l}
\text{there exists  an index } \alpha\in A,\\
\text{and a morphism } g: X_\alpha\rightarrow Y,\\
\text{such that } i_\alpha = f\circ g: \: X_\alpha\rightarrow X
\end{array} \right\}.\]
$(ii)$ The definition of the {\em \'etale topology} $J_{et}$  on $(\Af)$ is slightly more technical. For a given scheme $X \in\ob(\catc)$ one defines $\mathbf R\in J_{et}(X)$ if and only if there exists a family  $\{X_\alpha\}_{\alpha\in A}$ of Zariski open subsets $i_\alpha: X_\alpha\hookrightarrow X$ of $X$, such that $\bigcup_{\alpha\in A} X_\alpha = X$, and for each $\alpha\in A$ a family $\{f_{\alpha,\beta}: X_{\alpha,\beta} \rightarrow X_\alpha\}_{\beta\in B_\alpha}$, such that $B_\alpha$ is finite, and for each $\beta\in B_\alpha$ the morphism $f_{\alpha,\beta}$ is \'etale, and $\bigcup_{\beta\in B_\alpha} f_{\alpha,\beta}(X_{\alpha,\beta}) = X_\alpha$, with the property that 
for all schemes $Y\in\ob(\catc)$ holds
\[ 
\mathbf R(Y) = \left\{
f: Y\rightarrow X:   \begin{array}{l}
\text{there exists  an index } \alpha\in A, \\
\text{and there exists  an index } \beta\in B_\alpha, \\
\text{and a morphism } g:X_{\alpha,\beta}\rightarrow Y,\\
\text{such that } i_\alpha\circ f_{\alpha,\beta}  = f\circ g: \: X_{\alpha,\beta}\rightarrow X
\end{array} \right\}.\]
The family $\{i_\alpha\circ f_{\alpha,\beta}:X_{\alpha,\beta}\rightarrow X\}_{\alpha\in A,\beta\in B_\alpha}$ is called a {\em covering family} of $X$ with respect to the topology $J_{et}$. 

Details on Grothendieck topologies can for example be found in \cite{Mi}, and of course in \cite[IV]{SGA3}.
\end{abschnitt}

\begin{abschnitt}\em
There is a natural composition of $2$-functors. Thus they can be viewed as morphisms in the category of $2$-categories $(\TwoCat)$. This category can be endowed  with the structure of a $2$-category itself. If $\mathbf P$ and $\mathbf Q$ are two $2$-functors between two $2$-categories $\cata$ and $\catb$, then a $2$-morphism from $\mathbf P$ to  $\mathbf Q$,
or a {\em 2-transformation} between $\mathbf
P$ and $ \mathbf Q$, consists of a pair of families
\[ \{\eta_X\}_{X\in\ob(\catc)} \quad \Text{and} \quad \{u_f\}_{f\in
\mor_\catc(X,Y),\:  X,Y\in\ob(\catc)}, \]
such that for all objects $X \in\ob(\catc)$, one has a morphism 
\[ \eta_X : \:\: \mathbf P(X) \rightarrow \mathbf Q(X) ,\]
and for all $f : X \rightarrow Y$ in $\catc$, one has a $2$-morphism
\[ u_f : \:\: \eta_Y\circ \mathbf
P(f) \Rightarrow  \mathbf Q(f) \circ \eta_X ,\]
which is invertible. There are a few more compatibility conditions imposed, which are all natural, but lenghty when written explicitely. Again, they can be found in \cite{Hakim}. 

With $2$-transformations as its morphisms, the set $\Lax(\cata,\catb)$ of $2$-functors from $\cata$ to $\catb$ can be equipped with the structure of a category.  The category of $2$-categories $(\TwoCat)$ is then a $2$-category itself.
\end{abschnitt}

\begin{abschnitt}\em 
Let $\mathbf p: \catg \rightarrow \catc$ be a category fibred in
groupoids over $\catc$. Let a normalized lifting be chosen on $\catg$. One defines a $2$-functor 
\[ \mathbf G : \:\: \catc^{op} \rightarrow \Text{(Gp)}_2\]
 associated to $\catg$ as follows. 
For any $X \in \ob(\catc)$ put
\[ \mathbf G(X) := \catg(X), \]
and for any morphism $f: X' \rightarrow X$ in $\catc$ define
\[ \funktor{\mathbf G(f)}{\mathbf G(X)}{\mathbf G(X')}{G}{f^\ast
G}{g: \: G'\rightarrow G}{{f^\ast g}: \: f^\ast G' \rightarrow f^\ast
G.} \]
Note that $\mathbf G(f)$  is indeed a functor. However, $\mathbf G$ itself is not a functor. In general, for $f \in
\mor_\catc(X',X)$ and $f' \in
\mor_\catc(X'',X')$ the morphisms 
$ \mathbf G(f') \circ \mathbf G(f)$ and $\mathbf G(f\circ f')$ are not the same, 
since the choice made for some $(f\circ f')^\ast G$ needs not agree with the
choice of ${f'}^\ast(f^\ast G)$. But by axiom $(ii)$ of the definition \ref{appA3} of a fibred category,  there is a well-defined isomorphism of functors
\[ \zeta_{f,f'} : \:\: \mathbf G(f') \circ \mathbf G(f)\Rightarrow \mathbf G(f\circ f').  \] 
Indeed, together  with these  $2$-morphisms, $\mathbf G$ becomes a $2$-functor. \end{abschnitt}

\begin{abschnitt}\em
Just as the category of fibred categories has the structure of a $2$-category $\GrC_2$, the category $\Lax(\catc^{op},\Gp_2)$ of $2$-functors from $\catc^{op}$ to $\Gp_2$ can be given the structure of a $2$-category. This is done by defining for  a given pair of $2$-functors $\mathbf P,\mathbf Q :\catc^{op}\rightarrow \Gp_2 $, and for two given $2$-transformations  $\rho:= (\{\eta_X\}_{X},\{u_f\}_{f})$ and $\rho':= (\{\eta_X'\}_{X},\{u_f'\}_{f})$ between them, morphisms from $\rho$ to $\rho'$ as families $\{ \alpha_X\}_{X\in\ob(\catc)}$, where $\alpha_X: \eta_X\Rightarrow \eta_X'$ is a morphism in the category $\mor_{\Gp_2}(\mathbf P(X),\mathbf Q(X))$ for each $X\in\ob(\catc)$. As usual, certain natural compatibility conditions need to be satisfied. For details we refer again to \cite{Hakim}. 
\end{abschnitt}

\begin{rk}\em\label{correspondence}
There is a natural correspondence between the $2$-categories $\GrC_2$ and $\Lax(\catc^{op}, \Gp_2)_2$, which is formally expressed
by saying that there is  an equivalence of $2$-categories
\[ \GrC_2 \:\: \equiv_2 \:\:\Lax(\catc^{op}, \Gp_2)_2. \]
We indicated above how to obtain a $2$-functor from a fibred category. It is straightforward, even though lengthy, to write down a construction going in the inverse direction. Their composition, considered as a $2$-functor, is equivalent to the identity $2$-functor. This means that there is a pair of $2$-transformations $\rho$ and $\eta$ between them, going in opposite directions, such that $\rho\circ\eta$ and $\eta\circ\rho$ are isomorphic to the respective identity $2$-transformations, in the sense of the previous paragraph. 
\end{rk}

\begin{defi}\em\label{AppD16}
Let $\catc$ be a category, and let $J$ be a Grothendieck topology on it. 
A $2$-functor $ \mathbf F : \catc^{op}\rightarrow \Text{(Gp)}_2 $
is called a {\em 2-sheaf} on $\catc$,  if for all $X \in \ob(\catc)$, and
all sieves $\mathbf R \in J(X)$ the natural functor
\[ \mor_{\Lax( \scat{C}^{op},\Gp_2)}(X^\bullet, \mathbf F) \:\:  \longrightarrow \:\: 
\mor_{\Lax(\scat{C}^{op},\Gp_2)}(\mathbf R, \mathbf F) \]
is an equivalence of categories. 
\end{defi}

\begin{abschnitt}\em
The translation of the notion of a $2$-sheaf into the language of fibred categories, using the equivalence of remark \ref{correspondence}, leads to the notion of stacks, as we will see in proposition \ref{A0124} below.
\end{abschnitt}

\begin{prop}{\em (Yoneda's Lemma)}\label{YON}
Let $\mathbf F: \catc^{op} \rightarrow \Gp_2$ be a $2$-functor, and let $X \in \ob(\catc)$. Then there is a canonical  equivalence{\em
\[ \mor_{{\Lax(\scat{C}^{op},\Gp_2)}}(X^\bullet,\mathbf F) \:\: \equiv \:\: \mathbf F(X) \]}
of categories.
\end{prop}

\proof
This is a direct generalization of Yoneda's lemma as it is known for functors into the category of sets. For a proof see for example \cite[I.2.2]{Hakim}. 
\ebew

\begin{rk}\em\label{AppR114}
The $2$-category of groupoids over the category $\catc$ admits a notion of projective $2$-limits, see \cite{Hakim}. 
One can show that a $2$-morphism $\mathbf F: \catc^{op}\rightarrow \Gp_2$ is a $2$-sheaf if and only if for all $X\in\ob(\catc)$ and all sieves $\mathbf R\in J(X)$ the canonical functor
\[ \mathbf F(X) \:\rightarrow \:  \prlim \: \mathbf F \]
is an equivalence of categories. 
\end{rk}

\begin{defi}\em
Let $\mathbf p: \catg \rightarrow \catc$, $\mathbf p': \catg' \rightarrow \catc$ and $\mathbf p'': \catg'' \rightarrow \catc$ be fibred categories over $\catc$. Let $\mathbf F': \catg' \rightarrow \catg$ and $\mathbf F'': \catg'' \rightarrow \catg$ be morphisms of fibred categories. Then the {\em fibre product $\catg'\times_\scat{G} \catg''$ of $\catg'$ and $\catg''$ over $\catg$} is the category defined as follows. 
Its set of objects is given by
\[ \begin{array}{l}
\ob(\catg'\times_\scat{G} \catg'') := \\[2mm]
 := \left\{ (G',G'',\alpha) : 
\begin{array}{l}
G'\in \ob(\catg'), G''\in \ob(\catg'') \Text{ s.th. } \mathbf p'(G') = \mathbf p''(G''), \\
\alpha \in \mor_\scat{G}(\mathbf F'(G'),\mathbf F''(G'')) \Text{ s.th. } \mathbf p(\alpha) = \id_{\mathbf p'(G')} 
\end{array} \right\}. 
\end{array}\]
For any pair of objects $(G_1',G_1'',\alpha_1)$ and $(G_2',G_2'',\alpha_2)$, the set of morphisms is given by
\[\begin{array}{l}
\mor_{\catg'\times_\scat{G} \catg''}((G_1',G_1'',\alpha_1),(G_2',G_2'',\alpha_2)) :=\hfill\\[3mm]
 := \left\{(\phi',\phi''): 
\begin{array}{l}
\phi':G_1'\rightarrow G_2', \: \phi'':G_1''\rightarrow G_2''\Text{ s.th. } \mathbf p'(\phi') = \mathbf p''(\phi''), \\
\Text{and } \alpha_2 \circ \mathbf F'(\phi') = \mathbf F''(\phi'')\circ \alpha_1 
\end{array} 
\right\}.\hspace*{1cm}
\end{array} \]
There is a natural associative composition law for morphisms, so that the fibre product  $  \catg'\times_\scat{G} \catg''$ is a category, and in fact a category fibred in groupoides over $\catc$.  In those cases where we want to emphasize the defining morphisms of the fibre product, we will write $\catg'\times_{\mathbf F',\scat{G},\mathbf F''} \catg''$ instead of $\catg'\times_\scat{G} \catg''$. 
\end{defi}

\begin{abschnitt}\em\label{A0120}
$(i)$ Let $\catf$, $\catg$ and $\cath$ be categories fibred in groupoides over $\catc$. Let $\mathbf F : \catf \rightarrow \catg$, $\mathbf G : \catg\rightarrow \cath$ and $\mathbf H : \catf\rightarrow \cath$ be morphisms between them. We say that a diagram 
\[ \diagram 
\catf \rrto^{\mathbf H}\drto_{\mathbf F} && \cath\\
& \catg \urto_{\mathbf G}
\enddiagram\]
is {\em commutative}, if there exists a $2$-morphism $\eta: \mathbf G\circ\mathbf F \Rightarrow \mathbf H$, which is necessarily invertible. Using the same notation, a {\em commutative triangle} is by definition a tuple $(\mathbf F,\mathbf G,\mathbf H,\eta)$,  with  $\eta: \mathbf G\circ\mathbf F \Rightarrow \mathbf H$, . 

$(ii)$ Two fibred categories $\catf$ and $\catg$ are isomorphic if and only if there exist two morphisms of fibred categories $\mathbf F : \catf \rightarrow \catg$ and $\mathbf G : \catg \rightarrow \catf$, such that both diagrams
\[ {\diagram 
\catf \rrto^{\id_{\sscat{F}}} \xto[dr]_{\mathbf F}&& \catf\\
&\catg \urto_{\mathbf G}
\enddiagram }
\quad \Text{ and } \quad 
{ \diagram
& {\catf} \xto[dr]^{{\mathbf F}}\\
\catg\rrto_{\id_{\sscat{G}}} \urto^{\mathbf G} && \catg
\enddiagram}\]
commute. 
A diagram of fibred categories over $\catc$
\[\diagram
\cate \rto^{\mathbf E} \dto_{\mathbf G}&\catf\dto^{\mathbf F}\\
\catg \rto_{\mathbf H}&\cath
\enddiagram\]
is called {\em commutative}, if there exists a $2$-morphism $\eta: \mathbf F\circ \mathbf E \Rightarrow \mathbf H\circ\mathbf G$. A {\em commutative square} is a tuple $(\mathbf E,\mathbf F,\mathbf G,\mathbf H,\eta)$, explicitely specifying one $2$-morphism $\eta: \mathbf F\circ \mathbf E \Rightarrow \mathbf H\circ\mathbf G$, which makes the diagram commute. 

$(iii)$ The above commutative square is called {\em Cartesian}, if it satisfies the following universal property. For any fibred category $\catz$ over $\catc$, together with morphisms $\mathbf P : \catz \rightarrow \catf$ and $\mathbf Q:\catz\rightarrow\catg$, and a $2$-morphism $\rho: \mathbf F \circ \mathbf P \Rightarrow \mathbf H \circ \mathbf Q$, there exists a morphism $\mathbf U : \catz \rightarrow \cate$, together with $2$-morphisms $\mu : \mathbf E\circ \mathbf U \Rightarrow \mathbf P$ and $\nu:\mathbf G\circ\mathbf U \Rightarrow  \mathbf Q$, such that 
\[ (\id_{\mathbf F}\ast\mu)\circ(\eta\ast\id_{\mathbf U}) = \rho\circ(\id_{\mathbf H}\ast \nu) .\]
Furthermore, the morphism $\mathbf U$ must be unique in the following sense. If $\mathbf U': \catz\rightarrow \cate$ is another morphism, together with $2$-morphisms $\mu' : \mathbf E\circ \mathbf U' \Rightarrow \mathbf P$ and $\nu':\mathbf G\circ\mathbf U' \Rightarrow  \mathbf Q$, satisfying the analogous compatibility condition, then there exists a $2$-morphism $\delta:\mathbf U \Rightarrow \mathbf U'$, such that $\mu=\mu'\circ(\id_{\mathbf E}\ast\delta)$ and $\nu=\nu'\circ(\id_{\mathbf G}\ast\delta)$.
\end{abschnitt}

\begin{prop}\label{A0117}
Let $\mathbf p: \catg \rightarrow \catc$, $\mathbf p': \catg' \rightarrow \catc$ and $\mathbf p'': \catg'' \rightarrow \catc$ be fibred categories over $\catc$. Then the diagram
\[\diagram 
\catg'\times_\scat{G} \catg'' \rrto^{\bpr_1} \dto_{\bpr_2}&& \catg' \dto^{\mathbf F'} \\
\catg''\rrto_{\mathbf F''}&& \catg 
\enddiagram\]
is $2$-commutative with $2$-morphism $\eta: \mathbf F' \circ \bpr_1 \Rightarrow \mathbf F'' \circ \bpr_2$, such that for all objects $(G',G'',\alpha)\in \ob(\catg'\times_\catg\catg'')$ holds $\eta_{(G',G'',\alpha)}= \alpha$. The diagram is even Cartesian, satisfying the following stronger universal property.  For any fibred category $\catz$ over $\catc$, together with morphisms $\mathbf Q' : \catz \rightarrow \catg'$ and $\mathbf Q'':\catz\rightarrow\catg''$, and a $2$-morphism $\rho: 
\mathbf F' \circ \mathbf Q' \Rightarrow \mathbf F'' \circ \mathbf Q''$, there exists a unique morphism $\mathbf U : \catz \rightarrow \catg'\times_\scat{G} \catg''$, satisfying $\bpr_1\circ \mathbf U = \mathbf Q'$ and $\bpr_2\circ \mathbf U = \mathbf Q''$, and such that $\rho = \eta\ast\id_{\mathbf U}$.
\end{prop}

\proof
It is straightforward to verify that the diagram is $2$-commutative. For an object $(G',G''\alpha)\in \ob(\catg'\times_\scat{G}\catg'')$ one has by definition an isomorphism 
\[ \mathbf F' \circ \bpr_1(G',G'',\alpha) =\mathbf F'(G') \:\: \stackrel{\alpha}{\longrightarrow} \:\: \mathbf F''(G'')  = \mathbf F'' \circ \bpr_2(G',G'',\alpha).\]
Let $\catz$ be a fibred category, together with morphisms $\mathbf Q' : \catz \rightarrow \catg'$ and $\mathbf Q'':\catz\rightarrow\catg''$, and a $2$-morphism $\rho: 
\mathbf F' \circ \mathbf Q' \Rightarrow \mathbf F'' \circ \mathbf Q''$. For an object $Z \in \ob(\catz)$ define
\[ \mathbf U (Z) := (\mathbf Q'(Z),\mathbf Q''(Z), \rho_Z) .\]
Note that $\rho_Z$ is by assumption an isomorphism from $\mathbf F'(\mathbf Q'(Z))$ to $\mathbf F''(\mathbf Q''(Z))$. For a morphism $f : Z_1\rightarrow Z_2$ in $\catz$ define
\[ \mathbf U(f) := (\mathbf Q'(f), \mathbf Q''(f) ) ,\]
so that $\rho_{Z_1}\circ \mathbf F'(\mathbf Q'(f)) = \mathbf F''(\mathbf Q''(f))\circ \rho_{Z_1} $ is automatically satisfied. This assignment is functorial, and a morphism of fibred categories. 

Clearly, the two conditions $\bpr_1\circ \mathbf U = \mathbf Q'$ and $\bpr_2\circ \mathbf U = \mathbf Q''$ hold true, and these two conditions determine the first two coordinates of $\mathbf U(Z)$ uniquely, as well as the pair $\mathbf U(f)$.  The third condition $\rho=\eta\ast\id_{\mathbf U}$ translates into the equality
$\rho_Z = (\eta\ast\id_{\mathbf U})_Z$
for any object $Z\in\ob(\catz)$. By definition of $\eta$, one has  
\[ (\eta\ast\id_{\mathbf U})_Z = \eta_{\mathbf U(Z)} =\bpr_3(\mathbf U(Z)).\]
Hence the above definiton of $\mathbf U$ satisfies the third condition $\rho=\eta\ast\id_{\mathbf U}$. Conversely, this condition determines in combination with the first two the triple $\mathbf U(Z)$ uniquely.
\ebew

With the following proposition we give a useful criterion to decide whether a given fibred category is a fibre product. Essentially, what we do is to reduce the task of verifying the universal property for all stacks to a verification on objects and morphisms in the base category $\catc$ only. 

\begin{prop}\label{A0118}
Consider a commutative diagram of fibred categories
\[\diagram 
\catf \rrto^{\mathbf Q'} \dto_{\mathbf Q''}&& \catg' \dto^{\mathbf F'} \\
\catg''\rrto_{\mathbf F''}&& \catg, 
\enddiagram\]
with $2$-morphism $\rho:\mathbf F'\circ \mathbf Q' \Rightarrow \mathbf F''\circ \mathbf Q''$. Suppose that $\catf$ satisfies the following universal property.\\
$(i)$ For all objects $S\in\ob(\catc)$, together with morphisms $\mathbf s': S^\bullet \rightarrow \catg'$ and $\mathbf s'': S^\bullet \rightarrow \catg''$, and a $2$-morphism $\sigma: \mathbf F'\circ\mathbf s' \Rightarrow \mathbf F''\circ\mathbf s''$, there exists a unique morphism $\mathbf t: S^\bullet\rightarrow \catf$ such that $\mathbf Q'\circ \mathbf t = \mathbf s'$, $\mathbf Q''\circ \mathbf t = \mathbf s''$, and $\sigma = \rho\ast\id_{\mathbf t}$. \\
$(ii)$ For any pair of objects $S_1,S_2\in\ob(\catc)$, together with morphisms $\mathbf s_i':S_i^\bullet \rightarrow \catg'$ and $\mathbf s_i'':S_i^\bullet \rightarrow \catg''$, and a $2$-morphism $\sigma_i:\mathbf F'\circ \mathbf s_i' \Rightarrow \mathbf F''\circ \mathbf s_i''$, for $i=1,2$, as well as two morphisms $\phi': \mathbf s_1'\rightarrow \mathbf s_2'$ in $\catg'$ and $\phi'': \mathbf s_1''\rightarrow \mathbf s_2''$ in $\catg''$, satisfying
\[ \mathbf F'(\phi')\circ \sigma_1 = \sigma_2\circ\mathbf F''(\phi'') ,\]
there exists a unique morphism $\tau: \mathbf t_1\rightarrow \mathbf t_2$ in $\catf$, such that $\mathbf Q'(\tau)=\phi'$ and $\mathbf Q''(\tau)=\phi''$, where $\mathbf t_1$ and $\mathbf t_2$ are given as in $(i)$. 

Then there exists an isomorphism of fibred categories
\[ \catf \: \: \equiv \: \: \catg'\times_\scat{G} \catg'' ,\]
which is even an isomorphism of the underlying categories. 
\end{prop}

\proof
By the universal property of the fibre product, there exists a unique morphism $\mathbf U : \catf \rightarrow \catg'\times_\scat{G} \catg''$, satisfying $\bpr_1\circ\mathbf U=\mathbf Q'$, $\bpr_2\circ\mathbf U=\mathbf Q''$ and $\rho=\eta\ast\id_{\mathbf U}$. We construct an inverse functor $\mathbf V : \catg'\times_\scat{G} \catg'' \rightarrow \catf$ as follows. Let $(G',G'',\alpha)\in \catg'\times_\scat{G} \catg''(B)$ be given, for some $B\in\ob(\catc)$. By Yoneda's lemma, there are corresponding morphisms $G' : B^\bullet\rightarrow \catg'$ and $G'':B^\bullet\rightarrow \catg''$, and $\alpha$ is a $2$-morphism from $\mathbf F'\circ G'$ to $\mathbf F''\circ G''$. Hence there exists by assumption a distinguished morphism $\mathbf t: B \rightarrow \catf$, or equivalently, an object  $\mathbf t \in\catf(B)$. We put
\[ \mathbf V(G',G'',\alpha) := \mathbf t .\]
Note that  by construction $\mathbf U(\mathbf t)= (\mathbf Q'(\mathbf t),(\mathbf Q''(\mathbf t),\rho_{\mathbf t}) = (G', G'',\alpha)$, so we have $\mathbf U\circ \mathbf V = \id_{\scat{G}'\times_\sscat{G} \scat{G}''}$ on objects. Conversely, for an object $F\in\ob(\catf)$, we compute $\mathbf V\circ \mathbf U(F) = \mathbf V(\mathbf Q'\circ F, \mathbf Q''\circ F,\rho_F)=F$, where the last equality holds because of the uniqueness assumption of the universal property. 

Consider now a morphism $(\phi',\phi'') : (G_1',G_1'',\alpha_1)\rightarrow (G_2',G_2'',\alpha_2)$ in $\catg'\times_\scat{G} \catg''$. Again, we apply Yoneda's lemma to consider objects of fibre categories as morphisms into the fibred category. Thus we have morphisms $G_i': B_i \rightarrow \catg'$ and  $G_i'': B_i \rightarrow \catg''$, together with $2$-morphisms $\alpha_i: \mathbf F'\circ G_i' \Rightarrow \mathbf F''\circ G_i''$ for $i=1,2$, satisfying
$\mathbf F''(\phi'')\circ \alpha_1 = \alpha_2\circ\mathbf F'(\phi')$. Hence  by assumption, there exists a distinguished morphism $\tau: \mathbf V(G_1',G_1'',\alpha_1)\rightarrow \mathbf V(G_2',G_2'',\alpha_2)$ in $\catf$, and we put
\[ \mathbf V(\phi',\phi'') := \tau .\]
It is straightforward to verify that this assignment is functorial, and a morphism of fibred categories. Applying the functor $\mathbf U$ to the morphism $\tau$, we get $\mathbf U(\tau)=(\mathbf Q'(\tau),\mathbf Q''(\tau))= (\phi',\phi'')$ by construction. Conversely, for a morphism $f: F_1\rightarrow F_2$ in $\catf$, we find $\mathbf V\circ\mathbf U(f)=\mathbf V(\mathbf Q'(f),\mathbf Q''(f)) = f$, by the uniqueness assumption of the universal property.
\ebew

\begin{rk}\em
Note that the strong universal property of the fibre product of proposition \ref{A0117} is not stable with respect to isomorphisms of fibred categories. Recall that an isomorphism of fibred categories is only an equivalence on the underlying categories.

If $\catf$ is isomorphic to $ \catg'\times_\scat{G} \catg''$ as a fibred category, then $\catf$ satisfies in general only the weaker universal property of paragraph \ref{A0120}. The Cartesian property 
 characterizes all $\catf$ which are isomorphic to $ \catg'\times_\scat{G} \catg''$ as fibred categories. Compare also remark \ref{appC59} below, where we discuss an analogous situation in more detail.
\end{rk}

\vfill\pagebreak
%
%

\section{Stacks}

Troughout this section let $\catc$ be one of the following categories: the category $(\Sch)$ of schemes, the full subcategory $(\Af)$ of affine schemes, 
or the relative categories $(\Sch/S)$ or $(\Af/S)$, where $S$ is a scheme. 
To ensure\footnote{\cite[p. 1]{LM}} that there always exist  terminal objects and fibre products in any of the relative categories, we will in addition  assume that the scheme $S$ is affine if $\catc = (\Af/S)$. 

Let $J = J_{et}$ be the \'etale topology on $\catc$.

\subsection{Descent}

\begin{defi}\em
$(i)$ Let $\catf$ be a category fibred in groupoids over $\catc$ with a fixed normalized lifting.
Let $U \in \ob(\catc)$, and let ${\cal U} =
 \{ f_\alpha:U_\alpha\rightarrow U\}_{\alpha\in A}$ be a covering
family of $U$ with respect to the chosen topology $J$. A {\em descent datum of $\catf$ on $\cal U$} is a family of pairs
\[ \{(x_\alpha,\phi_{\beta,\alpha})\}_{\alpha,\beta\in A} \]
such that  $x_\alpha \in \catf(U_\alpha)$ for all $\alpha \in A$, and 
\[ \phi_{\beta, \alpha} : \:\: \overline{f}^\ast_\beta x_\alpha \simto  
\overline{f}^\ast_\alpha x_\beta \]
is an isomorphism for all $(\alpha,\beta)\in A\times A$, where
$\overline{f}_\alpha $ and $\overline{f}_\beta$ are defined by the
Cartesian diagram
\[ \diagram
U_{\alpha,\beta} := U_\alpha \times_U U_\beta
\rrto^{\overline{f}_\alpha} \dto_{\overline{f}_\beta} & & U_\beta \dto^{{f}_\beta}\\
U_\alpha \rrto_{f_\alpha}&& U .
\enddiagram \]
As a shorthand notation one often writes
\[ x_\alpha | U_{\alpha,\beta} := \overline{f}_\beta^\ast (x_\alpha), \]
and similarly for $x_\beta$. With this convention, we require
furthermore that for all $(\alpha,\beta,\gamma) \in A\times A\times A$
the cocycle condition  
\[ \phi_{\gamma,\alpha}|U_{\alpha,\beta,\gamma} = \phi_{\gamma,\beta} |
U_{\alpha,\beta,\gamma} \circ \phi_{\beta,\alpha}
|U_{\alpha,\beta,\gamma} \]
holds, where $U_{\alpha,\beta,\gamma} = U_\alpha \times_U U_\beta
\times_U U_\gamma$. \\

$(ii)$ A descent datum $
\{ (x_\alpha,\phi_{\beta,\alpha)} \}_{\alpha,\beta\in A}$ is called {\em
effective} if there exists an object $x \in \catf(U)$ and isomorphisms 
$\phi_\alpha: x|U_\alpha := f_\alpha^\ast x \simto  x_\alpha$ in
$\st{F}(U_\alpha),$  for all $\alpha \in A$,
such that 
\[ \phi_\beta | U_{\alpha, \beta} = \phi_{\beta,\alpha} \circ
\phi_\alpha|U_{\alpha,\beta} \]
for all $(\alpha,\beta)\in A\times A$. 
\end{defi}

\begin{rk}\em
Let us expand the cocycle condition a little bit. On $U_{\alpha,\beta}$ there is an isomorphism
\[ \phi_{\beta,\alpha} : \:\: x_\alpha|U_{\alpha,\beta} \simto x_\beta|U_{\alpha,\beta} .\]
It gets lifted to an isomorphism on the fibred product $U_{\alpha,\beta,\gamma} := U_{\alpha,\beta} \times_U U_\gamma$
\[ \phi_{\beta,\alpha}|U_{\alpha,\beta,\gamma} : \:\: x_\alpha|U_{\alpha,\beta,\gamma} \simto x_\beta|U_{\alpha,\beta,\gamma} .\] 
Similarly, the isomorphisms $\phi_{\gamma,\beta}$ on $U_{\beta,\gamma}$, and $\phi_{\gamma,\alpha}$ on $U_{\alpha,\gamma}$ induce isomorphisms on $U_{\alpha,\beta,\gamma}$, which are required to satisfy the cocycle condition as stated.  All squares of the following diagram are Cartesian. 

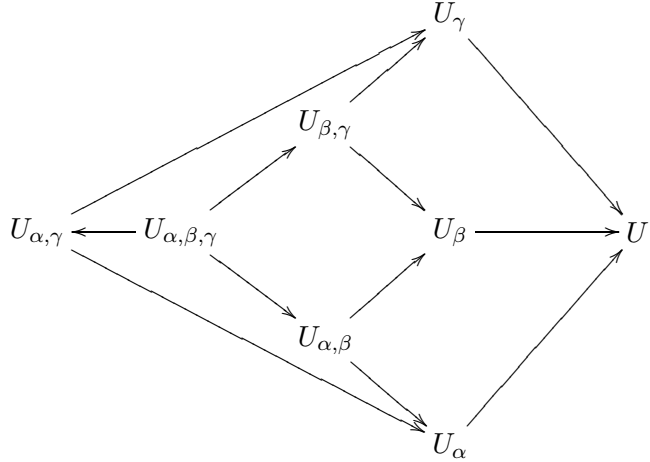
\[ \diagram
&&& U_\gamma \xto[rrdd] \\
&& U_{\beta,\gamma} \urto \drto\\
U_{\alpha,\gamma} \xto[rrruu] \xto[rrrdd] & U_{\alpha,\beta,\gamma} \lto \urto \drto& &U_\beta \rrto&& U \\
&& U_{\alpha,\beta} \urto \drto \\
&&& U_\alpha \xto[rruu]
\enddiagram \]

Note that here we tacidly assume that we have actually identities 
\[ U_{\alpha,\beta,\gamma} = U_{\alpha,\beta} \times_U U_\gamma = U_{\beta,\gamma} \times_U U_\alpha = U_{\alpha,\gamma} \times_U U_\beta ,\]
rather than isomorphisms, for purely typographical reasons. If one wants to be completely formal, one should replace $U_{\alpha,\beta,\gamma}$ in the diagram by the commutative triangle given by the three different fibre products involved, and the natural isomorphisms between them. 
\end{rk}

\begin{defi}\em\footnote{\cite[Definition 3.1]{LM}}
$(i)$ Let $\catf \in \ob\GrC$ be a category fibred in groupoids over
$\catc$ with a normalized lifting. The category $\catf$ is called a {\em prestack over \catc} if
for all $U \in \ob(\catc)$ and all $x,y \in \catf(U)$ the presheaf
\[ \funktor{\Isom_U(x,y)}{(\catc/U)^{op}}{\Set}{f: V \rightarrow
U}{\mor_{\catf(V)}(f^\ast x,f^\ast y)}{g: f\rightarrow f' 
}{c_{g,f}^y\circ g^\ast(.) \circ {c_{g,f}^x}^{-1}}
\]
is a sheaf on $\catc/U$. Here $g: V'\rightarrow V$ and $f': V'
\rightarrow U$ are morphisms in $\catc$ such that $f'= f \circ g$. The morphisms $c_{g,f}^y : g^\ast f^\ast y
\rightarrow (f\circ g)^\ast y = {f'}^\ast y$ are uniquely defined by the chosen normalized lifting, as in  remark
\ref{AD15}, together with the uniqueness property of definition \ref{appA3}$(ii)$, and analogously for $c_{g,f}^x$. 

$(ii)$ A prestack $\catf$ over $\catc$ is called a {\em stack over
$\catc$} if for all $U\in \ob(\catc)$, and for all $J$-covering families
$\cal U$ of $U$, all descent data of $\catf$ on $\cal U$ are effective. 
\end{defi}

\begin{rk}\em
$(i)$ Note that $\Isom_U(x,y)(f)$ is in fact the set of all isomorphisms between $f^\ast x$ and $f^\ast y$, since $\catf(V)$ is a groupoid. 

$(ii)$ Let $\{(x_\alpha,\phi_{\beta,\alpha})\}_{\alpha,\beta\in A}$ be an effective descent datum on $\cal U$, represented by two different elements $x,x' \in \catf(U)$, with $\phi_\alpha : f_\alpha^\ast x \simto x_\alpha$ and  $\phi_\alpha' : f_\alpha^\ast x' \simto x_\alpha$. Then of course
\[ \phi_\alpha^{-1} \circ \phi_\alpha'\in \Isom_U(x',x)(f_\alpha) .\]
If $\Isom_U(x',x)$ is a sheaf, then these composed morphisms glue to give  an isomorphism between $x'$ and $x$. 
\end{rk}

\begin{defi}\em
We define the {\em 2-category of stacks over $\catc$}, denoted by \linebreak 
$(\St/\catc)_2$,  as the full
2-subcategory  of $\GrC_2$, where the objects are the
stacks over $\catc$, and the morphisms are all morphisms of fibred categories in between
them. 
\end{defi}

Recall the notion of a $2$-sheaf from definition \ref{AppD16}. 

\begin{prop}\label{A0124}
Let $\catf$ be a category fibred in groupoids over $\catc$. This
determines, up to isomorphism, a $2$-functor  {\em $\mathbf F: \catc^{op}\rightarrow \Gp_2$}.  Then
$\catf$ is a stack over $\catc$ if and only if $\: \mathbf F$ is a 2-sheaf on
$\catc$. 
\end{prop}

\begin{remark}\em 
Note
that if $\cal U$ is a $J$-covering family of some $U\in\ob(\catc)$, generating the sieve
$\mathbf R := \mathbf R_{{\cal U}}$, then 
an object of the projective limit category $\prlim \: \mathbf F $ defines a descent datum of $\catf$ for $\cal U$, and vice
versa. 
\end{remark}

\proof
The correspondence between fibred categories over $\catc$ and lax functors from $\catc^{op}$  into $\Gp_2$ was given in remark \ref{correspondence}. 

Assume that $\catf$ is a stack over $\catc$. By remark  \ref{AppR114}, it is enough to show that, for any $U \in \ob(\catc)$ and any sieve $\mathbf R \in J(U)$, the canonically existing functor gives  an equivalence of categories
\[ \prlim \: \mathbf F \equiv \mathbf F(U) .\]
Objects of $\prlim \: \mathbf F $ are  nothing else but  descent data of $\catf$,
and since $\catf$ is a stack, each descent datum is effective. Hence
an object of the projective limit determines an object of $\mathbf F(U)$, which is  unique up to isomorphism.   
Let $x,y$ be two objects of $\prlim \: \mathbf F $, now viewed as objects of
$\mathbf F(U)$.  
Then a  morphism in $\prlim \: \mathbf F $ between objects $x$ and $y$ consists of a family of morphisms
\[ a_f : \:\: f^\ast x \rightarrow f^\ast y, \]
for any $f : V \rightarrow U \in \mathbf R(U)$. These morphisms are
necessarily isomorphisms. Since $\Isom_U(x,y)$ is a sheaf, they glue to give a morphism $a : x \rightarrow y$ in the fibre category $\mathbf F(U)$.  

This provides a functor from $\prlim \: \mathbf F $ to  $\mathbf F(U)$, which is in fact an equivalence of categories.  

Conversely, assume that $\mathbf F$ is a 2-sheaf on $\catc$. The first
thing we have to show is that for all $U \in \ob(\catc)$ and for all  $x,y \in \catf(U)$, the presheaf $\Isom_U(x,y)$ is a sheaf. In other words, for all $f : V \rightarrow U \in \ob(\catc/U)$  and all $\mathbf R \in J(f)$, there should be  a bijection of sets
\[ \Isom_U(x,y)(f) \:  \cong \: \lim\limits_{\stackrel{\longleftarrow}{\mathbf R}} \Isom_U(x,y). \]
Note that an element of
 $\lim\limits_{\stackrel{\longleftarrow}{\mathbf R}} \Isom_U(x,y)$
 determines a morphism in $\prlim \: \mathbf F $. 
 But since $\mathbf F$ is a 2-sheaf, such a morphism can be identified
 with a morphism between $f^\ast x $ and $f^\ast y$ in $\mathbf F(V)$,
 i.e. an element of $\Isom_U(x,y)(f)$.  

Analogously, the effectiveness of descent data on $J$-covering families
of $\, U$ follows immediately from the equivalence
\[ \mathbf F(U) \equiv \prlim \: \mathbf F,\]
for any $\mathbf R \in J(U)$. 
\ebew

\begin{rk}\em
\footnote{\cite[7.9]{Vi}}
Let $\mathbf r': \catf'\rightarrow \catf$ and $\mathbf r'': \catf''  \rightarrow \catf$ be morphisms of fibred categories over $\catc$. If $\catf'$, $\catf''$ and $\catf$ are stacks over $\catc$, then so is the fibre product $\catf'\times_{\mathbf r',\catf,\mathbf r''} \catf''$. 
This follows easily from the definition.
 \end{rk}

\subsection{Representability}

\begin{rk}\em
Recall that via the functor of points a scheme $X\in\ob(\catc)$ can be
viewed as a sheaf on $\catc$, so in particular as a 2-sheaf. Moreover,
it determines a fibred category $\cat{X}$ over $\catc$, with fibre
category $X^\bullet(B)$ for $B\in \ob(\catc)$, up to isomorphism of
fibred categories. Hence we can view $X$, or $X^\bullet$ rather, as a stack over
$\catc$. There is an inclusion
\[ \catc \:\: \subset \:\: (\St/\catc) \]
of $\catc$ as a full subcategory. 
\end{rk}

\begin{rk}\em
Let $\st{F}$ be a stack over $\catc$. As $(\St/\catc)_2$ is a
2-category, we have in particular that
$\mor_{\Lax(\scat{C}^{op},\Gp_2)}(X^\bullet,\st{F})$ is a category for any
$X\in\ob(\catc)$. By the generalized version of Yoneda's lemma
\ref{YON}, there is a natural equivalence of categories
\[\st{F}(X) \equiv  \mor_{\Lax(\scat{C}^{op},\Gp_2)}(X^\bullet,\st{F}) . \] 
By abuse of notation one therefore generally identifies morphisms 
 $\mathbf x : X^\bullet \rightarrow \st{F}$ with objects in the
category $\st{F}(X)$, and one even writes $\mathbf x \in \st{F}(X)$. 
\end{rk}

\begin{rk}\em
By definition, a morphism of stacks is a morphism of fibred
categories. Two stacks $\st{F}$ and $\st{G}$ are  called {\em
isomorphic}, if they are 
isomorphic as fibred categories (and hence in particular equivalent as
categories), or equivalently, if they correspond to isomorphic
2-sheaves on $\catc$. This shall be denoted by 
\[ \st{F} \cong \st{G} .\]
This is a weaker notion than isomorphism of categories, but stronger
than equivalence. 
\footnote{ \cite[7.5]{Vi}} 
\end{rk}

\begin{defi}\em
$(i)$ A stack $\st{F}$ is called {\em (strongly) representable}, if
there exists a scheme $X \in \ob(\catc)$ and an isomorphism of stacks,
such that 
\[  X^\bullet \cong \st{F} .\]
$(ii)$ A morphism of stacks $\mathbf r : \st{F} \rightarrow \st{G}$ is
called {\em representable}, if for all schemes $Y \in \ob(\catc)$ and all morphisms 
$\mathbf y : Y^\bullet  \rightarrow \st{G}$ the fibre product in the category
of stacks
\[ \st{F} \times_{\mathbf r, \scat{G},\mathbf y} Y^\bullet \]
can be represented by a scheme in $\catc$.
\end{defi}

\begin{rk}\em
Note that there exists also a more general notion of representability
by algebraic 
spaces instead of schemes. \footnote{\cite[Definition 2.3]{EHKV}} 
\end{rk}

\begin{rk}\em
$(i)$ A scheme representing a stack is by definition unique up to
isomorphism of stacks, and hence unique up to isomorphism of schemes. \\
$(ii)$
If a scheme representing $ \st{F} \times_{\mathbf r, \st{G},\mathbf y}
Y^\bullet $ exists, we will usually denote it by $ \st{F}
\times_{\mathbf r, \st{G},\mathbf y} Y$, so
\[  (\st{F} \times_{\mathbf r, \st{G},\mathbf y} Y)^\bullet \cong
\st{F} \times_{\mathbf r, \st{G},\mathbf y} Y^\bullet . \]
Since the fibre product is unique up to isomorphism of stacks, the
representing scheme is unique up to isomorphism. 
\end{rk}

\begin{ex}\em
$(i)$ If $X \in \ob(\catc)$ is a scheme, then obviously $X$ is
representable precisely by all schemes $X'\in \ob(\catc)$ isomorphic
to $X$.\\ 
$(ii)$ Let $f : X \rightarrow Y$ be a morphism of schemes in
$\catc$, and let $\mathbf f : X^\bullet \rightarrow Y^\bullet$ be the
induced morphism of stacks. Then $\mathbf f$ is representable. 
Indeed, let $Z \in \ob(\catc)$ be given, together with a morphism
$\mathbf g : Z^\bullet \rightarrow Y^\bullet$. By Yoneda's lemma, the morphism 
$\mathbf g$ corresponds to some morphism $g : Z \rightarrow Y$ in
$\catc$. Then one has  the Cartesian square
\[ \diagram
(X \times_{f,Y,g} Z)^\bullet \cong X^\bullet \times_{\mathbf f,
Y^\bullet,\mathbf g} Z^\bullet \rto \dto & Z^\bullet \dto^{\mathbf g} \\
X^\bullet \rto^{\mathbf f} & Y^\bullet 
\enddiagram \]
by the universal property of the fibre product.
\end{ex}

\begin{defi}\em
Let $\cal P$ be a property of morphisms in the category $\catc$. \\
$(i)$ Suppose that $\cal P$ holds for a morphism $f : X \rightarrow Y$
in $\catc$ if and only if $\cal P$ holds for the induced morphism
\[ f_\alpha : \:\: X\times_{f,Y,{u_\alpha}}U_\alpha  \rightarrow Y \]  
for all $\alpha \in A$, where ${\cal U} = \{ u_\alpha : U_\alpha
\rightarrow Y\}_{\alpha\in A}$ is a $J$-covering family of $Y$ in $\catc$. Then $\cal P$ is
called {\em local on the target} with respect to the topology $J$.\\
$(ii)$ Suppose that for any pair of morphisms $f : X\rightarrow Y$ and
$g : Z \rightarrow Y$ in $\catc$
the property $\cal P$ holds for the morphism $\overline{f}$ induced by
base change, if it holds for the 
morphism $f$. Here $\overline{f}$ is defined by the
Cartesian diagram 
\[ \diagram
Z \times_{g,Y,f} X \rto \dto_{\overline{f}} & X \dto^f\\
Z \rto_g & Y.
\enddiagram \]
Then $\cal P$ is said to be {\em stable under base change}.
\end{defi}

\begin{ex}\em
Examples of properties, which are local on the target with respect to the \'etale topology  and
stable under base change are morphisms being\\[2mm]
\begin{tabular}{lll}
- locally of finite type,&
- of finite type,&
- separated,\\
- locally of finite presentation,&
- proper,&
- affine,\\
- surjective,&
- flat,&
- smooth,\\
- open/closed embedding,&
- \'etale,&
- unramified,\\
- universally open/closed,&
- quasi-compact,&
- finite,\\
- integral, & 
- quasi-finite,&
- birational,\\
- dominant,&
- etc.
\end{tabular}\\
For more such  properties see for example \cite[3.10]{LM}.\footnote{\cite[4.10]{DM}}
\end{ex}

\begin{defi}\em
\label{St119}
Let ${\cal P}$ be a property of morphisms of schemes in $\catc$,
which is defined locally on the target and stable under base
change. The attribute ${\cal P}$ {\em holds for a 
representable morphism $\mathbf r : \st{F} \rightarrow \st{G}$} in
$(\St/\catc)$, if it holds for all schemes $S \in \ob(\catc)$ and all objects 
$\mathbf s\in \st{G}(S)$ for the morphism
\[ \overline{ r} : \:\: \st{F}\times_{\mathbf r,\st{G},\mathbf s} S
\: \rightarrow \: S, \]
which is induced by the base change diagram 
\[\diagram
(\st{F}\times_{\mathbf r,\st{G},\mathbf s} S)^\bullet  \rrto^{\overline{\mathbf
r}} \dto_{\overline{\mathbf s}} && S^\bullet \dto^{\mathbf s}\\
\st{F} \rrto_{\mathbf r} &&\st{G}.
\enddiagram \]
\end{defi}

\begin{rk}\em
Note that a property  ${\cal P}$ holds for a  
morphism $f : X \rightarrow Y$ in $\catc$ if and only if it holds 
for the induced morphism $\mathbf f: X^\bullet \rightarrow
Y^\bullet$ in $(\St/\catc)$.
\end{rk}
 
\begin{lemma}
Let $\mathbf r : \st{F}\rightarrow \st{G}$ be a representable morphism
of stacks, and let ${\cal P}$ be a property of $\mathbf r$ as in
definition \ref{St119}. Let
$\mathbf s: \st{G}' \rightarrow \st{G}$ be another morphism of
stacks. Then the morphism $\overline{\mathbf r}$ induced by base
change
\[ \diagram
\st{G}' \times_{\mathbf s,\st{G},\mathbf r} \st{F}
\rrto^{\overline{\mathbf s}} \dto_{\overline{\mathbf r}} && \st{F}
\dto^{\mathbf r}\\
\st{G}' \rrto_{\mathbf s} && \st{G}
\enddiagram \]
is again representable, and it has the property ${\cal P}$.
\end{lemma}

\proof 
Elementary, see \cite[Lemme 3.11]{LM}.
\ebew

\begin{lemma}
$(i)$ Let $S \in \ob(\catc)$ be a scheme, and let $\st{F} \in
\ob(\St/\catc)$ be a stack over $\catc$. If a morphism $\mathbf r:
\st{F}\rightarrow S^\bullet$ is representable, then $\st{F}$ is
representable.\\
$(ii)$ Let $\mathbf r:
\st{F}\rightarrow \st{G}$ and $\mathbf s: \st{G} \rightarrow \st{H}$
be representable  morphisms of stacks, both of them having a property
${\cal P}$ as above, which is furthermore stable under composition of morphisms in $\catc$. Then the composition $\mathbf s \circ
\mathbf r$ is representable as well, and has the property ${\cal P}$.
\end{lemma}

\proof
To prove
$(i)$ consider the Cartesian diagram induced by the identity morphism
on $S$.
\[ \diagram
\st{F} \times_{\mathbf r,S^\bullet,{\id_{S^\bullet}}} S^\bullet
\rrto^{\overline{\mathbf r}} \dto_{\overline{\id}_{S^\bullet}} && S^\bullet
\dto^{\id_{S^\bullet}}\\
\st{F} \rrto_{\mathbf r} && S^\bullet.
\enddiagram \]
Clearly, $\st{F} \times_{\mathbf r,S^\bullet,{\id_{S^\bullet}}}
S^\bullet$ is isomorphic to $\st{F}$, and representable since $\mathbf
r$ is representable. Again, the proof of $(ii)$ is elementary. \footnote{  \cite[Lemme 3.12]{LM}} 
\ebew

\subsection{Algebraic Stacks}

\begin{rk}\em
Let either $\catc = (\Sch/S)$ or $\catc = (\Af/S)$ as above, together with the \'etale topology $J_{et}$ on it\footnote{\cite[p. 1]{LM}}. Since $S$
is terminal in $\catc$, we have an isomorphism of categories $\catc
\cong S^\bullet$. If $\st{F}$ is a stack over $\catc$, then  the projection
functor $\mathbf p : \st{F} \rightarrow \catc$ can be viewed as a
morphism of stacks $\mathbf p : \st{F}\rightarrow S^\bullet$.
\end{rk}

\begin{defi}\em
$(i)$ Let $\st{F}$ be a stack over $\catc = S^\bullet$ with projection functor $\mathbf p: \st{F} \rightarrow S^\bullet$. The {\em diagonal morphism} is the unique morphism
\[ \Delta:  \quad \st{F} \rightarrow \st{F}\times_{\mathbf p,S^\bullet,\mathbf p} \catf  \]
induced by the universal property of the of the fibre  product from the identity morphisms $\id_{\scat{F}} : \st{F} \rightarrow  \st{F}$ on both factors, as in  proposition  \ref{A0117}. \\
$(ii)$ 
A stack \st{F} is called {\em quasi-separated} if the diagonal
morphism $\Delta: \st{F} \rightarrow \st{F}
\times_{\mathbf p, S^\bullet,\mathbf p} \st{F} $
is representable, quasi-compact and separated.
\end{defi}
 
\begin{prop}
\label{St123} 
Let $\st{F}$ be a stack over $\catc$. Then the following are
equivalent.\\
$(i)$ The diagonal morphism $\Delta: \st{F} \rightarrow \st{F}
\times_{\mathbf p, S^\bullet,\mathbf p} \st{F} $ is representable.\\
$(ii)$ For all $X,Y \in \ob(\catc)$, for all $\mathbf x\in \st{F}(X)$,
and for all $\mathbf y \in \st{F}(Y)$, the fibre product $X^\bullet
\times_{\mathbf x, \scat{F},\mathbf y} Y^\bullet$ is representable.\\
$(iii)$ For all $X \in \ob(\catc)$ and for all $\mathbf x, \mathbf x'
\in \st{F}(X)$, the sheaf $\, \Isom_X(x,x')$ is representable by a scheme
over $X$.\\
$(iv)$ For all $X\in \ob(\catc)$ and for all $\mathbf x\in \st{F}(X)$,
 the morphism $\mathbf x:   X^\bullet \rightarrow \st{F}$ is representable.
\end{prop}

\proof
See \cite[Cor. 3.13]{LM}.\footnote{Seminarnotizen}
\ebew

\begin{rk}\em\label{526}
\cite[7.12]{Vi}
Let $\st{F}$ be a stack over $\catc$, and consider the diagonal
morphism
$ \Delta: \st{F} \rightarrow \st{F}
\times_{\mathbf p, S^\bullet,\mathbf p} \st{F}$. 
A morphism $\mathbf x : X^\bullet \rightarrow \st{F}
\times_{\mathbf p, S^\bullet,\mathbf p} \st{F}$ is given by a pair of
objects $x_1,x_1 \in \st{F}(X)$, and an isomorphism between them. The fibre product $\st{F}\times_{\Delta, \scat{F}
\times_{\mathbf p, S^\bullet,\mathbf p} \scat{F},\mathbf x} X^\bullet$ is
isomorphic to $\Isom_X(x_1,x_2)$ as a stack over $X^\bullet$. For a given scheme $b : B\rightarrow X$ over $X$, an object of $\Isom_X(x_1,x_2)(B)$ is an isomorphism $\phi: b^\ast x_1 \rightarrow b^\ast x_2$. This corresponds to the triple $(b^\ast x_1,b,\phi)\in \st{F}\times_{\Delta, \scat{F}
\times_{\mathbf p, S^\bullet,\mathbf p} \scat{F},\mathbf x} X^\bullet(B)$. 

If $\st{F}$ is quasi-separated, then by definition there is a scheme over $X$,
representing  both $\st{F}\times_{\Delta, \scat{F}
\times_{\mathbf p, S^\bullet,\mathbf p} \scat{F},\mathbf x} X^\bullet$ and
$\Isom_X(x_1,x_2)$.
\end{rk}

\begin{defi}\em
A stack $\st{F}$ is called an {\em algebraic stack}, or an {\em Artin stack}, if $\st{F}$ is
quasi-separated, and if there exists a scheme $A \in \ob(\catc)$, and a
smooth surjective morphism 
\[ \mathbf a : \quad A^\bullet \rightarrow \st{F} .\]
Such a morphism $\mathbf a$ is then called an {\em atlas} of $\st{F}$.
\end{defi}

\begin{rk}\em
Note that the condition on $\mathbf a$ only makes sense if $\mathbf a$
is representable. But as we have seen in proposition \ref{St123},
representability of the diagonal morphism $\Delta: \st{F} \rightarrow \st{F}
\times_{\mathbf p,S^\bullet,\mathbf p} \st{F}$ is equivalent to the representability of
all morphisms $\mathbf x : X^\bullet \rightarrow \st{F}$ with $X \in \ob(\catc)$. 
\end{rk}

\begin{defi}\em
An algebraic stack $\st{F}$ is called a {\em Deligne-Mumford stack} if
there exists an atlas
\[ \mathbf a : \quad  A^\bullet \rightarrow \st{F} \]
such that $\mathbf a$ is an \'etale morphism.
\end{defi}

\begin{rk}\em
$(i)$ Algebraic stacks form a full subcategory (ASt/$\catc$) of the category
of stacks over $\catc$, and the same is true for the category of
Deligne-Mumford stacks (DM/$\catc$) over $\catc$. \\
$(ii)$ Let $\catc = (\Af/S)$. Then for all schemes $X \in \ob(\catc)$ we have
as well $X^\bullet \in\ob(\mbox{DM}/\catc)$, and
\[ (\Af/S) \subset (\mbox{DM}/\catc) \]
is a full subcategory.\\
$(iii)$ Each representable stack $\st{F}$ over $\catc = (\Af/S)$ is a Deligne-Mumford
stack. If $\st{F}$ is represented by a scheme $X \in \ob(\catc)$, then
the isomorphism $X^\bullet \rightarrow \st{F}$ is an \'etale atlas.\\
$(iii)$ A scheme $X\in \ob(\catc)$ is a Deligne-Mumford  stack if and only if $X$ is separated and quasi-compact.\\
$(iv)$ Both categories (ASt/$\catc$) and  (DM/$\catc$) contain their  fibre products, since
the properties separated, quasi-compact, smooth and \'etale are stable under
base change.
\end{rk}

\begin{lemma}
Let $\mathbf r : \st{F} \rightarrow \st{G}$ be a morphism of quasi-separated  stacks
over $\catc$. If $\mathbf r$ is representable, and $\st{G}$ is an algebric  stack or a 
Deligne-Mumford stack, then so is $\st{F}$.
\end{lemma}

\proof
The claim follows from the fact that the pullback of an atlas of $\st{G}$ via a representable morphism is an atlas of $\st{F}$. 
\ebew

\begin{prop}\label{532}
Let $\st{F}$ be a Deligne-Mumford stack. Then the diagonal morphism
\[ \Delta: \:\: \st{F} \rightarrow \st{F}
\times_{\mathbf p,S^\bullet,\mathbf p} \st{F} \]
is unramified.\footnote{\cite[Prop. 7.15]{Vi}} 
\end{prop}

\proof
Let $\mathbf a: A^\bullet \rightarrow\st{F}$ be an \'etale atlas. Consider a morphism $\mathbf x : X^\bullet \rightarrow \catf\times_{\mathbf p,S^\bullet,\mathbf p} \st{F}$ for some scheme $X\in\ob(\catc)$.   We need to show that the induced morphism 
\[ \overline{\Delta}: \quad Y^\bullet:= \st{F}\times_{\scat{F}\times_{\mathbf p,S^\bullet,\mathbf p} \scat{F}} X^\bullet  \: \longrightarrow \:  X^\bullet\]
is unramified, where $Y$ denotes  the representing scheme. Consider the commutative diagram
\[\diagram
Y^\bullet\times_\scat{F} A^\bullet \xto[rrr]^\phi \xto[dr]\xto[ddd]_{\mathbf p}&&& X^\bullet\times_{\scat{F}\times_{S^\bullet}\scat{F}}(A\times_S A)^\bullet
\xto[dl]\xto[ddd]^{\mathbf q}\\
& A^\bullet \xto[d]_{\mathbf a} \xto[r]^{\hspace*{-4mm}\tilde{\Delta}}  & (A\times_S A)^\bullet \xto[d]^{\mathbf a\times\mathbf a}\\
& \st{F}\xto[r]_\Delta&\st{F}\times_{S^\bullet}\st{F}\\
Y^\bullet\xto[ur]\xto[rrr]^{\overline{\Delta}}&&&X^\bullet \xto[ul]_{\mathbf x}
\enddiagram\]
The squares  to the right, to the left and at the bottom are Cartesian. The morphism $\phi$ exists and is uniquely determined by the universal property of $ X^\bullet\times_{\scat{F}\times_{S^\bullet}\scat{F}}(A\times_S A)^\bullet$. Some diagram chasing now shows that the upper diagram is Cartesian as well. 

By assumption, $\st{F}$ is a Deligne-Mumford stack, so $\Delta$ is quasi-separated, and hence $\tilde{\Delta}$ is  a closed embedding. Therefore $\phi$ is a closed embedding, too. Since $\mathbf a: A^\bullet \rightarrow \st{F}$ is an \'etale atlas, the morphisms $\mathbf p$ and $\mathbf q$ are \'etale, so in particular unramified. Hence $\overline{\Delta}\circ\mathbf p = \mathbf q\circ \phi$ is unramified, and therefore $\overline{\Delta}$ must be unramified, too.
\ebew

\begin{prop}\label{appB51}
Let $S$ be a Noetherian scheme, and  let $\st{F}$ be a quasi-separated  stack over
$\catc$. If the diagonal morphism $\Delta: \st{F} \rightarrow \st{F}
\times_{\mathbf p,S^\bullet,\mathbf p} \st{F}$ is unramified, and if there exists a scheme
$A \in \ob(\catc)$, which is of finite type over $S$, together with a
smooth and surjective morphism $\mathbf a: A^\bullet \rightarrow \st{F}$,
then $\st{F}$ is a Deligne-Mumford stack.
\end{prop}

\proof
See \cite[Thm. 4.12]{DM} and \cite{Ed}. \footnote{Note the additional assumption which was made there to prove the claim.}
\ebew

\begin{rk}\em
For a separated scheme $X\in\ob(\catc)$ the diagonal morphism $\Delta: X\rightarrow X\times_S X$ is an embedding by definiton. This is not true in general for an algebraic stack $\st{F}$. Here, the diagonal morphism $\Delta: \st{F}\rightarrow \st{F}\times_{\mathbf p,S^\bullet,\mathbf p} \st{F}$ is an embedding if for all schemes $Y\in\ob(\catc)$ and all objects $\mathbf y\in\st{F}\times_{\mathbf p,S^\bullet,\mathbf p}\st{F}(Y)$ the induced morphism 
\[\Delta_Y : \quad \st{F}\times_{\Delta,\scat{F}\times_{\mathbf p,S^\bullet,\mathbf p}\scat{F},\mathbf y}Y^\bullet\rightarrow Y^\bullet\]
is an embedding. By remark \ref{526} and proposition \ref{St123} there is an isomorphism of schemes
\[ \mbox{Isom}_Y(y_1,y_2) \cong \st{F}\times_{\Delta,\scat{F}\times_{\mathbf p,S^\bullet,\mathbf p}\scat{F},\mathbf y}Y \]
over $Y$, if $\mathbf y$ is given by the pair $(y_1,y_2)$ of objects in $\st{F}(Y)$. So $\Delta_Y$ is an embedding if and only if 
\[\delta_Y : \quad  \mbox{Isom}_Y(y_1,y_2) \rightarrow Y \]
is an embedding. In fact, the diagonal $\Delta_\scat{F}: \catc \rightarrow \catf\times_{S^\bullet}\catf$ of a Deligne-Mumford stack is an embedding if and only if $\catf$ is representable by an algebraic space.\footnote{\cite[Rk. 2.8]{Ed}} In general, by the definition of a quasi-separated stack, the morphism $\delta_Y$ is quasi-compact and separated, and unramified by proposition \ref{532}. 
\end{rk} 

\begin{prop}
Let $\st{F}$ be a Deligne-Mumford stack, let $X \in \ob(\catc)$ be
quasi-compact and let $\mathbf x \in \st{F}(X)$. Then the scheme {\em
$\mbox{Isom}_X(\mathbf x,\mathbf x)$}  of automorphisms of $\mathbf
x$ is finite. 
\end{prop}

\proof
See \cite[Cor. 2.1]{Ed}.
\ebew

\subsection{Properties of Deligne-Mumford stacks}

\begin{defi}\em
Let ${\cal P}$ be a property of schemes in $\catc$. Supose that ${\cal
P}$ holds for $X \in \ob(\catc)$ if and only if ${\cal P}$ holds for
all $U_\alpha$, where $\{ f_\alpha: U_\alpha \rightarrow
X\}_{\alpha\in A}$ is a $J_{et}$-covering family of $X$.
Then ${\cal P}$ is called  {\em local with respect to the \'etale
topology}. 
\end{defi} 

\begin{defi}\em
Let ${\cal P}$ be a property of schemes in $\catc$, which is local
with respect to the \'etale topology. We say that ${\cal P}$ holds
for a Deligne-Mumford stack $\st{F}$ if there is an \'etale atlas
$\mathbf a : A^\bullet \rightarrow \st{F}$ such that the property
${\cal P}$ holds for $A$.
\end{defi}

\begin{rk}\em
$(i)$ If a property ${\cal P}$ holds for one \'etale atlas $\mathbf a
: A^\bullet \rightarrow \st{F}$ of
$\st{F}$, then it holds for any other \'etale atlas $\mathbf b:
B^\bullet \rightarrow \st{F} $ as well. Indeed, consider the Cartesian diagram
\[ \diagram
A^\bullet \times_{\mathbf a, \scat{F},\mathbf b} B^\bullet
\rto^{\quad \overline{\mathbf a}} \dto_{\overline{\mathbf b}} & B^\bullet
\dto^{\mathbf b}\\
A^\bullet \rto_{\mathbf a}& \st{F}. 
\enddiagram \]
Since $\st{F}$ is a Deligne-Mumford stack, the fibre product 
$A^\bullet \times_{\mathbf a, \scat{F},\mathbf b} B^\bullet$ is
representable, and $\overline{\mathbf a}$ and $\overline{\mathbf b}$
are both surjective and \'etale. In particular $\mathbf b \circ
\overline{\mathbf a}$ and $\mathbf a \circ
\overline{\mathbf b}$ are \'etale atlases of $\st{F}$. Since ${\cal
P}$ is a local property, it holds for $A$ if and only if it holds for  
$A \times_{\mathbf a, \scat{F},\mathbf b} B$. Therefore it holds for $A$  if and only if holds
for $B$. 

$(ii)$ Let $\st{F}$ and $\st{G}$ be two Deligne-Mumford stacks, which are
isomorphic. If $\st{F}$ has a property ${\cal P}$, then so has
$\st{G}$. This can easily be seen with an argument similar to the one in $(i)$. 
\end{rk}

\begin{ex}\em
Examples of properties, which are local with respect to the
\'etale topology are schemes being\\
\begin{tabular}{ll}
- regular,  & - locally Noetherian,\\  
- normal, & - Cohen-Macaulay,\\
- reduced, & - of characteristic $p$.
\end{tabular}\\
Compare also  \cite[4.7]{LM} for the smooth topology.\footnote{\cite[2.7]{Ed}} 
\end{ex}

\begin{defi}\em
$(i)$ A Deligne-Mumford stack $\st{F}$ is called {\em quasi-compact},
if there is an \'etale atlas $\mathbf a : A \rightarrow \st{F}$ such
that $A$ is quasi-compact.\footnote{\cite[4.10]{DM},  \cite{Ed}}\\
$(ii)$ A Deligne-Mumford stack $\st{F}$ is called {\em Noetherian}, if
$\st{F}$ is quasi-compact and locally Noetherian. 
\end{defi}

\begin{rk}\em 
It follows immediately from the definition that a Deligne-Mumford stack $\catf$ is Noetherian if and only if there exists an atlas $\mathbf a: A^\bullet \rightarrow \catf$, where $A$ is Noetherian. 
Note that 
not every atlas of a quasi-compact Deligne-Mumford stack needs to be
quasi-compact, and not every atlas of a Noetherian Deligne-Mumford stack needs to be Noetherian.\footnote{\cite[7.22]{Vi}}
\end{rk}

\begin{defi}\em 
Let $\cal P$ be a property of morphisms in $\catc$. Suppose that $\cal P$ holds for a morphism $f : X \rightarrow Y$
in $\catc$ if and only if $\cal P$ holds for all morphisms
\[ f_\alpha :X_\alpha \rightarrow Y_\alpha, \]
where $\{g_\alpha: X_\alpha\rightarrow X\}_{\alpha\in A}$ and  $\{h_\alpha: Y_\alpha\rightarrow Y\}_{\alpha\in A}$ are covering families of $X$ and $Y$, respectively, with respect to the \'etale topology, such that for all $\alpha\in A$ the diagram 
\[ \diagram
X_\alpha \rto^{f_\alpha}\dto_{g_\alpha}&Y_\alpha\dto^{h_\alpha}\\
X\rto_f& Y 
\enddiagram\]
commutes. 
Then $\cal P$ is called {\em local on the source and on the target} with respect to the \'etale topology.
\end{defi}

\begin{defi}\em
Let $\mathbf r : \st{F} \rightarrow \st{G} $ be a morphism of
Deligne-Mumford stacks, and let ${\cal P}$ be a property of morphisms
of schemes, which is local with respect to the \'etale topology. We
say that the property ${\cal P}$ holds for $\mathbf r$, if there are \'etale
atlases $\mathbf a: A^\bullet \rightarrow \st{F} $ and $\mathbf b: B^\bullet
\rightarrow \st{G}$, together with a morphism $f: A \rightarrow B$,
such that the diagram
\[ \diagram
A^\bullet \rto^{\mathbf f} \dto_{\mathbf a}& B^\bullet \dto^{\mathbf b}\\
\st{F} \rto_{\mathbf r} & \st{G} 
\enddiagram \]
commutes, and $\cal P$ holds for the morphism $f$.  
\end{defi}

\begin{rk}\em
$(i)$ For any morphism $\mathbf r :\st{F}\rightarrow \st{G}$ of stacks,
there is always a pair of atlases $\mathbf a : A^\bullet \rightarrow
\st{F}$ and $\mathbf b: B^\bullet \rightarrow \st{G}$, together with
a morphism $f : A\rightarrow B$, such that $\mathbf b \circ \mathbf f
= \mathbf r \circ \mathbf a$. Indeed, take a pair of  atlases $\mathbf a_1: A_1^\bullet \rightarrow \catf$ and  $\mathbf b:
B^\bullet \rightarrow \st{G}$. Since $\catg$ is a  Deligne-Mumford stack, the composition $\mathbf r\circ\mathbf a_1$ is representable. Hence there is a scheme $A$ representing $A_1^\bullet\times_{\mathbf r\circ\mathbf a_1,\scat{G},\mathbf b} B^\bullet$. Since $\mathbf b$ is \'etale and surjective, the induced morphism $A^\bullet \rightarrow \catf$ is an atlas, and for the induced morphism $f: A\rightarrow B$  holds $\mathbf b\circ\mathbf f =\mathbf r\circ\mathbf a$.\\
$(ii)$ If $f' : A' \rightarrow B'$ is another morphism between atlases
of $\st{F}$ and $\st{G}$, making the corresponding diagram commutative, then ${\cal P}$ does not necessarily hold for $f'$ as well. \\
$(iii)$ If $\mathbf r$ is representable, and $\cal P$ is stable under base change, then this definition is compatible with definition \ref{St119}. 
\end{rk}

\begin{ex}\em
Examples of properties of morphisms of schemes, which are local on the source and on the target with
respect to the \'etale topology are being\footnote{\cite{Ed},\cite{DM}}\\
\begin{tabular}{ll}
- flat, & - locally of finite type,\\
- smooth,& - locally of finite presentation,\\
- \'etale, & - normal,\\
- unramified,& - etc.
\end{tabular}
\end{ex}

\begin{defi}\em
Let $\mathbf r: \st{F}\rightarrow \st{G}$ be a morphism of
Deligne-Mumford stacks. \\
$(i)$ We say that $\mathbf r$ is {\em quasi-compact}, if for all
schemes $X \in\ob(\catc)$, which are quasi-compact, and for all
$\mathbf x\in \st{G}(X)$ the fibre product $\st{F}\times_{\mathbf
r,\scat{G},\mathbf x} X^\bullet$ is again quasi-compact as a
Deligne-Mumford stack.\footnote{\cite[4.40]{DM}, unprecise in \cite[2.7]{Ed}}\\
$(ii)$  We say that $\mathbf r$ is {\em of finite type}, if $\mathbf
r$ is quasi-compact and locally of finite type.
\end{defi}

\begin{lemma}
Let $\mathbf r: \st{F}\rightarrow \st{G}$ be a morphism of Noetherian Deligne-Mumford stacks. If there exists an atlas $\mathbf a: A^\bullet \rightarrow \st{F}$ such that $\mathbf r\circ\mathbf a$ is of finite type, then $\mathbf r$ is of finite type.\footnote{\cite[7.23]{Vi}}
\end{lemma}

\proof
Note that $\mathbf r\circ\mathbf a$ is necessarily a representable morphism. Choose an atlas $\mathbf b: B^\bullet \rightarrow  \st{G}$ of $\st{G}$. Consider the Cartesian diagram
\[\diagram
A^\bullet\times_{\mathbf a,\scat{G},\mathbf b} B^\bullet \rrto^{\mathbf c} \dto_{\overline{\mathbf b}} && B^\bullet\dto^{\mathbf b}\\
A^\bullet\rto_{\mathbf a} & \st{F} \rto_{\mathbf r} & \st{G} .
\enddiagram \]
Clearly, the composition $\mathbf a\circ\overline{\mathbf b} : A^\bullet\times_{\mathbf a,\scat{G},\mathbf b} B^\bullet \rightarrow \st{F}$ is an \'etale atlas of $\st{F}$. By assumption on $\mathbf r\circ\mathbf a$,  the induced morphism $ c$ is quasi-compact and locally of finite type, and thus $\mathbf r$ is locally of finite type. If $X\in\ob(\catc)$ is a quasi-compact scheme, and $\mathbf x\in \catg(X)$, then the scheme $A\times_{\mathbf a,\scat{G},\mathbf x}X$ is a quasi-compact atlas of the Deligne-Mumford stack $\catf\times_{\mathbf r,\scat{G},\mathbf x} X^\bullet$. Thus $\mathbf r$ is quasi-compact, and hence of finite type. 
\ebew

\begin{defi}\em\label{appB55}
$(i)$ A Deligne-Mumford stack $\st{F}$ is called {\em separated}, if the diagonal morphism
$ \Delta:  \st{F} \rightarrow \st{F}\times_\scat{C}\st{F}$
is proper.\footnote{\cite[Def. 4.7]{DM}} \\
$(ii)$ A morphism $\mathbf r: \st{F}\rightarrow \st{B}$ between two Deligne-Mumford stacks is called {\em separated}, if for all separated schemes $B\in\ob(\catc)$ and for all object $\mathbf b\in\st{G}(B)$  the stack
$\st{F}\times_{\mathbf r,\scat{G},\mathbf b} B^\bullet$ is separated.
\end{defi}

\begin{rk}\em
$(i)$ In an equivalent way, one can define a separated Deligne-Mumford stack by requiring that the diagonal morphism $\Delta$  is finite.\\   
$(ii)$ The definition of a separated stack is compatible with equivalence of stacks. Separated morphism are preserved under base change and under composition.
\end{rk}

\begin{lemma}
Let $\mathbf r: \st{F}\rightarrow \st{G}$ be a morphism of Deligne-Mumford stacks. Then the induced relative diagonal morphism
\[\Delta_{\mathbf r}: \quad \st{F}\rightarrow\st{F}\times_{\mathbf r,\scat{G},\mathbf r}\st{F} \]
is representable, separated and of finite type. 
\end{lemma}

\proof
See \cite[Lemme 7.7]{LM}.
\ebew

\begin{lemma}\label{appB58}
Let  $\mathbf r: \st{F}\rightarrow \st{G}$ be a morphism of Deligne-Mumford stacks. The morphism $\mathbf r$ is separated if and only if the relative diagonal morphism $\Delta_{\mathbf r}:  \st{F}\rightarrow\st{F}\times_{\mathbf r,\scat{G},\mathbf r}\st{F}$ is universally closed. 
\end{lemma}

\proof
See \cite[Lemma 1.4]{Vi}. 
\ebew

\begin{rk}\em
The lemma shows in particular that definition \ref{appB55} is compatible with the earlier definition of separatedness  for representable morphisms.
\end{rk}

\begin{prop}
{\em (Valuative criterion for separatedness)} 
Let $\mathbf r: \st{F}\rightarrow \st{G}$ be a morphism of Deligne-Mumford stacks. The morphism $\mathbf r$ is separated if and only if for every complete discrete valuation ring $R$, with quotient field $K$, and inclusion morphism $i: \: V := \Spec(K) \hookrightarrow \Spec(R) =: U$, and for all objects $\mathbf g_1,\mathbf g_2\in \st{F}(U)$, for all morphisms $\xi \in\mor_{\st{F}(V)}(i^\ast \mathbf g_1, i^\ast \mathbf g_2)$ and $\eta\in\mor_{\st{G}(U)}(\mathbf r(\mathbf g_1),\mathbf r(\mathbf g_2))$ with $\mathbf r(\xi)=i^\ast\eta$, there exists a unique morphism
\[ \rho\in\mor_{\st{F}(U)}(\mathbf g_1,\mathbf g_2) \]
such that $i^\ast \rho = \xi$ and $\mathbf r(\rho) = \eta$.
\end{prop}

\begin{rk}\em
To visualize the requirements of the above proposition, consider the following commutative diagram.
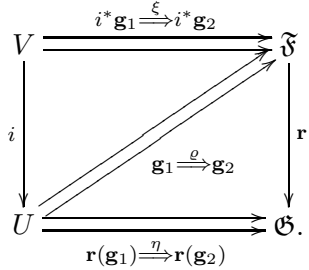
\[\diagram
V \xto[rrr]<-0.5ex> \xto[rrr]<0.5ex>^{i^\ast \mathbf g_1 \stackrel{\xi}{\Longrightarrow} i^\ast \mathbf g_2} \xto[dd]_i&&& \st{F} \xto[dd]^{\mathbf r}\\
\\
U \xto[rrr]<0.5ex> \xto[rrr]<-0.5ex>_{\mathbf r(\mathbf g_1)\stackrel{\eta}{\Longrightarrow} \mathbf    r(\mathbf g_2)} 
\xto[rrruu]<0.5ex>  \xto[rrruu]<-0.5ex>_{\mathbf g_1 \stackrel{\rho}{\Longrightarrow}  \mathbf g_2} &&& \st{G}.
\enddiagram \]
Roughly speaking, the proposition says that, up to $2$-morphisms, if $\mathbf r$ is separated, then for a morphism $V \rightarrow \st{F}$, whose composition with $\mathbf r$ extends to all of $U$,  there exists an extension to a morphism $ U \rightarrow \st{F}$, and this extension is unique.   
\end{rk}

\proof
See \cite[Prop. 3.19]{LM}. The proof uses the preceeding lemma \ref{appB58} and the valuative criterion for universally closed representable morphisms as in \cite[I.5.5.8]{EGA}.
\ebew

\begin{defi}\em 
A morphism $\mathbf r: \st{F}\rightarrow \st{G}$ of Noetherian Deligne-Mumford stacks is called {\em proper}, if 
there exists a scheme $X\in\ob(\catc)$, and a finite surjective morphism $\mathbf x : X^\bullet \rightarrow \st{F}$, such that the composed morphism $\mathbf r\circ\mathbf x: X^\bullet \rightarrow \st{G}$ is proper. 
\end{defi}

\begin{rk}\em
$(i)$ For any Noetherian Deligne-Mumford stack $\catf$ there exists a scheme $X\in\ob(\catc)$ and a finite surjective  morphism $\mathbf x : X^\bullet \rightarrow \st{F}$. The composition  $\mathbf r\circ\mathbf x$ is necessarily representable.\footnote{\cite[Def. 2.9]{Ed}}\\
$(ii)$ The above  definition of a proper morphism agrees with the definition for representable morphisms.  Representable morphisms are preserved under base change and composition. 
\end{rk}

The following proposition is taken from \cite[Thm. 4.19]{DM}. 

\begin{prop}
{\em (Valuative criterion for proper morphisms)}
Let  $\mathbf r: \st{F}\rightarrow \st{G}$ be a morphism of Deligne-Mumford stacks which  is separated and of finite type. Then $\mathbf r$ is proper if and only if for each discrete valuation ring $R$ with quotient field $K$, and inclusion morphism $i : V:= \Spec(K)\hookrightarrow \Spec(R) =: U$, and for all objects $\mathbf x\in \st{F}(V)$, for all $\mathbf y\in\st{G}(U)$ and all morphisms $\xi\in \mor_{\st{G}(V)}(i^\ast\mathbf y, \mathbf r(\mathbf x))$ there exists a triple $(K',\mathbf z, \eta)$ with the following properties. The first entry $K'$ of the triple is a  finite field extension of $K$, with $R'$ as the integral closure of the ring $R$ in $K'$,  and with morphisms $v: V':= \Spec(K') \rightarrow V$ and $u: U':= \Spec(R') \rightarrow U$, as well as $j : V'\rightarrow U'$. The second entry is an object $\mathbf z \in \st{F}(U')$, and the third is a morphism $\eta\in\mor_{\st{F}(V')}(\mathbf v^\ast \mathbf x,j^\ast \mathbf z)$.
\end{prop}

\begin{rk}\em
The valuative criterion for proper morphisms can be used as a definition for proper morphisms between Deligne-Mumford stacks, which are not Noetherian.\footnote{\cite[Def. 1.1]{Vi}}
\end{rk} 

\begin{rk}\em
Note that the separatedness of $\mathbf r$ implies that there exists also a morphism $\rho\in\mor_{\st{G}(V')}(g^\ast \mathbf y,\mathbf r(\mathbf z))$. All of the morphisms $\xi$, $\eta$ and $\rho$ are isomorphisms by the definition of a fibred category. We try to visualize the morphisms in the following diagram.
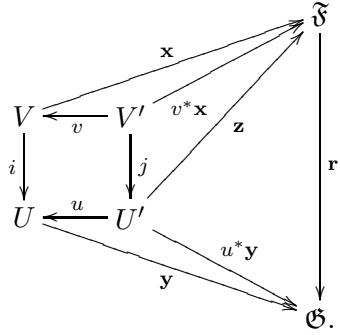
\[\diagram
&&&\st{F}\xto[ddd]^{\mathbf r}\\
V \dto_i\xto[urrr]^{\mathbf x}& V' \lto^{v} \dto^j \xto[urr]_<<<<{v^\ast \mathbf x}\\
U\xto[drrr]_{\mathbf y} & U' \lto_u \xto[drr]^{u^\ast\mathbf y}\xto[rruu]_{\mathbf z}\\
&&&\st{G}.
\enddiagram\]
\end{rk}

\begin{rk}\em
Assume under the conditions of the above proposition in addition that the stack $\st{F}$ is representable. Then the existence of $\mathbf z$ implies the existence of a morphism $\tilde{\mathbf z}: U^\bullet \rightarrow \st{F}$, such that \[ \mathbf r(\tilde{\mathbf z}) \cong \mathbf y, \]
compare \cite[Rk. 2.15]{Ed}.
\end{rk}

\begin{defi}\em
Let $\{\st{F}_i\}_{i\in I}$ be a family of categories fibred in groupoids over $\catc$, for some index set $I$. We define\footnote{\cite[3.3]{LM}} the {\em disjoint sum} of $\{\st{F}_i\}_{i\in I}$ as the fibred category associated to the $2$-functor
\[ \funktor{\coprod\limits_{i\in I}\st{F}_i}{\catc^{op}}{\Gp_2}{X}{\coprod_{i\in I}\st{F}_i(X)}{f: X\rightarrow Y}{\coprod_{i\in I}\st{F}_i(f).} \]
For a scheme $X\in\ob(\catc)$ the category $\coprod_{i\in I}\st{F}_i(X)$ is defined as folows. Its objects are
\[ 
\ob\left(\coprod_{i\in I}\st{F}_i(X)\right) := \left\{ 
\begin{array}{ll}
\{(X_i,F_i)\}_{i\in I} :& \{X_i\}_{i\in I} \text{ is a family of disjoint} \\
&\text{subschemes } X_i\in\ob(\catc) \text{ of } X, \\
&\text{such that} \coprod_{i\in I} X_i =X,  \text{ and}\\
&F_i\in \st{F}(X_i) \text{ for all } i\in I
\end{array} \right\}. \]
Morphism between two objects $\{(X_i,F_i)\}_{i\in I}$ and $\{(Y_i,G_i)\}_{i\in I}$ in $\coprod_{i\in I}\st{F}_i(X)$ are given by 
\[ \begin{array}{l}
\mor_{\coprod_{i\in I}\scat{F}_i(X))}(\{(X_i,F_i)\}_{i\in I},\{(Y_i,G_i)\}_{i\in I}) :=\\[3mm]
\hspace*{2cm} := \left\{ \{\zeta_i\}_{i\in I}: \:\: \zeta_i\in \mor_{\scat{F}_i(X_i)}(F_i,G_i) \text{ for all } i\in I \right\},
\end{array}\]
if $\{X_i\}_{i\in I} = \{Y_i\}_{i\in I}$, and by the empty set otherwise. 

For a morphism $f:X\rightarrow Y$ in $\catc$ one defines the functor $\coprod_{i\in I}\st{F}_i(f)$ as follows.
\[ \funktor{\coprod_{i\in I}\st{F}_i(f)}{\coprod_{i\in I}\st{F}_i(Y)}{\coprod_{i\in I}\st{F}_i(X)}{\{(Y_i,G_i)\}_{i\in I}}{\{(f^{-1}(Y_i),f^\ast G_i)\}_{i\in I}}{ \{\zeta_i\}_{i\in I}}{ \{f^\ast\zeta_i\}_{i\in I}.}\]

If $f:X\rightarrow Y$ anf $g: Y\rightarrow Z$ are two morphisms in $\catc$, then there exists a canonical natural transformation 
\[ \coprod_{i\in I}\st{F}_i(f)\circ 
\coprod_{i\in I}\st{F}_i(g) \Longrightarrow \coprod_{i\in I}\st{F}_i(g\circ f),\]
which is given by a family of maps
\[ \{c_i\}_{i\in I} : \quad \{(f^{-1}(g^{-1}(Z_i)),f^\ast g^\ast  H_i)\}_{i\in I} \rightarrow \{((g\circ f)^{-1}(Z_i),(g\circ f)^\ast  H_i)\}_{i\in I} \]
for objects $\{(Z_i,H_i)\}_{i\in I} \in  \coprod_{i\in I}\st{F}_i(Z)$. 
Clearly, $f^{-1}(g^{-1}(Z_i)) = (g\circ f)^{-1}(Z_i)$ for all $i\in I$, and on the second coordinate the morphism is induced by the canonical morphisms of the normalized liftings $c_{f,g} : f^\ast g^\ast H_i \rightarrow (g\circ f)^\ast H_i$. The compatibility of the $2$-morphisms with morphisms $\{\zeta_i\}_{i\in I}$ in   $\coprod_{i\in I}\st{F}_i(Z)$ follows automatically. 
\end{defi}

\begin{rk}\em
$(i)$ Note that by definition one has
\[ \mor_{\coprod_{i\in I}\scat{F}_i(X)}(\{(X_i,F_i)\}_{i\in I},\{(X_i,G_i)\}_{i\in I}) = \coprod_{i\in I} \mor_{\scat{F}_i(X)}(F_i,G_i) .\]
In particular, since $\st{F}_i(X_i)$ is a groupoid for all $i\in I$, all of the morphisms are isomorphisms, and hence $\coprod_{i\in I}\st{F}_i(X)$ is indeed a groupoid as well. \\
$(ii)$ If for all $i\in I$ the fibred category $\st{F}_i$ is a stack, then the disjoint sum is a stack as well. This is true because all necessary conditions can be verified componentwise. \\
$(iii)$ A disjoint sum $\coprod_{i\in I}\st{F}_i$  is a Deligne-Mumford stack if and only if for all $i\in I$ the stacks $\st{F}_i$ are Deligne-Mumford stacks. \end{rk}

\begin{defi}\em
$(i)$ A {\em substack} $\catg$ of a stack $\catf$ is a representable morphism of stacks $\mathbf i : \catg \rightarrow \catf$, such that for all schemes $B\in\ob(\catc)$ and all objects $\mathbf b\in \catf(B)$, the induced morphism of schemes $\catg\times_{\mathbf i,\scat{F},\mathbf b} B\rightarrow B$ is an embedding.\footnote{\cite[Def. 1.5]{Vi}} \\ 
$(ii)$ A fibred category is called a {\em union} of a family $\{\st{F}_i\}_{i\in I}$ of substacks of $\st{F}$, if the morphism of stacks, which is induced by the inclusion functors
\[ \coprod_{i\in I}\st{F}_i \longrightarrow \st{F}, \]
is surjective.\footnote{\cite[Def. 3.14]{LM}} We then write $\catf = \bigcup\limits_{i\in I} \catf_i$. 
\end{defi}

\begin{ex}\em
Let $B \in\ob(\catc)$ be a scheme, with connected components $\{B_i\}_{i\in I}$. Then one has for the associated stack the same decomposition as disjoint sum
\[ B^\bullet = \coprod_{i\in I}B_i^\bullet .\]
Indeed, let $X\in \ob(\catc)$. Then for the sets of objects we find
\[ \begin{array}{lcl}
\coprod_{i\in I}B_i^\bullet(X) &=& \left\{ \:\{(X_i,f_i)\}_{i\in I}: 
\: \coprod_{i\in I}X_i = X, \: f_i\in B_i^\bullet(X_i) \right\}\\
&=& \left\{ \: \{(X_i,f_i)\}_{i\in I}: \: f\in B^\bullet(X), \: X_i := f^{-1}(B_i), \: f_i:= f|X_i \right\}\\
&=& B^\bullet(X). 
\end{array}\]
A morphism $\{\zeta_i\}_{i\in I} \in \mor_{\coprod_{i\in I}B_i^\bullet}(\{(X_i,f_i)\}_{i\in I},\{(X_i,g_i)\}_{i\in I})$  is given by a family of morphisms
\[ \zeta_i : \:\: f_i \rightarrow g_i ,\]
i.e. morphisms $\zeta_i : X_i \rightarrow X_i$ such that $g_i= f_i\circ \zeta_i$. Note that $f_i = f|X_i$ and $g_i = g|X_i$ for  some morphisms $f,g\in B^\bullet(X)$. Since $\coprod_{i\in I}X_i = X$, the morphisms $\zeta_i$ glue to give a morphism $\zeta: X\rightarrow X$, which satisfies $g = f \circ \zeta$. This shows that 
\[ B^\bullet(X) \cong \coprod_{i\in I}B_i^\bullet(X) \]
as categories, and from this one deduces that $B^\bullet$ and $\coprod_{i\in I}B_i^\bullet$ are isomorphic as fibred categories. 
\end{ex}

\begin{defi}\em
$(i)$ A fibred category $\st{F}$ is called {\em empty} or {\em void}, if $\st{F}$  is represented by the empty scheme $\emptyset\in\ob(\catc)$.\\
$(ii)$  A fibred category $\st{F}$ is called {\em connected}, if $\st{F}$ is nonempty, and if an isomorphism $\st{F}\cong\coprod_{i=1,2}\st{F}_i$ for some fibred categories $\st{F}_1$ and $\st{F}_2$, implies that either $\st{F}_1$  or $ \st{F}_2$ is empty.
\end{defi}

\begin{prop}\label{appB73}
Let $\st{F}$ be a Deligne-Mumford stack, which is locally \linebreak Noetherian. Then there exists a family $\{\st{F}_i\}_{i\in I}$ of connnected Deligne-\linebreak Mumford stacks, such that $\st{F}\equiv \coprod_{i\in I}\st{F}_i $. This family is uniquely determined up order,  and the members of the family are uniquely determined up to equivalence of stacks. 
\end{prop} 

\proof See \cite[Prop. 4.9]{LM}. 
\ebew

\begin{defi}\em
The members of the family $\{\st{F}_i\}_{i\in I}$ of proposition \ref{appB73} are  called the  {\em connected components} of $\st{F}$. The set of connected components of $\st{F}$ is denoted by $\pi_0(\st{F})$.
\end{defi}

\begin{defi}\em\label{appB75}\footnote{\cite[p. 102]{DM}}
Let $\st{F}$ be a Deligne-Mumford stack. An {\em open subset} $\st{U}$ of $\st{F}$ is a full subcategory $\st{U}$ of $\st{F}$, which is again a Deligne-Mumford stack over $\catc$, with the following properties:\\
$(i)$ for all $U\in \ob(\st{U})$ and all $F\in\ob(\st{F})$ with $F\cong U$ holds $F\in \ob(\st{U})$;\\
$(ii)$ the inclusion functor $\iota: \st{U}\rightarrow\st{F}$ is representable by  open immersions, i.e. for all schemes $B\in\ob(\catc)$ and all objects $\mathbf b\in\st{F}(B)$, the morphism $\overline{\iota}$ induced by the Cartesian diagram 
\[\diagram
\st{U}\times_{\iota,\scat{F},\mathbf b} B^\bullet \rto^{\overline{\iota}}\dto& B^\bullet\dto^{\mathbf b}\\
\st{U} \rto_\iota&\st{F}
\enddiagram\]
is induced by a morphism of schemes $i: \st{U}\times_{\iota,\scat{F},\mathbf b} B
 \rightarrow B$, which is an open immersion. 
\end{defi}

\begin{prop}\label{appB76}
For each open subset $\st{U}$ of a Deligne-Mumford stack $\st{F}$  there exists a unique subcategory $\st{R}$, which is again a Deligne-Mumford stack over $\catc$, and which satisfies property $(i)$ of definition \ref{appB75}, such that \\
$(i)$ the stack $\st{R}$ is reduced;\\
$(ii)$ the natural morphism $\rho: \st{R}\rightarrow \st{F}$ is representable by closed immersions;\\
$(iii)$ for each \'etale atlas $\mathbf a: A^\bullet \rightarrow \st{F}$ holds that the scheme representing $\st{R}\times_{\rho,\scat{F},\mathbf a}A^\bullet$ is the complement of $ \st{U}\times_{\iota,\scat{F},\mathbf b} A$ in $A$. 
\end{prop}

\proof
See \cite[p. 102]{DM} and \cite[3.9]{LM}. The idea of the construction is to fix an atlas  $\mathbf a: A^\bullet \rightarrow \st{F}$, and then use property $(iii)$ to construct a full subcategory. 
\ebew

\begin{defi}\em
The substack $\st{R}$ of $\st{F}$ of proposition \ref{appB76} is called the {\em reduced complement} of $\st{U}$ in $\st{F}$. If one chooses $\st{U} = \emptyset^\bullet$, then the reduced complement is called the {\em underlying reduced Deligne-Mumford stack} $\st{F}_{red}$.
\end{defi}

\begin{defi}\em\footnote{\cite[p. 102]{DM}}
$(i)$ A full subcategory $\st{A}$ of a Deligne-Mumford stack $\st{F}$, which is a Deligne-Mumford stack itself and satisfies conditions $(i)$ and $(ii)$ of proposition \ref{appB76}, is called a {\em closed subset} of $\st{F}$.\\
$(ii)$ A Deligne-Mumford stack $\st{F}$ is called {\em irreducible}, if $\st{F}$ is not empty, and if for all open subsets $\st{F}_1$ and $\st{F}_2$ of $\st{F}$ , which are both nonempty, holds that their intersection 
\[ \st{F}_1\cap\st{F}_2 := \st{F}_1\times_{\scat{F}}\st{F}_2 \]
is nonempty.
\end{defi}

\begin{rk}\em
A Deligne-Mumford stack $\st{F}$ is irreducible if and only if for any two closed subsets  $\st{F}_1$ and $\st{F}_2$, such that $\st{F}=\st{F}_1 \cup \st{F}_2$,  follows that $\{\st{F}_1,\st{F}_2\} = \{ \emptyset^\bullet,\st{F}\}$. 
\end{rk}

\begin{prop}\label{appB80}
Let $\st{F}$ be a Noetherian Deligne-Mumford stack. Then there exists a unique family $\{\st{F}_i\}_{i\in I}$, up to permutation,  of closed (and hence reduced) subsets of $\st{F}$, such that\\
$(i)$  for all $i\in I$ the stack $\st{F}_i$ is irreducible;\\
$(ii)$ for all $i,j\in I$ the inclusion $\st{F}_i\subset \st{F}_j$ implies $i=j$;\\
$(iii)$ for any open and quasi-compact subset $\st{U}\subset \st{F}$ there exist only finitely many $i\in I$ such that $\st{U}\cap\st{F}_i := \st{U}\times_{\scat{F}}\st{F}_i$ is nonempty; \\
$(iv)$ for the underlying reduced stack holds $\st{F}_{red} =\coprod_{i\in I}\st{F}_i$.
\end{prop}

\begin{defi}\em
The irreducible substacks $\st{F}_i$ of $\st{F}$ of proposition \ref{appB80} are called  the {\em irreducible components} of $\st{F}$. 
\end{defi}

\begin{prop}
Let $\st{F}$ be a Deligne-Mumford stack, which is normal and Noetherian. Then the irreducible components of $\st{F}$ are exactly the connected components of $\st{F}$.
\end{prop} 

\proof See \cite[Prop. 4.13]{LM}.
\ebew

\begin{defi}\em\footnote{\cite[Def. 1.6]{Vi}}
A Deligne-Mumford stack $\catf$ is called {\em integral}, if it is both reduced and irreducible.
\end{defi}

\begin{defi}\em
Let $\st{F}$ and $\st{G}$ be integral Deligne-Mumford stacks, and let $\mathbf r: \st{F}\rightarrow \st{G}$ be a separated and dominant morphism of finite type between them. Let $ \mathbf a: A^\bullet \rightarrow \st{F}$ and  $ \mathbf b: B^\bullet \rightarrow \st{G}$ be  atlases, where $A$ and  $B$ are integral schemes. If $\mathbf r$ is representable, then we define the {\em degree} of $\mathbf r$ by
\[ \deg(\mathbf r) := \deg(\overline{ r}) ,\]
where $\overline{r} : \st{F}\times_{\mathbf r, \st{G},\mathbf b} B \rightarrow B $ denotes the induced morphism of schemes in the Cartesian diagram. If $\mathbf r$ is not representable, then   we define
\[ \deg(\mathbf r) := \frac{\deg(\mathbf r\circ\mathbf a)}{\deg(\mathbf a)} .\]
\end{defi}

\begin{rk}\em
Note that the degree of a morphism of stacks need not be an integral number. As usual, the degree is multiplicative with respect to the composition of morphisms, i.e. if $\mathbf s: \st{G} \rightarrow \st{H}$ is a second morphism, then
\[ \deg(\mathbf s\circ\mathbf r) = \deg(\mathbf r) \cdot \deg(\mathbf s) \] 
holds.
\end{rk}

\begin{defi}\em
Let $\st{F}$ be a Deligne-Mumford stack.  We define the {\em inertia group stack} $\st{I}_{\scat{F}}$ of $\st{F}$ as the following fibred category over $\catc$.\footnote{\cite[1.12]{Vi}} For a scheme $B \in \ob(\catc)$, put
\[ \st{I}_{\scat{F}}(B) := \{ (\mathbf b,\beta) : \:\: \mathbf  b\in\st{F}(B), \beta\in\mor_{\scat{F}(B)}(\mathbf  b,\mathbf b) \} .\]
For a pair of objects $(\mathbf b_1,\beta_1)\in\st{I}_{\scat{F}}(B_1)$ and $(\mathbf b_2,\beta_2)\in\st{I}_{\scat{F}}(B_2)$, morphisms are given by
\[ \mor_{\scat{I}_{\sscat{F}}}((\mathbf b_1,\beta_1),(\mathbf b_2,\beta_2)) := \{ \phi \in \mor_\scat{C}(\mathbf b_1,\mathbf b_2) : \:\: \phi\circ \beta_1 = \beta_2\circ \phi \} .\] 
The projection to $\catc$ is the obvious one. 
\end{defi}

\begin{rk}\em
Note that morphisms $\beta\in\mor_{\scat{F}(B)}(\mathbf  b,\mathbf b)$ are automorphisms by the definition of a fibred category. There is a natural morphism
\[ \mathbf p : \quad \st{I}_{\scat{F}} \rightarrow \st{F} \]
by projection. It follows directly from the construction of $\st{I}_{\scat{F}}$ and of the definition of the fibre  product that there is  an isomorphism of stacks
\[ \st{I}_{\scat{F}} \cong \st{F}\times_{\Delta,\scat{F}\times_\sscat{C}\scat{F},\Delta}\st{F} ,\]
where $\Delta: \st{F}\rightarrow \st{F}\times_\scat{C}\st{F}$ denotes the diagonal morphism. 

In fact, $\st{I}_{\scat{F}}$ is a Deligne-Mumford stack, and the morphism $\mathbf p$ is representable, separated, quasi-finite and unramified. If there exists a scheme $M$ and a separated morphism of finite type from $\st{F}$ to $M^\bullet$, then the morphism $\mathbf p: \st{I}_{\scat{F}} \rightarrow \st{F}$ is finite. 

For an integral stack $\st{F}$ one defines  the number
\[ \delta(\st{F}) := \deg(\mathbf p) .\]
This number equals the order of the group of automorphisms of a general geometric point of $\st{F}$. One has $\delta(\st{F})=1$ if and only if $\st{F}$ contains  an open subset isomorphic to a scheme. For a full account on the inertia group stack see \cite{Vi} and \cite{Gi}.
\end{rk}

\begin{rk}\em
Let $\catf$ be a Deligne-Mumford stack. For a scheme $X\in\ob(\catc)$ the {\em set of connected components of} $\catf(X)$ is defined as the set
\[ \pi_0(\catf(X)) := \catf(X) / \sim \]
where by definition the equivalence $\mathbf x_1 \sim \mathbf x_2$ holds for two elements $\mathbf x_1,\mathbf x_2\in\catf(X)$ if and only if there exists a morphism $\phi: \mathbf x_1 \rightarrow \mathbf x_2$. Note that $\phi$ is necessarily an isomorphism. If $\catf = Y^\bullet$ for some scheme $Y\in\ob(\catc)$, then $\pi_0(Y^\bullet(X)) = Y^\bullet(X)$. A morphism $\mathbf r : \catf\rightarrow \catg$ of Deligne-Mumford stacks induces a map of the sets of connected components
\[ \mathbf r_\ast : \:\: \pi_0(\catf(X)) \rightarrow \pi_0(\catg(X)) \]
for any scheme $X\in\ob(\catc)$.
\end{rk}

\begin{defi}\em
Let $\st{F}$ be a Deligne-Mumford stack. A {\em moduli space} for $\st{F}$ is a scheme $M \in \ob(\catc)$, together with   a proper morphism 
\[ \mathbf p : \:\: \st{F} \rightarrow M^\bullet ,\]
such that for all algebraically closed fields $k$, with $\Omega:= \Spec(k)\in\ob(\catc)$ the induced map
\[ p_\ast : \:\: \pi_0(\st{F}(\Omega)) \rightarrow M^\bullet(\Omega) \]
is a bijection  of the sets of connected components of the fibre categories.\footnote{\cite[Def. 3.10]{Gi}, \cite[Def 2.1]{Vi}} 
\end{defi}

\begin{defi}\em\label{appB84}
Let $\st{F}$ be a Deligne-Mumford stack.
A  {\em coarse moduli space} for $\catf$ is a scheme $M \in \ob(\catc)$, together with   a morphism 
$ \mathbf p :  \st{F} \rightarrow M^\bullet$, such that for all algebraically closed fields $k$, with $\Omega:= \Spec(k)\in\ob(\catc)$, the induced map
$ p_\ast : \pi_0(\st{F}(\Omega)) \rightarrow M^\bullet(\Omega) $ 
is a bijection, and satisfying the following universal property. For all schemes $N\in\ob(\catc)$, and all morphisms $\mathbf q: \st{F}\rightarrow N^\bullet $,  there exists a morphism $u: M\rightarrow N$ in $\catc$, such that the diagram
\[ \diagram
&\st{F} \dlto_{\mathbf p}\drto^{\mathbf q}\\
M^\bullet \rrto_u&& N^\bullet
\enddiagram\]
commutes.\footnote{\cite[Def. 4.1]{Ed}} 
\end{defi}

\begin{rk}\em
Let  $M$ be a moduli space of a Deligne-Mumford stack $\st{F}$, and let $f: N\rightarrow M$ be a morphism in $\catc$. Then $N$ is a moduli space for the stack $\st{F}\times_{\mathbf p, M^\bullet,\mathbf f} N^\bullet$. If $M$ is a coarse moduli space for $\st{F}$, then $N$ need not be a coarse moduli space for $\st{F}\times_{\mathbf p, M^\bullet,\mathbf f} N^\bullet$. If $M$ is a moduli space for $\catf$, satisfying in addition the  universal property of definition \ref{appB84}, then $N$ need not necessarily satisfy the universal property again. 

A coarse moduli space is uniquely determined by the universal property up to isomorphism. Two moduli spaces for a Deligne-Mumford stack $\catf$ may differ by a universal homeomorphism.\footnote{\cite[Rk. 4.1]{Ed}, \cite[Def. 2.1]{Vi}}
\end{rk}

\begin{ex}\em
Consider a (strict) functor $\mathbf F : \catc^{op}\rightarrow \Set$. Let $\st{F}$ be the associated stack. Since for any scheme $B\in\ob(\catc)$ the fibre category $\st{F}(B) = \mathbf F(B)$ is a set, the set of connected components equals
$\pi_0(\st{F}(B)) = \mathbf F(S)$. Thus a coarse moduli space for $\st{F}$ is also a coarse moduli space for $\mathbf F$ in the sense of \cite[Def. 5.6]{GIT}. 
\end{ex}

\begin{rk}\em
Let  $M$ be a moduli space of a Deligne-Mumford stack $\st{F}$. Then $\catf$ is irreducible if and only if $M$ is irreducible, and $\catf$ is connected if and only if $M$ is connected.\footnote{\cite[Lemma 2.3]{Vi}}
\end{rk}

\begin{prop}
$(i)$ Let $\catf$ be a Deligne-Mumford stack, together with a separated morphism $\mathbf q: \catf \rightarrow N^\bullet$ of finite type, for some scheme $N\in\ob(\catc)$. Then there exists a scheme $T\in\ob(\catc)$ and a morphism
\[ \mathbf t : \:\: T^\bullet \rightarrow \catf \]
which is proper and surjective.\\
$(ii)$  Let $\catf$ be a Deligne-Mumford stack, together with a moduli space  $\mathbf p: \catf \rightarrow M^\bullet$. Then there exists a scheme $T\in\ob(\catc)$ and a morphism
\[ \mathbf t : \:\: T^\bullet \rightarrow \catf \]
which is finite and surjective.
\end{prop}

\proof
For part $(i)$ of the proposition see \cite[Thm. 4.14]{DM}, and for part $(ii)$ see \cite[Prop. 2.6]{Vi}.
\ebew

\vfill\pagebreak
%
%

\section{Quotients}
Throughout this section let $\catc$ be one of the categories $(\Sch/S)$ or $(\Af/S)$, where $S$ is a scheme, or an affine scheme, respectively.

\subsection{Group actions}

\begin{defi}\em\label{Q11}
\footnote{\cite[Def. 0.1]{GIT}}
A {\em group scheme over $S$} is an object $p: G \rightarrow S \in \ob(\catc)$, together with morphisms
\[ \begin{array}{cccc}
\mu: & G \times_S G & \rightarrow & G,\\
\iota:& G & \rightarrow &G, \\
e : & S & \rightarrow & G ,
\end{array} \]
such that
\begin{enumerate}
\item[$(i)$] $\mu\circ(\id_G\times\mu) = \mu\circ(\mu\times\id_G) : G \times_S G \times_S G \rightarrow G $, \hfill ({\em associativity})

\item[$(ii)$] $\mu\circ(\id_G \times \iota) \circ \Delta = e \circ p : G \rightarrow G,$ \hfill ({\em inverse})\\
$\mu\circ(\iota \times \id_G) \circ \Delta = e \circ p : G \rightarrow G,$

\item[$(iii)$] $\mu(e\times \id_G) = \id_G : G \cong S \times_S G \rightarrow G$ \hfill ({\em identity})\\
$\mu(\id_G \times e) = \id_G : G \cong G \times_S S \rightarrow G$.

\end{enumerate}
Here, in $(ii)$ the morphism $\Delta: G \rightarrow G\times_S G$ is the diagonal morphism, and in $(iii)$ we identify 
$S \times_S G \cong G \cong G \times_S S$. 
\end{defi}

\begin{ex}\em
The {\em general linear group} $\GL_n(k)$ over some field $k$ can be viewed as a group scheme, even as an affine group variety, by
\[ \GL_n(k) = \Spec \left( k[x_{ij}][\det(x_{ij})^{-1}]\right), \]
for $1 \le i,j\le n$. Replacing $k$ by an arbitrary ring $R$ we obtain a group variety $\GL_n(R)$. 

Let $Z$ be a noetherian scheme. By definition, $Z$ can be covered by finitely many affine open subsets $\Spec( A_i)$, where all $A_i$ are Noetherian rings. Each of them defines an affine group variety $\GL_n(A_i)$, and the glueing of the $\Spec(A_i)$ induces a glueing of the group varieties, which gives a group scheme
\[ \GL_n(Z) = \bigcup_{i} \GL_n(A_i) .\]
This group scheme is also denoted by $\Aut_{\oka_Z}(\oka_Z^n)$.\footnote{\cite[Exp. I.4.5]{SGA3}}
\end{ex}

\begin{defi}\em
\footnote{\cite[Def. 0.3]{GIT}}
Let $X \in \ob(\catc)$ be a scheme, and let $G$ be a group scheme over
$S$. We say that {\em $G$ acts on $X$ from the right} if there is a morphism
\[ \sigma : \:\: X \times_S G \rightarrow X \]
in $\catc$, such that
\[ \sigma(\id_X \times \mu) = \sigma (\sigma \times \id_G) : \:\: X \times_S G  \times_S G \rightarrow X ,\]
as well as 
\[ \sigma(\id_X \times e) = \pr_1 : \:\: X \times_S S \rightarrow X .\]
\end{defi}

\begin{rk}\em
Suppose that $S = \Spec(k)$ for some field $k$, and  that $G$ acts on a scheme $X\in\ob(\catc)$ via $\sigma : X \times_S G \rightarrow X$. For any closed point $g\in G$ define
\[ \function{\sigma_g}{X}{X}{x}{\sigma(x,g).}\]
Then $\sigma_g$ is an automorphism of $X$ over $S$, and the map
\[ \function{\sigma_\bullet}{G}{\Aut(X)}{g}{\sigma_g}\]
is a homomorphism of groups.\footnote{\cite[Section 6.1]{LeP}}
\end{rk}

\begin{defi}\em
$(i)$ An action $\sigma : X \times_S G \rightarrow X$ is called {\em trivial}, if $\sigma = \pr_1$. \\
$(ii)$ An action $\sigma : X \times_S G \rightarrow X$ is called {\em free}, if \[ \Psi := (\pr_1, \sigma) : \:\: X \times_S G \rightarrow X \times_S X \]
is a closed immersion. The action $\sigma$ is called {\em proper}, if the morphism $\Psi$ is proper.\footnote{\cite[Def. 0.8]{GIT})}
\end{defi}

\begin{rk}\em
Assume that $S = \Spec(k)$ for some field $k$. Then the morphism $e: S \rightarrow G$ from definition \ref{Q11} determines a closed point $\id \in G$. 
If $G$ acts freely on $X$, then  the identity
\[ \sigma(x,g) = x, \]
for all closed points $x \in X$, implies $ g = \id$. 
\end{rk}

\begin{defi}\em
Let $G$ be a group scheme acting on $X \in \ob(\catc)$ with action $\sigma_X$. Let $Y \in \ob(\catc)$ be a scheme.\\
$(i)$ A morphism $f : X \rightarrow Y$ in $\catc$ is called {\em $G$-invariant}, if the diagram
\[ \diagram 
X \times_S G \rrto^{\pr_1} \dto_{\sigma_X} && X \dto^{f} \\
X \rrto_{f} & &Y 
\enddiagram
\]
commutes.\\
$(ii)$ Suppose that there is a $G$-action $\sigma_Y$ on $Y$. A morphism $f : X\rightarrow Y$ is called {\em $G$-equivariant}, if the diagram
\[ \diagram
X \times_S G \rrto^{f\times\id_G} \dto_{\sigma_X} && 
 Y \times_S G \dto^{\sigma_Y}\\
X \rrto_{f} && Y 
\enddiagram \]
commutes.\footnote{\cite[Def. 6.1.4]{LeP}}
\end{defi}

\begin{defi}\em
\footnote{\cite[Def. 6.1.1]{LeP}}
Let $S = \Spec (k)$, where $k$ is a field. An {\em (affine) linear algebraic group over $k$} is a group scheme $G$ over $S$, which is a closed subgroup scheme of $\GL_n(k)$ for some $n \in {\Bbb N}$. 
\end{defi}

\begin{rk}\em
\footnote{\cite[Remark 2.11]{EHKV}}
\label{Q18}
More generally, let $S$ be a Noetherian scheme. 
A {\em linear algebraic group over $S$} is a group scheme $G$ over $S$, which is a subgroup scheme of $\GL_n(S)$ for some $n \in {\Bbb N}$. 

Consider the special case $S = \Spec( k)$. Every affine group variety over a field $k$ is a subgroup variety of some $\GL_n(k)$, see \cite[Prop. I.1.10]{Borel}. In fact, any affine group scheme of finite type over $k$ is a subgroup scheme of $\GL_n(k)$. 
\end{rk}

\begin{defi}\em
\footnote{\cite[Def. 0.4]{GIT}} 
Let $G$ be a group scheme acting on a scheme $X\in\ob(\catc)$. Let $f: Z\rightarrow X$ be a morphism of schemes in $\catc$. The {\em stabilizer of} $f$ is the subgroup scheme $\Stab_G(f)$ of $Z\times_S G$ over $Z$, defined by the Cartesian diagram
\[ \diagram 
\Stab_G(f) \rrto \dto&& Z \dto^{(\id_Z,f)} \\
Z\times_S G \rrto_{(\pr_1,\sigma\circ(f,\id_G))} && Z\times_S X.
\enddiagram \]
If $Z$ is an irreducible subscheme of $X$, and $z_0$ is a generic point of $Z$, then the preimage of $z_0$ in $\Stab_G(f)$ can be interpreted as a subgroup of $G$, which is called the {\em stabilizer subgroup} $\Stab_G(Z)$ of $Z$.
\end{defi}

\subsection{Quotients of schemes}

\begin{defi}\em
\footnote{\cite[Def. 0.5]{GIT}}
Let $G$ be a group scheme over $S$ acting with action $\sigma$ on $X \in \ob(\catc)$. A {\em categorical quotient of $X$ by $G$} is a scheme $Y \in \ob(\catc)$ together with a morphism $q : X \rightarrow Y$, such that the diagram
\[ \diagram
X \times_S G \rto^{\sigma} \dto_{\pr_1} & X \dto^{q}\\
X \rto_{q} & Y 
\enddiagram \]
commutes, satisfying the following universal property. For all schemes $Z \in \ob(\catc)$ with morphisms $r : X \rightarrow Z$, such that $r \circ \sigma = r \circ \pr_1$, there is a unique morphism $u : Y \rightarrow Z$ in $\catc$, such that the following diagram commutes.
\[ \diagram
X \times_S G \rto^{\sigma} \dto_{\pr_1} & X \dto^{q} \xto[ddr]^{r} \\
X \rto_q \xto[drr]_{r}& Y \drto_{u} \\
&& Z .
\enddiagram \]
\end{defi}

\begin{rk}\em
If  a categorical quotient exists, then it is uniquely determined up to isomorphism by its universal property.
\end{rk}

The following definition is taken from \cite[2.10]{Vi}. Note that this notion of a quotient is weaker than that of a good or a geometric quotient. 

\begin{defi}\em
\footnote{\cite[2.10]{Vi}}
Let $G$ be a group scheme over $S$ acting with action $\sigma$ on $X \in \ob(\catc)$. A {\em quotient of $X$ by $G$} is a scheme $Y \in \ob(\catc)$ together with a $G$-invariant morphism $q: X \rightarrow Y$, satisfying the following  properties.
\begin{enumerate}
\item[$(i)$] The morphism $q$ is surjective
 and universally submersive.   
\item[$(ii)$] Consider the morphism $\Psi := (\pr_1,\sigma) : X \times_S G \rightarrow X \times_S X$. Then for its image holds
\[ \im(\Psi) \cong X \times_Y X .\]
\end{enumerate}
\end{defi}

\begin{rk}\em
A morphism $q: X \rightarrow Y$ is called {\em universally submersive}, if for every morphism $\overline{q}: X' := Y'\times_{f,Y,q}X \rightarrow Y'$, which is induced by a morphism $f: Y' \rightarrow Y$ in $\catc$, holds that a subset $U' \subset Y'$ is open if and only if $\overline{q}^{-1}(U') \subset X'$ is open. In particular, for $f=\id_Y$, it follows that $q$ itself is submersive, i.e. the topology on $Y$ is the topology induced via $q$ from $X$. 

Condition $(ii)$ of the definition says that the geometric fibres of $q$ are the orbits of the geometric points of $X$ with respect to the action of $G$. 

Note that a quotient in this sense need not be unique, so it is in general not a categorical quotient. 
\end{rk}

\begin{defi}\em
\footnote{\cite[Def. 0.6]{GIT}}
A morphism $q: X \rightarrow Y$ is called a {\em geometric quotient of $X$ by $G$}, if it is a quotient of $X$ by $G$, and if the natural morphism of sheaves
\[\oka_Y \:\: \rightarrow  \:\: (q_\ast\oka_X)^G \]
is an isomorphism. Here $(q_\ast\oka_X)^G$ denotes the subsheaf of $ q_\ast\oka_X$ of $G$-invariant sections. 
\end{defi}

\begin{defi}\em 
\footnote{\cite[Def. 0.7]{GIT}}
Let $G$ and $X$ be as above, and let $q: X \rightarrow Y$ be a categorical/geometric quotient. \\
$(i)$ It is called a {\em universal categorical/geometric quotient} if for all morphisms $f: Y'\rightarrow Y$ the morphism induced by base change $\overline{q} : X \times_Y Y'\rightarrow  Y'$ is a categorical/geometric quotient of $X \times_Y Y'$ by $G$. \\
$(ii)$ It is called a {\em uniform categorical/geometric quotient} if for all flat morphisms $f: Y'\rightarrow Y$ the morphism induced by base change $\overline{q} : X \times_Y Y' \rightarrow Y'$ is a categorical/geometric quotient of $X \times_Y Y'$ by $G$.
\end{defi}

\begin{rk}\em
A geometric quotient is also always  a categorical quotient, and hence in particular unique up to isomorphism. This is true in the universal case as well.\footnote{\cite[Prop. 0.1]{GIT}}
\end{rk}

\subsection{Principal bundles}

From now on we will assume that $S = \Spec(k)$ for some algebraically closed field $k$. In this subsection we follow the presentation in Sorger's paper \cite{So}.

\begin{abschnitt}\em\footnote{\cite[sect. 2]{So}}
Let $G$ be a linear algebraic group over a field $k$, and let $B \in \ob(\catc)$ be a scheme. A {\em $G$-fibration over $B$} is a scheme $E \in \ob(\catc)$, on which $G$ acts, together with a $G$-invariant morphism 
\[ p : \:\: E \rightarrow B .\]
A {\em morphism of $G$-fibrations} $p : E \rightarrow B$ and $p':E'\rightarrow B$ is a $G$-equivariant morphism
\[ \phi :\:\:  E \rightarrow E',\]
in $\catc$  such that $ p = p'\circ \phi$.
In other words, a $G$-fibration is a commutative diagram
\[ \diagram
 E \times_k G \rto^{\sigma} \dto_{\pr_1} & E \dto^{p} \\
E \rto_{p} & B 
\enddiagram \]
 in $\catc$, where $\sigma$ denotes the group action of $G$ on $E$.  
A $G$-fibration $p : E \rightarrow B$ is called {\em trivial}, if it is up isomorphism given by the diagram
 \[ \diagram
 B \times_k G \times_k G \rrto^{(\id_B, \mu)} \dto_{(\pr_1,\pr_2)} && B \times_k G \dto^{\pr_1} \\
B \times_k G  \rrto_{\pr_1} && B .
\enddiagram \]
\end{abschnitt}

\begin{defi}\em
Let $G$ and $B$ be as above. A {\em principal $G$-bundle over $B$} is a $G$-fibration over $B$, which is locally trivial with respect to the \'etale topology, i.e. if there exists an \'etale covering ${\cal U} = \{ f_\alpha: B_\alpha  \rightarrow B\}_{\alpha \in A}$ of $B$, such that for all $\alpha \in A$ 
\[ f_\alpha^\ast E : = E \times_B B_\alpha \]
is trivial over $B_\alpha$. A {\em morphism of principal $G$-bundles over $B$} is just a morphism of the underlying $G$-fibrations. 
\end{defi} 

\begin{ex}\em\label{E6118}
Let $p: V \rightarrow B$ be a $k$-vector bundle of rank $r$ over $B$. Then the {\em  associated frame bundle}
\[ F := \ISOM_{\oka_B}(\oka_B^r,V) \]
is a principal $\GL_r(k)$-bundle over $B$.

Note that the pullback $F \times_B V$ of $V$ to $F$ is a trivial vector bundle over $F$, but it is in general not trivial as a $\GL_r(k)$-bundle over $V$.
Indeed, consider the Cartesian diagram
\[ \diagram
F\times_B V \rto \dto& V\dto \\
F \rto & B.
\enddiagram
\]
There exists a section $s : F \rightarrow F\times_B V$, which is defined on the second factor by the natural map
\[ \begin{array}{ccc}
F = \ISOM_{\oka_B}(\oka_B^r,V) & \rightarrow &V \\
\sigma &\mapsto& \sigma(0),
\end{array}
\]
and by the identity on the first factor.
\end{ex}

\begin{ex}\em
As in example \ref{E6118}, let $p:V\rightarrow B$ be a $k$-vector bundle of rank $r$ over $B$. Let $\PP(V)$ denote the projective bundle over $B$ defined by $V$. Then the {\em $\PGL_k(r)$-bundle associated to} $\PP(V)$ 
\[ P := \ISOM_{\oka_B}(\PP(\oka_B^r),\PP(V))\]
is a principal $\PGL_k(r)$-bundle over $B$.  
\end{ex} 

\begin{rk}\em
Let $p : V\rightarrow B$ be a $k$-vector bundle of rank $r$ over some scheme $B$, and let $q: P\rightarrow B$ be the principal  bundle associated to $\PP(V)$. Then the pullback bundle $q^\ast \PP( V) \rightarrow P$ has a natural trivialization as a projective bundle over $P$. If $r: E\rightarrow B$ is a principal $\PGL_k(r)$-bundle over $B$, such that $r^\ast \PP(V)\rightarrow E$  is trivial, then $E$ is in a natural way isomorphic to $P$.  
\end{rk} 

\begin{rk}\em
A group $G$ is called {\em special}, or {\em locally isotrivial}, if local triviality with respect to the \'etale topology also implies local triviality with respect to the Zariski topology. 
An example of a special group is $G = \GL_n(k)$. Furthermore, for linear algebraic groups over a field $k$, all locally trivial principal bundles with respect to the \'etale topology over smooth curves are locally trivial with respect to the Zariski topology  as well.\footnote{ 
See  \cite[Remark 2.1.2]{So}, and his references to \cite{Serre}.}
\end{rk}

\begin{rk}\em 
Any morphism 
$\phi : E \rightarrow E'$ of principal $G$-bundles
 over $B$ is an isomorphism. Indeed, there exists an
 \'etale covering ${\cal U} = \{ f_\alpha: B_\alpha  \rightarrow
 B\}_{\alpha \in A}$ of $B$, such that for all $\alpha \in A$ the
 liftings $f^\ast_\alpha E$ and $f^\ast_\alpha E'$ are trivial,  and
 $f^\ast_\alpha \phi$ is given by multiplication with some morphism $g_\alpha : B_\alpha \rightarrow  G$. 
So over $B_\alpha$ we can define $(f^\ast_\alpha\phi)^{-1}$ by $g_\alpha^{-1}$, and the local definitions glue together to give a global inverse $\phi^{-1}$. 
\end{rk}

\begin{rk}\em
Let $p : E \rightarrow B$ be a principal $G$-bundle. If there exists a section $s : B \rightarrow E$, then $E$ is trivial. Indeed, an isomorphism is given by
\[ \function{\alpha}{B\times_k G}{E}{(b,g)}{s(b)g.} \]
\end{rk}

\begin{rk}\em
Let $p : E \rightarrow B$ be a principal $G$-bundle, and $f : B'\rightarrow B$ be a morphism in $\catc$. Then the Cartesian diagram
\[ \diagram
E'\rto^{\overline{f}} \dto_{p'} & E \dto^{p} \\
B'\rto_{f}& B 
\enddiagram \]
defines a principal $G$-bundle $p': E'\rightarrow B'$ up to  isomorphism. If $E$ is trivial, then so is $E'$.  
\end{rk}

\begin{rk}\em\label{appC27}
Let $G$ be a linear algebraic group, which acts freely on a scheme $X$, which is of finite type over $k$. Let $q : X\rightarrow Y$ be the geometric quotient morphism with respect to this action. Then $X$ is a principal $G$-bundle over $Y$, see \cite[Prop. 0.9]{GIT}.
\end{rk}

\begin{rk}\em\label{Q122}
Let $p : E \rightarrow B$ be a principal $G$-bundle, and let
$X\in\ob(\catc)$ be a quasi-projective scheme, on which $G$ acts from
the left. There is a $G$-action on the product $E \times_k X$ defined by
\[ \sigma_E(\pr_1,\pr_3) \times \sigma_X(\iota\circ\pr_3,\pr_2) : \:\: 
E \times_k X \times_k G \longrightarrow E \times_k X,\]
or simply, for $e \in E$, $x\in X$ and $g \in G$ by
\[ (e,x) \cdot g  := (e \cdot g, g^{-1} \cdot x) .\]
Then the quotient 
\[ E(X) :=( E \times_k X) /G \]
exists as a scheme in $\catc$.\footnote{ since $X$ is quasi-projective, see \cite[2.2]{So}} We have a natural morphism
\[ E(X) \rightarrow B, \]
where each  fibre  is isomorphic to $X$, which is locally trivial with respect to the \'etale topology.
\end{rk}

\begin{rk}\em
\footnote{\cite[2.2]{So}}
\label{Q121}
Let $\rho : H \rightarrow G$ be a morphism of linear algebraic groups. \\
$(i)$ Let $p : E \rightarrow B$ be a principal $H$-bundle. Via $\rho$, there is an action of $H$ on $G$. So we can form   
\[  E(G) = (E \times_k G) /H,  \]
which is in fact a principal $ G$-bundle 
on $B$, compare remark \ref{Q122}.  It is called the {\em extension of
$E$ by $G$}. The action of $G$ on $E(G)$ is on the second factor and from
the right. If $\rho$ is a closed immersion then  $E$ is in a natural way an $H$-subbundle of
$E(G)$. \\
$(ii)$ An {\em $H$-reduction} of a principal $G$-bundle $q:F \rightarrow B$ is a principal $H$-bundle $q': F' \rightarrow B$, such that $F'(G) \cong F$.
\end{rk}

\begin{prop}  
\label{Q123}
\footnote{\cite[2.2.3]{So}}
Let $\rho : H \rightarrow G$ be a closed immersion, and let $q : F
\rightarrow B$ be a principal $G$-bundle on $B$. Then there exists an
$H$-reduction of $F$ if and only if there exists a section 
\[ s : \:\: B \rightarrow F / H := F ( G/H) .\]
\end{prop} 

\proof
Assume that $p : E\rightarrow B$ is a reduction of $F$, so that  we have 
\[ F = (E\times_k G)/H .\]
In particular, $F$ is a locally trivial $G$-fibration, where the transition functions are given by elements of $H$. Thus $F/H$ is in fact trivial over $B$, and a zero-section $s : B \rightarrow F/H$ exists, which maps a point $b \in B$ to the equivalence class of the unit of $G$ in the fibre over $b$.

Conversely, if there is such a section $s : B \rightarrow F/H$, then we can form the Cartesian square
\[ \diagram
E'\rto^{\overline{s}} \dto_{\overline{\pi}} & F \dto^{\pi}\\
B \rto^{s} & F/H.
\enddiagram \]
Here $\pi : F \rightarrow F/H$ is the natural quotient map, and hence a principal $H$-fibration, and $E'$ is the subbundle of $F$ which is fibrewise the kernel of $\pi$. This implies that $\overline{\pi}: E'\rightarrow B$ is a principal $H$-fibration as well.

To see that $E'$ is indeed a reduction of $F$, we will write down
isomorphisms between $F$ and $E'(G)$.  Using the section $s : B
\rightarrow F/H$, we identify $F/H$ with $B \times_k G/H$, via
\[ \begin{array}{ccc}
B \times_k G/H & \rightarrow & F/H \\
(b,\overline{g}) & \mapsto& s(b)g .
\end{array}\]
Here $\overline{g}$ denotes the equivalence class of an element  $g\in G $ in
$G/H$. Conversely, any element $\overline{f} \in F/H$, represented by
$f\in F$ with $q(f)= b$, maps to $ (q(f), \overline{s(q(f))^{-1}\cdot
f}) \in  B \times_k G/H$. Note that even though $s(q(f))^{-1}$ is not
defined, the ``distance'' $s(q(f))^{-1}\cdot f$ of $f$ from the
section is well defined in $G$. Under this identification, 
we have $s(b) = (b ,\overline{\id})$, for $b \in B$, and the quotient map becomes 
\[ \function{\pi}{F}{B\times_k G/H}{f}{(q(f),\tilde{\pi}(f)).} \]
We now define 
\[ \function{\alpha}{F}{E'(G)=(E'\times_k G)/H=((F\times_{B\times_k G/H} B )\times_k G)/H}{f}{(f\cdot\tilde{\pi}(f)^{-1},q(f),\tilde{\pi}(f)).} \]
Note that indeed $ (f\cdot\tilde{\pi}(f)^{-1},q(f)) \in F\times_{B\times_k G/H} B$. The inverse morphism is simply given by
\[ \function{\beta}{((F\times_{B\times_k G/H} B )\times_k G)/H}{F}{(f,b,\overline{g})}{f \cdot g.} \]
Note here that $\tilde{\pi}(f) = \overline{\id}$, so that  $\tilde{\pi}(f\cdot g) = g$.
\ebew

\begin{rk}\em\label{rkB141}
A stronger version of the above proposition \ref{Q123} is given in \cite[lemma 2.2.3]{So}. 
In fact, there is a one-to-one correspondence between sections $s: B
\rightarrow F/H$ and $H$-subbundles $E'$ of $F$ with $F= E'(G)$. Note
that giving a section  $s: B
\rightarrow F/H$ is equivalent to specifying a trivialization $F/H
\cong B \times_k G/H$ such that $s$ becomes the zero-section, i.e. the
section given by $s(b) = (b, \overline{\id})$. The $H$-subbundle
$E'$ is defined 
fibrewise as the kernel of $ \pi_b : F_b \rightarrow (F/H)_b$, for any
$b\in B$. Conversely, if $E'(G)=F$, then there is a natural inclusion
 of $E'$ as an $H$-subbundle of $F$. Hence the composed map
 \[ B = E'/H \hookrightarrow F'/H \] 
defines a section $s: B \rightarrow F/H$. There is a unique
trivialization of $F/H\cong B \times_k G/H $, such that 
$s$ 
becomes the zero-section.
\end{rk}

\begin{rk}\em
 Note that any two $H$-reductions $E_1'$ and $E_2'$ are isomorphic to each other. Indeed, the corresponding sections $s_1,s_2 : B \rightarrow F/H$ define two trivializations of $F/H \cong B\times_k G/H$, such that $E_1'$ and $E_2'$ are the preimages of the respective zero-sections, and hence isomorphic. 
\end{rk}

\begin{rk}\em
\label{Q125}
Let $\rho : H \rightarrow G$ be a morphism of linear algebraic
groups, and let $\phi : E_1\rightarrow E_2$ be a morphism of principal
$H$-bundles over $B$. Then $\phi$ extends naturally to a morphism 
\[ \begin{array}{cccc}
\overline{\phi}: & E_1(G) &\rightarrow & E_2(G) \\
& (e,g) & \mapsto & (\overline{\phi}(e),g)
\end{array} \]
of principal $G$-bundles over $B$.\\ 

Conversely, let $E_1$ and $E_2$ be
$H$-reductions of principal $G$-bundles $F_1$ and $F_2$,
respectively.  
A morphism $\tilde{\phi}: F_1 \rightarrow F_2$
restricts to a morphism 
\[ \phi' : \:\: E_1 \rightarrow E_2 \]
if and only if for the induced map 
$ \tilde{\phi}' : \:\: F_1/H \rightarrow F_2/H $
holds
$ s_2 = \tilde{\phi}' \circ s_1$, 
where $s_i: B \rightarrow F_i/H$ are the sections corresponding to $E_i$, for $i=1,2$. 
If $E_i \subset F_i$ are subbundles, then this is equivalent to the condition that $\tilde{\phi}(E_1) \subset E_2$, and hence that $\tilde{\phi}(E_1) = E_2$.

If $\phi'$ exists, then its extension 
is equal to  $\tilde{\phi}$.  
\end {rk}

\begin{lemma}\label{appC34}
Let $G$ and $H$ be linear algebraic groups, and let $E$ be a principal $G\times_k H$-bundle over $B$. Then the quotient
$ E/G$ exists
as a  principal $H$-bundle, and the quotient $E/H$  exists as a principal $H$-bundle over $B$.
\end{lemma}

\proof
Consider an \'etale covering ${\cal U} =\{f_\alpha: B_\alpha \rightarrow B\}_{\alpha \in A}$ of $B$, such that for all $\alpha \in A$ holds
\[ f^\ast_\alpha E \cong B_\alpha\times_k G\times_k H .\]
The transition functions on intersections $B_{\alpha\beta}:= B_\alpha \times_B  B_\beta$, for $\alpha,\beta \in A$, are given by pairs of morphisms
\[ (g_{\alpha\beta},h_{\alpha\beta}): \:\: B_{\alpha\beta} \rightarrow G\times_k H\]
such that for all $\alpha,\beta,\gamma \in A$ holds
$ (g_{\alpha\gamma},h_{\alpha\gamma})=(g_{\beta\gamma},h_{\beta\gamma}) \circ (g_{\alpha\beta},h_{\alpha\beta})$. 
In particular, the transition functions $(g_{\alpha\beta})_{\alpha,\beta\in A}$ and $(h_{\alpha\beta})_{\alpha,\beta\in A}$ define principal bundles $E/H$ and $E/G$. 
\ebew

\subsection{Quotient fibred categories}

\begin{defi}\em
Let $G$ be a linear algebraic group over $S$, acting on a scheme $X\in\ob(\catc)$. The {\em quotient lax functor  $[X/G]$} is defined as the  following $2$-functor
\[ [X/G] : \:\: \catc^{op}  \rightarrow \Gp_2. \]
For all $B \in \ob(\catc)$ the category $[X/G](B)$ is defined by
\[ \ob([X/G](B)) := \left\{ \begin{array}{ll}
(E,p,\phi):\\
p:E \rightarrow B \Text{ is a principal $G$-bundle,}\\
\phi: E \rightarrow X \Text{ is $G$-equivariant}
\end{array} \right\} .\]
For two objects $(E,p,\phi)$ and $(E',p',\phi')$ in $[X/G](B)$ one defines 
\[ \mor_{[X/G](B)}((E,p,\phi),(E',p',\phi')) := \left\{ \begin{array}{l}
\alpha : E' \rightarrow E: \\
\Text{$\alpha$ is a morphism of} \\
\Text{principal $G$-bundles over $B$,}\\
\Text{such that }  \phi' = \phi\circ \alpha 
\end{array} \right\}. \]
Necessarily a morphism of principal $G$-bundles is an isomorphism, so 
$[X/G](B)$ is a groupoid. 

A morphism $f: B'\rightarrow B$ in $\catc$ gets mapped to a functor 
\[ \funktor{[X/G](f)}{[X/G](B)}{[X/G](B')}{(E,p,\phi)}{(f^\ast E,\overline{p}, \phi\circ\overline{f})}{\alpha}{f^\ast\alpha,} \]
where $\overline{p}$, $\overline{f}$ and $f^\ast E$ are defined by the Cartesian diagram 
\[ \diagram
f^\ast E \rto^{\overline{f}} \dto_{\overline{p}} & E \dto^{p}\\
B'\rto_{f} & B.
\enddiagram \]
Note that $[X/G]$ is not a functor. For a pair of morphisms $f': B''\rightarrow B'$ and $f: B'\rightarrow B$ there is an isomorphism
\[ c_{f,f'}: \:\: {f'}^\ast(f^\ast E) \simto (f\circ f')^\ast E \]
rather than equality. However, these $c_{f,f'}$ provide the
2-morphisms of a $2$-functor, which is thus well defined. (Recall that
for each morphism $f: B' \rightarrow B$ an object $f^\ast E$ has been
fixed a priory. Modifying the lifting changes the 2-functor by an
isomorphism.)
\end{defi}

\begin{defi}\em
Let $G$ be a linear algebraic group over $S$, acting on a scheme $X\in\ob(\catc)$. The {\em quotient fibred category} is defined up to isomorphism of fibred categories as the category fibred in groupoids over $\catc$ which corresponds to the quotient lax functor  $[X/G]$ under the natural equivalence of $2$-categories of remark \ref{correspondence}. It is denoted by $[X/G]$, too.
\end{defi}

\begin{rk}\em
In particular, objects of the fibred category $[X/G]$ are objects $(E,p,\phi) \in [X/G](B)$ for some $B \in \ob(\catc)$. A morphism from $(E,p,\phi) \in [X/G](B)$ to $(E',p',\phi') \in [X/G](B')$ is given by  a pair of morphisms $f : B'\rightarrow B$ and $\tilde{f} : E' \rightarrow f^\ast E$ in $\catc$, where $\tilde{f}$ is a morphism of principal $G$-bundles, such that the diagram
\[ \diagram
& &  & & X \\
E' \rto_{\tilde{f}}\drto_{p'}\xto[urrrr]^{\phi'}&f^\ast E\xto[rrru]|{\overline{\phi}}\rto_{\overline{f}} \dto^{\overline{p}} & E
\dto^{p}\xto[rru]_{\phi}& & \\
&B'\rto_f & B 
\enddiagram \]
is commutative, and such that 
 the square is Cartesian. 
\end{rk}

\begin{rk}\em\label{Q131}
Let $[X/G]$ be a quotient fibred category  as above, and let $B \in \ob(\catc)$. By Yoneda's lemma \ref{YON} there is an equivalence of categories
\[ \mor_{\Lax(\scat{C}^{op},\Gp_2)}(B^\bullet,[X/G]) \equiv [X/G](B) .\]
An object of $\mor_{\Lax(\scat{C}^{op},\Gp_2)}(B^\bullet,[X/G])$ is given by a pair of families
\[ (\{\eta_Y\}_Y,\{u_f\}_f), \]
where for all $Y\in \ob(\catc)$
\[ \eta_Y : \:\: B^\bullet(Y)\rightarrow [X/G](Y) \]
is a functor between the fibre categories, and for all morphisms $f:Y\rightarrow Z$ in $\catc$ 
\[ u_f : \:\:\:  \eta_Y\circ B^\bullet(f) \Rightarrow [X/G](f)\circ \eta_Z \]
is a $2$-isomorphism, subject to the usual compatibility conditions.

If a morphism ${\mathbf b} : B^\bullet \rightarrow [X/G]$ is given, then the corresponding object of $[X/G](B)$ with respect to  Yoneda's equivalence is given by $\eta_B(\id_B) \in [X/G](B)$. Hence 
\[ \eta_B(\id_B) = (E,p,\phi) \]
for some principal $G$-bundle $p:E\rightarrow B$, together with a $G$-equivariant morphism $\phi:E \rightarrow X$.

Let $Y\in \ob(\catc)$, and $f: Y \rightarrow B \in B^\bullet(Y)$. The $2$-isomorphism $u_f$, applied to $\id_B \in B^\bullet(B)$, gives an isomorphism
\[ \eta_Y\circ B^\bullet(f)(\id_B) \cong [X/G](f)\circ\eta_B(\id_B), \]
hence
\[ \begin{array}{lll}
\eta_Y(f) &\cong& [X/G](f)(E,p,\phi)\\[2mm]
&\cong& (f^\ast E,\overline{p},\phi\circ\overline{f}), 
\end{array} \]
where $\overline{p}$ and $\overline{f}$ are given by the Cartesian diagram
\[ \diagram
f^\ast E \rto^{\overline{f}} \dto_{\overline{p}} & E \rto^\phi \dto^p& X\\
Y \rto_f & B.
\enddiagram\]
\end{rk}

\begin{lemma}
Let $[X/G]$ be a fibred category as above, and let $A,B \in \ob(\catc)$, together with morphisms $f: A\rightarrow B$ in $\catc$  and ${\mathbf b} : B^\bullet \rightarrow [X/G]$ in {\em $\Lax(\catc^{op},\Gp_2)$}. If ${\mathbf b}$ is represented under Yoneda's equivalence by a triple $(E,p,\phi)\in[X/G](B)$, then the composed map
\[ {\mathbf a} := {\mathbf b}\circ {\mathbf f} : \:\:  A^\bullet \rightarrow [X/G] \]
is represented, up to isomorphism,  by the triple $(f^\ast E,\overline{p},\phi\circ\overline{f})$, 
where $\overline{p}$ and $\overline{f}$ are given by the Cartesian diagram
\[ \diagram
f^\ast E \rto^{\overline{f}} \dto_{\overline{p}} & E \rto^\phi \dto^p& X\\
A \rto_f & B.
\enddiagram\]
\end{lemma}

\proof
Let the notation be as in remark \ref{Q131}, and let the morphism ${\mathbf a}: A^\bullet \rightarrow [X/G]$ be represented by a pair of families
\[ (\{\tilde{\eta}_Y\}_Y,\{\tilde{u}_f\}_f). \]
Note that for any $g : Y \rightarrow A$ in $\catc$ holds
\[ {\mathbf b}\circ{\mathbf f}(g) = {\mathbf b}(f\circ g),\]
hence
\[ \tilde{\eta}_Y(g) = \eta_Y(f\circ g) .\]
In particular, under Yoneda's equivalence, ${\mathbf a}$ is represented by
\[ \tilde{\eta}_A(\id_A) = \eta_A(f),\]
hence by the triple $(f^\ast E, \overline{p},\phi\circ\overline{f})$ as claimed. 
\ebew

\begin{rk}\em
Let $[X/G]$ be as above. The triple $(X\times G, \pr_1, \sigma)$ defines  the {\em canonical quotient morphism}
\[ \pi: \:\:  X^\bullet \rightarrow [X/G], \]
using Yoneda's lemma. 
\end{rk}

\begin{lemma}\label{Q001}
Let $G$ be a linear algebraic group acting on a scheme $X\in\ob(\catc)$, with action $\sigma: X\times G \rightarrow X$. 
 Then the natural diagram 
\[ \diagram
X^\bullet\times_\scat{C} G^\bullet \rto^\sigma \dto_{\bpr_1} & X^\bullet\dto^\pi\\
X^\bullet\rto_\pi&[X/G] 
\enddiagram\]
is commutative, up to a $2$-morphism $\eta: \pi\circ\bpr_1 \Rightarrow \pi\circ\sigma$. It is  even a Cartesian diagram, i.e. $ X^\bullet\times_\scat{C} G^\bullet \equiv X^\bullet\times_{\pi,[X/G],\pi}X^\bullet$. 
\end{lemma}

\proof
We will show the $2$-commutativity of the diagram first. 
Let $f : B\rightarrow X\times G$ be an object of $X^\bullet\times_\catc G^\bullet$. The triple $\pi\circ\sigma(f)\in[X/G](B)$ is determined by the pullback diagram
\[\diagram
f^\ast((X\times G)\times G) \rto \dto & (X\times G)\times G \dto\rto & X\times G \rto^\sigma \dto^{\pr_1}&X\\
B\rto_f & X\times G \rto_\sigma& X.
\enddiagram \]
Hence we find $\pi\circ\sigma(f)= (B\times G,\pr_1,\sigma(\sigma\circ f,\id_G))$. Similarly, we obtain $\pi\circ \bpr_1(f)= (B\times G,\pr_1,\sigma(\pr_1\circ f,\id_G))$. There is an isomorphism $\eta_f : \pi\circ \bpr_1(f) \rightarrow \pi\circ\sigma(f)$ given by
\[ \eta_f := (\id_B, \mu(\pr_2(f),\pr_2)) ,\]
where $\mu: G\times G\rightarrow G$ denotes the multiplication in $G$. 
Note the direction of this morphism according to the definition of morphisms in $[X/G]$. 
In particular, on $B\times G$ there is an identity
\[ \sigma(\pr_1\circ f, \id_G) = \sigma(\sigma\circ f,\id_G) \circ \eta_f .\]
Because  the morphism $\eta_f$ is $G$-equivariant, it is indeed an isomorphism of principal $G$-bundles. It is easy to see that the cocycle condition for morphisms in 
$X^\bullet\times_\scat{C} G^\bullet$ is also satisfied, and we have thus constructed a $2$-morphism $\eta :\pi\circ \bpr_1\Rightarrow \pi\circ\sigma.$

To see that the diagram is Cartesian, we will verify the universal property. Here it is sufficient to restrict ourselves to schemes $Y\in \ob(\catc)$. Consider a commutative diagram
\[ \diagram
Y^\bullet \rto^{\mathbf g}\dto_{\mathbf f} & X^\bullet \dto^\pi\\
X^\bullet \rto_\pi &[X/G],
\enddiagram\]
together with a $2$-morphism $\rho: \pi\circ {\mathbf f} \Rightarrow \pi\circ{\mathbf g}$. We need to show that there exists a morphism $u: Y \rightarrow X\times G$, such that ${\mathbf g}= \sigma\circ {\mathbf u}$ and ${\mathbf f} = \bpr_1\circ{\mathbf u}$, as well as $\rho = \eta\ast\id_{\mathbf u}$. 

The composed morphism $\pi\circ{\mathbf f}$ corresponds to the triple $(Y\times G, \pr_1,\sigma\circ(f,\id_G)) \in [X/G](Y)$, while the morphism $\pi\circ{\mathbf g}$ corresponds to $(Y\times G, \pr_1,\sigma\circ(g,\id_G))$. The $2$-morphism $\rho$ defines a morphism 
\[ \rho_{\id_Y} : \:\: \pi\circ {\mathbf f}(\id_Y) \rightarrow \pi\circ{\mathbf g}(\id_Y) ,\]
i.e. a morphism from $(Y\times G, \pr_1,\sigma\circ(f,\id_G))$ to $(Y\times G, \pr_1,\sigma\circ(g,\id_G))$ in $ [X/G]$. On $Y\times G$,  this morphism is given by the pair $(\pr_1,\mu(\gamma\circ\pr_1,\pr_2))$, for some  $\gamma: Y \rightarrow G$. In particular, we have the identity
\[ g = \sigma(f,\gamma) .\]
Now define the morphism ${u}: Y \rightarrow X\times G$ by
\[ u := (f,\gamma) .\]
One sees immediately that $\pr_1\circ u = f$ and $\sigma\circ u = g$, as well as $\rho = \eta\ast\id_{\mathbf u}$. 

Note that these three conditions determine $u$ uniquely. Analogously, one verifies the condition on morphisms of schemes as in proposition \ref{A0118}. Therefore the above diagram is indeed Cartesian. 
\ebew

By the definition of the canonical quotient morphism $\pi:  X^\bullet \rightarrow [X/G]$ via Yoneda's lemma, the above lemma \ref{Q001} is just a special case of the more general proposition below.

\begin{prop}\label{P6137}
Let $[X/G]$ be a fibred category as above, and let $B\in \ob(\catc)$. By Yoneda's lemma,  a morphism of fibred categories
$ {\mathbf b} :  B^\bullet \rightarrow [ X/G] $ 
corresponds to a triple $(E,p,\phi) \in [X/G](B)$, where $p:E\rightarrow B$ 
is a principal $G$-bundle over $B$, together with an $G$-equivariant morphism $\phi: E\rightarrow X$. Then the diagram
\[ \diagram
E^\bullet \rto^{\mathbf p} \dto_{\phi^\bullet} & B^\bullet \dto^{\mathbf b}\\
X^\bullet \rto_\pi &[X/G]
\enddiagram \]
is Cartesian.
\end{prop}

\proof
The claim of the proposition is a generalization of lemma \ref{Q001}, so we may abbreviate some of the arguments. 
First we want to prove that the diagram is indeed commutative. 

The composed morphisms  ${\mathbf b}\circ {\mathbf p}$ and $\pi \circ \phi^\bullet$ from $E^\bullet$ to $[X/G]$ correspond to triples 
$(p^\ast E, \overline{p},\phi\circ p')$ and $(\phi^\ast (X\times G), \overline{\pr}_1, \sigma\circ \overline{\phi})$ in $[X/G](E)$, as defined by the following Cartesian diagrams
\[ \diagram
{{p^\ast} E} \rto^{{p'}} \dto_{\overline{p}} & E \rto^\phi \dto^p & X &\Text{and} & {\phi^\ast}(X\times G) \rto^{\overline{\phi}} \dto_{{\overline{\pr}_1}} & X\times G \rto^\sigma\dto^{{\pr_1}}&X \\
E \rto_p & B, &&& E \rto_\phi&X.
\enddiagram \]
Here $\sigma$ denotes the group action of $G$ on $X$. Clearly $p^\ast E = E \times_{p,B,p} E \rightarrow E$ admits a trivializing section, so
\[ p^\ast E \cong E\times G \cong \phi^\ast (X\times G) .\]
Under these isomorphisms, $\overline{p}$ is just the projection $\pr_1$, and thus identified with $\overline{\pr}_1$, and $p'$ equals the group action $\sigma_E: E\times G$ of $G$ on $E$. Because $\phi$ is $G$-equivariant, one has $\phi\circ \sigma_E = \sigma \circ \overline{\phi}$. In other words, the two triples in $[X/G](B)$ coincide, up to a unique isomorphism. Therefore the diagram commutes, up to a $2$-morphism $\eta$. 

It is also a Cartesian diagram. To see this we will verify the universal property. Note that it is enough to do this for schemes $Y\in\ob(\catc)$ only. 
Consider a  commutative diagram 
\[ \diagram
Y^\bullet \rto^{\mathbf g} \dto_{\mathbf f} & B^\bullet \dto^{\mathbf b} \\
X^\bullet \rto_\pi & [X/G] 
\enddiagram \]
for some scheme $Y \in \ob(\catc)$. The composed morphisms $\mathbf b\circ\mathbf g$ and $\pi\circ \mathbf f$ correspond to two  triples $(g^\ast E, \overline{p}, \phi\circ \overline{{g}})$ and $(f^\ast(X\times G), \overline{\pr}_1,\sigma\circ\overline{f})$ in $[X/G](Y)$, defined by Cartesian diagrams
\[ \diagram
{{g^\ast} E} \rto^{\overline{{ g}}} \dto_{\overline{p}} & E \rto^\phi \dto^p & X &\Text{and} & {f^\ast}(X\times G) \rto^{\overline{f}} \dto_{{\overline{\pr}_1}} & X\times G \rto^\sigma\dto^{{\pr_1}}&X \\
Y \rto_g & B, &&& Y \rto_f&X.
\enddiagram \]
By assumption, there exists a $2$-morphism $\eta: \mathbf b\circ\mathbf g \Rightarrow \pi\circ \mathbf f$, so in particular there is an isomorphism $\eta_{\id_Y}$ between the two triples. 
Since  $f^\ast(X\times G)$ is trivial over $Y$, and isomorphic to $g^\ast E$, there is a section $s : Y \rightarrow g^\ast E$. Consider the composed morphism $u := \overline{g}\circ s: Y\rightarrow E$  inducing  ${\mathbf u} : Y^\bullet \rightarrow E^\bullet$.  It clearly holds 
\[ {\mathbf p}\circ {\mathbf u} = {\mathbf g}, \quad \Text{and } \quad \phi\circ {u} = {f} . \]
One easily verifies $\rho=\eta\ast \id_{\mathbf u}$, 
and with   these three  conditions, ${\mathbf u}$ is uniquely determined.
\ebew

\begin{defi}\em
Let $G$ and $H$ be groups acting on a scheme $X\in\ob(\catc)$. Let $\sigma_G: X\times G \rightarrow X$ and  $\sigma_H: X\times H \rightarrow X$ denote the respective actions. The actions of $G$ and $H$ are called {\em independent}, if
\[ \sigma_H(\sigma_G(\pr_1,\pr_2),\pr_3) = \sigma_G(\sigma_H(\pr_1,\pr_3),\pr_2) \]
as morphisms from $X\times G \times H$ to $X$. 
\end{defi}

\begin{rk}\em
Let $G$ and $H$ be linear algebraic groups acting independently  on a scheme $X \in \ob(\catc)$. Then there exists a canonical morphism of fibred categories
\[ \mathbf q : \:\: [X/G] \rightarrow [X/G\times H] ,\]
sending a triple $(E,p,\phi) \in [X/G](B)$, for some $B \in \ob(\catc)$, to the triple $(E\times H, p\circ \pr_1, \sigma_H \circ (\phi,\id_H))\in [X/G\times H](B)$, where $\sigma_H : X\times H \rightarrow X$ denotes the group action of $H$ on $X$. The definition of $\mathbf q$ on morphisms in $[X/G]$ is the natural one. 
\end{rk}

Not much has been written on quotient stacks with respect to products of groups. The following generalization of lemma \ref{Q001} seems to be new. 

\begin{prop}\label{appC45}
Let $G$ and $H$ be linear algebraic groups acting independently on $X \in \ob(\catc)$. Then the diagram
\[ \diagram
X^\bullet\times_\scat{C} H^\bullet \rto^{\bpr_1} \dto_{\pi'} & X^\bullet \dto^{\pi}\\
[X/G] \rto_{\mathbf q} &[X/G\times H]
\enddiagram \]
is  commutative, with a $2$-morphism $\eta: \mathbf q\circ\pi' \Rightarrow \pi\circ\bpr_1$. It is even a Cartesian diagram.  Here,   ${\mathbf q}$ denotes the canonical morphism. The morphism $\pi'$ corresponds under Yoneda's equivalence to the triple $(X\times G \times H,(\pr_1,\pr_3),\sigma_{G\times H})\in [X/G](X\times H)$, where $\sigma_{G\times H} : X\times G \times H \rightarrow X$ is given by the group action of $G\times H$ on $X$. 
\end{prop} 

\begin{rk}\em
$(i)$
Note that in the special case $G = \{ \id\}$, the morphism $\pi'$ equals the group action $\sigma_H : X\times H \rightarrow X$ on $X$. \\
$(ii)$ If one replaces the upper arrow $\bpr_1$ by the group action $\sigma_H$, then the diagram is strictly commutative. 
\end{rk}

\proofof{the proposition}
Let us show the commutativity of the diagram first. Consider an object $f \in X^\bullet\times_\scat{C} H^\bullet(B)$ for some scheme   $B\in\ob(\catc)$, i.e. a morphism of schemes $f=(b,\beta) : B\rightarrow X\times H$. By definition, the images of $f$ under the composed morphism are the triples
\[ \pi\circ \bpr_1(f) =( B\times G\times H,\pr_1, \sigma_{G\times H}(b\circ\pr_1,\pr_2,\pr_3)), \]
and 
\[ \mathbf q\circ \pi'(f)  =( B\times G\times H,\pr_1, \sigma_H(\sigma_{G\times H}(b\circ\pr_1,\pr_2,\beta\circ\pr_1),\pr_3)). \]
We define an isomorphism of triples
\[ \eta_f : \:\: \mathbf q\circ \pi'(f)  \rightarrow \pi\circ\bpr_1(f) \]
through the isomorphism of principal $G\times H$-bundles 
\[ n_f := (\pr_1,\pr_2,\mu_H(\beta^\vee\circ\pr_1,\pr_3)) : \:\:  B\times G\times H \longrightarrow  B\times G\times H,\]
where $\mu_H: H\times H \rightarrow H$ denotes the multiplication in $H$, and $\beta^\vee := \iota\circ \beta$ denotes the composition of $\beta$ with the group inversion $\iota: H \rightarrow H$ of $H$. One verifies that 
\[ \sigma_H(\sigma_{G\times H}(b\circ \pr_1,\pr_2,\beta\circ\pr_1),\pr_3) \circ n_f = \sigma_{G\times H}(b\circ\pr_1,\pr_2,\pr_3),\]
so $n_f$ defines indeed a morphism of triples. Note that the construction of $\eta_f$ is compatible with compositions of morphisms in $X^\bullet\times_\scat{C} H^\bullet$, so there is indeed a $2$-morphism
\[\eta : \:\: \mathbf q\circ \pi' \Rightarrow \pi\circ\bpr_1 .\]

We now want to show that the diagram is Cartesian, i.e. that there is an isomorphism of stacks
\[ (X\times H)^\bullet \equiv [X/G]\times_{\mathbf q,[X/G\times H],\pi} X^\bullet.\]
By the universal property of the fibre product, there exists a morphism $\mathbf u$ from $(X\times H)^\bullet $ to $ [X/G]\times_{\mathbf q,[X/G\times H],\pi} X^\bullet$, as constructed in the proof of proposition \ref{A0117}. We want to define an inverse morphism $\mathbf v$. Consider an object  $\Omega=((E,p,\phi),b,\alpha) \in[X/G]\times_{\mathbf q,[X/G\times H],\pi} X^\bullet(B)$ for some scheme $B\in\ob(\catc)$. Here, $p:E\rightarrow B$ is a principal $G$-bundle over $B$, $\phi:E\rightarrow X$ is a  $G$-equivariant morphism, and $b: B\rightarrow X$ is a morphism in $\catc$, together with an isomorphism 
\[\alpha : \:\: \mathbf q(E,p,\phi) \rightarrow \pi(b).\]
Since the triples can be computed as
\[  \mathbf q(E,p,\phi) = (E\times H, p\circ\pr_1,\sigma_H(\phi\circ\pr_1,\pr_2)) ,\]
and
\[ \pi(b)= (B\times G\times H,\pr_1,\sigma_{G\times H}(b\circ\pr_1,\pr_2,\pr_3)),\]
the isomorphism $\alpha$ is given by a morphism of principal $G\times H$-bundles
\[ \lambda=(\lambda_E,\lambda_H) : \:\: B\times G\times H \rightarrow E\times H ,\]
satisfying the equality
\[ \sigma_{G\times H}(b\circ\pr_1,\pr_2,\pr_3) = \sigma_H(\phi\circ\pr_1,\pr_2)\circ \lambda .\]
Let $\nu: B\rightarrow B\times G\times H$ denote the zero-section. Define a morphism
$ \beta:  B \rightarrow H$ 
by
\[ \beta := \iota\circ\lambda_H\circ\nu,\]
where $\iota$ is the inversion morphism on $H$, 
and  define
\[ \mathbf v(\Omega) := (b,\beta)\in(X\times H)^\bullet(B).\]
Let another object $\Omega'=((E',p',\phi'),b',\alpha') \in[X/G]\times_{\mathbf q,[X/G\times H],\pi} X^\bullet(B')$ for some scheme $B'\in\catc$ be given, together with a morphism $\omega: \Omega\rightarrow \Omega'$. By definition, $\omega$ is given by a pair $(\omega_1,\omega_2)$, where 
\[\omega_1 : \:\: (E,p,\phi)\rightarrow (E',p',\phi')\]
is a morphism in $[X/G]$, and 
\[\omega_2 : \:\: b \rightarrow b'\]
is a morphism in $X^\bullet$. In particular, $\omega_2$ is given by a mophism $w: B'\rightarrow B$ in $\catc$, such that $b' = b\circ w$, i.e. $\omega_2 = w^\ast$. We define
\[\mathbf v(\omega) := w^\ast ,\]
considered as a morphism in $(X\times H)^\bullet$. It is clear from the construction that $\mathbf v$ is a functor from $[X/G]\times_{\mathbf q,[X/G\times H],\pi} X^\bullet$ to $(X\times H)^\bullet$, respecting the fibrations over $\catc$. 

For an object $f=(b,\beta):B\rightarrow X\times H$ in $(X\times H)^\bullet$ one computes
\[\begin{array}{rcl}
 \mathbf v(\mathbf u(f)) &=& \mathbf v(\pi'(f),\bpr_1(f),\eta_f)\\
&=& \mathbf v((B\times G,\pr_1,\sigma_{G\times H}(b\circ\pr_1,\pr_2,\beta\circ\pr_1)),b,\eta_f),
\end{array}\]
where $\eta_f$ is given by the automorphism $(\pr_1,\pr_2,\mu_H(\beta^\vee\circ\pr_1,\pr_3))$ on $B\times G\times H$. Therefore 
\[\begin{array}{rcl}
 \mathbf v(\mathbf u(f)) &=& 
 (b,\iota\circ\mu_H(\beta^\vee\circ\pr_1,\pr_3)   \circ\nu)\\
&=& (b,\beta) \\
&=& f. 
\end{array} \]
For a morphism $w: B' \rightarrow B$ in $\catc$, i.e. a morphism $w^\ast$ in $X^\bullet\times_\scat{C} H^\bullet$,  one finds
\[ \mathbf v(\mathbf u(w^\ast)) = \mathbf v(\pi'(w^\ast),w^\ast) = w^\ast .\]
Hence for the composition of functors holds
$ \mathbf v\circ\mathbf u = \id_{( X\times H)^\bullet}$. 

Conversely, look at the composition $\mathbf u\circ\mathbf v$. For a tuple $\Omega=((E,p,\phi),b,\alpha)$ as above, with $(b,\beta):= (b,\iota\circ\lambda_H\circ\nu)= \mathbf v(\Omega)$, one computes
\[ \mathbf u(\mathbf v(\Omega)) = ((B\times G,\pr_1,\sigma_{G\times H}(b\circ\pr_1,\pr_2,\beta\circ\pr_1)),b,\eta_{(b,\beta)}).\]
To construct a $2$-morphism 
\[ \gamma: \:\: \id_{[X/G]\times_{\mathbf q,[X/G\times H],\pi} X^\bullet} \Rightarrow \mathbf u \circ\mathbf v\]
we need to define an isomorphism for each $\Omega$ an isomorphism of tuples
\[ \gamma_\Omega: \:\: \Omega \rightarrow \mathbf u(\mathbf v(\Omega)),\]
i.e. a pair of morphisms $\gamma_\Omega=(\gamma_1,\gamma_2)$ in $[X/G]\times_\scat{C} X^\bullet$, such that 
\[ \eta_{(b,\beta)} \circ\mathbf q(\gamma_1) = \pi(\gamma_2)\circ\alpha .\]
Clearly, $\gamma_2 = \id_b$. To construct $ \gamma_1: (E,p,\phi)\rightarrow (B\times G,\pr_1, \sigma_{G\times H}(b\circ\pr_1,\pr_2,\beta\circ\pr_1))$, recall that by assumption there is an isomorphism 
\[ \lambda=(\lambda_E,\lambda_H) : \:\: B\times G\times H \rightarrow E\times H \]
of principal $G\times H$-bundles over $B$. Using the respective zero-sections, there are natural subbundles
\[ i_{B\times G}: \:\: B\times G \hookrightarrow B\times G\times H \quad \text{ and } \quad i_E: E \hookrightarrow E\times H.\]
It follows from the definition of the action of $H$ on $E\times H$, and the fact that $\lambda$ is $G$-equivariant, that the composed morphism 
\[ \lambda_H\circ i_{B\times G} : \:\: B\times G \rightarrow H \]
does not depend on $G$, i.e. 
\[ \lambda_H \circ i_{B\times G} = \lambda_H \circ i_{B\times G}\circ\nu\circ\pr_1,\]
where $\nu: B\rightarrow B\times G$ is the zero-section. Let $\sigma_H': (E\times H)\times H \rightarrow E\times H$ denote the action of $H$ on $E\times H$. The image of the morphism
\[ \delta := \sigma_H'(\lambda,\iota\circ\lambda_H)\circ i_{B\times H} : \:\: B\times G \rightarrow E\times H \]
is contained in $E\times\{\id\}$, i.e. $\delta:  B\times G \rightarrow E$ is a morphism of principal $G$-bundles. 

To show that  the morphism $\delta$   defines a morphism of triples
\[ \gamma_1 : \:\: (E,p,\phi) \rightarrow (B\times G,\pr_1,\sigma_{G\times H}(b\circ \pr_1,\pr_2,\beta\circ\pr_1)),\]
one needs to verify the identity 
\[ \phi\circ \delta = \sigma_{G\times H}(b\circ\pr_1,\pr_2,\beta\circ\pr_1) . \]
This follows easily from the fact that $\lambda$ satisfies the condition $ \sigma_{G\times H}(b\circ\pr_1,\pr_2,\pr_3) = \sigma_H(\phi\circ\pr_1,\pr_2)\circ\lambda$. 

Next we need to verify that the pair $\gamma_\Omega=(\gamma_1,\gamma_2)$ defines a morphism of tuples 
\[ \Omega \:  \Rightarrow \: \mathbf u(\mathbf v(\Omega))  ,\]
i.e. we need to verify that the equality $\eta_{(b,\beta)} \circ\mathbf q(\gamma_1)= \pi(\gamma_2)\circ\alpha$ holds, considered as morphisms from $\mathbf q(E,p,\phi)$ to $\pi(b)$. By construction, the morphism of principal $G\times H$-bundles defining $\pi(\gamma_2)\circ\alpha$ is just $\lambda$, and the morphism defining $\eta_{(b,\beta)} \circ\mathbf q(\gamma_1)$ is 
\[ (\delta,\id_H)\circ(\pr_1,\pr_2,\mu_H(\beta^\vee\circ\pr_1,\pr_3)) : \:\:  B\times G\times H\rightarrow E\times H. \]
It is straightforward to show that these two morphisms are the same.

Finally, one checks that the construction of the morphisms $\gamma_\Omega$ is compatible with morphisms in $[X/G]\times_{\mathbf q,[X/G\times H],\pi} X^\bullet$, so that they define indeed a natural transformation 
\[ \gamma: \:\: \id_{[X/G]\times_{\mathbf q,[X/G\times H],\pi} X^\bullet} \Rightarrow \mathbf u \circ\mathbf v\]
as claimed. 
\ebew

The above proposition \ref{appC45} is a special case of our much more general result below. 

\begin{prop}\label{appC45a}
Let $G$ and $H$ be linear algebraic groups acting independently on $X\in\ob(\catc)$. Let $\mathbf y : Y^\bullet \rightarrow [X/G\times H]$ be given by the triple $(T,t,\tau)$. Then the diagram
\[\diagram
(T/G)^\bullet \rto^{\mathbf t'}\dto_{\mathbf y'} & Y^\bullet\dto^{\mathbf y}\\
[X/G]\rto_{\mathbf q} & [X/G\times H]
\enddiagram\]
is Cartesian. Here $t':T/G\rightarrow Y$ denotes the morphism induced by the projection $t: T\rightarrow Y$, and $\mathbf y'$ is given by the triple $(T,\pi_H,\tau)\in[X/G](T/G)$, where $\pi_H: T\rightarrow T/G$ denotes the canonical quotient morphism.
\end{prop}

\proof
By lemma \ref{appC34}, there exist principal bundles $t_G: \: T_G := T/H \rightarrow Y$ and $t_H: \: T_H:= T/G \rightarrow Y$. We denote the canonical quotient morphisms by $\pi_G: T\rightarrow T_G$ and $\pi_H : T\rightarrow T_H$, respectively. To show $2$-commutativity of the above diagram it is enough to show that the triples representing $\mathbf y\circ \mathbf t'$ and $\mathbf q\circ\mathbf y'$ are isomorphic. It is easy to see that the diagram
\[\diagram
T\times H \rrto^{\sigma_H'}\dto_{\pi_H\circ\pr_1}&& T\dto^t\\
T_H =T/G\rrto_{t_H} &&Y 
\enddiagram \]
is Cartesian. Here $\sigma_H'$ denotes the group action of $H$ on $T$. Using this, one finds that the composed morphism $\mathbf y\circ \mathbf t'$ is represented by the triple 
\[ (T\times H, \pi_H\circ\pr_1,\tau\circ\sigma_H'),\]
up to isomorphism. Note that in the computation of this triple we used the  isomorphism $t^\ast_H T \cong T\times H$. This isomorphism may change $\mathbf y\circ \mathbf t'$ at most by a $2$-morphism. For the composed morphism $\mathbf q\circ \mathbf y'$ one obtains the triple
\[ (T\times H, \pi_H\circ\pr_1,\sigma_H(\tau\circ\pr_1,\pr_2)) .\]
Clearly these two triples are identical, so the two morphisms $\mathbf y\circ \mathbf t'$ and $\mathbf q\circ \mathbf y'$ differ at most by a $2$-morphism, i.e. the diagram is commutative, with a $2$-morphism $\eta: \mathbf q\circ\mathbf y' \Rightarrow \mathbf y\circ\mathbf t'$. 

To show that the diagram is Cartesian, we will show that there is an isomorphism of stacks
\[ (T/G)^\bullet \equiv [X/G]\times_{\mathbf q,[X/G\times H],\mathbf y} Y^\bullet.\]
By the universal property of the fibre product, there exists a morphism $\mathbf u$ from $(T/G)^\bullet $ to $ [X/G]\times_{\mathbf q,[X/G\times H],\mathbf y} Y^\bullet$, as constructed in the proof of proposition \ref{A0117}. We want to define an inverse morphism $\mathbf v$. Consider an object  $\Omega=((E,p,\phi),b,\alpha) \in[X/G]\times_{\mathbf q,[X/G\times H],\mathbf y} Y^\bullet(B)$ for some scheme $B\in\ob(\catc)$. Here, $p:E\rightarrow B$ is a principal $G$-bundle over $B$, $\phi:E\rightarrow X$ is a  $G$-equivariant morphism, and $b: B\rightarrow Y$ is a morphism in $\catc$, together with an isomorphism 
\[\alpha : \:\: \mathbf q(E,p,\phi) \rightarrow \mathbf y(b).\]
The triples can be computed as
\[  \mathbf q(E,p,\phi) = (E\times H, p\circ\pr_1,\sigma_H(\phi\circ\pr_1,\pr_2)) ,\]
and
\[ \mathbf y(b)= ( b^\ast T, \overline{t},\tau\circ \overline{b}) ,\]
where $\overline{t}$ and $\overline{b}$ are defined  by the Cartesian diagram
\[\diagram
b^\ast T \rto^{\overline{b}}\dto_{\overline{t}} & T\dto^{t}\\
B\rto_{b} & Y. 
\enddiagram\]
The isomorphism $\alpha$ is thus given by a morphism of principal $G\times H$-bundles
\[ \lambda : \:\: b^\ast T \rightarrow E\times H ,\]
satisfying the equality
\[ \tau\circ\overline{b}  = \sigma_H(\phi\circ\pr_1,\pr_2)\circ \lambda .\]
In particular there is an induced trivialization $\overline{\lambda}:b^\ast T/G \rightarrow B\times H$. 
Let $\nu: B\rightarrow B\times H$ denote the zero-section. Define a morphism
$ f:  B \rightarrow T/G$ 
by
\[ f := \hat{b}\circ\overline{\lambda}^{-1}\circ\nu,\]
where $\hat{b}:b^\ast T/G \rightarrow T/G$ is induced by $\overline{b}: b^\ast T \rightarrow T$. Define now
\[ \mathbf v(\Omega) := f \in (T/G)^\bullet(B).\]
Let another object $\Omega'=((E',p',\phi'),b',\alpha') \in[X/G]\times_{\mathbf q,[X/G\times H],\mathbf y} Y^\bullet(B')$ for some scheme $B'\in\catc$ be given, together with a morphism $\omega: \Omega\rightarrow \Omega'$. By definition, $\omega$ is given by a pair $(\omega_1,\omega_2)$, where 
\[\omega_1 : \:\: (E,p,\phi)\rightarrow (E',p',\phi')\]
is a morphism in $[X/G]$, and 
\[\omega_2 : \:\: b \rightarrow b'\]
is a morphism in $Y^\bullet$. In particular, $\omega_2$ is given by a mophism $w: B'\rightarrow B$ in $\catc$, such that $b' = b\circ w$, i.e. $\omega_2 = w^\ast$. We define
\[\mathbf v(\omega) := w^\ast ,\]
considered as a morphism in $(T/G)^\bullet$. It is clear from the construction that $\mathbf v$ is a functor from $[X/G]\times_{\mathbf q,[X/G\times H],\mathbf y} Y^\bullet$ to $(T/G)^\bullet$, respecting the fibrations over $\catc$. 

For an object $f:B\rightarrow T/G$ in $(T/G)^\bullet$  one computes
\[\begin{array}{rcl}
 \mathbf v(\mathbf u(f)) &=& \mathbf v(\mathbf y'(f),\mathbf t'(f),\eta_f)\\
&=& \mathbf v( (f^\ast T,\overline{\pi}_H,\tau\circ\overline{f}), t'\circ f, \eta_f),
\end{array}\]
where $\overline{\pi}_H$ and $\overline{f}$ are defined by the Cartesian diagram
\[\diagram
f^\ast T \rto^{\overline{f}} \dto_{\overline{\pi}_H} & T\dto^{\pi_H}\\
B \rto_f& T/G .
\enddiagram\]
Therefore 
\[\begin{array}{rcl}
 \mathbf v(\mathbf u(f))
&=& f. 
\end{array} \]
For a morphism $w: B' \rightarrow B$ in $\catc$, i.e. a morphism $w^\ast$ in $(T/G)^\bullet$,  one finds
\[ \mathbf v(\mathbf u(w^\ast)) = \mathbf v(\mathbf y'(w^\ast),w^\ast) = w^\ast .\]
Hence for the composition of functors holds
$ \mathbf v\circ\mathbf u = \id_{( T/G)^\bullet}$. 

Conversely, look at the composition $\mathbf u\circ\mathbf v$. For a tuple $\Omega=((E,p,\phi),b,\alpha)$ as above, with $f= \mathbf v(\Omega)$, one computes
\[ \mathbf u(\mathbf v(\Omega)) = ((f^\ast T,\overline{\pi}_H,\tau\circ \overline{f}),t'\circ f,\eta_f) .\]
Consider the commutative diagram
\[\diagram
E\rto \dto_p& E\times H \dto \rto^{\lambda^{-1}}& b^\ast T \dto\rto^{\overline{b}}& T\dto^{\pi_H}\rto^{\tau}& X \\
B \rto_\nu&B\times H\rto_{\overline{\lambda}^{-1}}& b^\ast T/G\dto\rto_{\hat{b}}& T/G \dto^{t'}\\
&& B \rto_b& Y. 
\enddiagram \] 
Recall that by definition $f= \hat{b}\circ\overline{\lambda}^{-1}\circ \nu$. From the above diagram one can read off  that there exists an isomorphism from $(f^\ast T,\overline{\pi}_H,\tau\circ\overline{f})$ to $(E,p,\phi)$, as well as an equality $t'\circ f= b$. In fact, we have an isomorphism of tuples $\gamma_\Omega: \Omega\rightarrow \mathbf u(\mathbf v(\Omega))$.  
Using this, one can conclude  that there exists a $2$-morphism $\gamma: \id_{[X/G]\times_{\mathbf q,[X/G\times H],\mathbf y} Y^\bullet} \Rightarrow \mathbf u\circ\mathbf v$, proving  that $\mathbf u$ is an isomorphism of stacks. The arguments are analogous to that of the proof of proposition \ref{appC45}, so we shall omit them here. 
\ebew

\begin{rk}\em
\label{St01}
Let $G$ be a linear algebraic group over $S$, acting on a scheme $X\in\ob(\catc)$. Let $\pi: X \rightarrow Y$ be a  quotient of $X$ by $G$. Then there exists a canonical morphism of fibred categories
\[ [ X /G] \:\: \rightarrow \:\: Y^\bullet .\]
Indeed, let $B \in \ob(\catc)$, and let $(E,p,\phi) \in [X/G](B)$. Since $Y$ is a quotient, and $\phi: E \rightarrow X$ is $G$-equivariant, the morphism $\phi$  descends to a morphism 
\[ \phi': \:\:  E/G \rightarrow Y .\]
The projection map $p : E \rightarrow B$ induces a natural  isomorphism  $E/G \cong B$, so we may assume $\phi'\in \mor_\scat{C}(B,Y) = Y^\bullet(B)$. We define a functor
\[ \funktor{\Sigma_B}{[X/G](B)}{Y^\bullet(B)}{(E,p,\phi)}{\phi'}{(f,\tilde{f})}{\id_{Y^\bullet(B)}.} \]

This functor induces itself a morphism of fibred categories 
\[ \funktor{\Sigma}{[X/G]}{Y^\bullet}{(E,p,\phi)}{\phi'}{(f,\tilde{f})}{f : \phi'\rightarrow \phi.} \]
Note that in general this functor is neither full nor faithfull.
\end{rk}

\begin{ex}\em
\label{Q130}
Let $S = \Spec(k)$, and let  $G$ be a linear algebraic group over $k$,
acting freely on a scheme $X\in\ob(\catc)$. Let $\pi : X
\rightarrow Y$ denote the  geometric quotient of $X$ by $G$. 
Then there is an isomorphism of fibred categories 
\[ [X/G] \:\: \equiv \:\: Y^\bullet .\]
Indeed, let $f : B \rightarrow Y \in Y^\bullet(B)$. By remark \ref{appC27}, the quotient map $\pi : X \rightarrow Y$ is a principal $G$-bundle. We define $E$ by the Cartesian square
\[ \diagram 
E := f^\ast X \rrto^{\phi:= \overline{f}} \dto_{p} && X\dto^{\pi} \\
B \rrto_f && Y .
\enddiagram \]
Then $p :E \rightarrow B$ is again a principal $G$-bundle, and
$(E,p,\phi) \in [X/G](B)$. Note that this triple is uniquely determined by
the pair $(f,\pi)$, provided normalized liftings have been fixed a priory. The 2-functor thus defined is in fact inverse, up to a $2$-transformation, to the functor $\Sigma$ defined above in remark \ref{St01}.
\end{ex}

\begin{ex}\em\label{appC49}
Let $G$ and $X$ be as in the previous example, but assume now that $G$ acts trivially on $X$.

Let $(E,p,\phi) \in [X/G](B)$. Clearly, $E/G \cong B$, and since $\phi : E \rightarrow X$ is invariant under the action of $G$, there is an induced morphism $\overline{\phi} : B \rightarrow X$. We thus get the  functor 
\[ \funktor{\Sigma_B}{[X/G](B)}{X^\bullet(B)}{(E,p,\phi)}{\overline{\phi}}{\alpha}{\id_{X^\bullet(B)},}  \]
which again induces a morphism of fibred categories
\[ \Sigma: \:\: [X/G] \rightarrow X^\bullet .\]
If $G \neq \{\id\}$, then $\Sigma_B$ is never an equivalence of categories. Indeed, any morphism $\alpha: E \rightarrow E'$ in $[X/G](B)$ will be mapped to the same morphism in $X^\bullet(B)$ as the morphism $\alpha\circ\sigma_g$, where $\sigma_g:E\rightarrow E$ is the automorphism given by multiplication with a fixed element $g\in G$. Hence $\Sigma_B$ is not faithfull, if $G \neq \{\id\}$. In particular, $\Sigma$ is not  an isomorphism  of fibred categories if $G\neq\{\id\}$.  

Of course, there is always a functor
\[ \funktor{T_B}{X^\bullet(B)}{[X/G](B)}{f:B\rightarrow X}{B \times G \quad ( \cong f^\ast(X\times G)) }{h : B'\rightarrow B}{h\times \id_G: B'\times G \rightarrow B\times G .} \]
For similar reasons as above, $T_B$ is not a full functor, unless $G$ is trivial. Hence the induced morphim of fibred categories $T : X^\bullet \rightarrow [X/G]$ will  only give an isomorphism if $G$ is trivial.

Obviously, for the composition with the morphism $\Sigma$ holds the identity
\[ \Sigma\circ T \: = \: \id_{X^{\bullet}} .\]
\end{ex}

We have the following apparently new result for a trivially acting group $H$.  

\begin{prop}\label{appC50}
Let $G$ and $H$ be linear algebraic groups acting independently on $X\in \ob(\catc)$, such that  $H$ acts trivially on $X$. Then there exists a natural morphism of fibred categories
\[ \omega : \:\: [X/G\times H] \rightarrow [X/G] \]
such that the diagram
\[ \diagram
[X/H] \rto^\Sigma \dto_{\mathbf q} & X \dto^\pi\\
[X/G\times H] \rto_\omega &[X/G] 
\enddiagram \]
is strictly commutative and Cartesian. Here $\mathbf q$ denotes the natural morphism. 
\end{prop}

\proof
The morphism $\omega$ is defined as follows. For a triple $(E,p,\phi)\in[X/G\times H](B)$, for some scheme $B\in\ob(\catc)$, consider the quotient bundle $\overline{p} : E/ H \rightarrow B$, which exists, since $E$ is a principal $G\times H$-bundle over $B$. By the $G\times H$-equivariance of $\phi$, and the fact that $H$ acts trivially on $X$, there is an induced $G$-equivariant morphism $\overline{\phi}: E/H \rightarrow X$. Now we define
\[ \omega(E,p,\phi) := (E/H,\overline{p},\overline{\phi}). \]
The definition of $\omega$ on morphisms of triples is the obvious one. 

To show the commutativity of the above diagram, consider a triple $(E,p,\phi)\in[X/H](B)$. Because of the $H$-invariance of $\phi$, there is an induced morphism $\nu : B = E/H \rightarrow X$, such that $\phi= \nu\circ p$. One computes for the image of $(E,p,\phi)$ the triple
\[ \pi\circ \Sigma(E,p,\phi) = \pi(\nu)  = (B\times G,\pr_1,\sigma_G(\nu,\pr_2)) \in[X/G](B) ,\]
and finds the same triple for $\omega\circ \mathbf q (E,p,\phi)$. Thus the diagram commutes. 

By the universal property of the fibre product there exists a canonical morphism $\mathbf u : [X/H]\rightarrow [X/G\times H]\times_{\omega,[X/G],\pi} X^\bullet$. We want to define an inverse morphism $\mathbf v$. Consider an object $\Omega= ((E,p,\phi),b,\alpha)\in[X/G\times H]\times_{\omega,[X/G],\pi} X^\bullet(B)$ for some scheme $B\in\ob(\catc)$. Here $p: E\rightarrow B$ is a principal $G\times H$-bundle, $\phi: E\rightarrow X$ is a $G\times H$-equivariant morphism, $b: B\rightarrow X$ is a morphism in $\catc$, and $\alpha: \omega(E,p,\phi)\rightarrow \pi(b)$ is an isomorphism. In fact, $\alpha$ is a morphism from $(E/H,\overline{p},\overline{\phi})$ to $(B\times G,\pr_1,\sigma_G(b\circ\pr_1,\pr_2))$, and is thus given by a morphism of principal $G$-bundles
\[ \lambda : \:\: B\times G \rightarrow E/H ,\]
which satisfies the equality
\[ \sigma_G(b\circ\pr_1,\pr_2) = \overline{\phi}\circ\lambda .\]
The morphism $p$ induces a morphism $p': E/G \rightarrow B$, and the composition $b\circ p': E/G \rightarrow X$ is clearly $H$-invariant. We define
\[  \mathbf v(\Omega) := (E/G,p',b\circ p') \in [X/G](B). \]
Consider objects $\Omega_i = ((E_i,p_i,\phi_i),b_i,\alpha_i)\in[X/G\times H]\times_{\omega,[X/G],\pi} X^\bullet(B_i)$ for some scheme $B_i\in\ob(\catc)$, for $i=1,2$, together with a morphism $\omega: \Omega_1 \rightarrow \Omega_2$ between them. By definition, $\omega$ is given by a pair of morphisms $(\omega_1,\omega_2)$, where $\omega_1:(E_1,p_1,\phi_1)\rightarrow (E_2,p_2,\phi_2)$, and $\omega_2: b_1\rightarrow b_2$. In particular, $\omega_2$ is given by a morphism $w: B_2\rightarrow B_1$ in $\catc$, such that $\omega_2 = w^\ast$, and $\omega_1$ is defined by a morphism of principal $G\times H$-bundles $\rho: E_2\rightarrow E_1$, such that $\phi_2= \phi_1\circ\rho$. Note that under the fibre projection from $[X/G\times H]$ to $\catc$, the morphism $\omega_1$ maps to $w$. 

Since $\rho$ is $G$-equivariant, there is an induced morphism of principal $H$-bundles $\rho': E_2/G \rightarrow E_1/G$, satisfying $b_1\circ p_1' = b_2\circ p_2'\circ\rho'$. Thus $\rho'$ defines a morphism of triples
\[ r : \:\: \mathbf v(\Omega_1)= (E_1/G,p_1',b_1\circ p_1') \rightarrow (E_2/G,p_2',b_2\circ p_2')=\mathbf v(\Omega_2).\]
We define
\[ \mathbf v(\omega) := r ,\]
as a morphism in $[X/H]$. Clearly we have thus defined a functor $\mathbf v : 
[X/G\times H]\times_{\omega,[X/G],\pi} X^\bullet \rightarrow [X/H]$. 

For an object $(E,p,\phi)\in[X/H](B)$ one computes
\[\begin{array}{l}
\mathbf v(\mathbf u(E,p,\phi))  =  \mathbf v(\mathbf q(E,p,\phi),\Sigma(E,p,\phi), \id_{(B\times G,\pr_1,\sigma_G(\phi',\pr_2))})\\[1mm]
\hspace*{1cm} = \mathbf v((E\times G,p\circ\pr_1,\sigma_G(\phi\circ\pr_1,\pr_2)),\phi',\id_{(B\times G,\pr_1,\sigma_G(\phi'\circ\pr_1,\pr_2))} )\\[1mm]
\hspace*{1cm} = (E,p,\phi),
\end{array}\]
where $\phi' : B=E/H \rightarrow X$ denotes the morphism induced by the $H$-invariant morphism $\phi: E\rightarrow X$. 

For a morphism of triples $\tau: (E_1,p_1,\phi_1)\rightarrow (E_2,p_2,\phi_2)$, given by a morphism of principal $H$-bundles $t: E_2\rightarrow E_1$, one finds
\[\begin{array}{rcl}
\mathbf v(\mathbf u(\tau)) &=& \mathbf v(\mathbf q(\tau),\Sigma(\tau))\\
&=& \mathbf v (\tilde{\tau},t') \\
&=& \tau. 
\end{array}\]
Here $t': E_2/H= B_2\rightarrow B_1 =E_1/H$ is induced by $t$, or equivalently, $t'$ is the image of $\tau$ under the projection from $[X/H]$ to $\catc$. The morphism $\tilde{\tau}$ is the morphism of triples defined by 
\[ (t,\id_G) : \:\: E_2\times G \rightarrow E_1\times G.\]
Thus  we obtain the identity of functors
\[ \mathbf v\circ\mathbf u = \id_{[X/H]}.\]
Consider now an object $\Omega= ((E,p,\phi),b,\alpha)\in[X/G\times H]\times_{\omega,[X/G],\pi} X^\bullet(B)$ as  above. Then
\[\begin{array}{l}
\mathbf u (\mathbf v(\Omega)) =  \mathbf u(E/G,p',b\circ p')\\[1mm]
\hspace*{5mm} = ((E/G\times G ,p'\circ \pr_1,\sigma_G(b\circ p'\circ\pr_1,\pr_2)),b,\id_{(B\times G,\pr_1,\sigma_G(b\circ\pr_1,\pr_2))}) .
\end{array}\]
To construct a $2$-morphism 
\[ \gamma: \:\:  \mathbf u \circ \mathbf v \Rightarrow \id_{[X/G\times H]\times_{\omega,[X/G],\pi} X^\bullet(B)}, \]
we need to define for each $\Omega$ an isomorphism of tuples 
\[ \gamma_\Omega  : \:\: ((E/G\times G ,p'\circ \pr_1,\sigma_G(b\circ\pr_1,\pr_2)),b,\id_{(B\times G,\pr_1,\sigma_G(b\circ\pr_1,\pr_2))}) \rightarrow \Omega,\]
such that 
\[\alpha\circ\omega(\gamma_1) = \pi(\gamma_2)\circ \id_{(B\times G,\pr_1,\sigma_G(b\circ\pr_1,\pr_2))} .\]
Note that $\gamma_\Omega$ is given by a pair $(\gamma_1,\gamma_2)$. 
For the second coordinate we obviously have $\gamma_2 = \id_b$. To construct $\gamma_1$ we define a morphism of principal $G\times H$-bundles 
\[\theta: \:\: E\rightarrow E/G \times G \]
by
\[ \theta:= (q_G,\pr_2\circ\lambda^{-1}\circ q_H),\]
where $q_G : E\rightarrow E/G$ and $q_H: E\rightarrow E/H$ denote the canonical quotient morphisms. 

Note that $\theta$ is $G\times H$-equivariant, so it is indeed a morphism of principal $G\times H$-bundles over $B$. Using the identity $\overline{\phi}\circ\lambda= \sigma_G(b\circ\pr_1,\pr_2)$ from above, one can verify the equality
\[ \phi = \sigma_G(b\circ \pr_1,\pr_2)\circ\theta .\]
Thus $\theta$ defines a morphism of triples 
\[\gamma_1 : \:\:  (E/G\times G,p'\circ\pr_1,\sigma_G(b\circ p'\circ\pr_1,\pr_2)) \rightarrow   (E,p,\phi).\]
It remains to show that the pair $(\gamma_1,\gamma_2)$ defines a morphism of tuples, i.e. that the equality $\alpha\circ\omega(\gamma_1)= \pi(\gamma_2)\circ\id_{(B\times G,\pr_1,\sigma_G(b\circ\pr_1,\pr_2))} $ is satisfied, where both sides are  considered as morphisms from $\omega(E/G\times G,p'\circ\pr_1,\sigma_G(b\circ p'\circ\pr_1,\pr_2))= (B\times G,\pr_1,\sigma_G(b\circ\pr_1,\pr_2))$  to $\pi(b)=(B\times G,\pr_1,\sigma_{G}(b\circ\pr_1,\pr_2))$. The morphism of principal $G$-bundles representing $\alpha$ is by definition $\lambda$. Hence $\pi(\gamma_2)\circ  \id_{(B\times G,\pr_1,\sigma_G(b\circ\pr_1,\pr_2))} $ is given by the morphism of principal $G$-bundles 
\[ \id_{B\times G} : \:\: B\times G \rightarrow  B\times G ,\]
and similarly, the morphism $\alpha\circ\omega(\gamma_1)$ is given by
\[ \lambda^{-1}\circ\lambda: \:\: B\times G \rightarrow B\times G  .\]
Finally, one verifies that the definition of $\gamma_\Omega$ is compatible with morphisms in $[X/G]$, so that there exists indeed a $2$-morphism
\[ \gamma: \:\:  \mathbf u\circ\mathbf v \Rightarrow \id_{ [X/G\times H]\times_{\omega,[X/G],\pi} X^\bullet}.\]
This concludes the proof.
\ebew

\begin{rk}\em
Let the notation be as above, with $H$ acting trivially on $X$.  Consider compositions of the canonical   morphisms 
\[ \mathbf r : \:\: [X/G] \rightarrow [X/G\times H ] \]
and
\[ \omega : \:\: [X/G\times H] \rightarrow [X/G]. \]
By construction, one obtains
\[ \omega\circ \mathbf r = \id_{[X/G]} .\]
However, if $(E,p,\phi) \in[X/G\times H](B)$, for some scheme $B\in\ob(\catc)$, then $\mathbf q\circ \omega(E,p,\phi) = (E/H\times H,p',\phi')$ needs not be isomorphic to $(E,p,\phi)$. Hence there is in general no $2$-morphism from $\mathbf r\circ \omega$ to $\id_{[X/G\times H]}$.
\end{rk}

\begin{lemma}
Let $H$ be a linear algebraic group, which is a subgroup of a linear algebraic group $G$  over $k$, acting on some  scheme $X\in\ob(\catc)$. 
Suppose that a scheme 
\[ X':= (X \times G)/H \in\ob(\catc)\]
exists as a  quotient of $X\times G$ by $H$, where the action is defined by $(x,g) \cdot h := (xh,h^{-1}g)$, for
$x\in X$, $g \in G$ and $h \in H$. Then there is a canonical  isomorphism of
fibred categories
\[ [X/H] \:\: \equiv \:\: [X'/G] .\]
\end{lemma}

\proof
Let $(E,p,\phi)\in [X/H](B)$ for some $B \in \ob(\catc)$. Then there exists an extension bundle $p': E(G) = (E\times G)/H \rightarrow B$, see remark \ref{Q121}. We define a $G$-equivariant morphism
\[ \overline{\phi}: \:\: (E\times G)/H \rightarrow (X \times G) /H \]
in the obvious way, that is, as the morphism induced  by $(\phi,\id_G): E\times G \rightarrow X\times G$.  Put
\[ \mathbf v_B (E,p,\phi) := ((E\times G)/H,p',\overline{\phi}) .\] 
Similarly, it is possible to extend morphisms of
principal $H$-bundles, see remark \ref{Q125}.
So we obtain a functor 
\[ \mathbf v_B  :\:\: [X/H](B) \rightarrow [X'/G](B) \]
between the fibre categories over $B$, 
which in turn leads to a morphism of fibred categories
\[ \mathbf v : \:\: [X/H] \rightarrow [X'/G], \]
appropriately defined. 

To construct the inverse morphism, let $(E',p',\phi')\in [X'/G](B)$ be given, 
for some $B \in \ob(\catc)$. Consider the canonical inclusion map $\nu: X \rightarrow (X\times G)/H = X'$, and form the Cartesian square
\[ \diagram
E'' \rto^{\phi''} \dto_{{\nu}'} & X \dto^\nu\\
E'\rto_{\phi'}\dto_{p'} & X'\\
B.
\enddiagram \]
Let $E''$ be the restriction  of $E'$ to $X$. We claim that the composed map $p'':= p'\circ \nu': E'' \rightarrow B$
is a principal $H$-bundle, and an $H$-reduction of $E'$. 
Choose an  \'etale covering ${\cal U} = \{ f_\alpha: B_\alpha  \rightarrow B\}_{\alpha \in A}$ of $B$, such that $f_\alpha^\ast E' \cong B_\alpha \times G \  \rightarrow B_\alpha$ is a local trivialization of $E'$. 

Composition of the zero-section of $ B_\alpha \times G \rightarrow B_\alpha$ with $\phi'$ gives a morphism $\xi'_\alpha: B_\alpha \rightarrow X'$. By the $G$-equivariance of $\phi'$, the morphism $\phi'$ is locally given by 
\[ \function{\phi'|\,f_\alpha^\ast E'}{B_\alpha\times G}{(X\times G)/H}{(b,g)}{[\xi'_\alpha(b),g], } \]
up to an automorphism of $B_\alpha\times G$.  
Here the square brackets indicate equivalence classes modulo the action of $H$.  The inclusion $\nu'$ is given by the map $x \mapsto [x,1]$, hence $\phi'(g,b) \in X$ if and only if $g \in H$.  Consider the preimage
\[ \begin{array}{rcl}
f_\alpha^\ast E'' &=& \{ e\in f_\alpha^\ast E': \phi'(e) \in X\}\\
&\cong& \{ (b,g) \in B_\alpha\times G: \: g\in H \}\\
&=& B_\alpha\times H .
\end{array} \]
This shows  that $E''$ is locally trivial with respect to the \'etale topology, and, after verifying the appropriate cocycle condition, that $E''$ is  a principal $H$-bundle over $B$. Thus $(E'',p'',\phi'') \in [X/H](B)$. 

There exists a well-defined morphism of principal $G$-bundles 
\[ \function{\alpha}{(E''\times G)/H}{E'}{[e,g]}{\nu'(e) \cdot g,} \]
which of course is an isomorphism. So $E'$ is indeed an extension of
$E''$ by $G$. 

Note that under this construction morphisms of $G$-bundles over $X$ get as well  
restricted to morphisms of $H$-bundles over $X$. Hence this defines a functor, inverting
the morphism of fibred categories from above. 
\ebew

\subsection{Quotient stacks}

Recall that a category fibred in groupoids is a stack if and only if
its associated $2$-functor is a 2-sheaf. We claim that for quotient
fibred categories this is always the case. The key to this is the
fact, that both 2-sheaves and principal bundles are defined with
respect to the \'etale topology. 

\begin{prop}\footnote{\cite[Prop. 2.2]{Ed}} 
Let $G$ be a smooth linear algebraic group over $S$, acting on a scheme $X
\in \ob(\catc)$. Then the quotient fibred category $[X/G]$ is a stack.
\end{prop}

\proof
By proposition \ref{A0124} and remark \ref{AppR114}  we need to show that for all $B\in
\ob(\catc)$ and all sieves $\mathbf R\in J(B)$ the natural functor
\[ [X/G](B) \:\: \rightarrow \:\: \Prlim{\mathbf R}{[X/G]} \]
is an equivalence of categories.
An object of $ \Prlim{\mathbf R}{[X/G]}$ consists of data
$ \{(E_u,v_f)_{u,f}\},$ 
where $u: U \rightarrow B $ is an \'etale morphism contained in $\mathbf R(U)$, and $f : U'
\rightarrow U$ is a morphism such that $u \circ f = u' \in \mathbf R(U')$. $E_u$ is
a principal $G$-bundle over $U$, together with a $G$-equivariant
morphism $\phi_u$ into $X$ as in the diagram
\[ \diagram
E_u \rto^{\phi_u} \dto_{p_u}& X\\
U \dto_{u}\\
B,
\enddiagram \]
and $v_f$ is an isomorphism between the diagrams
\begin{center}
$\diagram 
E_{u'} \rto^{\phi_{u'}} \dto_{p_{u'}}& X\\
U' \dto_{u'}\\
B,
\enddiagram $ 
\quad and \quad  
$\diagram 
f^\ast E_{u} \rto  \dto& X\\
U' \dto \\
B.
\enddiagram $ 
\end{center}
Here both diagrams are considered as objects in $ \Prlim{\mathbf
R}{[X/G]}$. Clearly the data $\{(E_u,v_f)_{u,f} \}$ glue together to
define a principal $G$-bundle $E$ over $B$, together with a
$G$-equivariant morphism to $X$, i.e. an element 
\begin{center}
$
\diagram
E \rto^{\phi} \dto_{p}& X\\
B
\enddiagram$
$\in \ob([X/G](B))$.
\end{center}
Analogously, morphisms in $\Prlim{\mathbf
R}{[X/G]}$ define morphisms in $[X/G](B)$, and one checks that the
resulting functor is inverse, up to isomorphism, to the natural
functor above.
\ebew

\begin{rk}\em
Example \ref{Q130} implies that if $\pi: X \rightarrow Y$ is a
geometric quotient of $X$ under the action of a group $G$, which acts
freely on $X$, then the quotient stack $[X/G]$ is representable by the
scheme $Y$.
\end{rk}

\begin{rk}\em
There is the notion of the {\em dimension} of a locally Noetherian algebraic stack $\catf$. This can be found in full detail in \cite[chap. 11]{LM}. In general, $\dim \catf$ is an integral number, which is not necessarily non-negative. 

In the case where $\catf = [X/G]$ is a quotient stack with $X$ and $G$ as above, the dimension can be computed by
\[ \dim\, [X/G] = \dim X - \dim G ,\]
using the dimensions of $X$ and $G$ as schemes. 
\end{rk}

Before we proceed we note a simple technical lemma. 

\begin{lemma}\label{L167}
Let $[X/G]$ be a quotient stack, and let $\catf$ be a stack, together with a morphism $\bq: X^\bullet \rightarrow  {\catf}$. Suppose that there exists a $2$-morphism $\rho: \bq\circ  \pr_1 \Rightarrow \bq\circ \sigma$. Suppose that $\eta$ is compatible with the group action of $G$, in the sense that for all $\in X^\bullet(B)$ holds $\rho_{(x,e)}= \id_{\mathbf q \circ \bpr_1(x,e)}$, where $e: B\rightarrow G$ is the constant morphism of unity.
Then for all schemes $B$, and all $x\in X^\bullet(B)$ and all $g_1,g_2 \in G^\bullet(B)$ holds
\[ \rho_{(x,\mu(g_2,g_1))} = \rho_{(\sigma(x,g_2),g_1)} \circ \rho_{(x,g_2)} .\]
\end{lemma}

\proof
By definition of a $2$-morphism, there is a  commutative diagram 
\[ \diagram
\bq\circ \bpr_1 (x,g_2) \xto[rr]^{\rho_{(x,g_2)}} \dto_{\bq\circ\bpr_1(g_2,g_2^{-1})} && 
\bq\circ \sigma(x,g_2) = 
\bq\circ \bpr_1(\sigma(x,g_2),e) \dto_{\bq\circ\sigma(g_2,g_2^{-1})= \bq(e) = \id_{\bq\circ \bpr_1(\sigma(x,g_2),e)}} \\
\bq\circ \bpr_1(\sigma(x,g_2),e) \xto[rr]_{\id_{\bq\circ \bpr_1(\sigma(x,g_2),e)}} & &
\bq\circ \sigma(\sigma(x,g_2),e) = \bq\circ \bpr_1(\sigma(x,g_2),e) ,
\enddiagram \]
where $e: B\rightarrow G$ is the constant morphism of unity.  
From this follows the equality
\[ \rho_{(x,g_2)} = \bq\circ \bpr_1(g_2,g_2^{-1}) .\]
Now consider the commutative diagram
\[\diagram
\bq\circ\bpr_1(x,g_2) = \bq\circ\bpr_1(x,\mu(g_2,g_1)) \dto_{\bq\circ\bpr_1(g_2,g_2^{-1})} \xto[rrr]^{\rho_{(x,\mu(g_2,g_1))}} &&&
\bq\circ \sigma(x,\mu(g_2,g_1)) \dto_{\bq\circ\sigma(g_2,g_2^{-1})= \bq(e) = \id_{
\bq\circ \sigma(x,\mu(g_2,g_1))}} \\
\bq\circ\bpr_1(\sigma(x,g_2),e)\dto_= && & \bq\circ\sigma(x,\mu(g_2,g_1))\dto_=\\
 \bq\circ\bpr_1(\sigma(x,g_2),g_1)\xto[rrr]^{\rho_{(\sigma(x,g_2),g_1)}}&&&  \bq\circ\sigma(\sigma(x,g_2),g_1) .
\enddiagram\]
Together with the identity from above this shows the claim. 
\ebew

\begin{prop}\label{appC3}
Let $X$ and $G$ be as above. Then $[X/G]$ is a categorical quotient of $X$ with respect to the action of $G$ in the category of stacks, in the following sense. The quotient stack satisfies the universal property, that for each stack ${\catf}$, together with a morphism $\bq : X^\bullet \rightarrow {\catf}$ and a $2$-morphism $\rho: \bq\circ\bpr_1 \Rightarrow \bq\circ\sigma$, which is compatible with the group action as in lemma \ref{L167}, there exists a unique morphism of stacks ${\bf u}: [X/G] \rightarrow {\catf}$ such that the diagram
\[ \diagram
X^\bullet\times G^\bullet \rto^\sigma\dto_{\bpr_1}& X^\bullet\dto^\pi\xto[ddr]^\bq\\
X^\bullet\rto_\pi\xto[rrd]_\bq& [X/G]\drto_{\mathbf u}\\
&& {\catf} 
\enddiagram \]
commutes with ${\mathbf u}\circ \pi = \bq$, and such that $\rho=\id_{\mathbf u}\ast  \eta$. 
\end{prop}

\proof
Consider an object $(E,p,\phi)\in[X/G](B)$ for some scheme $B$. Suppose that it is given by a descent datum $\{((E_\alpha,p_\alpha,\phi_\alpha),t_{\beta,\alpha})\}_{\alpha,\beta\in A}$, with respect to some \'etale covering family $\{b_\alpha: B_\alpha\rightarrow B\}_{\alpha\in A}$. By the definition of a principal $G$-bundle, we may assume that $E_\alpha = B_\alpha\times G$, and that $p_\alpha = \pr_1$. There are isomorphisms 
\[ \tau_\alpha :\:\:\: B_\alpha\times G\rightarrow b_\alpha^\ast E ,\]
for all $\alpha\in A$, and for the transition morphisms on $B_{\alpha,\beta} := B_\alpha \times_B B_\beta$ holds
\[ t_{\beta,\alpha} =\overline{ \tau}_\beta\circ\overline{\tau}_\alpha^{-1} \]
for all $\alpha,\beta \in A$, where $\overline{ \tau}_\alpha$ and $ \overline{ \tau}_\beta$ denote the respective liftings of $\tau_\alpha$ and $\tau_\beta$ to $B_{\alpha,\beta}$. For $\alpha,\beta,\gamma\in A$ one has  obviously the cocycle identity $t_{\gamma,\alpha}=t_{\gamma,\beta}\circ t_{\beta,\alpha}$. 

By the sheaf property of the stack $\catf$ it thus suffices to define a functor $\mathbf u$ on objects $(E,p,\phi)$ in $[X/G]$, where $E$ is a trivial $G$-bundle, and on morphisms between such triples only. Once the condition of  functoriality 
\[ {\mathbf u}(\phi_2\circ\phi_1) = {\mathbf u}(\phi_2)\circ {\mathbf u}(\phi_1) \]
is verified for compatible morphisms $\phi_1$ and $\phi_2$ between such objects,   a functor from $[X/G]$ to ${\catf}$ is then well-defined by glueing. 

Consider an object $(B\times G,\pr_1,\phi)\in[X/G](B)$. Let $\nu: B \rightarrow B\times G$ denote the section of unity. Then $\phi\circ \nu \in X^\bullet(B)$, and hence
\[ {\mathbf u}(B\times G,\pr_1,\phi) := \bq(\phi\circ \nu) \in \catf(B) \]
is defined. 
$e: B\rightarrow G$ is constant morphism of unity.
A morphism between two objects $(B_2\times G,p_2:=\pr_1,\phi_2)$ and $(B_1\times G,p_1:=\pr_1,\phi_1)$ consists of a pair of morphisms $(f,\lambda)$, such that the diagram
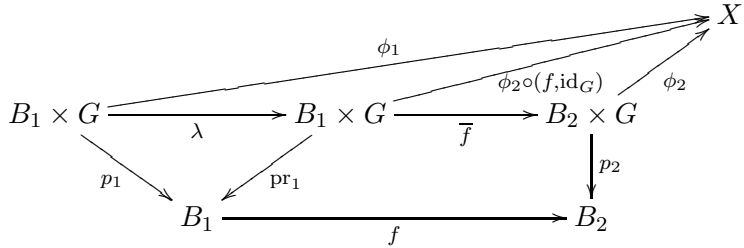
\[\diagram
&&&&&X\\
B_1\times G \rrto_\lambda\drto_{p_1}\xto[rrrrru]^{\phi_1} &&
B_1\times G \dlto^{\pr_1} \xto[urrr]_{\phi_2\circ(f,\id_G)} \rrto_{\overline{f}}& &
B_2\times G \urto_{\phi_2}\dto^{p_2}\\
&B_1\xto[rrr]_f&&& B_2
\enddiagram\]
commutes. Note that by its $G$-equivariance, the morphism $\lambda$ is given by a pair $(\pr_1,\mu(\gamma\circ\pr_1,\pr_2))$, where 
$\gamma: B_1 \rightarrow G$. We have the equality
\[ \phi_1= \phi_2\circ(f\circ\pr_1,\pr_2)\circ \lambda = \phi_2\circ(f\circ\pr_1,\mu(\gamma\circ\pr_1,\pr_2)),\]
and thus, for the sections of unity $\nu_i: B_i\rightarrow B_i\times G$, for $i=1,2$, 
\[ \phi_1\circ\nu_1 = \phi_2\circ(f\circ\pr_1,\gamma\circ\pr_1) =
\sigma(\phi_2\circ\nu_2\circ f\circ\pr_1,\gamma\circ\pr_1),\]
using the $G$-equivariance of $\phi_2$. The pair
$(\phi_2\circ\nu_2\circ f,\gamma)$ is in particular an object of $(X\times G)^\bullet(B_1)$, whereas $f$ can be considered as a morphism from $\phi_2\circ\nu_2$ to $\phi_2\circ\nu_2\circ f$ in $X^\bullet$. We define
\[ {\mathbf u}(f,\lambda) := \rho_{(\phi_2\circ\nu_2\circ f,\gamma)}\circ \bq(f) ,\]
by composing the morphisms
\[ 
{\mathbf u}(B_2\times G,p_2,\phi_2) = \bq(\phi_2\circ \nu_2) \:\:  \stackrel{\mathbf q(f)}{\longrightarrow } \:\:  \bq(\phi_2\circ \nu_2\circ f)
\]
and 
\[ \diagram
 \bq(\phi_2\circ \nu_2\circ f) = \bq\circ\bpr_1(\phi_2\circ\nu_2\circ f,\gamma)
\dto<-1.5cm>^{ \rho_{(\phi_2\circ\nu_2\circ f,\gamma)}} \hspace*{4.8cm} \\
   \hspace*{2.4cm} \bq\circ\sigma(\phi_2\circ\nu_2\circ f,\gamma) = \bq(\phi_1\circ \nu_1)  = \mathbf u(B_1\times G,p_1,\phi_1).
\enddiagram \]

To see that this defines indeed a functor, we need to verify that for a pair of morphisms 
\[ (f_1,\lambda_1) : \:\: (B_2\times G ,p_2,\phi_2) \rightarrow  (B_1\times G ,p_1,\phi_1)\]
and
\[ (f_2,\lambda_2) : \:\: (B_3\times G ,p_3,\phi_3) \rightarrow  (B_2\times G ,p_2,\phi_2)\]
in $[X/G]$ holds 
\[ {\mathbf u}((f_1,\lambda_1)\circ(f_2,\lambda_2)) = {\mathbf u}(f_1,\lambda_1)\circ{\mathbf u}(f_2,\lambda_2).\]
Note that $(f_1,\lambda_1)\circ(f_2,\lambda_2) = (f_2\circ f_1,\lambda_2\circ\lambda_1)$. Here $\lambda_2\circ\lambda_1$ is actually a lax notation for  the automorphism of $B_1\times G$, which is represented by the pair $
(\pr_1,\mu(\mu(\gamma_2\circ f_1\circ \pr_1,\gamma_1\circ\pr_1),\pr_2))$.
By definition, the right hand side of the above identity equals
\[ \mbox{r.h.s.} = \rho_{(\phi_2\circ\nu_2\circ f_1,\gamma_1)}\circ \bq(f_1)\circ \rho_{(\phi_3\circ\nu_3\circ f_2,\gamma_2)}\circ \bq(f_2),\]
if $\lambda_1=(\pr_1,\mu(\gamma_1\circ \pr_1,\pr_2))$ and $\lambda_2=(\pr_1,\mu(\gamma_2\circ \pr_1,\pr_2))$ are the automorphisms of $B_1\times G$ and of $B_2\times G$, respectively. Consider the commutative diagram
\[ \diagram
\bq\circ \bpr_1(\phi_3\circ\nu_3\circ f_2,\gamma_2) 
 \xto[d]^{\rho_{(\phi_3\circ\nu_3\circ f_2,\gamma_2)}} 
 \xto[rrr]^{\bq\circ\bpr_1(f_1,\id_G)}_{=\bq(f_1)} &&&
\bq\circ \bpr_1(\phi_3\circ\nu_3\circ f_2\circ f_1,\gamma_2\circ f_1)
 \xto[d]_{\rho_{(\phi_3\circ\nu_3\circ f_2\circ f_1,\gamma_2\circ f_1)}}\\
\bq\circ \sigma(\phi_3\circ\nu_3\circ f_2,\gamma_2)
 \xto[rrr]^{\bq\circ\sigma(f_1,\id_G)}_{=\bq(f_1)}&&& 
\bq\circ\sigma(\phi_3\circ\nu_3\circ f_2\circ f_1,\gamma_2\circ f_1) .
\enddiagram\]
This implies the identity
\[ \bq(f_1)\circ \rho_{(\phi_3\circ\nu_3\circ f_2,\gamma_2)} =\rho_{(\phi_3\circ\nu_3\circ f_2\circ f_1,\gamma_2\circ f_1)}\circ \bq(f_1) .\]
Since $\phi_2\circ\nu_2 = \sigma(\phi_3\circ\nu_3\circ f_2,\gamma_2)$, we have also
\[ \phi_2\circ\nu_2 \circ f_1= \sigma(\phi_3\circ\nu_3\circ f_2\circ f_1,\gamma_2\circ f_1) .\]
Thus, from lemma \ref{L167}, we obtain $\mbox{r.h.s.}=$ 
\[ \begin{array}{l}
 \rho_{(\sigma(\phi_3\circ\nu_3\circ f_2\circ f_1,\gamma_2\circ f_1),\gamma_1)} \circ \rho_{(\phi_3\circ\nu_3\circ f_2\circ f_1,\gamma_2\circ f_1)} \circ \bq(f_1)\circ \bq(f_2)= \\
\qquad = \rho_{(\phi_3\circ\nu_3\circ f_2\circ f_1,\mu(\gamma_2\circ f_1, \gamma_1))} \circ \bq(f_2\circ f_1) \\
\qquad =  {\mathbf u}(f_2\circ f_1,\lambda_2\circ \lambda_1) ,
\end{array}\]
as claimed. 

The morphism ${\mathbf u} :[X/G]\rightarrow {\catf}$ is uniquely determined by the condition ${\mathbf u}\circ \pi = \bq$. Indeed, this condition determines $\mathbf u$ uniquely on triples $(E,p,\phi)$, where $E$ is a trivial bundle, and $p$ its projection. A principal bundle is uniquely determined by its descent datum, so $\mathbf u$ is uniquely determined on objects. We still have to show that $\mathbf u$  it also uniquely determined on morphisms, if it is determined on morphisms between triples with trivial bundles.  Suppose that 
\[(f,\lambda) : \quad (B_2\times G,p_2:=\pr_1,\phi_2)\rightarrow (B_1\times G,p_1:=\pr_1,\phi_1) \]
is such a morphism. Let $\lambda$ be given by the pair $(\pr_1,\mu(\gamma\circ\pr_1,\pr_2))$, with $\gamma: B_1\rightarrow G$. If $\gamma$ is the constant morphism of unity $e$,   then $(f,\lambda)=(f,\id_{B_1\times G})= \pi(f)$. Then the condition ${\mathbf q}= {\mathbf u}\circ \pi$ implies that we must necessarily have ${\mathbf u}(f,\lambda)= {\mathbf q}(f)$. Note that this agrees with our definition from above, as 
\[ \rho_{(\phi_2\circ\nu_2\circ f, e)}= \id_{{\mathbf q}(\phi_2\circ \nu_2\circ f)} = \id_{{\mathbf u}(B_1\times G,p_1,\phi_1)} \]
holds by assumption. 

Suppose now that $\gamma: B_1\rightarrow G$ is nontrivial. By the definition of the $2$-morphism $\eta: \pi\circ \bpr_1\Rightarrow \pi\circ \sigma$ as in the proof of  lemma \ref{Q001},  we have the equality
\[ (f,\lambda) =\eta_{(\phi_2\circ\nu_2\circ f,\gamma)}\circ\pi(f) .\]
Note that the pair $(f,\gamma)$  can be viewed as a morphism
\[ (f,\gamma) : \quad (\phi_2\circ\nu_2,\id_G) \rightarrow (\phi_2\circ\nu_2\circ f,\gamma) \]
in the category $(X\times G)^\bullet(B_1)$. Thus, by the definition of a $2$-morphism, there is a commutative diagram
\[ \diagram
\pi\circ\bpr_1(\phi_2\circ\nu_2,\id_G) \dto_{\pi(f)=\pi\circ\bpr_1(f,\gamma)} \rrto^{\eta_{(\theta_2\circ\nu_2,\id_G)}=\id_{\pi(\theta_2\circ\nu_2)}} &&
\pi\circ\sigma(\theta_2\circ\nu_2,\id_G)= \pi(\theta_2\circ\nu_2) \dto^{\pi\circ\sigma(f,\gamma)}\\
\pi\circ\bpr_1(\phi_2\circ\nu_2\circ f,\gamma) \rrto^{\eta_{(\phi_2\circ\nu_2\circ f,\gamma)}}&& \pi\circ\sigma(\phi_2\circ\nu_2\circ f,\gamma) = \pi(\phi_1\circ \nu_1). 
\enddiagram \]
From this we obtain the identity of morphisms
\[ \pi\circ\sigma(f,\gamma) = \eta_{(\phi_2\circ\nu_2\circ f,\gamma)}\circ \pi(f) .\]
Together with the first equality we find
\[ (f,\lambda) = \pi\circ\sigma(f,\gamma) .\]
The condition $\mathbf u\circ \pi = \mathbf q$ now implies that if $\mathbf u$ exists, then it must satisfy $\mathbf u(f,\lambda) = \mathbf q\circ \sigma(f,\gamma)$. 
An analogous diagram as above shows that
\[ {\mathbf q}\circ\sigma(f,\gamma) = \rho_{(\phi_2\circ\nu_2\circ f,\gamma)}\circ{\mathbf q}(f),\]
and hence
\[ {\mathbf u}(f,\lambda) = \rho_{(\phi_2\circ\nu_2\circ f,\gamma)}\circ{\mathbf q}(f), \]
which confirms our definition of $\mathbf u$, and shows that this is the only possible choice of definition. Since any morphism in $[X/G]$ is uniquely determined by its restictions to triples with trivial bundles, we are done. 

Note that the equality $\rho=\id_{\mathbf u}\ast\eta$ is automatically satisfied. Indeed, for any $f:B\rightarrow X\times G$ in $(X\times G)^\bullet(B)$, one verifies the identity
\[ (\id_{\mathbf u}\ast\eta)_f = \mathbf u(\eta_f) = \rho_f\]
directly from the definitions.  
\ebew

\begin{rk}\em
Note that even if a categorical quotient $X/G$ of $X$ under $G$ exists in the category of schemes, then $(X/G)^\bullet$ is in general not the categorical quotient in the category of stacks. In fact, there is usually not even a morphism from $(X/G)^\bullet$ to $[X/G]$ which makes the respective diagram commutative. As an example think of the case where $G$ is a nontrivial group, and $X =G$. 
\end{rk}

\begin{rk}\em\label{appC59}
The universal property of the quotient stack from proposition \ref{appC3} is not preserved under isomorphisms of stacks. Let $\catq$ be a stack isomorphic to $[X/G]$, with two morphisms $\mathbf a : \catq\rightarrow [X/G]$ and $\mathbf b: [X/G] \rightarrow \catq$, as well as two $2$-morphisms $\alpha: \mathbf b\circ \mathbf a \Rightarrow \id_{\catq}$ and $\beta: \mathbf a\circ \mathbf b \Rightarrow \id_{[X/G]}$. If $\pi' := \mathbf b\circ\pi$, then there is a commutative square with a $2$-morphism $\eta' := \id_{\mathbf b}\ast\eta: \:\: \pi'\circ\pr_1\Rightarrow \pi' \circ\sigma$. 

Let $\catf$ be a stack as in proposition \ref{appC3}. Obviously there exists a morphism $\mathbf u' := \mathbf u\circ \mathbf a: \:\: \catq\rightarrow \catf$, such that the diagram
\[ \diagram
X^\bullet\times G^\bullet \rto^\sigma\dto_{\bpr_1}& X^\bullet\dto^{\pi'}\xto[ddrr]^\bq\\
X^\bullet\rto^{\pi'}\xto[rrrd]_\bq& \catq\xto[rrd]^{\mathbf u'}\\
&&& {\catf} 
\enddiagram \]
commutes, with a $2$-morphism $\mathbf u'\circ\pi' \Rightarrow \mathbf q$, and such that $\rho=\id_{\mathbf u'}\ast \eta'$. 

A stack $\catq$, which is isomorphic to $[X/G]$, satisfies the following universal property. Let ${\catf}$ be a stack, together with a morphism $\bq : X^\bullet \rightarrow {\catf}$ and a $2$-morphism $\rho: \bq\circ\bpr_1 \Rightarrow \bq\circ\sigma$, which is compatible with the action of $G$ as in lemma \ref{L167}. Then there exists a  morphism of stacks ${\bf u'}: \catq \rightarrow {\catf}$ and a $2$-morphism $\tau' : \mathbf u'\circ \pi' \Rightarrow \mathbf q$, such that the above diagram
commutes, with  $\rho=\id_{\mathbf u'}\ast  \eta'$. The morphism $\mathbf u'$ is unique in the following sense. If there exists another morphism $\mathbf u'' : \catq \rightarrow \catf$, together with a $2$-morphism $\tau'' : \mathbf u''\circ \pi' \Rightarrow \mathbf q$, such that $\rho=\id_{\mathbf u''}\ast  \eta'$, then there exists a unique $2$-morphism $\lambda: \mathbf u'' \Rightarrow \mathbf u'$ such that $\tau'\circ(\lambda\ast\id_{\pi'}) = \tau''$. 

Indeed, the existence of such a morphism $\mathbf u'$ is clear. One simply defines $\mathbf u' := \mathbf u\circ\mathbf a$, and finds immediately the identity $\tau'=\id_{\mathbf u}\ast\beta\ast\id_{\pi}$. Suppose now that there is a second morphism $\mathbf u'':\catq\rightarrow \catf$ together with a $2$-morphism  $\tau'' : \mathbf u''\circ \pi' \Rightarrow \mathbf q$. Since $\mathbf q=\mathbf u\circ \pi$, and $\pi'= \mathbf b\circ \pi$, and since $\pi$ is surjective, there exists a $2$-morphism $\tau: \mathbf u''\circ \mathbf b\Rightarrow \mathbf u$, and therefore a $2$-morphism
\[ \mathbf u''\circ \mathbf b\circ\mathbf a\Rightarrow \mathbf u\circ\mathbf a
 \Rightarrow \mathbf u' .\]
Using the $2$-morphism $\alpha: \mathbf b\circ\mathbf a\Rightarrow \id_\catq$, one can construct a $2$-morphism $\mathbf u'' \Rightarrow \mathbf u'$ as claimed. 
\end{rk}

\begin{rk}\em\label{A0060}
The canonical quotient morphism $\pi:  X^\bullet \rightarrow [X/G]$ is
representable, surjective and smooth. To see this, we first have to prove that for any scheme $B
\in \ob(\catc)$, and all morphisms $\mathbf b : B^\bullet \rightarrow
[X/G]$ the stack
$ X^\bullet \times_{\pi,[X/G],\mathbf b} B^\bullet $
defined by the Cartesian diagram
\[ \diagram
 X^\bullet \times_{\pi,[X/G],\mathbf b} B^\bullet \rto \dto &
 B^\bullet  \dto^{\mathbf b}\\
X^\bullet \rto_\pi & [X/G] 
\enddiagram 
\]
is representable. But by Yoneda's lemma, the morphism $\mathbf b$
corresponds to an object  $(E,p,\phi)$ in $[X/G](B)$, 
where $p : E \rightarrow B$ is a principal $G$-bundle, and $\phi$ is
a $G$-equivariant morphism. By proposition \ref{P6137}, the diagram
\[ \diagram
E^\bullet \rto^{\mathbf p}\dto_{\phi^\bullet} & B^\bullet\dto^{\mathbf
b}\\
X^\bullet \rto_\pi & [X/G]  
\enddiagram\]
is Cartesian.  Hence there is an isomorphism 
\[ E^\bullet \:\: \equiv \:\: X^\bullet \times_{\pi,[X/G],\mathbf b} B^\bullet \]
of stacks, and thus  $
X^\bullet \times_{\pi,[X/G],\mathbf b} B^\bullet$ is representable by
 the scheme $E$. 

Note that the morphism $ p$ is surjective by definition, so the
quotient morphism $\pi$ is surjective, and $\pi$ is smooth since the projection $p: E\rightarrow B$ is smooth. 
\end{rk}

Using the same techniques as in the proof of proposition \ref{appC3}, we can generalize the Cartesian diagram of remark \ref{A0060}. With the following observation we provide a useful tool for studying quotient stacks, of which proposition \ref{P6137} is just  one special case.

\begin{prop}\label{A0062}
Let $[X/G]$ be a quotient stack as above, and let $U,V\in\ob(\catc)$ be schemes. Consider morphisms $\mathbf u: U^\bullet \rightarrow [X/G]$ and $\mathbf v: V^\bullet \rightarrow [X/G]$, represented by triples $(E_U,p_U,\phi_U)\in[X/G](U)$ and $(E_V,p_V,\phi_V)\in[X/G](V)$, respectively. There is a natural diagonal action of $G$ on $E_U\times_X E_V$, and the diagram
\[\diagram
(E_U\times_X E_V)/G^\bullet \rto \dto & V^\bullet\dto^{\mathbf v}\\
U^\bullet \rto_{\mathbf u} & [X/G]
\enddiagram\]
is Cartesian.
\end{prop}

\begin{rk}\em
There is a commutative diagram
\[\diagram
(E_U\times_X E_V)^\bullet \rto\dto& E_V^\bullet\dto^{\phi_V}\rto^{\mathbf p_V}&V^\bullet\dto^{\mathbf v}\\
E_U^\bullet \rto^{\phi_U}\dto_{\mathbf p_U}& X^\bullet \rto^\pi\dto^\pi& [X/G]\\
U^\bullet \rto_{\mathbf u} & [X/G],
\enddiagram \]
where all squares are Cartesian. In particular, one has natural isomorphisms
\[ U\times_{[X/G]} E_V \cong E_U\times_X E_V \cong E_U \times_{[X/G]} V .\]
\end{rk}

\proofof{proposition \ref{A0062}}
Note that the quotient $(E_U\times_X E_V)/G$ exists as a scheme. 
Indeed, $E_U\times_X E_V$ is a $G$-invariant closed subscheme of $E_U\times E_V$, and the quotient $(E_U\times E_V)/G$ exists. 
We need to construct an isomorphism of stacks 
\[ U^\bullet \times_{[X/G]} V^\bullet \cong (E_U\times_X E_V)/G^\bullet .\]
For a scheme $B\in \ob(\catc)$, an object of $U^\bullet \times_{[X/G]} V^\bullet(B)$ is given by a triple $(f,g,\lambda)$, where $f:B\rightarrow U$ and $g: B\rightarrow V$ are morphisms of schemes in $\catc$, and $\lambda :g^\ast E_V \rightarrow f^\ast E_U$ is an isomorphism of principal $G$-bundles, satisfying the usual compatibility condition with $\phi_U$ and $\phi_V$.

In analogy to the proof of proposition \ref{appC3}, it is enough to construct a functor on such triples, where $g^\ast E_V$ is trivial, and on morphisms between such triples. Once the functoriality is established, the functor can be extended to all of $U^\bullet \times_{[X/G]} V^\bullet $ by glueing. 

Let $(f,g,\lambda)\in U^\bullet \times_{[X/G]} V^\bullet(B)$ be given for some scheme $B$, and such that $g^\ast E_V$ is trivial. Then there exists the zero-section $s: B\rightarrow g^\ast E_V$. Composition with the canonical  morphism $\overline{g}: g^\ast E_V\rightarrow E_V$ defines a morphism 
\[ \tilde{g}:= \overline{g}\circ s : \:\:\: B \rightarrow E_V.\]
Using the isomorphism $\lambda:g^\ast E_V \rightarrow f^\ast E_U$ and the canonical morphism $\overline{f}: f^\ast E_U \rightarrow E_U$, we define
\[ \tilde{f} := \overline{f}\circ\lambda\circ s: \:\:\: B \rightarrow E_U.\] 
Note that by assumption the isomorphism $\lambda$ satisfies the equality
\[ \phi_U \circ \overline{f}\circ \lambda = \phi_V \circ \overline{g} .\]
From this one immediately verifies the identity
\[ \phi_U \circ\tilde{f} = \phi_V \circ \tilde{g} .\]
Therefore, the pair $(\tilde{f},\tilde{g})$ defines an object of $(E_U\times_X E_V)^\bullet(B)$, and thus  an object of $(E_U\times_X E_V)/G^\bullet(B)$. 

A morphism between triples $(f,g,\lambda)\in U^\bullet \times_{[X/G]} V^\bullet(B)$ and $(f',g',\lambda')\in U^\bullet \times_{[X/G]} V^\bullet(B')$ is given by a morphism $h : B' \rightarrow B$ of schemes, such that $f'=f\circ h$ and $g' = g\circ h$, and which is compatible with $\lambda$ and $\lambda'$. Thus it clearly induces a morphism from  $(\tilde{f},\tilde{g})$ to  
 $(\tilde{f}',\tilde{g'})$, and this assignment is functorial. 

For non-trivial bundles $g^\ast E_V$, the functor is constructed by glueing local trivializations. Note that two trivializations over some \'etale open subscheme $B_0 \rightarrow B$ may differ by the action of a section $\gamma: B_0\rightarrow G$. Let  $(f,g,\lambda)\in U^\bullet \times_{[X/G]} V^\bullet(B_0)$ be given, and let $(\tilde{f}_1,\tilde{g}_1)$ and $(\tilde{f}_2,\tilde{g}_2)$ be two pairs obtained using  two different trivializations of $g^\ast E_V|B_0$. By construction, and by the $G$-equivariance of $\lambda$, one has the identities
$ \tilde{f}_1= \sigma_{E_U}(\tilde{f}_2,\gamma\circ p_U\circ\tilde{f}_2 )$ and $ \tilde{g}_1= \sigma_{E_V}(\tilde{g}_2,\gamma\circ p_V\circ\tilde{g}_2 )$, where $\sigma_{E_U}$ and $\sigma_{E_V}$ denote the actions of $G$ on $E_U$ and $E_V$, respectively. Thus by definition, the pairs $(\tilde{f}_1,\tilde{g}_1)$ and $(\tilde{f}_2,\tilde{g}_2)$ define the same object in the quotient $(E_U\times_X E_V)/G^\bullet(B_0)$, i.e. the assignment is independent of the chosen trivialization.  

Now the local definitions can be glued together because cocycles are preserved due to the functoriality of the construction. Here we use the descent property of algebraic stacks. 

To define the inverse functor from   $(E_U\times_X E_V)/G^\bullet$ to $U^\bullet \times_{[X/G]} V^\bullet$ consider a pair $(\tilde{f},\tilde{g}) \in (E_U\times_X E_V)^\bullet(B)$ for some scheme $B$. Suppose that it is represented by a pair of morphisms $\tilde{f}_0:B\rightarrow E_U$ and $\tilde{g}_0:B\rightarrow E_V$. Then for the pair
\[ (p_U\circ \tilde{f}_0,p_V\circ\tilde{g}_0) : \:\: B \rightarrow U\times V \]
holds $(p_U\circ \tilde{f}_0)^\ast E_U \cong B\times G \cong (p_V\circ \tilde{g}_0)^\ast E_V$. Thus there is an isomorphism $\lambda: (p_V\circ\tilde{g}_0)^\ast E_V \rightarrow (p_U\circ\tilde{f}_0)^\ast E_U$ of principal $G$-bundles over $B$. Consider the diagram
\[\diagram
&& X \\
E_V \dto_{p_V} \xto[urr]^{\phi_V}& (p_V\circ \tilde{g}_0)^\ast E_V \lto_{q_V}\drto^{\pi_V} \rrto^\lambda && (p_U\circ\tilde{f}_0)^\ast E_U \dlto_{\pi_U} \rto^{q_U} & E_U \xto[ull]_{\phi_U}\dto^{p_U}\\
V && B \llto^{p_V\circ\tilde{g}_0} \xto[ull]^{\tilde{g}_0} \xto[urr]_{\tilde{f}_0} \rrto_{p_U\circ \tilde{f}_0} && U .
\enddiagram \] 
By assumption, we have the identity
\[ \phi_V \circ \tilde{g}_0 = \phi_U \circ \tilde{f}_0 ,\]
since $(\tilde{f}_0,\tilde{g}_0)\in (E_U\times_X E_V)^\bullet(B)$. The isomorphism $\lambda$ satisfies $\pi_V = \pi_U\circ \lambda$. Let $s_V : B\rightarrow (p_V\circ \tilde{g}_0)^\ast E_V$ denote the canonical trivializing section of $\pi_V$. Then we have $q_V \circ s_V = \tilde{g}_0$. The corresponding  trivializing section of $\pi_U$ is $s_U = \lambda\circ s_V : B\rightarrow (p_U\circ \tilde{f}_0)^\ast E_U$, satisfying $q_U\circ s_U= \tilde{f}_0$. Thus we obtain the equality
\[ \begin{array}{ccl}
\phi_V\circ q_V\circ s_V &=& \phi_U\circ q_U \circ s_U \\
&=& \phi_U\circ q_U \circ \lambda\circ s_V .
\end{array}\]
All morphisms except the sections are $G$-equivariant. This implies that there is an identity
\[ \phi_V\circ q_V = \phi_U\circ q_U \circ \lambda,\]
showing that the triple 
 $(p_U\circ \tilde{f}_0,p_V\circ\tilde{g}_0,\lambda)$ is an object of $U^\bullet \times_{[X/G]} V^\bullet(B)$. 

Since the construction involves the application of the projections $p_U$ and $p_V$, this triple is indeed independent of the chosen representation of $(\tilde{f},\tilde{g})$. It is straightforward to make the appropriate construction  for morphisms in $(E_U\times_X E_V)/G^\bullet$, and to show that this defines a functor. This functor  is inverse to the functor defined in the first part of the proof, in the sense of morphisms of stacks. 
\ebew 

\begin{rk}\em\label{A0174}
Let $G$ and $H$ be linear algebraic groups acting independently on a scheme $X\in\ob(\catc)$.  There is a natural morphism of stacks
\[ \sigma: \:\: [X/G] \times_\scat{C} H^\bullet \rightarrow [X/G], \]
which is defined as follows. For a tuple $((E,p,\phi),h) \in [X/G] \times_\scat{C} H^\bullet(B)$, for some scheme $B\in\ob(\catc)$, we put
\[ \sigma((E,p,\phi),h) := (E,p,\sigma_H(\phi,h\circ p)) \in [X/G](B),\]
where $\sigma_H$ denotes the action of $H$ on $X$ as usual. Note that $\sigma_H(\phi,h\circ p)$ is indeed $G$-equivariant since the actions of $G$ and $H$ are independent. A morphism from $((E_2,p_2,\phi_2),h_2) \in [X/G] \times_\scat{C} H^\bullet(B_2)$ to  $((E_1,p_1,\phi_1),h_1) \in [X/G] \times_\scat{C} H^\bullet(B_1)$ is a pair of morphisms $((f,\lambda),\chi)$, where $(f,\lambda)$ is a morphism in $[X/G]$, and $\chi$ is a morphism in $H^\bullet$, such that both morphisms project to the same morphism in $\catc$. In particular, 
$f: B_1\rightarrow B_2$ is a morphism of schemes, $\lambda: E_1\rightarrow f^\ast E_2$ is a morphism of principal $G$-bundles, and $\chi: B_1\rightarrow B_2$ is a morphism in $\catc$, satisfying $h_1= h_2\circ\chi$. By definition of the projection $[X/G]\rightarrow \catc$ we must have $\chi = f$. We now define
\[ \sigma((f,\lambda),\chi) := (f,\lambda) .\]
One sees easily that this defines indeed  a functor from $[X/G] \times_\scat{C} H^\bullet$ to $[X/G]$. 

This morphism $\sigma$ can be interpreted as an action of the group $H$ on the stack $[X/G]$. Indeed, if $\mu_H: H\times H \rightarrow H $ denotes the group multiplication of $H$, one immediately verifies the identity
\[ \sigma(\id_{[X/G]}\times \mu_H) = \sigma(\sigma\times \id_{H^\bullet}) \]
as morphisms from $ [X/G] \times_\scat{C} H^\bullet \times_\scat{C} H^\bullet $ to $[X/G]$, as well as the identity
\[ \sigma(\id_{[X/G]}\times \mathbf e) =\bpr_1 \]
as morphisms from $[X/G]\times_\scat{C} S^\bullet $ to $[X/G]$, where $e: S \rightarrow H$ denotes the constant morphism of unity. 
\end{rk}

\begin{lemma}
Let $G$ and $H$ be linear algebraic groups acting independently on a scheme $X\in\ob(\catc)$. Then the diagram
\[ \diagram
[X/G] \times_\scat{C} H^\bullet \rto^{\sigma}\dto_{\bpr_1} & [X/G]\dto^\pi\\
[X/G]\rto_\pi & [X/G\times H]
\enddiagram\]
is commutative, with a $2$-morphism $\eta: \pi\circ\bpr_1\Rightarrow \pi\circ\sigma$, where $\pi$ denotes the canonical morphism. 
\end{lemma}

\proof
Let $((E,p,\phi),h) \in [X/G] \times_\scat{C} H^\bullet(B)$ be given for some scheme $B\in\ob(\catc)$. One computes 
\[ \pi\circ\bpr_1 ((E,p,\phi),h) = (E\times H ,p\circ\pr_1, \sigma_H(\phi\circ\pr_1,\pr_2)) ,\]
as well as
\[ \pi\circ\sigma ((E,p,\phi),h) = (E\times H ,p\circ\pr_1, \sigma_H(\sigma_H(\phi,h\circ p)\circ\pr_1,\pr_2)) .\]
There is an isomorphism $\eta_{((E,p,\phi),h)} :  \pi\circ\bpr_1 ((E,p,\phi),h) \rightarrow  \pi\circ\sigma ((E,p,\phi),h)$, which is  given on $E\times H$ by
\[ \eta_{((E,p,\phi),h)} := (\pr_1, \mu(h\circ p\circ\pr_1,\pr_2)).\]
Since $\eta_{((E,p,\phi),h)}$ is $G\times H$-equivariant, this is indeed an isomorphism of principal $G\times H$-bundles. One easily verifies that the cocycle condition with respect to morphism in $[X/G]\times_\scat{C} H^\bullet$ is satisfied, so that there is indeed a $2$-morphism  $\eta: \pi\circ\bpr_1\Rightarrow \pi\circ\sigma$. 
\ebew
     
\begin{prop}
Let $G$ and $H$ be linear algebraic groups acting independently on a scheme $X\in\ob(\catc)$. Then $[X/G\times H]$ is a categorical quotient of $[X/G]$ with respect to the action of $H$ in the category of stacks. 
\end{prop} 

\proof
Let $\catf$ be an algebraic stack, together with a morphism $\mathbf q: [X/G]\rightarrow \catf$, and a $2$-morphism $\rho: \mathbf q\circ\bpr_1\Rightarrow \mathbf q\circ \sigma$. We need to show that there exists a unique morphism $\mathbf u : [X/G\times H]\rightarrow \catf$ such that $\mathbf u\circ\pi= \mathbf q$. 

From the sheaf property of the stack $\catf$ is follows that it is enough to construct such a functor on triples $(E,p,\phi)\in[X/G\times H](B)$, where $E = E' \times H$ for some principal $G$-bundle $p': E'\rightarrow B$, with $p=p'\circ \pr_1$, and morphisms between such triples. Compare this to our strategy in the proof of proposition \ref{appC3}. 

Let $(E'\times H ,p,\phi)\in [X/G\times H](B)$ be given.  Let $\nu: B\rightarrow H$ denote the constant morphism of unity. Then $p' : E' \rightarrow B$ is a subbundle of $p: E'\times H\rightarrow B$ via $j :=(\id_{E'},\nu\circ p')$, and $(E',p',\phi\circ j)$ is an object of $[X/G](B)$. 
We now define
\[ \mathbf u(E'\times H , p, \phi) := \mathbf q(E',p',\phi\circ j).\]
Consider a morphism from $(E_2'\times H,p_2,\phi_2)\in[X/G\times H](B_2)$ to $(E_1'\times H,p_1,\phi_1)\in[X/G\times H](B_1)$. It is given by a pair $(f,\lambda)$, where $f: B_1\rightarrow B_2$ is a morphism of schemes, and $\lambda: E_1'\times H \rightarrow f^\ast (E_2'\times H)$ is an isomorphism of principal $G\times H$-bundles. Thus $\lambda$ is given by a pair $(\lambda',\mu_H(\chi\circ p_1,\pr_2))$, where $\lambda' : E_1'\rightarrow f^\ast E_2'$ is a morphism of principal $G$-bundles over $B_1$, and $\chi : B_1\rightarrow H$ is a morphism of schemes. We write $\overline{f}:f^\ast E_2'\rightarrow E_2'$ for the induced morphism. In particular, $(f,\lambda')$ is a morphism from the triple $(E_2',p_2',\phi_2\circ j_2)$ to $(E_1',p_1',\phi_2\circ j_2\circ \overline{f}\circ \lambda')$ in $[X/G]$. Therefore there is a morphism
\[ \mathbf q(f,\lambda') :  \mathbf u(E_2'\times H,p_2,\phi_2) =  \mathbf q(E_2',p_2',\phi_2\circ j_2) \rightarrow \mathbf q(E_1',p_1',\phi_2\circ j_2\circ \overline{f}\circ\lambda'). \]
Applying the $2$-morphism $\rho$, we obtain a morphism 
\[ \begin{array}{c}
\mathbf q(E_1',p_1',\phi_2\circ j_2\circ \overline{f}\circ \lambda')= \mathbf q\circ\bpr_1((E_1',p_1',\phi_2\circ j_2\circ\overline{f}\circ \lambda'),\chi) \\[3mm]
\downarrow \:\: 
{\rho_{((E_1',p_1',\phi_2\circ j_2\circ\overline{f}\circ\lambda'),\chi)} }\\[3mm]
\mathbf q\circ\sigma((E_1',p_1',\phi_2\circ j_2\circ\overline{f}\circ \lambda'),\chi) =\mathbf q(E_1',p_1',\phi_1\circ j_1) = \mathbf u(E_1'\times H, p_1,\phi_1).
\end{array} \]
We define
\[ \mathbf u(f,\lambda) := \rho_{((E_1',p_1',\phi_2\circ j_2\circ \overline{f}\circ \lambda'),\chi)} \circ \mathbf q(f,\lambda') .\]
As in the proof of proposition \ref{appC3}, one can verify that this definition is functorial on such morphisms as considered above. Finally, using the functoriality, we can glue the local definitions of the functor to define a global morphism $\mathbf u$ from $[X/G\times H]$ to $\catf$, as we did in the proof of proposition \ref{appC3}. The verification of the uniqueness of the morphism $\mathbf u$ is also completely analogous.
\ebew

\begin{rk}\em
Note that the above universal property of a categorical quotient is not compatible with isomorphisms of stacks. The situation here is analogous to that of remark \ref{appC59}.
\end{rk}

\begin{prop}
Let $G$ and $H$ be linear algebraic groups acting independently on a scheme $X\in\ob(\catc)$. Suppose that $G$ acts freely on $X$, and that $X'\in\ob(\catc)$ is a quotient of $X$ by the action of $G$. Then there is a natural isomorphism of stacks
\[ [X/ G\times H] \equiv [X'/H] .\]
\end{prop}

\proof
A morphism $\Theta$ from $[X/ G\times H]$ to $[X'/H]$ is constructed by sending a triple $(E,p,\phi)\in  [X/ G\times H](B)$, for some scheme $B\in\ob(\catc)$, to the triple $(E/G,\overline{p},\overline{\phi})\in[X'/H](B)$, where $\overline{p}: E/G \rightarrow B$ is induced by  $p: E\rightarrow B$, and $\overline{\phi}:E/G \rightarrow X'=X/G$ is induced by $\phi: E\rightarrow X$. The definition of $\Theta$ on morphisms is the obvious one. 

Conversely, let $(E',p',\phi')\in[X'/H](B)$ be given. Since the action of $G$ on $X$ is free, the Cartesian diagram
\[\diagram 
\tilde{E} \rto^{\tilde{\phi}}\dto_{\tilde{p}'} & X \dto\\
E' \rto_{p'}& X'=X/G
\enddiagram \]
defines a principal $G$-bundle $\tilde{p}' : \tilde{E}\rightarrow E'$, compare remark \ref{appC27}.  In particular, the composition $\tilde{p}:= p'\circ\tilde{p}': \tilde{E}\rightarrow B$ is a principal $G\times H$-bundle over $B$. Thus the triple $(\tilde{E},\tilde{p},\tilde{\phi})$ is an object of $[X/ G\times H](B)$. 

Up to a $2$-morphism, this construction is inverse to the construction of $\Theta$. Indeed, one has obviously $\tilde{E}/G = E'$. Conversely, a morphism $\tau : E \rightarrow \tilde{E}$ exists by the universal property of the fibre product. Since $\tau$ is in fact $G\times H$-equivariant, it is necessarily an isomorphism of principal $G\times H$-bundles over $B$. 
\ebew

\begin{prop}
Let $X\in\ob(\catc)$ be a Noetherian scheme of finite type, and let $G$ be a smooth linear algebraic group acting on $X$. Suppose that the stabilizers of geometric points of $X$ are finite and reduced. Then $[X/G]$ is a Deligne-Mumford stack. The quotient stack $[X/G]$ is separated if and only if the action of $G$ on $X$ is proper.
\end{prop}

\proof
See \cite[Cor. 2.2]{Ed}. To show the first part of the claim, it suffices to verify that the conditions of proposition \ref{appB51} are satisfied. 
\ebew

\begin{rk}\em\footnote{\cite[Ex. 7.17]{Vi}}
$(i)$ Suppose that $[X/G]$ is a Deligne-Mumford stack. Then there exists an \'etale atlas $\mathbf a: A^\bullet \rightarrow [X/G]$, for some scheme $A\in\ob(\catc)$. By Yoneda's lemma, this translates to the existence of a principal $G$-bundle $p: E\rightarrow A$, together with a $G$-equivariant morphism $\phi: E\rightarrow X$, which is \'etale and surjective.\\
$(ii)$  If $G$ is a smooth group, then the canonical morphism $\pi: X^\bullet\rightarrow [X/G]$ is an atlas of $[X/G]$, considered as an Artin stack.\footnote{\cite[Ex. 2.25]{Go}}. If the group $G$ is \'etale over $S$, then this morphism is even an atlas for the Deligne-Mumford stack $[X/G]$, see \cite[Ex. 4.8]{DM}. 
\end{rk}

\begin{lemma}\label{appC62}
Let $G$ and $H$ be linear algebraic groups acting independently on a scheme $X\in\ob(\catc)$. Suppose that both $[X/G]$ and $[X/G\times H]$ are integral Deligne-Mumford stacks. If $H$ is a finite group, then the canonical morphism
\[ {\mathbf q}: \:\: 
[X/G] \rightarrow [X/G\times H] \]
is unramified and finite of degree 
equal to the order of $H$.
\end{lemma}

\proof
Choose an integral atlas $\mathbf a: A^\bullet \rightarrow [X/G\times H]$. Let $\mathbf a$ be represented by the triple $(E,p,\phi)\in[X/G\times H](A)$. Then proposition \ref{appC45a} implies that the diagram
\[\diagram
(E/G)^\bullet \rto^{\mathbf p'} \dto& A^\bullet\dto^{\mathbf a}\\
[X/G]\rto_{\mathbf q}& [X/G\times H]
\enddiagram\]
is Cartesian, where $p'$ is induced by the projection $p:E\rightarrow A$. Thus $p'$ is finite and unramified, and its degree is  equal to the order of $H$. 
\ebew

\begin{rk}\em
In the special case of $G = \{\id\}$ and $H$ finite, this implies that the canonical morphism
\[ \pi : \:\: X \rightarrow [X/H] \]
is unramified and finite of degree equal to the order of $H$.
\end{rk}

\begin{cor}\label{appC64}
Let $G$ be a finite linear algebraic group over $k$ of order $n$, acting trivially on $X \in \ob(\catc)$. Then the natural morphism
\[ \Sigma  : \:\: [X/G] \rightarrow X^\bullet \]
is finite of degree $\frac{1}{n}$.
\end{cor}

\proof
The identity morphism $\id_{X^\bullet} : X^\bullet \rightarrow X^\bullet$ factors as
\[ \diagram
X^\bullet \rto^\pi & [X/G] \rto^\Sigma & X^\bullet ,
\enddiagram\]
where $\Sigma$ is defined as in example \ref{appC49}. 
From  lemma \ref{appC62} it follows that the degree of the quotient morphism $\pi$ is equal to $n$.
\ebew

\begin{rk}\em
Let $G$ and $H$ be linear algebraic groups acting independently on $X\in \ob(\catc)$, and suppose that both $[X/G]$ and $[X/G\times H]$ are integral Deligne-Mumford stacks. If $H$ is finite of order $n$, and  acts trivially on $X$, then 
it follows from lemma \ref{appC62} that the natural morphism 
\[ \omega : \:\: [X/G\times H] \rightarrow [X/G] \]
 is finite of degree $\frac{1}{n}$, since $\id_{[X/G]} = \omega\circ\mathbf q$.  If $H \neq \{\id\}$ then $\omega$ is not representable.
\end{rk} 

\begin{prop}
Let $[X/G]$ be a Deligne-Mumford stack, and let $M\in\ob(\catc)$ be a scheme. The scheme $M$ is a moduli space for the stack $[X/G]$ if and only if $M$ is a quotient of $X$ by $G$.
\end{prop}

\proof
See \cite[Prop. 2.11]{Vi}.
\ebew 

\begin{rk}\em
Let $[X/G]$ be a Deligne-Mumford stack with moduli space $\mathbf q: [X/G]\rightarrow M^\bullet$, such that $M$ is a quotient of $X$ by $G$. \\
$(i)$ If $(E,p,\phi)\in[X/G](B)$ for some scheme $B\in\ob(\catc)$, then $\mathbf q(E,p,\phi)=\overline{\phi} \in M^\bullet(B)$, where $\overline{\phi}: B=E/G \rightarrow M$ is the morphism induced by $\phi: E\rightarrow X$.\\
$(ii)$ Let $M_0 \subset M$ be a subscheme, and let $X_0$ denote its preimage under the composed morphism $X^\bullet \rightarrow [X/G]\rightarrow M^\bullet$. Then both squares of the diagram
\[\diagram
X_0^\bullet \rto\dto & X^\bullet\dto\\
[X_0/G] \rto\dto &[X/G] \dto \\
M_0^\bullet\rto & M^\bullet
\enddiagram \]
are Cartesian.
\end{rk}

%
%


\chapter*{List of Notations}\label{symbole}
\addtocontents{toc}{\protect \contentsline {chapter}{\protect \numberline
{}List of Notations}{\pageref{symbole}} }

\pagestyle{plain}

For the comfort of the reader we collect below the most frequently used notations,  and indicate on which page they were introduced. Symbols, which occur with and without overlines, are usually only listed in their closed (i.e. overlined) version. 

\begin{tabular}{llr}
$\msbar_g$  & moduli space of stable curves of genus $g$ & \pageref{n1} \\
$\msbar_g^{(k)}$ & locus of curves with at least $k$ nodes in $\msbar_g$ &\pageref{n2}\\
$\Sigma_m$ & symmetric group on $m$ elements &\pageref{n3}\\
$\overline{H}_{g,n,m}$
& locus of $m$-pointed stable curves of genus $g$,&\pageref{n3a}\\
& pluricanonically embedded in $\PP^N$\\ 
$[C_0]$& point in $\Hgn$  representing a curve $C_0\rightarrow \Spec(k)$&\pageref{n3b}\\
$\msbar_{g,m}$& moduli space of $m$-pointed stable curves of genus $g$ &\pageref{n4}\\
$\mgstack{g,m}$&moduli stack of $m$-pointed stable curves of genus $g$ &\pageref{n5}\\
$\mgamma$ &$\mgamma$-pointed stable curve of some genus &\pageref{n6}\\
$ \msbar_{g,\mgamma}$&  moduli space of $\mgamma$-pointed stable curves of genus $g$&\pageref{n7}\\
$\Hgmp$& locus of embedded $\mgamma$-pointed curves of genus $g$&\pageref{n8}\\
$\mgstack{g,\mgamma}$&moduli stack of $\mgamma$-pointed stable curves of genus $g$ &\pageref{n9}\\
$ C_0 $& fixed stable curve $f_0: C_0 \rightarrow \Spec(k)$ of genus $g$ &\pageref{n10}\\
&with labelled nodes $P_1,\ldots , P_{3g-3}$ \\
$\overline{D}(C_0;P_1)$& component of $\msbar_g^{(3g-4)}$ distinguished by $C_0$ and $P_1$ &\pageref{n11}\\ 
$C_0^+, \; C_0^-$ & $m$-pointed subcurves of $C_0$ of genus $g^+$ and $g^-$, &\pageref{n12}\\
& with set of marked points $\calq = \{ Q_1,\ldots,Q_m\}$ \\
$C_i^-$& connected component of $C_0^-$ of genus $g_i^-$, &\pageref{n13}\\
&with $m_i$ marked points, for $i=1,\ldots, r$
\end{tabular}

\begin{tabular}{llr}
$\Aut(C_0^-)$& automorphism group of $C_0^-$, respecting $\calq$ pointwise& \pageref{n15}\\
$\Aut(C_0^-;\calq)$& automorphism group of $C_0^-$, respecting $\calq$ as a set&\pageref{n16}\\
$\calq$& set of marked points of $C_0^-$ &\pageref{n16a}\\
$\Gamma(C_0^-;\CQ)$& subgroup of $\Sigma_m$ induced by $\Aut(C_0^-;\calq)$&\pageref{n17}\\
$\Gamma^\circ(C_0^-;\CQ)$& subgroup induced by $\Aut(C_0^-;\calq)$ and $\Aut(C_0^+;\calq)$&\pageref{n18}\\
$\Aut_\Gamma(C_0^-;\calq)$& subgroup of $\Aut(C_0^-;\calq)$ induced by $\Gamma$&\pageref{n19}\\
$\Aut_\Gamma(C^+;\calq)$ & subgroup of $\Aut(C_0^+;\calq)$ induced by $\Gamma$&\pageref{m1}\\
$ \hxbar$&  locus of curves with splitting automorphism group &\pageref{m2}\\&  inside $\hilb{g^+,n,m}$  having the right number of nodes\\
$\thx$& fixed morphism  $\thx :  \hxbar \rightarrow \hilb{g,n,0}$ &\pageref{n14}\\
$\mxbar$& locus of curves with splitting automorphism group&\pageref{m3}\\&  inside $\msbar_{g^+,m}$  having the right number of nodes\\[1mm]
$\mgprimebar$& subscheme of $\msbar_{g^+,m}$ defined by  $\dpobar$ via $\overline{\Theta}$&\pageref{n21}\\[1mm]
$\overline{\Theta}$& natural morphism $\overline{\Theta} : \: \: \msbar_{g^+,m} \rightarrow \msbar_g$&\pageref{n20}\\
$\hgprimebar$ & subscheme of $\hilb{g^+,n,m}$ defined by $\mgprimebar$&\pageref{n22}\\[1mm]
$\mgammaxbar$&  locus of curves with splitting automorphism group&\pageref{m4}\\&  inside $\mgbarquot$  having the right number of nodes\\[1mm]
$\hgammaxbar$&  subscheme of $\hilb{g^+,n,m}$ defined by $\mgammaxbar$ &\pageref{m4}\\[1mm]
$\mgprimebarquot$& subscheme of $\mgbarquot$ defined by  $\dpobar$ via ${\Theta}$&\pageref{m4}\\[1mm]
$\hgprimebarquot$ & subscheme of $\hgbarquot$ defined by $\mgprimebarquot$&\pageref{m4}\\[1mm]
$\dxbar$& locus of curves with splitting automorphism group&\pageref{m5}\\&  inside $\dpobar$ \\
$ \Theta$& induced embedding  $ \Theta : \hxbar / \Gamma(C_0^-;\CQ)   \rightarrow  \hilb{g,n,0}$&\pageref{n23}\\
$\overline{K}(C_0;P_1)$& reduced preimage of  $\dpobar$ in $\hilb{g,n,0}$&\pageref{n24}\\
$\kxbar$& reduced preimage of  $\dxbar$ in $\hilb{g,n,0}$  &\pageref{m6}\\
$\FDPtilde$& preimage substack  in $\mgstack{g}$ defined by $\dpobar$&\pageref{m7}\\
$\FDP$& reduced preimage substack of $\dpobar$ in $\mgstack{g}$ &\pageref{n25}\\
$\FDPc$& reduced preimage substack of $\dxbar$ in $\mgstack{g}$ &\pageref{m8}\\
$\cmgprimebarquot$ &reduced preimage substack of $\mgprimebarquot$ in $\cmgbarquot$ &\pageref{k1}\\[1mm]
$\cmgammaxbar$& reduced preimage substack of $\mgammaxbar$ in $\cmgbarquot$&\pageref{k2}\\[1mm]
$C_0^{++}, \; C_0^{--} $&$\mu$-pointed subcurves of $C_0$ of genus $g^{++}$ and $g^{--}$, &\pageref{n28}\\ 
& with set of marked points $\calr = \{ R_1,\ldots,R_\mu\}$ \\
$C_i^{++}$& connected component of $C_0^{++}$ of genus $g_i^{++}$, &\pageref{n29}\\
&with $\mu_i^{++}$ marked points, for $i=1,\ldots, r^{++}$\\
$C_j^{--}$& connected component of $C_0^{--} $ of genus $ g_j^{--}$, &\pageref{n30}\\
&with $\mu_j^{--}$ marked points, for  $j=1,\ldots, r^{--}$
\end{tabular}

\begin{tabular}{llr}
$\Aut(C_0^{--};{\cal R})$&automorphism group of $C_0^{--}$, respecting $\calr$ as a set&\pageref{n31}\\ 
$\Gamma(C_0^{--};\calr)$&subgroup of $\Sigma_\mu$ induced by $\Aut(C_0^{--};\calr)$&\pageref{n32}\\
$\overline{H}_{g',n',\mu'}$& locus of $\mu'$-pointed prestable curves of genus $g'$&\pageref{n33}\\
& pluricanonically embedded in $\PP^{N'}$\\ 
$ \hxxbar$&  locus of curves with splitting automorphism group\\[-1mm]&  inside $\hilb{g^{++},n,\mu}$  having the right number of nodes &\pageref{k3}\\[1mm]
$\thxx$& fixed morphism $ \thxx : \hxxbar \longrightarrow \overline{H}_{g,n,0}$&\pageref{n34}\\
$\hgammaxxbar$& locus of curves with splitting automorphism group\\[-1mm]&  inside $\hilb{g^{++},n,\mu/\Gamma}$  having the right number of nodes &\pageref{n35}\\
$\hat{D}(C_0;P_1)$& partial compactifiaction of $D(P_1)$ at the point $[C_0]$&\pageref{n36}\\
$\dxxbar$& locus of curves with splitting automorphism group&\pageref{k4}\\
&  inside $\hat{D}(C_0;P_1)$ \\
$\hat{K}(C_0;P_1)$& reduced preimage of  $\hat{D}(C_0;P_1)$ in $\hilb{g,n,0}$&\pageref{n37}\\
$\kxxbar$& reduced preimage of $\dxxbar$ in $\hilb{g,n,0}$ &\pageref{k5}\\[1mm]
$\msbar^{\,ps}_{g^{++},\mu}$&moduli space of $\mu$-pointed prestable curves &\pageref{n40}\\
&of genus $g^{++}$\\
$\mxxbar$&reduced image of $\hxxbar$ in $\msbar^{\,ps}_{g^{++},\mu}$&\pageref{n41}\\[1mm]
$\mgammaxxbar$&quotient  of  $\mxxbar$ by $\hat{\Gamma}$&\pageref{n42}\\
$\hat{\Theta}$& induced embedding  $ \hat{\Theta} :  \hgammaxxbar \rightarrow \kxxbar$&\pageref{n43}\\[1mm]
$\FDPxx$& reduced preimage substack of $\dxxbar$ in $\mgstack{g}$ &\pageref{n45}\\
$\cmgammaxxbar$& stack of $\mu/\hat{\Gamma}$-pointed prestable curves of genus $g^{++}$ &\pageref{n46}\\[-1mm]
& with moduli space $\mgammaxxbar$   \\
$\doo$ & intersection of $\dx$  and $\dxxbar$ & \pageref{cc1}\\
$\cdoo$ & substack of $\FDPopen$ corresponding to $\doo$ & \pageref{cc2}\\
$\koo$ & intersection of $\kx$  and $\kxxbar$ & \pageref{cc3}\\[1mm]
$\hoogamma$ & reduced preimage of $\koo$  in  $H_{g^+,n,\mgamma} $& \pageref{cc5a}\\
$\hoogammahat$ & reduced preimage of $\koo$  in  ${H}_{g^{++},n,\mu/\hat{\Gamma}} $& \pageref{cc5a}\\[1mm]
$ {\hat{\cal D}(C_i;P_{i,1})}$ & reduced preimage of $ {\hat{D}(C_i;P_{i,1})}$  in  $\FDP$& \pageref{lb1}\\[1mm]
$C_0^{+-}$& subcurve of $C_0$ with $\mu^{+-}$ marked points in $\cal N$&\pageref{n47}\\ 
$\CQ^{+-}, \, \CQ^{--}$&distribution of nodes of $C_0$ in $\calq$ &\pageref{n48}\\
$\CR^{+-}, \, \CQ^{--}$&distribution of nodes of $C_0$ in $\calr$ &\pageref{n49}\\
$\Aut(C;Y_i,\underline{Z}_j)$&automorphisms of $C$, stabilizing $Y_i$ as sets,&\pageref{n50}\\
& and $Z_j$ pointwise\\
$\Gamma(C;Y_i,\underline{Z}_j)$&permutation subgroup induced by  $\Aut(C;Y_i,\underline{Z}_j)$ &\pageref{n51}
\end{tabular}

\begin{tabular}{llr}
$\dii$ & locus of curves with threefold splitting groups  & \pageref{cd1}\\
& of automorphisms in $\dpo$\\
$\cdii$ & substack of $\FDPopen$ corresponding to $\dii$ & \pageref{cd2}\\
$\kii$ & reduced preimage of $\dii$ in $K(C_0;P_1)$ & \pageref{cd3}\\[1mm]
$\hiigamma$ & reduced preimage of $\kii$  in  $H_{g^+,n,\mgamma} $& \pageref{cd5b}\\
$\hiigammahat$ & reduced preimage of $\kii$  in  ${H}_{g^{++},n,\mu/\hat{\Gamma}} $& \pageref{cd5b}\\[1mm]
$\Theta^{+-}$& fixed morphism $\Theta^{+-} :\hii \longrightarrow \hiihat$&\pageref{n52}\\
${\Theta}^\diamond$&induced embedding of the quotient  &\pageref{n53}\\
&$\hii/\Gamma(C_0^{+-};\CQ^{+-},\underline{\CR}^{+-}) $ into $\hiihat$\\
$\overline{\Theta}{}^\diamond$& induced morphism from the scheme &\pageref{n54}\\
&$ {H}^\diamond_{g^+,n,m/\Gamma(C_0^-;\CQ)}$ to $ \hat{H}{}^\diamond_{g^{++},n,\mu/\Gamma(C_0^{--};\CR)}$\\[2mm]
${\Omega}$&morphism of stacks $ \cmgammadetailx \rightarrow \cdx$&\pageref{n55}\\[1mm]
$\hat{M}_{g^+,m/\Gamma(C_0^-;\CQ)}^{\: \dagger}$&partial compactification of $\hat{M}^\times_{g^{+},m/\Gamma(C_0^{-};\CQ)}$&\pageref{n57}\\[2mm]
$\hat{H}^{\: \dagger}_{g^+,m/\Gamma(C_0^-;\CQ)}$&reduced preimage of $\hat{M}_{g^+,m/\Gamma(C_0^-;\CQ)}^{\: \dagger}$ in the &\pageref{n56}\\
& scheme $\overline{H}_{g^+,m/\Gamma(C_0^-;\CQ)}$\\[2mm]
$d_C$&distance of points on a curve $C$&\pageref{n58}
\end{tabular}


%
%

\thispagestyle{plain}

\pagestyle{headings}

\end{document}